\documentclass[12pt]{amsart}

\usepackage{amsmath}
\usepackage{amsfonts}
\usepackage{amssymb}
\usepackage{amscd}
\usepackage{graphicx}
\usepackage[abbrev,alphabetic]{amsrefs}
\RequirePackage[dvipsnames,usenames]{color}
\usepackage{soul,xcolor}
\setstcolor{red}
\usepackage{stmaryrd}
\usepackage{mathtools}
\usepackage{booktabs}
\usepackage{multirow}
\newtagform{tiny}{\tiny(}{)}

\usepackage{mathtools}
\usepackage{hyperref}
\usepackage[margin=1.25in]{geometry}

\usepackage{amsthm}
\usepackage{comment}
\usepackage[all,cmtip]{xy}
\usepackage{tikz-cd}
\usetikzlibrary{cd}

\usepackage[all]{xy}

\makeatletter
\def\@tocline#1#2#3#4#5#6#7{\relax
  \ifnum #1>\c@tocdepth 
  \else
    \par \addpenalty\@secpenalty\addvspace{#2}%
    \begingroup \hyphenpenalty\@M
    \@ifempty{#4}{%
      \@tempdima\csname r@tocindent\number#1\endcsname\relax
    }{%
      \@tempdima#4\relax
    }%
    \parindent\z@ \leftskip#3\relax \advance\leftskip\@tempdima\relax
    \rightskip\@pnumwidth plus4em \parfillskip-\@pnumwidth
    #5\leavevmode\hskip-\@tempdima
      \ifcase #1
       \or\or \hskip 1em \or \hskip 2em \else \hskip 3em \fi%
      #6\nobreak\relax
    \hfill\hbox to\@pnumwidth{\@tocpagenum{#7}}\par
    \nobreak
    \endgroup
  \fi}
\makeatother

\newsavebox{\pullback}
\sbox\pullback{%
\begin{tikzpicture}%
\draw (0,0) -- (1ex,0ex);%
\draw (1ex,0ex) -- (1ex,1ex);%
\end{tikzpicture}}

\newsavebox{\pullbackdl}
\sbox\pullbackdl{%
\begin{tikzpicture}%
\draw (-1ex,0ex) -- (0ex,0ex);%
\draw (0ex,-1ex) -- (0ex,0ex);%
\end{tikzpicture}}

\newsavebox{\pushoutdr}
\sbox\pushoutdr{%
\begin{tikzpicture}%
\draw (-1ex,-1ex) -- (-1ex,0ex);%
\draw (-1ex,0ex) -- (0ex,0ex);%
\end{tikzpicture}}

\newcommand{\cyan}{\color{cyan}}

\newcommand{\cred}{\color{black}}
\newcommand{\ccred}{\color{black}}
\newcommand{\cccred}{\color{black}}

\newcommand{\rup}[1]{\lceil #1 \rceil}
\newcommand{\rdown}[1]{\lfloor #1 \rfloor}

\renewcommand{\mod}{\ \textrm{mod}\ }

\newcommand{\Z}{\mathbb{Z}}
\newcommand{\Q}{\mathbb{Q}}

\newcommand{\F}{\mathbb{F}}

\newcommand{\cHom}{\mathcal{H}om}

\newcommand{\bF}{\mathbb{F}}

\newcommand{\bQ}{\mathbb{Q}}

\newcommand{\bZ}{\mathbb{Z}}

\newcommand{\cO}{\mathcal{O}}

\newcommand{\MO}{\mathcal{O}}
\newcommand{\sO}{\mathcal{O}}
\newcommand{\PsiT}{T}

\newcommand{\m}{\mathfrak{m}}

\newcommand{\fp}{\mathfrak{p}}

\newcommand{\Proj}{\mathrm{Proj}}
\newcommand{\wt}{\widetilde}

\DeclareMathOperator{\Supp}{Supp}
\DeclareMathOperator{\Spec}{Spec}

\DeclareMathOperator{\Coker}{Coker}
\DeclareMathOperator{\Hom}{Hom}

\DeclareMathOperator{\Ext}{Ext}

\DeclareMathOperator{\Exc}{Exc}
\DeclareMathOperator{\Ex}{Exc}

\DeclareMathOperator{\Ker}{Ker}
\DeclareMathOperator{\reg}{reg}
\DeclareMathOperator{\Ann}{Ann}

\renewcommand{\div}{{\rm div}}

\newcommand{\mydot}{{{\,\begin{picture}(1,1)(-1,-2)\circle*{2}\end{picture}\ }}}

\theoremstyle{plain}
\newtheorem{theorem}{Theorem}[section]

\newtheorem{proposition}[theorem]{Proposition}
\newtheorem{prop}[theorem]{Proposition}
\newtheorem{lemma}[theorem]{Lemma}
\newtheorem{lem}[theorem]{Lemma}
\newtheorem{corollary}[theorem]{Corollary}
\newtheorem{cor}[theorem]{Corollary}

\newtheorem{claim}[theorem]{Claim}
\newtheorem*{claim*}{Claim}
\newtheorem{step}{Step}

\newtheorem{theoremA}{Theorem}

\theoremstyle{definition}
\newtheorem{definition}[theorem]{Definition}
\newtheorem{dfn}[theorem]{Definition}

\newtheorem{example}[theorem]{Example}
\newtheorem{notation}[theorem]{Notation}
\newtheorem{nothing}[theorem]{}
\newtheorem*{setup*}{Setup}

\theoremstyle{remark}
\newtheorem{remark}[theorem]{Remark}

\theoremstyle{plain}

\def\commentbox#1{\textcolor{Mahogany}%
{\footnotesize\newline{\color{Mahogany}\fbox{\parbox{\textwidth-15pt}{\textbf{comment: } #1}}}\newline}}



\numberwithin{equation}{theorem}

\title[Quasi-F-splittings in birational geometry III]
{Quasi-F-splittings in birational geometry III}

\author{Tatsuro Kawakami}
\address{Graduate School of Mathematical Sciences, University of Tokyo, 3-8-1 Komaba,
Meguro-ku, Tokyo 153-8914, Japan}
\email{kawakami@ms.u-tokyo.ac.jp}
\author{Teppei Takamatsu}
\address{Department of Mathematics, Faculty of Science,
Saitama University,
255 Shimo-Okubo, Sakura-ku,
Saitama-shi, Saitama 338-8570,
Japan}
\email{teppeitakamatsu.math@gmail.com}
\author{Hiromu Tanaka} 
\address{Department of Mathematics, 
Graduate School of Science, 
Kyoto University, 
Kyoto 606-8502, Japan} 
\email{tanaka.hiromu.7z@kyoto-u.ac.jp}

\author{Jakub Witaszek} 
\address{Northwestern University, Department of Mathematics, Lunt Hall, 2033 Sheridan Road, Evanston, IL 60208, USA}
\email{jakub.witaszek@northwestern.edu}
\author{Fuetaro Yobuko}
\address{Tokyo University of Science, Faculty of Science and Technology, Department of Mathematics}
\email{soratobumusasabidesu@gmail.com}
\author{Shou Yoshikawa}
\address{Institute of Science Tokyo, Tokyo 152-8551, Japan}
\email{yoshikawa.s.al@m.titech.ac.jp}

\begin{document}

\begin{abstract}
We prove that $\mathbb Q$-Gorenstein quasi-$F$-regular singularities are klt. 
To this end, we 
introduce quasi-test ideals. 
\end{abstract}


\maketitle

\setcounter{tocdepth}{2}

\tableofcontents


\section{Introduction}


One of the fundamental topics of research in both positive characteristic birational geometry and commutative algebra is the theory of Frobenius splittings and Frobenius regularity. Recently, the fifth author introduced a new notion in \cite{Yob23}, called \emph{quasi-$F$-splitting}, motivated by the theory of crystalline cohomology. It shares many properties with the usual $F$-splitting but is much less restrictive. In turn, in \cite{KTTWYY1} and \cite{KTTWYY2}, we engaged in the development of the theory of quasi-$F$-splittings in the context of birational geometry, from which we deduced new results on the liftability of singularities and extensions of differential forms. We refer to the introductions of \cite{Yob23} and \cite{KTTWYY1} for more information on quasi-$F$-splittings, and to \cite{KTY22} for Fedder's criterion for quasi-$F$-splittings. 

Having gained a basic understanding of quasi-$F$-splittings in the context of birational geometry, \cite{TWY} introduced and studied, amongst other things, the notion of quasi-$F$-regularity, generalising strong-$F$-regularity. Our paper continues the work of \cite{KTTWYY1}, \cite{KTTWYY2}, and \cite{TWY}, with a focus on the local theory of singularities.

The discovery of the connection between $F$-singularities and birational-geometric singularities, such as rational or Kawamata log terminal, triggered impactful collaborations between commutative algebraists and birational geometers. In our main result, we establish such connections for some quasi-$F$-singularities. In particular, we generalise the fact that $\bQ$-Gorenstein strongly $F$-regular singularities are Kawamata log terminal.


\begin{theoremA}
[Theorem \ref{thm:qFR to klt pair}]\label{theo:intro: test ideal}
Let $R$ be an $F$-finite normal Noetherian $\bQ$-Gorenstein domain over $\F_p$. 
If $R$ is quasi-$F$-regular, then it is klt.
\end{theoremA}


In order to prove Theorem \ref{theo:intro: test ideal}, 
we introduce an analogue $\tau^q(R)$ of the test ideal $\tau(R)$, which we call the \textit{quasi-test ideal}. 
This ideal measures how far the ring is from being quasi-$F$-regular.  A big part of this article is devoted to establishing foundational results on $\tau^q(R)$. 
For example, we prove the following. 
\begin{enumerate}
    \item[(a)] $\tau^q(R)$ commutes with localisation and completion 
    (Remark \ref{r q test ideal clear}). 
    \item[(b)] $\tau^q(R)=R$ if and only if $R$ is quasi-$F$-regular (Proposition \ref{prop: test ideal pair non-tilda}).
    \item[(c)] $\tau^q(R) \subseteq \mathcal{J}(R)$ 
    for the multiplier ideal $\mathcal{J}(R)$ (Corollary \ref{cor:trace under bir}). 
\end{enumerate}
 Theorem \ref{theo:intro: test ideal} is  a consequence of (b) and (c). 
 

The opposite implication to that in Theorem \ref{theo:intro: test ideal} does not hold in general, that is, there exist
klt singularities which are not quasi-$F$-regular 
(see \cite{KTTWYY2}*{Theorem 8.3}). 
On the other hand, it is tempting to expect that such an implication holds in low dimensions and for high characteristic. In \cite{KTTWYY2}, we proved that three-dimensional $\bQ$-Gorenstein klt singularities of characteristic $p>42$ are quasi-$F$-split. We can now generalise this result to prove their quasi-$F$-regularity.

\begin{theoremA}[Theorem \ref{t-3-dim-klt}]\label{intro-3-dim-klt}
Let $k$ be a perfect field of characteristic $p>42$ and 
let $R$ be a three-dimensional $\Q$-factorial klt ring of finite type over $k$. 
Then $R$ is quasi-$F$-regular. 
\end{theoremA}

In fact, quasi-$F$-regularity and quasi-test ideals $\tau^q(R,\Delta)$ can be defined for arbitrary pairs $(R,\Delta)$. Similarly to Theorem \ref{intro-3-dim-klt}, we can now generalise the fact that two-dimensional klt pairs are quasi-$F$-split to the quasi-$F$-regular case.
\begin{theoremA}[Theorem \ref{thm:QFR for klt}]\label{intro-2-dim-klt}
Let $k$ be a perfect field of characteristic $p>0$ and 
let $(R, \Delta)$ be a two-dimensional klt pair of finite type over $k$. 
Then $(R, \Delta)$ is quasi-$F$-regular.
\end{theoremA}

Finally, one of the fundamental properties of strongly $F$-regular singularities is that they are Cohen-Macaulay, and so it is natural to wonder if the same holds for quasi-$F$-regular singularities.

\begin{theoremA}
Let $R$ be an $F$-finite normal Noetherian $\bQ$-Gorenstein domain over $\F_p$. Assume that $R$ is quasi-$F$-regular and $\dim(R) \leq 3$. Then $R$ is Cohen-Macaulay.
\end{theoremA}
Alas, we do not know how to prove this statement in higher dimensions. Let us point out that the restriction on the dimension is not a consequence of any classification or the minimal model program, but arises from surprising cohomological constraints.


\subsection{Various classes of singularities} \label{ss:various-classes-intro}
Recall that 
{\cred a Noetherian $F$-finite $\F_p$-algebra} 
$R$ is \emph{$n$-quasi-$F$-split} for an integer $n>0$ if given a pushout diagram {of $W_nR$-modules}
\[
\begin{tikzcd}
W_nR \ar{r}{F} \ar{d}{{\cccred \mathbf{R}}^{n-1}} & F_*W_nR \ar{d}\\
R \ar{r}{\Phi_{R,n}} & Q_{R,n},
\end{tikzcd}
\]
the $R$-module homomorphism $\Phi_{R,n} \colon R \to Q_{R,n}$ splits. This is equivalent to the surjectivity of 
\[
\Hom_R(\Phi_{R,n},R) \colon \Hom_R(Q_{R,n}, R) \to R.
\]
In order to incorporate higher powers of Frobenius, \cite{TWY} considers the pushout
\begin{equation} \label{eq:intro-fe}
\begin{tikzcd}
W_nR \ar{r}{F^e} \ar{d}{{\cccred \mathbf{R}}^{n-1}} & F^e_*W_nR \ar{d}\\
R \ar{r}{\Phi^e_{R,n}} & Q^e_{R,n},
\end{tikzcd}
\end{equation}
and defines $R$ to be \emph{$n$-quasi-$F^e$-split} if 
\begin{equation} \label{eq:intro-fe-map}
\Hom_{W_nR}(\Phi^e_{R,n},W_n\omega_R(-K_R)) \colon \Hom_{W_nR}(Q^e_{R,n},W_n\omega_R(-K_R)) \to R
\end{equation}
is surjective, 
{\cred where $W_n\omega_R(-K_R) := \Hom_{W_nR}(W_nR(K_R), W_n\omega_R)$, 
$W_n\omega_R$ is the dualising $W_nR$-module (cf.\ Theorem \ref{thm:X-exist dualizing complex}(1)), 
and $W_nR(K_R) := \Gamma( \Spec R, W_n\MO_{\Spec R}(K_R))$ 
for the Witt divisorial sheaf $W_n\MO_{\Spec R}(K_R)$ 
(Subsection \ref{ss-notation})}. 
The reader should be warned that $Q^e_{R,n}$ is not an $R$-module anymore if $e>1$. We have the following two implications:
\begin{align*}
\text{ $n$-quasi-$F^e$-split } &\implies \text{ $(n+1)$-quasi-$F^e$-split.}\\
\text{ $n$-quasi-$F^e$-split } &\implies \text{ $n$-quasi-$F^{e-1}$-split.}
\end{align*}
One can formulate similar definitions for pairs $(R,\Delta)$. 
\begin{definition}
{\cred 
Let $n$ be a positive integer and 
let $R$ be an $F$-finite normal Noetherian domain over $\F_p$. 
We say that} $R$ is \emph{n-quasi-$F$-regular} if for every {\cred effective} Weil divisor $E$, there exists $\epsilon>0$ such that $(R,\epsilon E)$ is $n$-quasi-$F^e$-split for every $e \gg 0$ (depending on $E$ and $\epsilon$).

Then $R$ is \emph{quasi-$F$-regular} if it is $n$-quasi-$F$-regular for some integer $n>0$. Similarly, one defines the notion of \emph{quasi-$+$-regularity}
 (``quasi-splinter'') by replacing $F^e$ by every possible finite surjective morphism in (\ref{eq:intro-fe}) and (\ref{eq:intro-fe-map}).
\end{definition} 

\begin{remark} \label{rem:how-to-prove-klt3fold-qfsplit} Before proceeding, let us briefly sketch how three-dimensional klt singularities in large characteristic were proven to be quasi-$F$-split in \cites{KTTWYY1,KTTWYY2} so that we can later sketch the proof of Theorem \ref{intro-3-dim-klt}. First, one considers a special projective birational morphism $\pi \colon Y \to \Spec R$, called a plt blow-up, which has a property that the exceptional locus is a (log) del Pezzo surface. Then one shows that this (log) del Pezzo surface is globally quasi-$F$-split if $p > 42$. By using a special case of global inversion of adjunction, one deduces that $Y$ is globally quasi-$F$-split, which, by pushing forward, shows that $\Spec R$ is quasi-$F$-split.

Since necessary variants of inversion of adjunction for global $+$-regularity have been already established in \cite[Corollary 6.8]{TWY}, in order to conclude Theorem \ref{intro-3-dim-klt}, we need to do two things. First, we need to compare quasi-$F$-regularity with quasi-$+$-regularity (see Theorem \ref{theo:intro: q+R to qFR}). Second we need to show that (log) del Pezzo surfaces are globally quasi-$+$-regular (or quasi-$F$-regular) for $p > 42$ and not just quasi-$F$-split (see Theorem \ref{intro-Fano}).
\end{remark}

Singh has shown in \cite{Singh99} that a $\bQ$-Gorenstein ring is 
strongly $F$-regular 
if and only if it is $+$-regular 
(see also {\cred \cite{HH94-2} and} \cite{BST15}), and 
the same statement for non-$\bQ$-Gorenstein rings is a big open problem in the theory of $F$-singularities. We generalise Singh's result to the quasi setting. Note that the easier implication, from left to right, has already been established in \cite[Proposition 4.9]{TWY}. 
\begin{theoremA}[Theorem \ref{thm: rel q+R, qFrat and qFR}]\label{theo:intro: q+R to qFR}
Let $R$ be an $F$-finite normal Noetherian $\bQ$-Gorenstein domain over $\bF_p$. Then $R$ is quasi-$F$-regular if and only if $R$ is quasi-$+$-regular.
\end{theoremA}

The interaction between the parameters $n$ and $e$ in the definition of $n$-quasi-$F^e$-splitting, as well as with the divisor $E$ and $\epsilon>0$ in the definition of quasi-$F$-regularity, is very subtle (see the discussion below the sketch of the proof of Theorem \ref{intro-Fano}). The following fundamental result streamlines the definition of quasi-$F$-regularity significantly.\footnote{Let us briefly explain the key idea of the proof of this and many similar results in this paper. In the usual theory of $F$-splittings and $F$-regularity, various definitions become often equivalent up to multiplication by $t$ for a test element $t$. In the quasi-setting for height $n$, we show by induction on $n$ that the same is often true up to multiplication by $t^2$. Here one $t$ allows to reduce the induction from height $n$ to height $n-1$, in which case the other $t$ turns into $t^p$, where clearly $p \geq 2$, so the induction can continue to be run. Note that for technical reasons, we sometimes consider $t^4$ as opposed to $t^2$.} 

\begin{theoremA}[Theorem \ref{thm:best-def-quasi-F-regular-Section4}]\label{intro:thm:best-def-quasi-F-regular}
{\cred 
Let $R$ be an $F$-finite normal Noetherian $\bQ$-Gorenstein domain over $\F_p$. Let $D$ be an {effective} {Cartier} divisor on $\Spec R$ such that 
the inclusion $I_D \subseteq \tau(R)$ holds 
for the test ideal $\tau(R)$ of $R$ and the defining ideal $I_D$ of $D$. 
Then $R$ is quasi-$F$-regular 
if and only if there exist integers $e>0$ and $n>0$ 
such that 
 \[
 \Big( R, \frac{4}{p^e}D\Big)
 \]
is $n$-quasi-$F^e$-split.} 
\end{theoremA}

\noindent {\cred 
For an effective Cartier divisor $E$ containing the non-regular locus of $\Spec R$, 
Theorem \ref{intro:thm:best-def-quasi-F-regular} is applicable for $D := m E$ with $m \gg 0$. 
In this case, after enlarging $m$ further (specifically, replacing $m$ by $4m$), 
it is enough to check that $(R,\frac{1}{p^e}D)$ is $n$-quasi-$F^e$-split for some integers $n>0$ and $e > 0$.}

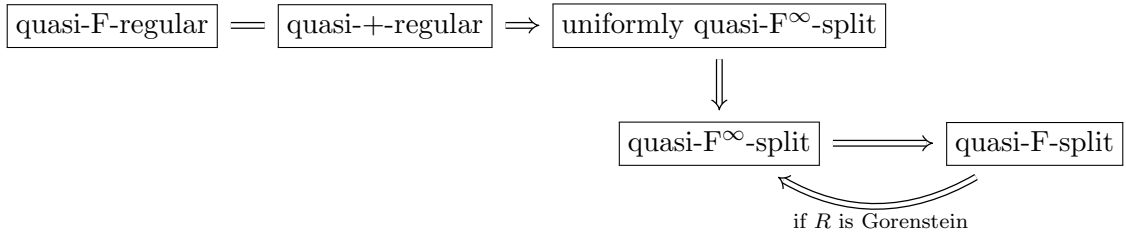
\begin{figure}[h]  \label{figure-diagram-notions}
\[
{\small \begin{tikzcd}[column sep = small]
\fbox{\text{quasi-F-regular}} \ar[r, equal] & \fbox{\text{quasi-+-regular}}  \ar[r, Rightarrow]  & \fbox{$\text{{uniformly} quasi-F}^\infty\text{-split}$} \ar[d, Rightarrow]  & \\
&  & \fbox{$\text{quasi-F}^\infty\text{-split}$} \ar[Rightarrow]{r}  & \ar[Rightarrow, bend left = 30]{l}{\text{ if $R$ {is} {\cccred Gorenstein}}} \fbox{\text{quasi-F-split}}
\end{tikzcd}}
\]
\caption{Relation between singularities when $R$ is $\bQ$-Gorenstein.}
\end{figure}

Finally, we address how 
to deduce quasi-$F$-regularity from quasi-$F$-splitting 
 in the case of (log) Fano varieties 
(cf.\ Remark \ref{rem:how-to-prove-klt3fold-qfsplit}). This is particularly meaningful, 
because quasi-$F$-splitting is the weakest notion and quasi-$F$-regularity is the strongest one amongst what we have introduced (see Figure \ref{figure-diagram-notions}). The result is interesting even 
for {\cred smooth Fano varieties.} 
Also note that, in fact, we show that the same result holds for strongly $F$-regular pairs with standard coefficients. 

\begin{theoremA}[Theorem \ref{thm:qFs to qFr}, 
{Corollary \ref{qFs to qFr for am Fano}}]\label{intro-Fano}
Let $k$ be an $F$-finite field of characteristic $p>0$. 
Let $X$ be a projective {normal} $\bQ$-factorial 
{strongly $F$-regular} variety {over $k$ such that $-K_X$ is ample}.
Then $X$ is 
quasi-$F$-split if and only if it is globally quasi-$F$-regular.
\end{theoremA}

{\cred \noindent 
Here we warn the reader that strong $F$-regularity is a local notion (e.g., a smooth variety is strongly $F$-regular), whilst quasi-$F$-splitting and global quasi-$F$-regularity are global ones.}

\begin{proof}[Sketch of the proof]
The proof consists of two parts. First, using Cartier operators, we show that $X$ is $n$-quasi-$F^e$-split for some fixed $n>0$ and all $e>0$. In fact, we show a bit more: that $(X,\frac{1}{p^e}H)$ is $n$-quasi-$F^e$-split for a fixed, carefully chosen {\cred effective} divisor $H$ and all $e \gg 0$.

Second, in the context of (log) Fano varieties, we establish a connection between quasi-$F$-regularity of a projective variety and its cone (see Corollary \ref{cor: ample cartier case}). Thus by taking a cone of $X$ and applying Theorem \ref{intro:thm:best-def-quasi-F-regular}, we can deduce that $X$ is globally quasi-$F$-regular.
\end{proof}



In 
{\cred \cite{TWY},}
two variants of quasi-$F^\infty$ splittings were defined:
\begin{enumerate}
\item $R$ is \emph{quasi-$F^\infty$-split} if  for every integer $e>0$, there exists an integer $n>0$ such that $R$ is $n$-quasi-$F^e$-split.
\item $R$ is \emph{uniformly quasi-$F^\infty$-split} if there exists an integer $n>0$ such that $R$ is $n$-quasi-$F^e$-split for every integer $e>0$.
\end{enumerate}
Note that quasi-$F$-regular rings are uniformly quasi-$F^\infty$-split. On the other hand, elliptic curves are $n$-quasi-$F^e$-split if and only if $e \leq n+1$; in particular, they are quasi-$F^\infty$-split but not uniformly quasi-$F^\infty$-split. 

What was not known before however is whether quasi-$F$-split rings are automatically quasi-$F^\infty$-split.

\begin{theoremA}[{Corollary \ref{cor:quasiF=quasiFinfty}}]\label{intro-qF-qFinfty}
Let $R$ be an $F$-finite normal Noetherian 
{\cred Gorenstein} domain over $\bF_p$. 
Then $R$ is quasi-$F$-split if and only if it is quasi-$F^\infty$-split.
\end{theoremA}
\noindent In fact, if $R$ is $n$-quasi-$F$-split, then it is automatically $(ne -e +1)$-quasi-$F^{e+1}$-split for every $e>0$. We do not know however if the assumption on the {\cccred Gorensteinness} of $R$ can be dropped in the above theorem.

Combining the above work, we can now significantly extend the diagram of singularities from {\cred \cite[p.\ 5]{TWY}.}


One can also extend this diagram by adding notions of feeble quasi-$F$-regularity or non-uniform quasi-$F$-regularity, but these notions are now known to agree with the usual quasi-$F$-regularity by Theorem \ref{thm:chara qFr pair2} and Remark \ref{remark:non-uniform-quasi-F-regularity}. 

Also note that in a forthcoming paper, Sato-Takagi-Yoshikawa will prove that quasi-$F^\infty$-split  singularities are log canonical. In particular, combined with Theorem \ref{intro-qF-qFinfty}, this addresses a long-standing question, at least in the Gorenstein case, whether the usual quasi-$F$-split singularities are log canonical.

\subsection{Construction of {\cred quasi-test submodules and} quasi-test ideals}\label{ss intro q F rat} 
For simplicity, we will start by explaining how to construct the quasi-test submodule $\tau_n(\omega_R)$ of $\omega_R$. When $R$ is Gorenstein, this object agrees with the quasi-test ideal $\tau_n(R)$.

{\cred Fix an integer $n>0$. 
Let $R$ be an $F$-finite Noetherian integral domain $R$ of characteristic $p>0$. 
Set $R^\circ = R  \backslash \{0\}$.}  
    In the spirit of 
    {\cred \cite{TWY}}
    and the discussion in Subsection \ref{ss:various-classes-intro}, one can define $\tau_n(\omega_R)$ as follows:
\begin{equation} \label{eq:intro-taun}
\tau_n(\omega_R) := \bigcap_{c \in R^\circ} \sum_{e_0>0} \bigcap_{e \geq e_0} {\cred {\rm Im}}(\Hom_{W_nR}(Q^e_{R,n},W_n\omega_R(-K_R)) \xrightarrow{(\star)} R),
\end{equation}
where $(\star)$ is the composition of $\Hom_{W_nR}(\Phi^e_{R,n},W_n\omega_R(-K_R))$ from (\ref{eq:intro-fe-map}) and multiplication by the Teichm\"uller lift $[c^{p^{e-e_0}}]$. In this subsection, we shall present a different perspective on $\tau_n(\omega_R)$, coming from commutative algebra, and which has some additional benefits. We eventually show in {Corollary \ref{cor:Qdefs-of-quasi-tight-submodule-PAIR}} that this different approach agrees with (\ref{eq:intro-taun}). 
{The key idea that allows us to show Theorem \ref{theo:intro: test ideal} is to redefine $\tau_n(\omega_R)$ in terms of a new object $\tau(W_n\omega_R)$ which behaves more naturally in terms of birational morphisms.}  For more details, we refer to Section \ref{s q-F-rat}.


We then define 
$\tau(W_n\omega_R)$ 
as the smallest co-small ${\cccred T_n}$-stable $W_nR$-submodule of $W_n\omega_R$, {\cred although its existence is not trivial (Proposition \ref{prop: test submod}(1))}. 
Here we say that a $W_nR$-submodule $M \subseteq W_n\omega_R$ is 
\begin{itemize}
\item {\em co-small} if there exists $c \in R^{\circ}$ such that 
$[c] \cdot W_n\omega_R \subseteq M$, where $[c] \in W_nR$ denotes the Teichm\"uller lift of $c$; 
\item {\em ${\cccred T_n}$-stable} if $T_n(F_*M) \subseteq M$, 
where $T_n : F_*W_n\omega_X \to W_n\omega_X$ is the trace map, that is, 
the $W_nR$-module homomorphism obtained by applying $\Hom_{W_nR}(-, W_n\omega_R)$ to Frobenius $F: W_nR \to F_*W_nR$. 
\end{itemize}
{Note that one can also construct $\tau(W_n\omega_R)$ via appropriate images of Frobenius trace maps (see Proposition \ref{prop: test submodule gen pair} and Proposition \ref{p-tau-WnX}).

Using $\tau(W_n\omega_R)$,} we introduce the {\em quasi-test $R$-submodule} $\tau_n(\omega_R)$ of $\omega_R$, which is given as the inverse image of $\tau(W_n\omega_R)$ by $({\cccred \mathbf{R}}^{n-1})^* : \omega_R \to W_n\omega_R$: 
\begin{equation}\label{e intro tau omega}
\begin{tikzcd}
W_n\omega_R \arrow[r, phantom, sloped, "\supseteq"]  &  \tau(W_n\omega_R)\\
\omega_R \ar[hook]{u}{({\cccred \mathbf{R}}^{n-1})^*} \arrow[r, phantom, sloped, "\supseteq"] & \tau_n(\omega_R) \ar[hook]{u}
\end{tikzcd}
\quad \text{where} \quad 
\tau_n(\omega_R) := (({\cccred \mathbf{R}}^{n-1})^*)^{-1}(\tau(W_n\omega_R)). 
\end{equation}
We say that $R$ is {\em quasi-$F$-rational} if 
there exists an integer $n>0$ such that 
$R$ is {\em $n$-quasi-$F$-rational}, i.e., 
$\tau_n(\omega_R) =\omega_R$ and $R$ is Cohen-Macaulay. 


When $(R, \m)$ is a local ring, it is often more convenient to work with the dual statements invol{\cred v}ing local cohomology 
 {\cred via Matlis duality}. In our setup, we  
 introduce two modules 
${\wt{0^*_n}}$ and $0^*_n$, 
which are summarised in the following square: 
\begin{equation}\label{e intro 0^*}
\begin{tikzcd}
H_{\m}^d(W_nR) \arrow[r, phantom, sloped, "\supseteq"] \arrow[d, two heads, "{\cccred \mathbf{R}}^{n-1}"'] & \wt{0^*_n} \arrow[d, rightarrow, two heads]
\\
H_{\m}^d(R) \arrow[r, phantom, sloped, "\supseteq"] & 0^*_n
\end{tikzcd}
\quad 
{\cred \text{where} \quad 
0^*_n := {\cccred \mathbf{R}}^{n-1}( \wt{0^*_n}).}
\end{equation}
This square (\ref{e intro 0^*}) corresponds to the square 
(\ref{e intro tau omega})  via Matlis duality, for example, $\Hom_R(0^*_n, E) \simeq \omega_R/\tau_n(\omega_R)$. 
As their notations indicate, 
$\wt{0^*_n}$ and $0^*_n$ are obtained as analogues of the usual tight closure $0^*_M$ of $0$ in an $R$-module $M$. 
For example, $0^*_n$ is defined as follows: 
\[
0^*_n := \bigcup_{\substack{c \in R^{\circ},\\ e_0 \in \Z_{>0}}} \bigcap_{e \geq e_0} K^{e, c}_{n} \subseteq H^d_\m(R), 
\]
where 
    \[
    K^{e, c}_{n} := \Ker\left( H^d_\m(R) \xrightarrow{\Phi^{e, c}_{R,n}} H^d_\m(Q^{e,c}_{R,n})\right)
    \]
and $\Phi^{e, c}_{R,n}$ and $Q^{e,c}_{R,n}$ are given by the following pushout diagram: 
\[
\begin{tikzcd}
W_nR \arrow[r, "F^e"] \arrow[d, "{\cccred \mathbf{R}}^{n-1}"'] & F_*^eW_nR \arrow[d] \arrow[r, "{ \cdot F_*^e[c]}"] & F_*^eW_nR \arrow[d]\\
R \arrow[r, "\Phi^e_{R, n}"] \arrow[rr, "\Phi^{e, c}_{R, n}"', bend right] & Q^e_{R, n} \arrow[r] & Q^{e, c}_{R, n}. 
\end{tikzcd}
\]





\subsubsection{Definition of $\tau^q({\cred R}, \Delta)$}\label{ss intro q test ideal}

We now overview how to define the quasi-test ideal $\tau^q(R, \Delta)$, which is the main object of interest of this article. 
Let $R$ be an $F$-finite normal integral domain of characteristic $p>0$ 
and let  $\Delta$ be an effective $\Q$-divisor on $\Spec R$ such that $K_R+\Delta$ is $\Q$-Cartier. 
{\cred For $D := K_R + \Delta$,} we introduce the log versions of $0^*_n$ and $\tau_n(\omega_R)$, specifically
the {\em $n$-quasi-tight closure} $0^*_{D, n} \subseteq H^d_{\m}(R(D))$ and 
the {\em $n$-quasi-test submodule} $\tau_n(\omega_R, D) \subseteq \omega_R(-D) := \cHom_R(R(D), \omega_R)$, respectively.  
Then the {\em $n$-quasi-test ideal} $\tau_n(R, \Delta)$ is defined as 
\[
\tau_n(R, \Delta) := \tau_n(\omega_R, K_R+\Delta) \subseteq R(\rup{-\Delta}) \simeq \omega_R(\rup{-(K_R+\Delta)}). 
\]
It is easy to see that $\{\tau_n(R, \Delta)\}_{n \geq 1}$ 
forms an ascending chain of ideals 
\[
\tau_1(R, \Delta) \subseteq \tau_2(R, \Delta) \subseteq \cdots \subseteq R. 
\]
Since $R$ is a Noetherian ring, we have 
\[
\tau^q(R, \Delta) := \bigcup_{n=1}^{\infty}\tau_n(R, \Delta) = \tau_n(R, \Delta)
\]
for $n \gg 0$.

In constrast to the usual test ideal $\tau(R, \Delta)$, we define the $n$-quasi-test ideal $\tau_n(R, \Delta)$ only through the $n$-quasi-test submodule $\tau_n(\omega_R, D)$. This definition is advantageous due to its simplicity, but it requires the assumption that $K_R + \Delta$ is $\Q$-Cartier.

{\cccred
\subsection{Remark}\label{ss remark}
In this short remark, we clarify the difference between $F$-regularity and quasi-$F$-regularity by collecting several known examples.
Let $k$ be a perfect field of characteristic $p>0$.
\begin{example}
\label{eg:rdp}
By Theorem \ref{intro-2-dim-klt}, any two-dimensional klt pair of finite type over $k$ is quasi-$F$-regular.
On the other hand, if $p\in \{2,3,5\}$, a rational double point 
\[
\Spec k[x,y,z]/(x^2 + y^3 +z^5)
\]
is not $F$-pure (and hence not $F$-regular) (cf.\ \cite[(4.4)]{hara98}).
\end{example}

\begin{example}
Smooth Fano threefolds over $k$ with Picard number $>1$ are globally quasi-$F$-regular by \cite[Theorem~E]{Kawakami-Tanaka4} and Theorem~\ref{thm:qFs to qFr}, whereas they are not necessarily globally  $F$-regular in characteristics $p = 2,3,5$ (see \cite[Theorem~G and Remark~1.2]{Kawakami-Tanaka4}).

Similarly, smooth del Pezzo varieties over $k$ are globally quasi-$F$-regular by \cite[Theorem~A]{Kawakami-Tanaka} and Theorem~\ref{thm:qFs to qFr}, while they are not necessarily globally $F$-regular in characteristics $p = 2,3,5$ (cf.\ \cite[Theorem~B]{Kawakami-Tanaka}).


\end{example}

\begin{example}
By Theorem~\ref{intro-3-dim-klt}, any $\Q$-factorial klt threefold singularity over $k$ is quasi-$F$-regular if $p>42$.
Moreover, by Corollaries~\ref{c klt dP QFR} and \ref{cor:QFR for log dP}, any klt del Pezzo surface (resp.\ log del Pezzo surface with standard coefficients) is globally quasi-$F$-regular if $p>5$ (resp.\ $p>42$).
On the other hand, if $k$ is algebraically closed of characteristic $p>0$, \cite[Theorem 1.1]{CTW15a}, Example~\ref{eg:rdp}, and Proposition \ref{prop:non-F-pure canonical Gorenstein} (that refines \cite[Theorem 1.2]{CTW15a}) show that
there is a klt del Pezzo surface (resp.\ $\Q$-factorial canonical Gorenstein threefold singularity over $k$) that is not globally $F$-regular (resp.\ not $F$-regular).
\end{example}
}

\medskip
\noindent {\bf Acknowledgements.}

\begin{itemize}
    \item The first author was supported by JSPS KAKENHI Grant number JP22KJ1771 and JP24K16897.
    \item The second author was supported by JSPS KAKENHI Grant number JP22J00962 and JP25K17228.
    \item The third author was supported by JSPS KAKENHI Grant numbers JP22H01112 and JP23K03028. 
    \item The fourth author was supported by NSF research grant DMS-2101897.
    \item The fifth author was supported by JSPS KAKENHI Grant number JP19K14501.
    \item The sixth author was supported by JSPS KAKENHI Grant number JP20J11886 and RIKEN iTHEMS Program. 
    \item 
The authors thank the referees for their constructive suggestions and careful reading of the manuscript.
\end{itemize}

\section{Preliminaries}

\subsection{Notation}\label{ss-notation}
In this subsection, we summarise notation and 
terminologies used in this {article}. 

\begin{enumerate}
\item Throughout the paper, $p$ denotes a prime number and we set $\F_p := \Z/p\Z$. 
\item 
We say that a ring $R$ is {\em of characteristic $p>0$} if 
$R$ is an $\F_p$-algebra, i.e., the ring homomorphism $\Z \to R$ factors through $\F_p$. 
We say that a scheme $X$ is {\em of characteristic $p>0$} if 
$X$ is an $\F_p$-scheme, i.e., the morphism $X \to \Spec \Z$ factors through $\Spec \F_p$. 
\item For a ring $R$ of characteristic $p>0$, $F: R\to R$ denotes the absolute Frobenius ring homomorphism, i.e., $F(r)=r^p$ for every $r \in R$. 
Given a scheme $X$ of characteristic $p>0$, $F: X \to X$ denotes the absolute Frobenius morphism. 
\item Given a ring $R$, $R^{\circ}$ denotes the subset of $R$ consisting of  
non-zero-divisors. In other words, for $x \in R$, we have $x \in R^{\circ}$ if and only if 
the map $R \to R$ given by  $r \mapsto xr$ is injective. 
In particular, $R^{\circ}$ is a multiplicatively closed subset of $R$. 
\begin{itemize}
\item If $R$ is an integral domain, then $R^{\circ} = R \setminus \{0\}$. 
\item 
If $R$ is a Noetherian reduced ring, then 
we have $R^{\circ} = R \setminus (\mathfrak p_1 \cup \cdots \cup \mathfrak p_r)$, 
where $\mathfrak p_1, ...,  \mathfrak p_r$ are all the minimal prime ideals of $R$. 
\end{itemize}
\item 
{\cred For $n \in \Z_{>0}$ and an $\F_p$-algebra $R$, 
$W_nR$ denotes the ring of Witt vectors of length $n$ 
(for details, see \cite[Subsection 2.2]{KTTWYY1}).}
Given $r \in R$, 
we define 
\[
[r] := (r, 0, ..., 0) \in W_nR
\]
and call it the {\em Teichm\"{u}ller lift} of $r \in R$. 
\item We say that a Noetherian scheme $X$ is {\em equi-dimensional} if 
$\dim Y_1  = \dim Y_2 = \cdots = \dim Y_r$  for all the irreducible 
components $Y_1, Y_2, ..., Y_r$ of $X$. 
We say that a Noetherian ring $R$ is {\em equi-dimensional} if 
so is $\Spec R$. 
    \item \label{sss-Ffinite} {We say that an $\F_p$-scheme $X$ is {\em $F$-finite} if 
$F \colon X \to X$ is a finite morphism,  
 and  we say that an $\F_p$-algebra $R$ is {\em $F$-finite} if $\Spec R$ is $F$-finite. Such schemes admit many good properties.
\begin{enumerate}
\renewcommand{\labelenumi}{(\roman{enumi})}
    \item If $R$ is an $F$-finite Noetherian $\F_p$-algebra, 
    then it is a homomorphic image of a regular ring of finite Krull dimension \cite[Remark 13.6]{Gabber}; 
    in particular, $R$ is excellent, it admits a dualising complex, and $\dim R < \infty$.  
    \item $F$-finite {\cred Noetherian $\F_p$-algebras} 
    are stable under localisation and ideal-adic completions \cite[Example 9]{Has15}.
    \item If a scheme $X$ is of finite type over an $F$-finite Noetherian $\F_p$-scheme $Y$, then it is also $F$-finite. 
\end{enumerate}}
\item 
Let $X$ be an integral normal $F$-finite Noetherian scheme $X$ of characteristic $p>0$. 
For a $\Q$-divisor $D$ on $X$, 
we define $\MO_X(D)$ as the coherent $\MO_X$-submodule 
of the constant sheaf $\underline{K(X)}$ given by the formula 
\[
\Gamma(U, \MO_X(D)) := \{ f \in K(X) \,|\, (\div(f) + D)|_{U} \geq 0\}. 
\]
We define $W_n\MO_X(D)$ as the subsheaf  
of the constant sheaf $\underline{W_n(K(X))}$ given by the formula 
\begin{align*}
\Gamma(U, W_n\MO_X(D)) &:= 
\prod_{m=0}^{n-1} \Gamma(U, \MO_X(p^mD)) \\
&= \big\{ (f_0, ..., f_{n-1}) \in W_n(K(X)) \,\big|\, \\
&\qquad \quad (\div(f_i) + p^iD)|_{U} \geq 0\text{ for all }0 \leq i \leq n-1\big\}. 
\end{align*}
It is known that $W_n\MO_X(D)$ is a coherent $W_n\MO_X$-submodule 
of the constant sheaf $\underline{W_n(K(X))}$ \cite[Proposition 3.8]{tanaka22}. 
{\cred We then have $\MO_X( D) = \MO_X(\rdown{D})$, 
whilst $W_n\MO_X(D) \neq W_n\MO_X(\rdown{D})$ in general. 
When $X$ is affine: $X = \Spec R$, we set $R(D) := \Gamma(X, \MO_X(D))$ and $W_nR(D) := \Gamma(X, W_n\MO_X(D))$.} 
\item 
We say that $f : Y \to X$ is an {\em alteration} if 
$f$ is a proper surjective morphism of Noetherian integral schemes such that 
induced field extension $K(X) \subseteq K(Y)$ is of finite degree. 
\item 
{\cred An integral normal excellent scheme $X$ admitting a dualising complex is 
{\em $\Q$-Gorenstein} if $nK_X$ is Cartier for some integer $n>0$. 
We say that a normal Noetherian domain $R$ is {\em $\Q$-Gorenstein} if $\Spec R$ is $\Q$-Gorenstein.} 
\item 
{\cred For a Noetherian local $\F_p$-algebra $(R, \m)$, an integer $n>0$, and a $W_nR$-module $M$, 
we set 
\[
H^i_{\m}(M) := H^i_{W_n\m}(M) 
\]
by abuse of notation. 
Here there should arise no confusion, because the induced morphism $\Spec R \to \Spec W_nR$ is a surjective closed immersion, and hence a homeomorphism.}
\item 
{\ccred 
Let $X$ be a Noetherian integral scheme admitting dualising complex. 
We say that $X$ is {\em pseudo-rational} if $f_*\omega_Y = \omega_X$ for every proper birational morphism $f: Y \to X$ from an integral normal scheme $Y$. 
We warn the reader that we do not require the condition that $X$ is Cohen-Macaulay, although 
it is imposed in some literature} 

\end{enumerate}

{\cred 
\subsubsection{Dualising sheaves and dualising modules}\label{subsubsection dualising}

Let $X$ be a Notherian scheme admitting a dualising complex. 
Take the decomposition $X = \coprod_{i \in I} X_i$ into the connected components. 
By applying suitable shifts on each connected component $X_i$, 
we can find an integer $e$ and a dualising complex $\omega_X^{\bullet}$ such that,  for every $i \in I$, 
$e$ is the smallest integer satisfying 
$\mathcal H^e(\omega_X^{\bullet}|_{X_i}) \neq 0$. 
In this case, $\omega_X := \mathcal  H^e(\omega_X^{\bullet})$ 
is called the {\em dualising sheaf} on $X$. 
It is easy to see that $\omega_X$ does not depend on the choice of $\omega_X^{\bullet}$. 
Moreover, it holds that $\omega_X|_U \simeq \omega_U$ 
if $\Supp(\omega_X) = X$ and $U$ is an open subset of $X$.


For a Noetherian ring $R$ such that $\Spec\,R$ admits a dualising complex, we set $\omega_R := \Gamma(\Spec R, \omega_{\Spec R})$, which is called the {\em dualising module} of $R$. 
Furthermore, for a multiplicative subset $S$ of $R$, it holds that $\omega_R \otimes_R S^{-1}R \simeq \omega_{S^{-1}R}$ if $\Supp(\omega_R)=\Spec R$.}

\subsubsection{Frobenius push-forward modules $F^e_*M$}\label{sss-F_*M}

Let $R$ be a ring of characteristic $p>0$. 
Fix integers $e \geq 0$ and $n>0$. 
For a $W_nR$-module $M$, 
we define a $W_nR$-module $F_*^eM$ as follows. 
\begin{enumerate}
\item As an additive group, we set $F_*^eM := M$. 
\item For $m \in M$, we set  $F_*^eM \ni F_*^em :=m$ (in order to distinguish $F_*^eM$ from  $M$). 
\item For $x \in W_nR$ and $m \in M$, the {$W_nR$}-linear structure on $F_*^eM$ is defined by 
\[
x \cdot (F_*^em) := F_*^e( F^e(x)m),
\]
where the product $F^e(x)m$ is given by the one in $M$. 
\end{enumerate}
In particular, we have $F_*^{e_1}F_*^{e_2}M = F_*^{e_1+e_2}M$. 
Similarly, if $0 \leq e' \leq e$, $x \in W_nR$, and $m \in M$, then we set 
\begin{equation}\label{e1-compu-F*M}
(F_*^{e'}x) \cdot (F_*^em) := F_*^{e'}( x \cdot (F_*^{e-e'}m)) \in F_*^eM. 
\end{equation}
For example, if $x = F^{e'}(y)$ for some $y \in W_nR$, then we have 
\begin{equation}\label{e2-compu-F*M}
(F_*^{e'}(F^{e'}(y))) \cdot (F_*^em) = y\cdot (F_*^em)\quad \text{in}\quad F_*^eM, 
\end{equation}
because 
\[
(F_*^{e'}(F^{e'}(y))) \cdot (F_*^em) = F_*^{e'}( (F^{e'}(y)) \cdot (F_*^{e-e'}m)) 
= y \cdot 
F_*^{e'}(F_*^{e-e'}m) =y\cdot (F_*^em). 
\]

\begin{remark}\label{r-compu-F_*M}
Given a scheme $X$ of characteristic $p>0$ and a quasi-coherent $W_n\MO_X$-module $M$ (i.e., a quasi-coherent sheaf on the scheme $(X, W_n\MO_X)$), 
we have the $e$-th iterated Frobenius morphism $F_*^e : W_nX \to W_nX$. 
The above definition coincides with the usual scheme-theoretic push-forward $F_*^eM$ 
for the case when $X$ is affine via the identification 
of quasi-coherent $W_n\MO_X$-modules and $\Gamma(X,  W_n\MO_X)$-modules. 
\end{remark}
    

\subsection{Recollection on the theory of $F$-singularities via closures}
The goal of this subsection is to provide a brief introduction to the theory of F-singularities via closures and to justify some of the definitions we will use later in the context of quasi-$F$-singularities. The results, proofs, and definitions in this subsection will not be used anywhere else in our paper and are only provided for the convenience of the reader. For this reason, we do not strive for the highest level of generality. Most general definitions may be found in Subsection \ref{ss:tight-closure-and-test-ideals}. For more on this topic, we refer the reader to \cite[Section 3]{schwedetucker12}, \cite[Section 2]{Tak21}, as well as the forthcoming book \cite{schwedesmithbook}. 

Let $(R,\m)$ be a $d$-dimensional Noetherian $F$-finite {\cred complete} local domain  of 
characteristic $p>0$, {\cred so that} $R^\circ = R\, \backslash \{0\}$. 
By applying $\Hom_R(-, \omega_R)$ to $F^e : R \to F^e_*R$, we obtain the trace of Frobenius $T^e \colon F^e_*\omega_R \to \omega_R$. We will often denote $T^1 \colon F_* \omega_R \to \omega_R$ by $T$.

\begin{definition}[{cf.\ \cite[Definition 8.7 and Exercise 8.8]{schwedetucker12}}]
{For simplicity, assume that $R$ is Cohen-Macaulay.} In this case, we say that $R$ is \emph{$F$-injective} if one of the following equivalent conditions hold:
\begin{enumerate}
\item $T^e \colon F^e_*\omega_R \to \omega_R$ is surjective for every $e>0$,\\[-1em]
\item $H^d_\m(R) \to H^d_\m(F^e_*R)$ is injective for every $e>0$. 
\end{enumerate}
\end{definition}
The equivalence of these two conditions follows by applying Matlis duality $\Hom_R(-, E)$, where $E = H^d_\m(\omega_R)$ is the injective hull of the $R$-module $R/\m$. Moreover, it is enough to check these conditions for a single $e>0$, for example, $e=1$. Last, we point out that if $\omega_R \simeq R$, then $R$ is $F$-injective if and only if $R$ is $F$-split.

Next, one constructs 
an $R$-submodule $\sigma(\omega_R)$ of $\omega_R$
that measures how far a ring is from being $F$-injective. To this end, we define the Frobenius closure of $0$ in the $R$-module $H^d_\m(R)$ as follows:
\begin{align*}
0^{F,e}_{H^d_\m(R)} &:= {\rm Ker}\Big(F^e \colon H^d_{\m}(R) \to H^d_{\m}(F^e_*R)\Big), \text{ and }\\
0^{F}_{H^d_\m(R)} &:= \bigcup_{e>0} 0^{F,e}_{H^d_\m(R)}.
\end{align*}

\begin{definition} \label{def:nonFpure-submodule-recap}
We define the 
$R$-submodule 
$\sigma(\omega_R)$ of $\omega_R$ by one of the following four equivalent assertions:
\begin{enumerate}
    \item $\sigma(\omega_R) := \bigcap_{e>0} {\rm Im}(T^e \colon F^e_*\omega_R \to \omega_R)$,\\[-0.5em]
    \item $\sigma(\omega_R) := \big(\varinjlim_{e>0}{\rm Im}( F^e \colon H^d_{\m}(R) \to H^d_{\m}(F^e_*R)) \big)^\vee$, \\[-0.5em]
    \item $\sigma(\omega_R) := \Big(H^d_{\m}(R)/{0^F_{H^d_\m(R)}}\Big)^\vee$, \\[-0.5em]
    \item $\sigma(\omega_R) := {\rm Ann}_{\omega_R}(0^F_{H^d_\m(R)})$.\footnote{the annihilator is taken with respect to the multiplicative pairing $\eta \colon \omega_R \otimes H^d_\m(R) \to H^d_\m(\omega_R)=E$; {specifically: ${\rm Ann}_{\omega_R}(0^F_{H^d_\m(R)}) = \{ x \in \omega_R\,|\, \eta( x \otimes y)=0 \text{ for every $y \in 0^F_{H^d_\m(R)}$}\}$.} This pairing is perfect by Matlis duality. }
\end{enumerate} 
Here $(-)^{\vee} := \Hom_R(-, E)$ 
{\cred and the inductive limit $\varinjlim_{e>0}{\rm Im}( F^e \colon H^d_{\m}(R) \to H^d_{\m}(F^e_*R))$ in (2) is defined by using $F^{e+1} \colon H^d_{\m}(R) \xrightarrow{F^e} H^d_{\m}(F^e_*R) \xrightarrow{F}  H^d_{\m}(F^{e+1}_*R)$}. 
\end{definition}
One can show that the intersection and colimit in these definitions stabilise for $e \gg 0$ {(see 
{\cred \cite[Lemma 8.1]{TWY}})}. 


\begin{proof}[Brief sketch of the proof {\cred of the equivalences}]
By stabilisation, the equivalence of (1) and (2) follows from Matlis duality. Further, the equivalence of (2) and (3) is ensured by the short exact sequence:
\begin{equation} \label{eq:basic-frobenius-closure-sequence}
0 \to 0^{F,e}_{H^d_\m(R)} \to H^d_\m(R) \to {\rm Im}( F^e \colon H^d_{\m}(R) \to H^d_{\m}(F^e_*R)) \to 0.
\end{equation}
As for the equivalence of (3) and (4), since $\omega_R \otimes H^d_\m(R) \to H^d_\m(\omega_R)$ is a pe{\cred r}fect pairing 
{\cred by Matlis duality}, we can identify the annihilator ${\rm Ann}_{\omega_R}(0^F_{H^d_\m(R)})$ with the $R$-submodule
\[
\Hom_R(H^d_\m(R)/0^F_{H^d_\m(R)}, H^d_\m(\omega_R)) \subseteq \Hom_R(H^d_\m(R), H^d_\m(\omega_R)).
\]
But the left hand side is exactly $\Big(H^d_{\m}(R)/{0^F_{H^d_\m(R)}}\Big)^\vee$ {\cite[Corollary 3.5.9]{BH93}}.
\end{proof}


\begin{remark} \label{remark:frobenius-quotient-identity}
Note that (3) is equivalent to 
$\big(\frac{\omega_R}{\sigma(\omega_R)}\big)^\vee = 0^F_{H^d_\m(R)}$.
Indeed, by applying Matlis duality to (\ref{eq:basic-frobenius-closure-sequence}), we get the following short exact sequence
\[
0 \to \sigma(\omega_R) \to \omega_R \to \Big(0^F_{H^d_\m(R)}\Big)^{\vee} \to 0,
\]
thanks to the fact that Matlis duality is exact and $H^d_\m(R)^\vee \simeq \omega_R$. 
\end{remark}

Now, we move to the definition of $F$-rationality. To this end, given $c \in R$, we  define $T^{e,c}$ as the composition of the following $R$-module homomorphisms:
\[
T^{e,c} \colon F^e_*\omega_R \xrightarrow{\cdot F^e_*c} F^e_* \omega_R \xrightarrow{T^e} \omega_R.
\]

\begin{definition}[{cf.\ \cite[Definition 8.1]{ST08} and  \cite[Definition 2.31]{BST15}}]
We say that $R$ is \emph{$F$-rational} if {$R$ is Cohen-Macaulay and}, for every $c \in R^\circ$, there exists an integer $e>0$ such that the following equivalent conditions hold:\\[-0.7em]
\begin{enumerate}
\item the $R$-module homomorphism $T^{e,c} \colon F^e_*\omega_R \to \omega_R$
is surjective,\\[-0.7em]
\item the composition $H^d_\m(R) \xrightarrow{F^e} H^d_\m(F^e_*R) \xrightarrow{\cdot F^e_*c} H^d_\m(F^e_*R)$
is injective.\\[-1em]
\end{enumerate}
\end{definition}

\noindent
{\cred We point out that if $\omega_R \simeq R$, then $R$ is $F$-rational if and only if $R$ is strongly $F$-regular.} 

To construct 
{the $R$-submodule $\tau(\omega_R)$ of $\omega_R$} 
measuring $F$-rationality, define 
\[
K^{e,c} := {\rm Ker}\Big(H^d_{\m}(R) \xrightarrow{F^e} H^d_{\m}(F^e_*R) \xrightarrow{\cdot F^e_*c} H^d_{\m}(F^e_*R) \Big) 
\]
for any integer $e>0$ and $c \in R^\circ$. 


\begin{definition}[{\cite[Definition 5.19]{ST08}}] The \emph{tight closure} of $0$ in $H^d_\m(R)$ is given by
\[
0^*_{H^d_\m(R)} := \bigcup_{\substack{c \in R^{\circ},\\ e_0 \in \Z_{>0}}} \bigcap_{e \geq e_0} K^{e,c} \subseteq H^d_\m(R).
\]
Equivalently, given $z \in H^d_\m(R)$, 
    we have $z \in 0^*_{H^d_\m(R)}$ if and only if 
    there exist $c \in R^{\circ}$ and an integer $e_0>0$ such that 
    $z$ is contained in $K^{e, c}$ for every integer $e \geq e_0$.
\end{definition}

\begin{remark}A key result in the theory of $F$-singularities is the existence of a test element (cf.\ Definition \ref{d-test-ideal}, {\cred \cite{schwedesmithbook}}). In particular, this implies the existence of $t \in R^\circ$ such that $t0^*_{H^d_\m(R)} = 0$. We will see by equality (\ref{eq:dual-test-quotient})  that this condition is the same as $t\omega_R \subseteq \tau(\omega_R)$ for the test 
{submodule}
$\tau(\omega_R)$ defined below.
\end{remark}


Finally, we say that an $R$-submodule $M \subseteq \omega_R$ is 
\emph{$T$-stable} 
if $T(F_*M) \subseteq M$ for $T \colon F_*\omega_R \to \omega_R$.


Before proceeding, we need the following lemma.

\begin{lem}\label{l recollection stab}
Take  $t \in {\rm Ann}_{R}(0^*_{H^d_\m(R)}) \cap R^\circ$ and $c \in (t^2) \cap R^{\circ}$. 
Set 
\[
I^{e, c} := {\rm Im}(T^{e,c} \colon F^e_*\omega_R \to \omega_R).
\]
Then there exists an integer $e_0>0$ such that 
\begin{enumerate}
\item 
$K^{0, c} \supseteq K^{1, c} \supseteq  \cdots \supseteq K^{e_0, c} 
= K^{e_0+1, c} = \cdots$, and 
\item 
$I^{0, c} \subseteq I^{1, c} \subseteq \cdots  \subseteq  I^{e_0, c} = I^{e_0+1, c} = \cdots$.
\end{enumerate}
\end{lem}

\begin{proof}
By $0^F_{H^d_{\m}(R)} \subseteq 0^*_{H^d_{\m}(R)}$, 
we get $t \in {\rm Ann}_R(0^*_{H^d_{\m}(R)}) \subseteq {\rm Ann}_R(0^F_{H^d_{\m}(R)})$. 
This, together with $(0^F_{H^d_{\m}(R)})^{\vee} = \omega_R/\sigma(\omega_R)$, 
implies $t\omega_R \subseteq \sigma(\omega_R)$. 
For $c = t^2c'$, assertion (2) follows from 
\begin{align*}
I^{e, c} &= {\rm Im}(T^{e, t^2c'} \colon F^e_*\omega_R \to \omega_R) \\
&\subseteq 
{\rm Im}(T^{e,tc'} \colon F^e_*\sigma(\omega_R) \to \omega_R) \\
&\subseteq {\rm Im}(T^{e+1,t^pc'^p} \colon F^{e+1}_*\omega_R \to \omega_R) \\
&\subseteq {\rm Im}(T^{e+1,t^2c'} \colon F^{e+1}_*\omega_R \to \omega_R) =I^{e+1, c}.
\end{align*}
{\cred Here the sequence $I^{0, c} \subseteq I^{1, c} \subseteq \cdots$ stablises, because 
$\omega_R$ is a Noetherian $R$-module.} 
By Matlis duality, we have $(I^{e, c})^{\vee} \simeq H^d_{\m}(R)/K^{e, c}$, 
and hence (2) implies (1). 
\end{proof}

\begin{definition}[{cf.\ \cite[Definition 2.33]{BST15}}] \label{def:test-submodule-recap}
We define the \emph{test submodule} $\tau(\omega_R)$ of $\omega_R$ by one of the following equivalent assertions:
\begin{enumerate}
    \item $\tau(\omega_R) := \bigcap_{c \in R^\circ} \bigcap_{e_0>0} \sum_{e \geq e_0} {\rm Im}(T^{e,c} \colon F^e_*\omega_R \to \omega_R)$,\\[-0.5em]
    \addtocounter{enumi}{1}
    \item $\tau(\omega_R) := \Big(H^d_\m(R)/0^*_{H^d_\m(R)} \Big)^\vee$, \\[-0.5em] 
    \item $\tau(\omega_R) := {\rm Ann}_{\omega_R}(0^*_{H^d_\m(R)})$,\\[-0.5em]
    \item $\tau(\omega_R) := \bigcap_{c \in R^\circ} \sum_{e_0>0} {\rm Im}(T^{e_0,c} \colon F^{e_0}_*\omega_R \to \omega_R)$,\\[-0.5em]
    \item $\tau(\omega_R)$ is the smallest {nonzero} $R$-submodule of $\omega_R$ which is $T$-stable, \\[-0.5em]
     \item $\tau(\omega_R) := {\rm Im}(T^{e,t^2} \colon F^e_*\omega_R \to \omega_R)$, where $t \in {\rm Ann}_{{R}}(0^*_{H^d_\m(R)}) \cap R^\circ$ and $e \gg 0$ (depending on $t$).
     \item $\tau(\omega_R) := \bigcap_{c \in R^\circ} \sum_{e_0>0} \bigcap_{e \geq e_0} {\rm Im}(T^{e,c^{p^{e-e_0}}} \colon F^e_*\omega_R \to \omega_R)$.
\end{enumerate}
\end{definition}
\noindent 
It 
{\cred follows from (5) and (7)} 
that $\tau(\omega_R)= {\cred \omega_R}$ if and only if $R$ is $F$-rational. Also note that we skipped the analogue of (2) from Definition \ref{def:nonFpure-submodule-recap} as formulating a correct statement is a bit subtle. 

\begin{proof}[Sketch of the proof of the equivalences]For the sake of the proof, we define $\tau(\omega_R)$ by assertion (1): \[
\tau(\omega_R) := \bigcap_{c \in R^\circ} \bigcap_{e_0>0} \sum_{e \geq e_0} {\rm Im}(T^{e,c} \colon F^e_*\omega_R \to \omega_R)
\]
and show that it is equivalent to (3)--(7). 
Fix $t \in {\rm Ann}_{{R}}(0^*_{H^d_\m(R)}) \cap R^\circ$.
{\cred Note that, for all $c, c' \in R^{\circ}$, we have 
\[
 {\rm Im}(T^{e,c} \colon F^e_*\omega_R \to \omega_R)
\supseteq  {\rm Im}(T^{e,cc'} \colon F^e_*\omega_R \to \omega_R), 
\]
which implies 
\[
\tau(\omega_R) = 
\bigcap_{c \in R^\circ} \bigcap_{e_0>0} \sum_{e \geq e_0} {\rm Im}(T^{e,c} \colon F^e_*\omega_R \to \omega_R)
= 
\bigcap_{c \in c' R^\circ} \bigcap_{e_0>0} \sum_{e \geq e_0} {\rm Im}(T^{e,c} \colon F^e_*\omega_R \to \omega_R). 
\]
Similarly, for every $c' \in R^{\circ}$, we get 
\[
0^*_{H^d_\m(R)}  = \bigcup_{\substack{c \in R^{\circ},\\ e_0 \in \Z_{>0}}} \bigcap_{e \geq e_0} K^{e,c} 
= \bigcup_{\substack{c \in c'R^{\circ},\\ e_0 \in \Z_{>0}}} \bigcap_{e \geq e_0} K^{e,c}. 
\]
}
\\[0.2em]

\noindent (1) $\iff$ (3): By Matlis duality, ${\cred I^{e,c}} := {\rm Im}(T^{e,c} \colon F^e_*\omega_R \to \omega_R) = (H^d_\m(R)/K^{e,c})^\vee.$ 
{\cred 
It follows from  Lemma \ref{l recollection stab} that, for every $c \in t^2 R^{\circ}$, 
we can find an integer $e_1 >0$ satisfying 
\[
I^{0, c} \subseteq I^{1, c} \subseteq \cdots \subseteq I^{e_1, c} = I^{e_1+1, c} = \cdots =:I^{\infty, c}, 
\]
\[
K^{0,c } \supseteq K^{1, c} \supseteq \cdots \supseteq 
K^{e_1, c} = K^{e_1+1, c} = \cdots  =: K^{\infty, c}. 
\]
Then it holds that 
\[
0^*_{H^d_\m(R)}  = \bigcup_{c \in R^{\circ}}
 \bigcup_{{e_0 \in \Z_{>0}}}
\bigcap_{e \geq e_0} K^{e,c} 
=\bigcup_{c \in R^{\circ}}
 \bigcup_{{e_0 \in \Z_{>0}}} K^{\infty, c}
=\bigcup_{c \in R^{\circ}} K^{\infty, c}. 
\]}
Since {\cred the} Matlis duality {\cred functor $(-)^{\vee}$  
turns colimits into limits}, 
we get that 
{\cred 
\begin{align*}
\tau(\omega_R) &= 
\bigcap_{c \in R^\circ} \bigcap_{e_0>0} \sum_{e \geq e_0}  (H^d_\m(R)/K^{e,c})^\vee\\
&= 
\bigcap_{c \in t^2R^\circ} \bigcap_{e_0>0} \sum_{e \geq e_0}  I^{e, c}\\
&= \bigcap_{c \in t^2R^\circ}  I^{\infty, c}\\
&= \varprojlim_{c \in t^2R^\circ}  (H^d_\m(R)/K^{{\cred \infty},c})^\vee \\
&= \Big(\varinjlim_{c \in t^2R^\circ} H^d_\m(R)/K^{{\cred \infty},c}\Big)^\vee \\ 
&= \Big(H^d_\m(R)/
\big(\varinjlim_{c \in t^2R^\circ}  K^{{\cred \infty},c}\big) \Big)^\vee \\ 
&=\Big(H^d_\m(R)/ 0^*_{H^d_\m(R)}\Big)^\vee.
\end{align*}
}


 
\noindent (3) $\iff$ (4): This is exactly the same as in Definition \ref{def:nonFpure-submodule-recap}.\\

\noindent (1) $\iff$ (5): Define $\tau'(\omega_R) := \bigcap_{c \in R^\circ} \sum_{e_0>0} {\rm Im}(T^{e_0,c} \colon F^{e_0}_*\omega_R \to \omega_R)$. Clearly $\tau(\omega_R) \subseteq \tau'(\omega_R)$, 
so it suffices 
to prove the opposite inclusion. Take any $t \in R^\circ$ such that $t\omega_R \subseteq \sigma(\omega_R)$. Then 
\[
{\rm Im}(T^{e_0,ct}) \subseteq {\rm Im}(T^{e,c^{p^{e-e_0}}}) \subseteq {\rm Im}(T^{e,c})   
\]
for all integers $e \geq e_0 \geq 0$. This implies inclusion ($\dagger$) in 
\[
\tau'(\omega_R) = \bigcap_{c \in R^\circ} \sum_{e_0>0} {\rm Im}(T^{e_0,ct}) \overset{(\dagger)}{\subseteq}
\bigcap_{c \in R^\circ} \bigcap_{e_0>0} \sum_{e\geq e_0} {\rm Im}(T^{e,c}) = \tau(\omega_R).
\]

\noindent (5) $\iff$ (6): 
It is enough to prove that $\tau(\omega_R) = \bigcap_{c \in R^\circ} \sum_{e_0>0} {\rm Im}(T^{e_0,c} \colon F^{e_0}_*\omega_R \to \omega_R)$ is $T$-stable and is contained in any nonzero $T$-stable $R$-submodule $M$ of $\omega_R$.  Clearly, 
\begin{align*}
T(F_*\tau(\omega_R)) &\subseteq \bigcap_{c \in R^\circ} \sum_{e_0>0} {\rm Im}(T^{e_0+1,c^p} \colon F^{e_0+1}_*\omega_R \to \omega_R) \\
&\subseteq  \bigcap_{c \in R^\circ} \sum_{e_0>0} {\rm Im}(T^{e_0,c^p} \colon F^{e_0}_*\omega_R \to \omega_R) \\
&= \tau(\omega_R).
\end{align*} 
Hence $\tau(\omega_R)$ is $T$-stable. To show that it is contained in $M$, pick $d \in R^\circ$ such that $d\omega_R \subseteq M$ 
{\cred (we can find such an element $d \in R^\circ$, because $\omega_U \simeq \MO_U$ 
for some non-empty open subset $U$ of $\Spec R$)}. Then
\[
\tau(\omega_R) = \bigcap_{c \in R^\circ} \sum_{e_0>0} {\rm Im}(T^{e_0,cd}) \subseteq \bigcap_{c \in R^\circ} \sum_{e_0>0} T^{e_0,c}(F^{e_0}_*M) \subseteq \sum_{e_0>0} T^{e_0}(F^{e_0}_*M)  \subseteq M,
\]
where the last inclusion is clear as $M$ is $T$-stable.\\

\noindent (6) $\iff$ (7): We will implicitly use that the equivalence of (1), {\cred (3),} and (6) is already proven. 
Define 
\[
\tau''(\omega_R) := \sum_{e>0} {\rm Im}(T^{e,t^2} \colon F^e_*\omega_R \to \omega_R) 
= {\rm Im}(T^{e,t^2} \colon F^e_*\omega_R \to \omega_R), 
\]
{where the latter equality holds for $e \gg 0$ (Lemma \ref{l recollection stab}).} 



To finish the proof, 
it suffices to prove that $\tau''(\omega_R)$ is $T$-stable and is contained in $\tau(\omega_R)$. Clearly
\[
T(F_*\tau''(\omega_R)) \subseteq {\rm Im}(T^{e+1,t^{2}} \colon F^{e+1}_*\omega_R \to \omega_R)  = \tau''(\omega_R).
\]
Moreover,
\[
\tau''(\omega_R) = {\rm Im}(T^{e,t^2} \colon F^e_*\omega_R \to \omega_R) 
\overset{{\cred {\rm (3)}}}{\subseteq} {\rm Im}(T^{e,t} \colon F^e_*\tau(\omega_R) \to \omega_R) \overset{{{(6)}}}{\subseteq} \tau(\omega_R).
\]
This concludes the proof.\\

\noindent (5) $\iff$ (8): 
Let us denote  the right hand side of (8) by $\tau'''(\omega_R)$. 
It holds that 
\begin{align*}
\tau'''(\omega_R) &= \bigcap_{c \in R^\circ} \sum_{e_0>0} \bigcap_{e \geq e_0} {\rm Im}(T^{e,c^{p^{e-e_0}}} \colon F^{e}_*\omega_R \to \omega_R) \\
&\subseteq \bigcap_{c \in R^\circ} \sum_{e_0>0}  {\rm Im}(T^{e_0,c} \colon F^{e_0}_*\omega_R \to \omega_R) \\
&\overset{{\rm (5)}}{=} \tau(\omega_R). 
\end{align*}

It suffices to show the opposite inclusion 
$\tau'''(\omega_R) \supseteq \tau(\omega_R)$. 
We now compute  
$\bigcap_{e \geq e_0} {\rm Im}(T^{e,c^{p^{e-e_0}}} \colon F^{e}_*\omega_R \to \omega_R)$. 
For $c \in R^{\circ}$, $\zeta \in \omega_R$, and integers $e \geq e_0$, 
we have  that 
\begin{align*}
T^e((F_*^ec^{p^{e-e_0}}) \cdot (F^e_*\zeta))
&= 
T^{e_0}(F^{e_0}_* T^{e-e_0} ((F_*^{e-e_0}c^{p^{e-e_0}}) \cdot (F^{e-e_0}_*\zeta)))\\
&= 
T^{e_0}(F^{e_0}_* T^{e-e_0} ( c \cdot (F^{e-e_0}_*\zeta)))\\
&= 
T^{e_0}(F^{e_0}_*( c\cdot  T^{e-e_0} (F^{e-e_0}_*\zeta))).
\end{align*}
Therefore, we get 
\begin{align*}
\bigcap_{e \geq e_0} {\rm Im}(T^{e,c^{p^{e-e_0}}} \colon F^{e}_*\omega_R \to \omega_R) &= 
\bigcap_{e \geq e_0} T^{e_0}(F^{e_0}_*( c\cdot  {\rm Im}(T^{e-e_0})))\\
&\supseteq 
 T^{e_0}(F^{e_0}_*( c\cdot (\bigcap_{e \geq e_0}  {\rm Im}(T^{e-e_0}))))\\
 & =  T^{e_0}(F^{e_0}_*( c\cdot \sigma(\omega_R))). 
\end{align*}
{\cred 
We have that $\sigma(\omega_R) \neq 0$, 
because the trace of Frobenius $T: F_*\omega_U \to \omega_U$ is surjective on some non-empty open subset $U \subseteq \Spec R$.  
Hence we can find} an element $d \in R^{\circ}$ satisfying 
$\sigma(\omega_R) \supseteq d \omega_R$. 
Then 
it holds that 
\[
\tau'''(\omega_R) \supseteq \bigcap_{c \in R^{\circ}} \sum_{e_0>0} T^{e_0}(F^{e_0}_*( c\cdot \sigma(\omega_R))) 
\supseteq \bigcap_{c \in R^{\circ}} \sum_{e_0>0} T^{e_0}(F^{e_0}_*( cd\cdot \omega_R)) 
\]
\[
=
\bigcap_{c \in R^{\circ}} \sum_{e_0>0} T^{e_0}(F^{e_0}_*( c\cdot \omega_R)) \overset{{\rm (5)}}{=} \tau(\omega_R).
\]
We are done.
\end{proof}


\begin{remark} 
As in Remark \ref{remark:frobenius-quotient-identity}, the third assertion is equivalent to
\begin{equation} \label{eq:dual-test-quotient}
\Big(\frac{\omega_R}{\tau(\omega_R)}\Big)^\vee = 0^*_{H^d_\m(R)}.
\end{equation}
\end{remark}




One way to define the test ideal $\tau(R)$ is to first define a test submodule for pairs, specifically, $\tau(\omega_R,D)$ for an effective $\bQ$-divisor $D$. Then, given a $\bQ$-divisor $\Delta$ such that $K_R + \Delta$ is $\bQ$-Cartier, one sets 
\[
\tau(R,\Delta) := \tau(\omega_R, K_R+\Delta).
\]
If $K_R$ is $\bQ$-Cartier, then $\tau(R) := \tau(R,0)$. This object is discussed in more detail in Subsection \ref{ss:tight-closure-and-test-ideals}.

\subsection{Tight closure and test ideals} \label{ss:tight-closure-and-test-ideals}
In this subsection, we further summarise the definitions and some known results on tight closure and test ideals. 
For more details, we refer to 
\cite[Section 3]{schwedetucker12}, \cite[Section 2]{Tak21}, and references therein.

\begin{definition}[\cite{HH90}*{Definition 8.2}, \cite{Tak04}*{Definition 2.1}]\label{d-tight-closure}
Let $R$ be an $F$-finite Noetherian reduced ring of characteristic $p>0$. 
Let $M$ be an $R$-module. 
\begin{enumerate}
    \item     
    We define the {\em tight closure} $0^*_M$ (of $0$ in $M$)
    {\cred as the $R$-submodule of $M$ that satisfies $(\star)$ below}. 
    \begin{enumerate}
    \item[($\star$)] Take $z \in M$. Then  $z \in 0^*_M$ if and only if 
    there exist an integer $e_0\geq 0$  and $c \in R^{\circ}$ 
    such that, for every $e \geq e_0$, 
    the equality $z \otimes F_*^ec=0$ holds in $M \otimes_{{\cred R}} F^e_*R$. 
    \end{enumerate}
    \item 
    Assume that $R$ is a normal integral domain. 
    Take an effective $\Q$-divisor $\Delta$ on $\Spec R$. 
    We define the {\em $\Delta$-tight closure} $0^{*\Delta}_M$ (of $0$ in $M$) 
    {\cred as the $R$-submodule of $M$ that satisfies $(\star\star)$ below}. 
\begin{enumerate}
\item[($\star\star$)] 
    Take $z \in M$. Then  $z \in 0^{*\Delta}_M$ if and only if 
    there exist an integer $e_0\geq 0$  and $c \in R^{\circ}$ 
    such that, for every $e \geq e_0$, 
    the equality  $z \otimes F^e_*c =0$ holds in $M \otimes_{{\cred R}} F^e_*R((p^e-1)\Delta)$. 
    \end{enumerate}   
\end{enumerate}
\end{definition}

\begin{lemma}\label{l-tight-M/N-vs-0}
Let $(R, \m)$ be a $d$-dimensional $F$-finite Noetherian normal local ring of characteristic $p>0$. 
Take a Weil divisor $D$ and a $\Q$-divisor $\Gamma$ on $\Spec R$. 
Then we have
\[
H^d_{\m}(R(D)) \otimes_R F^e_*R(\Gamma) \simeq H^d_{\m}(F^e_*R(p^eD+\Gamma))
\]
for every integer $e \geq 0$.
\end{lemma}

\begin{proof}
We consider the natural $R$-module homomorphism 
\[
\eta \colon R(D) \otimes_R F^e_*R(\Gamma) \to F^e_*R(p^eD+\Gamma)
\]
and the exact sequences
\[
0 \to K:=\mathrm{Ker}(\eta) \to R(D) \otimes_R F^e_*R(\Gamma) \to \mathrm{Im}(\eta) \to 0,
\]\[
0 \to \mathrm{Im}(\eta) \to F^e_*R(p^eD +\Gamma) \to 
{\rm Coker}(\eta)=:C \to 0.
\]
Since $\eta$ is an isomorphism at every point of codimension $\leq 1$, 
we get 
\[
\dim \Supp K \leq d-2\qquad \text{and}\qquad 
\dim \Supp C \leq d-2, 
\]
which imply $H^j_{\m}(K)=H^j_{\m}(C)= 0$ for $j \in \{d-1, d\}$. 
Therefore, we get the induced $R$-module isomorphism: 
\[
H^d_{\m}(\eta) : 
H^d_{\m}(R(D) \otimes_R F^e_*R(\Gamma)) \xrightarrow{\simeq} 
H^d_{\m}(F^e_*R(p^eD+\Gamma)). 
\]

Then it is enough to show that 
\begin{equation}\label{e-H^d-tensor}
H^d_{\m}(M) \otimes_R N \simeq
H^d_{\m}(M \otimes_R N)
\end{equation}
{\cred for all $R$-modules $M$ and $N$.} 
To this end, let $x_1,\ldots,x_d$ be a system of parameters of $R$.
Set $R_d:=R[x_1^{-1},\ldots,x_d^{-1}]$ and 
$R_{d-1}:=\prod_{1 \leq i \leq d} R[x_1^{-1},\ldots, x_{i-1}^{-1},x_{i+1}^{-1},\ldots,x_d^{-1}]$. 
In view of  \cite[\href{https://stacks.math.columbia.edu/tag/0A6R}{Tag 0A6R}]{stacks-project} or \cite[Theorem 3.5.6]{BH93}, 
we have an exact sequence
\begin{equation}\label{ex1}
    M \otimes_R R_{d-1} \to M \otimes_R R_d \to H^d_{\m}(M) \to 0.
\end{equation}
Then, by tensoring the above exact sequence (\ref{ex1}) with $N$ we obtain
\begin{equation}\label{ex2}
   M \otimes_R N \otimes_R R_{d-1} \to M \otimes_R N \otimes_{{\cred R}} R_d \to H^d_{\m}(M) \otimes_R N \to 0. 
\end{equation}
Moreover, we also have an exact sequence
\begin{equation}\label{ex3}
    M \otimes_R N \otimes_R R_{d-1} \to M \otimes_R N \otimes_R R_d \to H^d_{\m}(M \otimes_R N) \to 0.
\end{equation}
Comparing exact sequences (\ref{ex2}) and (\ref{ex3}), 
we obtain (\ref{e-H^d-tensor}), as required. 
\qedhere

\end{proof}

\begin{remark}\label{r-tight-M/N-vs-0-log}
Let $R$ be an $F$-finite Noetherian normal integral domain of characteristic $p>0$. 
Take an effective $\Q$-divisor $\Delta$ on $\Spec R$.
    We now prove that we may replace \lq\lq $p^e-1$" by \lq\lq $p^e$" 
    in Definition \ref{d-tight-closure}. 
    Specifically, we show 
    \begin{equation}\label{e1-tight-M/N-vs-0-log}
0^{*\Delta}_M = (0^{*\Delta}_M)'
\end{equation}
    for {\cred the $R$-submodule} $(0^{*\Delta}_M)'$ {\cred of $M$ that satisfies 
    $(\star\star)'$.}
    \begin{enumerate}
    \item[$(\star\star)'$] Take $z \in M$. Then  $z \in (0^{*\Delta}_M)'$ if and only if 
    there exist an integer $e_0\geq 0$  and $c \in R^{\circ}$ 
    such that, for every $e \geq e_0$, 
    the equality $z \otimes F^e_*c =0$ holds in $M \otimes_{{\cred R}} F^e_*R(p^e\Delta)$. 
    \end{enumerate}   

\end{remark}

\begin{proof}[Proof of (\ref{e1-tight-M/N-vs-0-log})]
    By the natural inclusion 
    \[
    R((p^e-1)\Delta) \hookrightarrow R(p^e\Delta),
    \]
    we get $0^{*\Delta}_M \subseteq (0^{*\Delta}_M)'$.
Let us prove the opposite inclusion $0^{*\Delta}_M \supseteq (0^{*\Delta}_M)'$. 
Fix $z \in (0^{*\Delta}_M)'$. 
    Then there exist $e_0 \geq 0$ and $c \in R^{\circ}$ 
    such that
    $z \otimes F^e_*c =0$ in  $M \otimes F^e_*R(p^e\Delta)$ for every $e \geq e_0$. 
    Take $c' \in R^{\circ}$ satisfying $\mathrm{div}(c') \geq \Delta$. 
    By $c' \in R(-\Delta)$, we get 
    \[
\cdot c' :     R(p^e\Delta) \to R( (p^e-1)\Delta), 
    \]
    which induces 
\[
M \otimes F^e_*R(p^e\Delta) \xrightarrow{\cdot (1 \otimes F^e_*c')} M \otimes F^e_*R((p^e-1)\Delta). 
\]
    Therefore, $z \otimes F^e_*(cc')=0$ in $M \otimes F^e_*\sO_X((p^e-1)\Delta)$. 
\end{proof}

\begin{prop}\label{p wt0 vs tight}
Let $(R, \m)$ be a $d$-dimensional $F$-finite Noetherian normal local ring of characteristic $p>0$. 
Take a $\Q$-divisor $D$ on $\Spec R$ and $z \in H^d_{\m}(R(D))$. 
Then the following are equivalent. 
\begin{enumerate}
\item $z \in 0^{*\{D\}}_{H^d_{\m}(R(D))}$. 
\item There exist 
an integer $e_0\geq 0$  and $c \in R^{\circ}$ 
    such that
\[
z \in \bigcap_{e \geq e_0} 
\Ker(H^d_{\m}(R(D)) 
\xrightarrow{F^e} H^d_{\m}(R(p^eD)) \xrightarrow{\cdot F_*^ec} H^d_{\m}(R(p^eD))). 
\]
\end{enumerate}
\end{prop}

\begin{proof}
It follows from  Remark \ref{r-tight-M/N-vs-0-log} that (1) $\iff$ (2)'. 
\begin{enumerate}
\item[(2)'] There exist an integer $e_0\geq 0$  and $c \in R^{\circ}$ 
    such that, for every $e \geq e_0$, 
    the equality $z \otimes F^e_*c =0$ holds in $H^d_{\m}(R(p^eD)) \otimes F^e_*R(p^e\{D\})$. 
\end{enumerate}
By Lemma \ref{l-tight-M/N-vs-0}, it holds that 
\begin{eqnarray*}
H^d_{\m}(R(D)) \otimes F_*^eR(p^e\{D\}) 
&\simeq& H^d_{\m}(F_*^eR(p^e\rdown{D} + p^e\{D\})) 
=H^d_{\m}(F_*^eR(p^eD))\\
z \otimes F_*^ec &\leftrightarrow& (F_*^ec)\cdot F^e(z). 
\end{eqnarray*}
Therefore, we get (2) $\iff$ (2)'. 
\end{proof}


\begin{definition}[\cite{BSTZ10}*{Definition-Proposition 3.3}]\label{d-test-ideal}
Let $R$ be a Noetherian $F$-finite reduced ring of characteristic $p>0$. 
Set $E := \bigoplus_{\m} E_R(R/\m)$, 
where 
$\m$ runs over all the maximal ideals of $R$, and 
$E_R(R/\m)$ denotes the injective hull of the $R$-module $R/\m$. 
\begin{enumerate}
\item 
We set 
\[
\tau(R) := \Ann_R(0_E^{*}), 
\]
which is called the {\em test ideal} of $R$. 
\item 
We say that $t \in R$ is a {\em test element} of $R$ if 
$t \in \tau(R)$. In other words, $t \in R$ is a test element if and only if 
$t \cdot 0_M^{*}=0$ for every $R$-module $M$ (Remark \ref{r-test-ideal-summary}(1)). 
\item 
Assume that $R$ is a normal integral domain. 
{\cred Take an effective $\Q$-divisor $\Delta$ on $\Spec R$.} 
We set 
\[
\tau(R, \Delta) := \Ann_R(0_E^{*\Delta}).
\]
This ideal is called the {\em test ideal} of $(R, \Delta)$. 
\end{enumerate}




\end{definition}

\begin{remark}
In some literature, our test ideal $\tau(R)$ is denoted by $\wt{\tau}(R)$ or $\tau_b(R)$ 
(cf.\ \cite[Remark 2.6]{Tak21}). 
\end{remark}

\begin{remark}\label{r-test-ideal-summary}
Let $R$ be a Noetherian $F$-finite reduced ring of characteristic $p>0$. It is known that the following hold. 
\begin{enumerate}
\item If $t$ is a test element of $R$ (i.e., $t \in \tau(R)$) and $M$ is an $R$-module, then 
$t \cdot 0^*_M=0$ \cite{BSTZ10}*{the same proof as in Definition-Proposition 3.3}. 
\item If $R$ is a normal integral domain, $\Delta$ is an effective $\Q$-divisor on $\Spec R$, 
and $t \in \tau(R, \Delta)$, 
then 
$t \cdot 0^{*\Delta}_M=0$ {\cred for every $R$-module $M$} \cite{BSTZ10}*{Definition-Proposition 3.3}. 
\item If $t$ is a test element of $R$, then $t$ is a test element of $R_{\mathfrak p}$ 
{\cred and $\widehat{R_{\mathfrak p}}$}
for every prime ideal $\mathfrak p$ of $R$ 
{\cred and the $\mathfrak p R_{\mathfrak p}$-adic completion $\widehat{R_{\mathfrak p}}$ of 
$R_{\mathfrak p}$} \cite{LS01}*{Theorem 7.1(2){\cred (3)}}. 
\item $(\omega_R/\tau(\omega_R))^{\vee} \simeq 0^*_{H^d_{\m}(R)}$ by 
a similar argument to \cite[the proof of Lemma 3.5]{bst} 
{(cf.\  (\ref{eq:dual-test-quotient}))}. 
\end{enumerate}
\end{remark}

\subsection{Witt dualising modules}




\begin{notation}\label{n-global-non-normal}
Let $X$ be an $F$-finite Noetherian reduced scheme of characteristic $p>0$  
{\cred such that $\Supp(\omega_X) = X$.}
    For $e \geq 0$ and $n  \geq m \geq 1$, 
we define the finite morphism $\iota^e_{X, m, n} \colon W_mX \to W_nX$ as the composition 
\[
\iota^e_{X, m, n} \colon W_mX \xrightarrow{F^e} W_mX \hookrightarrow W_nX, 
\]
where $W_nX := (X, W_n\MO_X), W_mX := (X, W_m\MO_X)$, 
and 
$W_mX \hookrightarrow W_nX$ denotes the closed immersion corresponding to $R^{n-m} \colon W_n\MO_X \to W_m\MO_X$. 
Assume that  there exist objects 
\[
W_1\omega_X^{\mydot},\, W_2\omega_X^{\mydot},\, ...
\]
of the derived category of $W_n\cO_X$-modules and a set of isomorphisms
\[
\{\rho^e_{X,m,n} \colon W_m\omega_X^{\mydot} \xrightarrow{\simeq} (\iota^e_{X,m,n})^{!}W_n\omega_X^{\mydot} 
\,|\, e \geq 0,  n  \geq m \geq 1\}
\] 
in the derived category of $W_n\cO_X$-modules which satisfy the following.
\begin{enumerate}
\item $W_n\omega_X^{\mydot}$ is a dualising complex on 
$W_nX = (X, W_n\cO_X)$ for every integer $n>0$. 
\item 
The equality $(\iota^{e'}_{X,l,m})^{!}\rho^{e}_{X,m,n} \circ \rho^{e'}_{X,l,m}=\rho^{e+e'}_{X,l,n}$ holds for all $e,e',l,m,n \in \Z$ satisfying $e,e' \geq 0$ and $n \geq m \geq l \geq 1$.
\end{enumerate}
{\cred We define the {\em dualising $W_n\MO_X$-module} $W_n\omega_X$ as the dualising sheaf on the Noetherian scheme $(X, W_n\MO_X)$ 
in the sense of Subsection \ref{subsubsection dualising}.} 
\end{notation}


\begin{notation}\label{n-global}
We use the same notation as in Notation \ref{n-global-non-normal}. 
Moreover, assume that $X$ is an integral normal scheme. 
\end{notation}

\begin{remark}
If $R$ is an $F$-finite Noetherian ring of characteristic $p>0$ and 
$X$ is a separated scheme of finite type over $R$, 
then there exist dualising complexes $\{W_n\omega_X^{\mydot}\}_{n\geq 1}$ 
and a set of isomorphisms $\{ \rho^e_{X, m, n}\}_{e\geq 0, n\geq m \geq 1}$ satisfying 
the above properties (1) and (2) ({Theorem \ref{thm:X-exist dualizing complex}}). 
\end{remark}

{\cred
\begin{remark}\label{support-omega}
Let $X$ be an $F$-finite Noetherian reduced scheme of characteristic $p>0$ admitting a dualising complex $\omega_X^{\mydot}$. 
\begin{enumerate}
    \item If $X=\Spec R$ for an equi-dimensional local ring $R$, then $\Supp(\omega_R)=\Spec R$.
    Indeed, let $\delta$ be the dimension function associated with $\omega_R^{\mydot}$ as defined in \cite[Tag 0AWF]{stacks-project}.
    Then, for any generic points $\eta_1,\eta_2 \in \Spec R$, we have
    \[
    \delta(\eta_1)-\delta(\eta_2)=(\delta(\eta_1)-\delta(P))-(\delta(\eta_2)-\delta(P))=\dim R-\dim R=0,
    \]
    where $P \in \Spec R$ is the closed point of $R$.
    Therefore, by \cite[Tag 0AWK]{stacks-project}, we conclude that $\Supp(\omega_R)=\Spec R$.
    \item If $X$ is equi-dimensional and of finite type over a field $k$, then $\Supp(\omega_X)=X$.
    Indeed, let $f \colon X \to \Spec k$ and assume $\omega_X^{\mydot}=f^{!}k$.
    Let $\delta$ be the dimension function associated with $\omega_X^{\mydot}$.
    Then, for every generic point $\eta \in X$, we have
    \[
    \delta(\eta)=\mathrm{trdeg}_k \kappa(\eta)=\dim V(\eta)
    \]
    by \cite[Tag 0AWM]{stacks-project}.
    Since $X$ is equi-dimensional, it follows from \cite[Tag 0AWK]{stacks-project} that $\Supp(\omega_X)=X$.
\end{enumerate}
\end{remark}
}

\begin{remark}\label{r global non normal FR compatible}
We use Notation \ref{n-global-non-normal}. 
\begin{enumerate}
\item 
For all $e,m,n \in \Z$ satisfying $e \geq 0$ and $n \geq m \geq 1$, 
we define the morphism $T^e_{X,m,n}$ by
\[
T^e_{X,m,n} \colon (\iota^e_{X,m,n})_*W_m\omega_X^{\mydot} \xrightarrow{\rho^e_{X,m,n}} (\iota^e_{X,m,n})_*(\iota^e_{X,m,n})^{!}W_n\omega_X^{\mydot} \xrightarrow{{\rm adj}} W_n\omega_X^{\mydot},
\]
where ${\rm adj}$ is the natural morphism induced by the fact that $(\iota^e_{X,m,n})^{!}$ is the right adjoint of $(\iota^e_{X,m,n})_*$ 
\cite[Tag 0A9Y]{stacks-project}. 
Then the equality 
\[
T^{e}_{X,m,n} \circ (\iota^{e}_{X,m,n})_*T^{e'}_{X,l,m}=T^{e+e'}_{X,l,n}
\]
holds for all $e,e',l,m,n \in \Z$ satisfying $e,e' \geq 0$ and $n \geq m \geq l \geq 1$ (cf.\ Step \ref{s1:X-exist dualizing complex} in the proof of Theorem \ref{thm:X-exist dualizing complex}).
\item {There exists an $F_*^eW_n\MO_X$-module isomorphism} 
\[
R\cHom_{W_n\MO_X}(F_*^eW_m\MO_X, W_n\omega_X^{\mydot}) \simeq F_*^eW_m\omega_X^{\mydot}
\] 
for every triple $(e, m, n) \in \Z^3$ satisfying $e \geq 0$ and $1 \leq m \leq n$. 
\end{enumerate}
\end{remark}

{\cred 
\begin{remark}\label{r dualising open}
We use Notation \ref{n-global-non-normal}.
\begin{enumerate}
    \item Let $j \colon U \hookrightarrow X$ be an open immersion and $j_n \colon W_nU \hookrightarrow W_nX$ the induced morphism. 
    Since $j^* = j^{!}$ by \cite[{Tag 0AU0}]{stacks-project}, we set $W_n\omega_U^{\mydot}:=j_n^{*}W_n\omega_X^{\mydot}$, which satisfies the conditions in Notation \ref{n-global-non-normal} by Step~2 of Theorem \ref{thm:X-exist dualizing complex}.
    Furthermore, since $\Supp(\omega_X)=X$, we have $\omega_X|_U=\omega_U$ and $\Supp(\omega_U)=U$.
    \item Assume that $X$ is affine: $X=\Spec R$. 
    Let $S$ be a multiplicatively closed subset of $R$.
    Then 
    \[
    W_n\omega_{S^{-1}R}^{\mydot}:=W_n\omega_{R}^{\mydot} \otimes_{W_n(R)} W_n(S^{-1}R)
    \]
    satisfies the conditions in Notation \ref{n-global-non-normal}.
    Indeed, we have
    \[
    W_n\omega_{R}^{\mydot} \otimes_{W_n(R)} W_n(S^{-1}R) \simeq \varinjlim_{f \in S} W_n\omega^{\mydot}_{R_f},
    \]
    where $W_n\omega^{\mydot}_{R_f}$ satisfies the conditions in Notation \ref{n-global-non-normal} for every $f \in S$ by (1).
    Since the conditions in Notation \ref{n-global-non-normal} are preserved under filtered colimits, it follows that $W_n\omega_{S^{-1}R}^{\mydot}$ satisfies the conditions in Notation \ref{n-global-non-normal}.
    Furthermore, since $W_n\omega_{S^{-1}R}=\varinjlim_{f \in S}W_n\omega_{R_f}$, we obtain $\Supp(\omega_{S^{-1}R})=\Spec(S^{-1}R)$.    
\end{enumerate}
\end{remark}
}



\begin{nothing}[${\cccred \mathbf{V}}^*, {\cccred \mathbf{R}}^*, T^e_n$]\label{n V^* R^* T}
We use Notation \ref{n-global-non-normal}. 
In this article, {we shall often consider the $W_n\omega_X$-dual denoted by} $(-)^* := {\cred \mathcal Hom}_{W_n\MO_X}(-, W_n\omega_X)$. 
For example, the $W_n\MO_X$-module homomorphisms 
\[
{\cccred \mathbf{V}} : F_*W_{n-1}\MO_X \to W_{n}\MO_X\qquad\text{and}\qquad {\cccred \mathbf{R}}: W_n\MO_X \to W_{n-1}\MO_X
\]
induce 
\[
{\cccred \mathbf{V}}^* : W_n\omega_X \to F_*W_{n-1}\omega_X\qquad\text{and}\qquad 
{\cccred \mathbf{R}}^* : W_{n-1}\omega_X \to W_n\omega_X. 
\]
Here, in order to get ${\cccred \mathbf{V}}^*$, we used  the following isomorphisms 
({\cred Remark \ref{r global non normal FR compatible}}): 
\[
\cHom_{W_n\MO_X}(W_n\MO_X,  W_{n}\omega_X) \simeq W_n\omega_X, \qquad 
\cHom_{W_n\MO_X}(F_*W_{n-1}\MO_X,  W_{n}\omega_X) \simeq F_*W_{n-1}\omega_X. 
\]
For integers $n>0$ and  $e \geq 0$, we define
\[
\PsiT^e_n : F_*^eW_n\omega_X \to W_n\omega_X
\]
as the trace map of $F^e : W_n\MO_X \to F_*^eW_n\MO_X$, that is, 
$T^e_n$ is the $W_n\MO_X$-module homomorphism 
obtained by applying $\cHom_{W_nR}(-, W_n\omega_R)$ 
to the $e$-the iterated Frobenius homomorphism $F^e : W_n\MO_X \to F_*^eW_n\MO_X$, 
which is $W_n\MO_X$-linear. 
Set $\PsiT_n := \PsiT^1_n$. 
\end{nothing}

\begin{remark} \label{remark:dual-omega-sequence}
By applying $\mathcal Hom_{W_n\MO_X}(-, W_n\omega_X)$ to the short exact sequence
\[
0 \to F_*W_{n-1}\MO_X \xrightarrow{{\cccred \mathbf{V}}}  W_n\MO_X \xrightarrow{{\cccred \mathbf{R}}^{n-1}}  \MO_X \to 0,
\]
we get the following exact sequence
\[
0 \to \omega_X \xrightarrow{({\cccred \mathbf{R}}^{n-1})^*} W_n\omega_X \xrightarrow{{\cccred \mathbf{V}}^*} F_*W_{n-1}\omega_X \to {\cred \mathcal Ext^1_{W_n\MO_X}(\MO_X, W_n\omega_X)}. 
\]
If $X$ is Cohen-Macaulay, then ${\cred \mathcal Ext^1_{W_n\MO_X}(\MO_X, W_n\omega_X)} \simeq \mathcal Ext^1_{\MO_X}(\MO_X, \omega_X)=0$ 
{\cred by Grothendieck duality}, and hence ${\cccred \mathbf{V}}^* \colon W_n\omega_X \to F_*W_{n-1}\omega_X$ is surjective. 

{\cred Similarly, if $X$ is regular, then we see that $T_n^e$ is surjective for all $e\geq 0$ and $n>0$ by induction on $n$.}
\end{remark}

\begin{dfn}\label{d-W_n omega D}
We use Notation \ref{n-global}. 
For an integer $n>0$ and a $\Q$-divisor $D$ on $X$, we define 
the coherent $W_n\MO_X$-module $W_n\omega_X(D)$ as follows 
\[
W_n\omega_X(D) := \cHom_{W_n\MO_X}(W_n\MO_X(-D), W_n\omega_X).
\]
Set $\omega_X(D) := W_1\omega_X(D)$. 
\end{dfn}

\begin{remark}\label{r-W_n omega D}
We use the same notation as in Definition \ref{d-W_n omega D}. 
\begin{enumerate}
\item 
If $n=1$, then 
\begin{align}\label{r-e1 W_n omega D}
\omega_X(D) &= W_1\omega_X(D) = \cHom_{\MO_X}(\MO_X(-D), \omega_X) \\
&\simeq \MO_X(K_X-\rdown{-D}) = \MO_X(K_X+\rup{D}). \nonumber
\end{align}
In particular, if $D$ is a Weil divisor, then $\omega_X(D) \simeq \MO_X(K_X+D)$. 
On the other hand, we have 
\[
\omega_X(D) \not\simeq \MO_X(K_X+D)\, (=\MO_X(\rdown{K_X+D}))
\]
in general. 
\item 
We now consider the case when $n$ is an arbitrary positive integer. 
In contrast to (\ref{r-e1 W_n omega D}), we have $W_n\MO_X(-D) \not\simeq W_n\MO_X(\rdown{-D})$ 
and  
\[
W_n\omega_X(D) \not\simeq W_n\MO_X(K_X+\rup{D}) 
\]
in general. 
\end{enumerate}
\end{remark}

\begin{remark}\label{r-W_nomega-j^*}
We use Notation \ref{n-global}. 
Take $\Q$-divisors $D_1$ and $D_2$ on $X$ satisfying $D_1 \leq D_2$. 
Then the natural inclusion 
$j : W_n\MO_X(-D_2) \hookrightarrow W_n\MO_X(-D_1)$ 
induces 
\[
j^* : W_n\omega_X(D_1) \to W_n\omega_X(D_2),
\]
which can be explicitly written as 
\[
j^* : \cHom_{W_n\MO_X}(W_n\MO_X(-D_1), W_n\omega_X) \to \cHom_{W_n\MO_X}(W_n\MO_X(-D_2), W_n\omega_X). 
\]

Let us show that $j^* : W_n\omega_X(D_1) \to W_n\omega_X(D_2)$ is injective. 
Both $W_n\omega_X(D_1)$ and $W_n\omega_X(D_2)$ are $S_2$, and hence we may assume that $X$ is regular. 
It is easy to see that we have the following commutative diagram in which each horizontal sequence is exact: 
\[
\begin{tikzcd}
0 \arrow[r] & \omega_X(D_1) \arrow[r, "({\cccred \mathbf{R}}^{n-1})^*"] \arrow[d, "j^*"] & W_n\omega_X(D_1) \arrow[r, "{\cccred \mathbf{V}}^*"] \arrow[d, "j^*"] & 
F_*W_{n-1}\omega_X(pD_1) \arrow[r]\arrow[d, "j^*"]  & 0\\
0 \arrow[r] & \omega_X(D_2) \arrow[r, "({\cccred \mathbf{R}}^{n-1})^*"] & W_n\omega_X(D_2) \arrow[r, "{\cccred \mathbf{V}}^*"] & F_*W_{n-1}\omega_X(pD_2) \arrow[r] & 0. 
\end{tikzcd}
\]
Then the injectivity of $j^*$ follows from the snake lemma and induction on $n$. 
\end{remark}

\begin{nothing}[The log versions of ${\cccred \mathbf{V}}^*, {\cccred \mathbf{R}}^*, T^e_n$]\label{n V^* R^* T log}
We now introduce the log version of (\ref{n V^* R^* T}). 
We use Notation \ref{n-global}. 
Take an integer $e \geq 0$  and a $\Q$-divisor $D$ on $X$. 
By applying $(-)^* := \mathcal Hom_{W_n\MO_X}(-, W_n\omega_X)$ to the $W_n\MO_X$-module homomorphisms 
\begin{eqnarray*}
&F^e:& W_n\MO_X(-D) \to F^e_*W_n\MO_X(-p^eD)\\
&{\cccred \mathbf{V}} :& F_*W_{n-1}\MO_X(-pD) \to W_{n}\MO_X(-D)\\
&{\cccred \mathbf{R}}:& W_n\MO_X(-D) \to W_{n-1}\MO_X(-D),    
\end{eqnarray*}
we get the following $W_n\MO_X$-module homomorphisms 
\begin{eqnarray*}
&T^e_{n}:& F^e_*W_n\omega_X(p^eD) \to W_n\omega_X(D)\\
&{\cccred \mathbf{V}}^* :& W_{n}\omega_X(D) \to F_*W_{n-1}\omega_X(pD)\\
&{\cccred \mathbf{R}}^*:& W_{n-1}\omega_X(D) \to W_{n}\omega_X(D)    
\end{eqnarray*}
Set $T_n := T^1_n$. 
\end{nothing}

\begin{lem}\label{l-W_nR-local-cplt}
Let $R$ be a Noetherian $F$-finite 
reduced ring of characteristic $p>0$ {\cred such that $\Supp(\omega_R) =\Spec R$}. 
\begin{enumerate}
\item 
If $S$ is a multiplicatively closed subset of $R$, 
then 
\[
W_n\omega_R \otimes_{W_nR} W_n(S^{-1}R) \simeq W_n\omega_{S^{-1}R}.
\]
\item 
If $(R, \m)$ is a local ring and $\widehat{R}$ denotes the $\m$-adic completion of $R$, 
then
\[
W_n\omega_R \otimes_{W_nR} W_n\widehat{R} \simeq W_n\omega_{\widehat{R}}.
\]
\end{enumerate}
\end{lem}

\begin{proof}
{\cred The assertion (1) follows from Remark \ref{r dualising open}.} 
Let us show (2). 
Recall that we have 
\[
E:=E_{W_nR} = E_{W_nR} \otimes_{W_nR} W_n\widehat{R} = E_{W_n\widehat{R}}. 
\]
for the injective hulls $E_{W_nR}$ and $E_{W_n\widehat{R}}$ of the residue field. 
By Matlis duality, we have 
 \[
 H^d_{\m}(W_nR) = \Hom_{W_nR}(W_n\omega_{R}, E) \quad \text{and} \quad
 H^d_{\m}(W_n\widehat{R}) = \Hom_{W_n\widehat{R}}(W_n\omega_{\widehat{R}}, E) 
 =(W_n\omega_{\widehat{R}})^{\vee}, 
 \]
 where $(-)^{\vee} :=\Hom_{W_n\widehat{R}}(-, E)$. 
 Therefore, 
\begin{eqnarray*}
 H^d_{\m}(W_n\widehat{R}) 
 &=&H^d_{\m}(W_nR) \otimes_{W_nR} W_n\widehat{R}\\
 &=& \Hom_{W_nR}(W_n\omega_{R}, E) \otimes_{W_nR} W_n\widehat{R}\\
   &=& \Hom_{W_n\widehat{R}}(W_n\omega_{R} \otimes_{W_nR} W_n\widehat{R}, E \otimes_{W_nR} W_n\widehat{R}) \\
   &=& \Hom_{W_n\widehat{R}}(W_n\omega_{R} \otimes_{W_nR} W_n\widehat{R}, E)\\
   &=& (W_n\omega_{R} \otimes_{W_nR} W_n\widehat{R})^{\vee}, 
\end{eqnarray*}
{\cred where the third equality holds, because $W_nR \to W_n\widehat{R}$ is flat and $W_n\omega_R$ is a finitely generated $W_nR$-module \cite[Theorem 7.11]{Matsumura}}. 
Applying $(-)^{\vee} = \Hom_{W_n\widehat{R}}(-, E)$ again, we obtain 
\[
W_n\omega_{R} \otimes_{W_nR} W_n\widehat{R} \simeq 
(W_n\omega_{R} \otimes_{W_nR} W_n\widehat{R})^{\vee\vee} \simeq 
H^d_{\m}(W_n\widehat{R})^{\vee} \simeq (W_n\omega_{\widehat{R}})^{\vee\vee} \simeq W_n\omega_{\widehat{R}}.  
\]
Thus (2) holds.  
\end{proof}




\subsection{Quasi-$F^e$-splitting and quasi-$F$-regularity}\label{ss QF^eS}

In this subsection, we recall the definition of a quasi-$F^e$-splitting and its variants. 
For their foundational properties, we refer to 
\cite{TWY}.

We use Notation \ref{n-global}. 
Take $e \in \Z_{>0}$ and let $\Delta$ be a (non-necessarily effective) $\bQ$-divisor on $X$. We define a $W_n\MO_X$-module $Q^e_{X,\Delta,n}$ and 
a $W_n\MO_X$-module homomorphism $\Phi^e_{X, \Delta, n}$ by the following pushout diagram of $W_n\MO_X$-module homomorphisms: 
\begin{equation} \label{diagram:quasi-F-split-definition}
\begin{tikzcd}
W_n\cO_X(\Delta) \arrow{r}{F^e} \arrow{d}{{\cccred \mathbf{R}}^{n-1}} & F^e_* W_n \cO_X(p^e\Delta)  \arrow{d}\\
\cO_X(\Delta) \arrow{r}{\Phi^e_{X, \Delta, n}}&
\arrow[lu, phantom, "\usebox\pushoutdr" , very near start, yshift=0em, xshift=0.6em, color=black] Q^e_{X, \Delta, n}.
\end{tikzcd}
\end{equation}


\begin{dfn}\label{d-IQFS}
We use Notation \ref{n-global}. 
Let $\Delta$ be a $\bQ$-divisor on $X$. 
Take integers $n>0$ and $e>0$. 
We say that $(X, \Delta)$ is {\em $n$-quasi-$F^e$-split} 
if $\rdown{\Delta}=0$ and the map  
\begin{multline*}
\Hom_{W_n\MO_X}(\Phi_{X, \Delta, n}^e, W_n\omega_X(-K_X)) : \Hom_{W_n\MO_X}(Q^e_{X, \Delta, n}, W_n\omega_X(-K_X)) \\ \longrightarrow 
\Hom_{W_n\MO_X}(\MO_X, W_n\omega_X(-K_X))
\end{multline*}
is surjective. 

We say that  $(X, \Delta)$ {\em quasi-$F^e$-split} if it is $n$-quasi-$F^e$-split for some $n \in \Z_{>0}$. 
We say that  $(X, \Delta)$ 
 {\em $n$-quasi-$F^e$-pure} (resp.\ {\em quasi-$F^e$-pure}) 
if 
 there exists an open cover $X = \bigcup_{i \in I} X_i$ such that 
 $(X_i, \Delta|_{X_i})$ is  {$n$-quasi-$F^e$-split} (resp.\ 
 {quasi-$F^e$-split}) for all $i \in I$. 
\end{dfn}

{The following characterisation of $n$-quasi-$F^e$-splitting in terms of local cohomology will be used extensively in this paper. 
{\cred Recall that we have $R(D) = \Gamma(\Spec R, \MO_{\Spec R}(D))$ 
for a $\Q$-divisor $D$ on $\Spec R$.} 

\begin{lemma}\label{lem:char-quasi-Fe-local-coh} 
{Let $(R, \m)$ be a $d$-dimensional Noetherian $F$-finite normal local ring of characteristic $p>0$.} 
Let $\Delta$ be a $\bQ$-divisor on {$\Spec R$}. 
{Take integers $n>0$ and $e>0$.} 
Then $(R,\Delta)$ is $n$-quasi-$F^e$-split if and only if $\rdown{\Delta}=0$ and the map 
\[
H_{\mathrm{m}}^d\left(\Phi_{R, K_R+\Delta, n}^e\right): H_{\mathrm{m}}^d\left(R\left(K_R+\Delta\right)\right) \rightarrow H_{\mathrm{m}}^d\left(Q_{R, K_R+\Delta, n}^e\right)
\]
is injective.
\end{lemma}}

\begin{proof}
See {\cite[Lemma 3.10]{TWY}}. 
\end{proof}

\begin{dfn}\label{d-QFR}
We use Notation \ref{n-global}. 
Let $\Delta$ be a $\bQ$-divisor on $X$. 
Take $n \in \Z_{>0}$. 
\begin{enumerate}
\item 
We say that $(X, \Delta)$ is {\em globally $n$-quasi-$F$-regular} 
if 
\begin{enumerate}
\item $\rdown{\Delta}=0$, and 
\item  for every effective $\Q$-divisor $E$ on $X$, 
there exists $\epsilon \in \Q_{>0}$  
such that $(X, \Delta + \epsilon E)$ is $n$-quasi-$F^e$-split 
for all $e \in \Z_{>0}$. 
\end{enumerate}
\item 
We say that $(X, \Delta)$ is {\em globally quasi-$F$-regular} 
if $(X, \Delta)$ is globally $m$-quasi-$F$-regular for some $m \in \Z_{>0}$. 
\end{enumerate}
\end{dfn}

\begin{dfn}\label{d-weak-QFR}
We use Notation \ref{n-global}.  
Let $\Delta$ be a $\bQ$-divisor on $X$. 
We say that $(X, \Delta)$ is {\em {feebly} globally quasi-$F$-regular} 
if 
\begin{enumerate}
\item $\rdown{\Delta}=0$, and 
\item 
for every effective $\Q$-divisor $E$ on $X$, 
there exist 
$n \in \Z_{>0}$ and $\epsilon \in \Q_{>0}$  
such that $(X, \Delta + \epsilon E)$ is $n$-quasi-$F^e$-split 
for all $e \in \Z_{>0}$. 
\end{enumerate}
\end{dfn}

In this article, we mainly treat the case when $X$ is affine. 
{Under this assumption,} we drop \lq\lq globally". 

\begin{dfn}\label{d QFR prelim}
    We use Notation \ref{n-global}. {\cred 
Let $\Delta$ be a $\bQ$-divisor on $X$.} Further, assume that $X$ is affine. 
    \begin{enumerate}
    \item 
    We say that $(X, \Delta)$ is ($n$-){\em quasi-$F$-regular} 
if $(X, \Delta)$ is globally ($n$-)quasi-$F$-regular. 
    \item 
We say that $(X, \Delta)$ is {\em feebly quasi-$F$-regular} 
if $(X, \Delta)$ is feebly  globally quasi-$F$-regular. 
\item Let $R$ be an $F$-finite Noetherian normal integral domain of characteristic $p>0$. 
We say that $(R, \Delta)$ is  ($n$-){\em quasi-$F$-regular} (resp.\ {\em feebly quasi-$F$-regular}) if so is $(\Spec R, \Delta)$. 
    \end{enumerate}
\end{dfn}

By definition, it holds that 
\[
\text{quasi-$F$-regular} \Rightarrow 
\text{feebly quasi-$F$-regular}.  
\]
We shall prove that the opposite implication holds when $X$ is affine  and $K_X$ is $\Q$-Cartier (Theorem \ref{thm:chara qFr pair2}).

\section{Quasi-test submodules and quasi-\texorpdfstring{$F$}--rationality}\label{s q-F-rat}

Let $R$ be a {\cred Noetherian} $F$-finite equi-dimensional reduced ring of characteristic $p>0$. 
The purpose of this section is to introduce quasi-$F$-rational singularities and establish some foundational results about them. 
The key objects we consider are 
\begin{enumerate}
\item the {\em $n$-quasi-test $R$-submodule} $\tau_n(\omega_R) \subseteq \omega_R$, and
\item the {\em $n$-quasi-tight closure} $0^*_n \subseteq H^{\dim R}_{\m}(R)$, 
\end{enumerate}
where {(2)} is defined only when $(R, \m)$ is a local ring, whilst 
we do not need such an additional assumption for {(1)}. 
For an overview of  this section, we refer to Subsection \ref{ss intro q F rat}. 

{We will work with two natural settings.}
{\cred By Remark \ref{support-omega}, for a ring $R$ in Notation \ref{n-local} satisfies the condition in Notation \ref{n-non-local}.
}
\begin{notation}\label{n-local}
Let $(R, \m)$ be a Noetherian $F$-finite equi-dimensional reduced local ring of characteristic $p>0$. 
Set $d := \dim R$. 
\end{notation}

\begin{notation}\label{n-non-local}
{
Let $R$ be a Noetherian $F$-finite reduced  ring of characteristic $p>0$ {\cred such that $\Supp(\omega_R)=\Spec R$}. 
Set $d := \dim R$.
}
\end{notation}







\subsection{Quasi-test submodules (no boundary) {and quasi-$F$-rationality} }\label{ss q test submod}
Recall that we have 
the trace map 
\[
T_n: F_*W_n\omega_R \to W_n\omega_R,
\]
which is the $W_nR$-module homomorphism 
obtained by applying  $\Hom_{W_nR}(-, W_n\omega_R)$ to 
Frobenius $F: W_nR \to F_*W_nR$ (\ref{n V^* R^* T}).

\begin{definition}
We work in the general setting (Notation \ref{n-non-local}).  
Let $M$ be a $W_nR$-submodule of $W_n\omega_R$. 
\begin{enumerate}
    \item We say that $M$ is {\em ${\cccred T_n}$-stable} if $\PsiT_n(F_*M) \subseteq M$.
    \item We say that $M$ is {\em co-small} if there exists $c \in R^{\circ}$ such that $[c] \cdot W_n\omega_R \subseteq M$.
    In other words, there exists an open dense subset $U$ of $X$ such that $\widetilde{M}|_U=W_n\omega_U$.
\end{enumerate}
\end{definition}

\begin{remark}\label{r nonzero cosmall example}
\begin{enumerate}
\item 
If $M$ is co-small, then $[c] \cdot M$ is also co-small for $c \in R^{\circ}$.
On the other hand, even if $M$ is ${\cccred T_n}$-stable, 
$[c] \cdot M$ is not necessarily ${\cccred T_n}$-stable.
For example, 
if $R$ is $F$-rational, 
then $R$ has no non-trivial ${\cccred T_1}$-stable $R$-submodules of $\omega_R$ \cite[Lemma 2.34]{bst}. 
In particular, $[c] \cdot \omega_R$ is not ${\cccred T_1}$-stable for every $c \in R^{\circ}$ 
satisfying $[c] \cdot \omega_R \neq \omega_R$. 
\item 
{\cred 
Opposed to the classical case (i.e., the case when $n=1$), 
a nonzero $W_n\omega_R$-submodule of $W_n\omega_R$ is not necessarily co-small. 
For example, if $n=2$, $R = \F_p$, and $M:=p W_2\omega_R$, 
then we get $M \neq 0$ and $[c] \cdot W_2\omega_R \not\subseteq M$ 
for every $c \in R^{\circ} = \F_p^{\times}$.

Indeed, $W_2(\F_p)  = \Z/p^2\Z$ is a Gorenstein artinian local ring, and hence $W_2\omega_{\F_p} = W_2(\F_p)$. 
In particular, $M = p W_2\omega_R = p\Z / p^2\Z \neq 0$. 
For $c \in \F_p^{\times}$, 
we have $[c] \in W_2(\F_p)^{\times} = ( \Z/p^2\Z)^{\times}$, 
which implies 
\[
[c] \cdot W_2\omega_R  = [c] \cdot \Z/p^2\Z = \Z/p^2\Z \not\subseteq 
p\Z/p^2\Z = M. 
\]
}


\end{enumerate}
\end{remark}

\medskip

Recall that we have the following $W_nR$-module homomorphisms (\ref{n V^* R^* T}): 
\[
{\cccred \mathbf{V}}^* : W_n\omega_R \to F_*W_{n-1}\omega_R
\]
\[
{\cccred \mathbf{R}}^* : W_{n-1}\omega_R \to W_{n}\omega_R. 
\]

\begin{proposition}\label{prop: first prop of F-stable}
We work in the general setting (Notation \ref{n-non-local}).  
Fix a test element $t \in R^{\circ}$ {(see Definition \ref{d-test-ideal}(2))}.
Let $M$ be a ${\cccred T_n}$-stable co-small $W_nR$-submodule of $W_n\omega_R$.
Then the following hold. 
\begin{enumerate}
    \item $({\cccred \mathbf{R}}^{*})^{-1}(M)$
    and $\PsiT_{n-1} \circ {\cccred \mathbf{V}}^*(M)$ {\cccred(resp.\ $T_n(F_*M)$)} are {\cccred(resp.\ is)} co-small and ${\cccred T_{n-1}}$-stable {\cccred(resp.\ $T_n$-stable)}. 
    \item $[t^2] \cdot W_n\omega_R \subseteq M$.
\end{enumerate}
\end{proposition}

\begin{proof}
Let us show (1). 
First, we prove that $({\cccred \mathbf{R}}^{*})^{-1}(M)$, {\cred $T_n(F_*M)$}, and $\PsiT_{n-1} \circ {\cccred \mathbf{V}}^*(M)$ are co-small.
We take an affine open dense subset $U$ of $X := \Spec R$ 
such that 
$\widetilde{M}|_{U}=W_n\omega_U$, and 
each of 
$\PsiT_{n-1}$, {\cred $T_n$}, and ${\cccred \mathbf{V}}^*$ is surjective {\cred (Remark \ref{remark:dual-omega-sequence})}.
Then we have
\begin{align*}
    ({\cccred \mathbf{R}}^*)^{-1}(\wt{M})|_U&=({\cccred \mathbf{R}}^*)^{-1}(W_n\omega_U)=W_{n-1}\omega_U, \\
    {\cred \PsiT_{n}(\wt{M})|_U} 
    &{\cred =\PsiT_{n}(W_n\omega_U)=
    W_{n}\omega_U,} \\
    \PsiT_{n-1}\circ {\cccred \mathbf{V}}^*(\wt{M})|_U&=\PsiT_{n-1}\circ {\cccred \mathbf{V}}^*(W_n\omega_U)=
    T_{n-1}(F_*W_{n-1}\omega_U)= W_{n-1}\omega_U. 
\end{align*}
Thus $({\cccred \mathbf{R}}^*)^{-1}(M)$, {\cred $T_n(F_*M)$}, and $\PsiT_{n-1} \circ {\cccred \mathbf{V}}^*(M)$ are co-small.

{\cred 
We see that $T_n(F_*M)$ is {\cccred $T_{n}$-stable}, because 
the inclusion $T_n(F_*M) \subseteq M$ implies $T_n(F_*(T_n(F_*M))) \subseteq T_n(F_*M)$.}  
Now, the fact that $({\cccred \mathbf{R}}^{*})^{-1}(M)$ and $\PsiT_{n-1} \circ {\cccred \mathbf{V}}^*(M)$ are ${\cccred T_{n-1}}$-stable follows from the commutative diagrams
\begin{equation*}
\begin{tikzcd}
F_*W_{n-1}\omega_R \arrow{r}{F_*{\cccred \mathbf{R}}^*} \arrow{d}{\PsiT_{n-1}} & 
F_*W_n\omega_R \arrow{d}{\PsiT_{n}}\\
W_{n-1}\omega_R \arrow{r}{{\cccred \mathbf{R}}^*} & 
W_n\omega_R
\end{tikzcd}
\qquad
\begin{tikzcd}
F_*W_{n}\omega_R \arrow{r}{F_*{\cccred \mathbf{V}}^*} \arrow[phantom, from=1-1, to=2-2, "\tiny (\dagger)"] \arrow{d}{\PsiT_{n}} & 
F^2_*W_{n-1}\omega_R \arrow{r}{F_*T_{n-1}} \arrow{d}{F_*\PsiT_{n-1}} & 
F_*W_{n-1}\omega_R \arrow{d}{\PsiT_{n-1}} \\
W_{n}\omega_R \arrow{r}{{\cccred \mathbf{V}}^*} & F_*W_{n-1}\omega_R \arrow{r}{T_{n-1}} & W_{n-1}\omega_R.
\end{tikzcd}
\end{equation*}
Thus (1) holds. 

Let us show (2) by induction on $n$. 
We may assume that $R$ is {\cred a complete} local {\cred ring} with maximal ideal $\m$ (Remark \ref{r-test-ideal-summary}(3)). 
The base case $n=1$ of this induction follows  from 
$\tau(\omega_R) \subseteq M$ \cite[Definition 2.33]{bst} and 
Matlis duality: $\omega_R/\tau(\omega_R) \simeq \Hom_R(0^*_{H^d_\m(R)}, E)$ (Remark \ref{r-test-ideal-summary}). 
In particular, we have $t \cdot \omega_R  \subseteq \tau(\omega_R) \subseteq M$.
We assume $n \geq 2$.
By the following exact sequence
\[
0 \to \frac{\omega_R}{(({\cccred \mathbf{R}}^{n-1})^*)^{-1}(M)} \xrightarrow{({\cccred \mathbf{R}}^{n-1})^*} \frac{W_n\omega_R}{M} \xrightarrow{{\cccred \mathbf{V}}^*} \frac{F_*W_{n-1}\omega_R}{{\cccred \mathbf{V}}^*(M)},
\]
it is enough to show that 
\[
[t] \cdot F_*W_{n-1}\omega_R \subseteq {\cccred \mathbf{V}}^*(M).
\]
This is a consequence of
\[
[t] \cdot F_*W_{n-1}\omega_R = F_*([t^p]\cdot W_{n-1}\omega_R) \subseteq F_*(\PsiT_{n-1} \circ {\cccred \mathbf{V}}^*(M)) ={\cccred \mathbf{V}}^* \circ \PsiT_n(F_*M) \subseteq {\cccred \mathbf{V}}^*(M),
\]
where the first inclusion follows from 
 the induction hypothesis and the fact that $\PsiT_{n-1}\circ {\cccred \mathbf{V}}^*(M)$ is co-small and ${\cccred T_{n-1}}$-stable, {whilst the subsequent equality is the commutativity of the diagram ($\dagger$) above}. 
 Thus (2) holds. 
\end{proof}


\begin{proposition}\label{prop: test submod}
We work in the general setting (Notation \ref{n-non-local}).   
Fix an integer $n>0$. 
Let $t \in R^{\circ}$ be a test element and let $c \in (t^2) \cap R^{\circ}$.
\begin{enumerate}
    \item There exists the smallest co-small ${\cccred T_{n}}$-stable $W_nR$-submodule $\tau(W_n\omega_R)$ of $W_n\omega_R$.
    \item For every co-small $W_nR$-submodule $M$ of $W_n\omega_R$, it holds that 
    \[
    \tau(W_n\omega_R)=\sum_{e \geq 0} \PsiT^e_n(F^e_*([c]\cdot M)). 
    \]
    \item 
    For every integer $e \geq 0$, 
    the inclusion  
   \[
     T^e_n(F^e_*([c]\cdot \tau(W_n\omega_R))) \subseteq 
   T^{e+1}_n(F^{e+1}_*([c]\cdot \tau(W_n\omega_R)))
   \]
   holds. 
    \item 
    There exists an integer $e_0 \geq 0$ such that 
the equality 
    \[
    \tau(W_n\omega_R)=T^{e}_n(F^e_*([c]\cdot M))
    \]
    holds for every  integer $e \geq e_0$ and every co-small ${\cccred T_{n}}$-stable $W_nR$-submodule $M$ of $W_n\omega_R$. 
    \item Let $L$ be a co-small $W_nR$-submodule of $W_n\omega_R$. 
    Then there exists an integer $e_L \geq 0$ such that 
    \[
    \tau(W_n\omega_R) = T^e_n(F^e_*([c] \cdot L))
    \]
    for every integer $e \geq e_L$. 
\end{enumerate}
\end{proposition}
We call $\tau(W_n\omega_R)$ the {\em $n$-quasi-test $(W_nR$-$)$submodule} of $W_n\omega_R$. {By {(2) and} (3), we have $\tau(W_n\omega_R) = T^e_n(F^e_*([c] \cdot W_n\omega_R))$ for $e \gg 0$, where $T^e_n \colon F^e_*W_n\omega_R \to W_n\omega_R$.}

\begin{proof}
Let us show (1). 
Set 
\[
\tau(W_n\omega_R) := \bigcap_{M} M,
\]
where $M$ runs over all the ${\cccred T_{n}}$-stable co-small $W_nR$-submodules $M$ of $W_n\omega_R$. 
Then it is clear that $\tau(W_n\omega_R)$ is ${\cccred T_{n}}$-stable. 
Moreover, $[t^2] \cdot W_n\omega_R \subseteq \tau(W_n\omega_R)$ 
(Proposition \ref{prop: first prop of F-stable}(2)), and hence 
$\tau(W_n\omega_R)$ is co-small. 
Thus (1) holds. 

Let us show $(2)$.
Set $N$ to be the right hand side of $(2)$. 
Since $\PsiT^0_n(F^0_*[c]\cdot M) =[c] \cdot M$ is co-small, also 
$N$ is co-small. 
Furthermore,  we have
\begin{align*}
     \PsiT_n(F_*N) 
=    T_n \left( F_* \sum_{e \geq 0} \PsiT^e_n(F^e_*([c]\cdot M))\right)
    =\sum_{e\geq 1} \PsiT_n^e(F^e_*([c] \cdot M)) 
    \subseteq N, 
\end{align*}
and hence  $N$ is ${\cccred T_{n}}$-stable. 
Since $N$ is co-small and  ${\cccred T_{n}}$-stable, 
(1) implies  
$\tau(W_n\omega_R) \subseteq N$. 
On the other hand, we have
\begin{align*}
    N
    &= \sum_{e \geq 0} \PsiT^e_n(F^e_*([c] \cdot M)) \\
    & \subseteq \sum_{e \geq 0} \PsiT^e_n(F^e_*([c] \cdot W_n\omega_R)) \\
    & \overset{(\star)}{\subseteq} \sum_{e \geq 0} \PsiT^e_n(F^e_*\tau(W_n\omega_R))\\ 
    & \subseteq \tau(W_n\omega_R), 
\end{align*}
where $(\star)$ follows from (1) and Proposition \ref{prop: first prop of F-stable}(2). 
Therefore, we get $\tau(W_n\omega_R) =N$. Thus (2) holds. 

Let us show (3). 
Since $T_n(F_*\tau(W_n\omega_R))$ is co-small and ${\cccred T_{n}}$-stable 
{\cred (Proposition \ref{prop: first prop of F-stable}(1))}, we have $\tau(W_n\omega_R) \subseteq T_n(F_*\tau(W_n\omega_R))$.
Furthermore, since $\tau(W_n\omega_R)$ is ${\cccred T_{n}}$-stable, we obtain 
the opposite inclusion $T_n(F_*\tau(W_n\omega_R)) \subseteq \tau(W_n\omega_R)$, and hence $T_n(F_*\tau(W_n\omega_R)) =\tau(W_n\omega_R)$. 
Therefore, it holds that 
\begin{align*}
T^{e+1}_n(F^{e+1}_*([c]\cdot \tau(W_n\omega_R)) )
&\supseteq T^{e+1}_n(F^{e+1}_*([c^p]\cdot\tau(W_n\omega_R))) \\
&=T^e_n(F^e_*([c]\cdot T_n(F_*\tau(W_n\omega_R)))) \\
&=T^e_n(F^e_*([c]\cdot \tau(W_n\omega_R))).
\end{align*}
Thus (3) holds.

Let us show (4). 
Since $W_nR$ is a Noetherian ring and $W_n\omega_R$ is a finitely generated $W_nR$-module, 
the assertion (3) enables us to find $e_0>0$ such that 
\[
T^{e_0}_n(F^{e_0}_*([c]\cdot \tau(W_n\omega_R))) = 
T^{e_0+1}_n(F^{e_0+1}_*([c]\cdot \tau(W_n\omega_R))) = \cdots. 
\]
Fix an integer $e \geq e_0$. 
Then (2) and (3) imply 
\[
\tau(W_n\omega_R) = T^e_n(F_*^e([c] \cdot \tau(W_n\omega_R))). 
\]
Therefore, 
\[
T^e_n(F_*^e([c] \cdot M)) 
\overset{(2)}{\subseteq} 
\tau(W_n\omega_R) = T^e_n(F_*^e([c] \cdot \tau(W_n\omega_R)))
\overset{(1)}{\subseteq} 
T^e_n(F_*^e([c] \cdot M)). 
\]
Thus (4) holds. 

Let us show (5). 
Since $L$ is co-small, there exists $c_L \in R^{\circ}$
such that $[c_L] W_n\omega_R \subseteq L \subseteq W_n\omega_R$, 
which implies 
\[
[c c_L] W_n\omega_R \subseteq [c] L \subseteq [c] W_n\omega_R. 
\]
By (4), we can find an integer $e_L \geq 0$ such that 
\[
\tau(W_n\omega_R) = T_n^e(F_*^e([cc_L]W_n\omega_R)) \quad\text{and}\quad 
\tau(W_n\omega_R) = T_n^e(F_*^e([c]W_n\omega_R))
\]
for every integer $e \geq e_L$, 
where the first equality is obtained by applying (4) after replacing $cc_L$ by $c$ (note that $e_L$ depends on $c_L$, and hence on $L$). 
Then $\tau(W_n\omega_R) = T_n^e(F_*^e([c] \cdot L))$. 
This concludes the proof of  (5). 
\end{proof}

\begin{proposition}\label{prop: test submod compat}
We work in the general setting (Notation \ref{n-non-local}).   
Then the following hold. 
\begin{enumerate}
    \item 
    If $S$ is a multiplicatively closed subset of $R$, then 
    \[
    \tau(W_n\omega_R) \otimes_{W_nR} W_n(S^{-1}R) \simeq \tau(W_n\omega_{S^{-1}R}).
    \]
    \item If $(R,\m)$ is a local ring, then
    \[
    \tau(W_n\omega_R) \otimes_{W_nR} W_n\widehat{R} \simeq \tau(W_n\omega_{\widehat{R}}), 
    \]
    where $\widehat{R}$ denotes the $\m$-adic completion of $R$. 
\end{enumerate}
\end{proposition}

\begin{proof}
Let $t \in R^{\circ}$ be a test element.
Recall that
\[
W_{n} (S^{-1}R) = [S]^{-1} W_{n} R, 
\]
where $[S] := \{ [s] \in W_nR \,|\, s \in S\}$, which is a multiplicatively closed subset of $W_nR$ (indeed, it is easy to construct a ring homomorphism 
$\theta : [S]^{-1} W_{n} R \to W_{n} (S^{-1}R)$; then we can directly check that $\theta$ is injective and surjective). 
We have  the trace maps 
\[
T^e_{R, n} : F_*^eW_n\omega_R \to W_n\omega_R \qquad\text{and}\qquad 
T^e_{S^{-1}R, n} : F_*^eW_n\omega_{S^{-1}R} \to W_n\omega_{S^{-1}R}. 
\]
Then we get 
\begin{align*}
    \tau(W_n\omega_R) \otimes_{W_nR} W_n(S^{-1}R)
    &\overset{{\rm (i)}}{=} \left(\sum_{e \geq 0}\PsiT^e_{R,n}(F^e_*[t^2] \cdot W_n\omega_R)\right) \otimes_{W_nR} W_n(S^{-1}R) \\
    &\overset{{\rm (ii)}}{=} \sum_{e \geq 0}\left(\PsiT^e_{R,n}(F^e_*[t^2] \cdot W_n\omega_R) \otimes_{W_nR} W_n(S^{-1}R)\right)  \\
    &\overset{{\rm (iii)}}{=}  \sum_{e \geq 0}  \PsiT^e_{S^{-1}R,n}(F^e_*[t^2] \cdot W_n\omega_{S^{-1}R}) \\
    &\overset{{\rm (iv)}}{=} \tau(W_n\omega_{S^{-1}R}), 
\end{align*}
where (i) and (iv) hold by Proposition \ref{prop: test submod}(2), 
(ii) follows from the fact that $\sum_{e \geq 0}$ commutes with the flat base change $(-) \otimes_{W_nR} W_n(S^{-1}R)$, 
and (iii) is assured by 
Lemma \ref{l-W_nR-local-cplt}. 
Thus we obtain the assertion $(1)$.
The assertion $(2)$ follows from a similar argument as above 
by using $W_n\omega_R \otimes_{W_nR} W_n\widehat{R} \simeq W_n\omega_{\widehat{R}}$ (Lemma \ref{l-W_nR-local-cplt}). 
\end{proof}

\begin{definition}\label{d-q-tau-omega}
We work in the general setting (Notation \ref{n-non-local}).   
\begin{enumerate}
    \item For an integer $n >0$, we set $\tau_n(\omega_R) := (({\cccred \mathbf{R}}^{n-1})^*)^{-1}(\tau(W_n\omega_R))$, 
    that is, $\tau_n(\omega_R)$ is the $R$-submodule  of $\omega_R$ which is 
    the inverse image of $\tau(W_n\omega_R)$ by {\cred the} $W_nR$-module homomorphism 
    \[
    ({\cccred \mathbf{R}}^{n-1})^* \colon \omega_R \to W_n\omega_R.
    \]
    We call $\tau_n(\omega_R)$ the {\em $n$-quasi-test $(R$-$)$submodule} of $\omega_R$. 
We have the following inclusions (Remark \ref{r-tau-omega-n-vs-n+1}): 
\begin{equation}\label{e1-q-tau-omega}
\tau_1(\omega_R) \subseteq \tau_2(\omega_R) \subseteq \cdots \subseteq 
\tau_n(\omega_R) \subseteq \tau_{n+1}(\omega_R) \subseteq \cdots \subseteq \omega_R. 
\end{equation}
    \item We set 
    $\tau^q(\omega_R) := \bigcup_{n =1}^{\infty} \tau_n(\omega_R)$, 
    which we call  the {\em quasi-test $(R$-$)$submodule} of $\omega_R$. 
    Since (\ref{e1-q-tau-omega}) is an ascending chain of $R$-submodules of $\omega_R$, 
    there exists $n_0>0$ such that 
\[
\tau^q(\omega_R) = \bigcup_{n =1}^{\infty} \tau_n(\omega_R) = \tau_n(\omega_R) \subseteq \omega_R
\]
for every integer $n \geq n_0$. 
In particular, $\tau^q(\omega_R)$ is an $R$-submodule of $\omega_R$. 
\end{enumerate}
\end{definition}

\begin{remark}\label{r-tau-omega-n-vs-n+1}
By applying $\Hom_{W_nR}(-, W_n\omega_R)$ to 
\[
0 \to F_*W_{n-1}R \xrightarrow{{\cccred \mathbf{V}}} W_nR \xrightarrow{{\cccred \mathbf{R}}^{n-1}} R \to 0, 
\]
we get the following commutative diagram 
in which each horizontal sequence is exact: 
\[
\begin{tikzcd}
0 \arrow[r] & \omega_R \arrow[r, "({\cccred \mathbf{R}}^{n-1})^*"] & W_n\omega_R \arrow[r, "{{\cccred \mathbf{V}}^*}"] & F_*W_{n-1}\omega_R\\
0 \arrow[r] & \tau_n(\omega_R) \arrow[r] \arrow[u, hook] & \tau(W_n\omega_R). \arrow[u, hook]
\end{tikzcd}
\]
By diagram chase, $({\cccred \mathbf{R}}^{n-1})^*$ induces an isomorphism 
\[
\tau_n(\omega_R) \xrightarrow{\simeq} 
\Ker\Big( \tau(W_n\omega_R)\hookrightarrow W_n\omega_R \xrightarrow{{\cccred \mathbf{V}}^*}F_*W_{n-1}\omega_R\Big). 
\]
By the injectivity of ${\cccred \mathbf{R}}^*$ and the following diagram of exact sequences
    \begin{equation*}
    \begin{tikzcd}
        0 \arrow{r} &\tau_n(\omega_R) \arrow{r} & \tau(W_n\omega_R) \arrow{r}{{\cccred \mathbf{V}}^*} \arrow{d}{{\cccred \mathbf{R}}^*} & F_*W_{n-1}\omega_R \arrow{d}{{\cccred \mathbf{R}}^*} \\
    0 \arrow{r} & \tau_{n+1}(\omega_R) \arrow{r} & \tau(W_{n+1}\omega_R) \arrow{r}{{\cccred \mathbf{V}}^*} & F_*W_n\omega_R, 
    \end{tikzcd}
    \end{equation*}
    we obtain an inclusion $\tau_n(\omega_R) \subseteq \tau_{n+1}(\omega_R)$.
\end{remark}

\begin{theorem}\label{thm:test submodule locali}
We work in the general setting (Notation \ref{n-non-local}).    
Then the following hold. 
\begin{enumerate}
    \item 
    If $S$ be a multiplicatively closed subset of $R$, 
    then 
    \[
    \tau_n(\omega_R) \otimes_R S^{-1}R \simeq \tau_n(\omega_{S^{-1}R})
    \quad\text{and}\quad \tau^q(\omega_R) \otimes_R S^{-1}R \simeq \tau^q(\omega_{S^{-1}R}).
    \]
    \item 
    {\cred If $(R, \m)$ is a local ring, then}
    \[
    \tau_n(\omega_R) \otimes_R \widehat{R} \simeq \tau_n(\omega_{\widehat{R}})\quad\text{and}\quad \tau^q(\omega_{R}) \otimes_R \widehat{R} \simeq \tau^q(\omega_{\widehat{R}}), 
    \]
    {\cred where $\widehat{R}$ denotes the $\m$-adic completion of $R$.}
\end{enumerate}
\end{theorem}

\begin{proof}
Both assertions follow from Proposition \ref{prop: test submod compat}. 
 {\ccred Here, we implicitly use 
{$R \otimes_{W_nR} 
W_n \widehat{R} \simeq \widehat R$ 
(by induction on $m$, we can prove that $(F_*^eW_mR) \otimes_{W_nR} 
W_n \widehat{R} \simeq F_*^eW_m\widehat R$ holds for $e \geq 0$ and $1 \leq m \leq n$}).} 
\end{proof}

{The following proposition describes the behaviour of $\tau_n$ under alterations and birational morphisms.} 
\begin{proposition}\label{prop: rat submod}
We work in the general setting (Notation \ref{n-non-local}).   
Fix an integer $n>0$. 
Let $f \colon Y \to X :=\Spec R$ be an alteration 
between integral schemes. 
Then the following hold. 
\begin{enumerate}
\item We have that 
\[
\tau_n(\omega_X) \subseteq 
(({\cccred \mathbf{R}}^{n-1})^*)^{-1}(\mathrm{Im}(T_{n}^f : f_*W_n\omega_Y \to W_n\omega_X)), 
\]
where $T_{n}^f : f_*W_n\omega_Y \to W_n\omega_X$ denotes the $W_n\MO_X$-module homomorphism obtained by applying $\mathcal Hom_{W_n\MO_X}(-, W_n\omega_X)$ to 
the induced $W_n\MO_X$-module homomorphism $f^* : W_n\MO_X \to f_*W_n\MO_Y$. 
\item {If $f : Y \to X$ is a proper birational morphism, then $\tau_n(\omega_X) \subseteq f_*\omega_Y \subseteq \omega_X$.}
\item 
If $\tau^q(\omega_X)=\omega_X$, then $X$ is pseudo-rational, that is, 
$f_*\omega_Y = \omega_X$ for every proper birational morphism $f: Y \to X$ from an integral normal scheme $Y$. 
\end{enumerate}
\end{proposition}

\begin{proof}
Let us show (1). 
Consider the  following commutative diagram in which each horizontal sequence is exact: 
\begin{equation} \label{eq:diagram-trace-wnomega}
\xymatrix{
0 \ar[r] & f_*\omega_Y \ar[r]^-{({\cccred \mathbf{R}}^{n-1})^*} \ar[d]^{T^f_1} & f_*W_n\omega_Y \ar[r]^-{{\cccred \mathbf{V}}^*} \ar[d]^{T^f_{n}} & f_*F_*W_{n-1}\omega_Y \ar[d]^{F_*T^f_{n-1}} \\
0 \ar[r] & \omega_X \ar[r]^-{({\cccred \mathbf{R}}^{n-1})^*} & W_n\omega_X \ar[r]^-{{\cccred \mathbf{V}}^*} & F_*W_{n-1}\omega_X. 
}
\end{equation}
{\cred On an open dense subset {\cred $X'$} of $X$ 
such that $X'$ is regular and 
the induced morphism $f|_{f^{-1}(X')} : f^{-1}(X') \to X'$ is finite}, 
the homomorphism $T^f_{1}: f_*\omega_Y \to \omega_X$ is surjective and so are the horizontal homomorphisms ${\cccred \mathbf{V}}^*$ in the above diagram. 
Hence $T^f_n : f_*W_n\omega_Y \to W_n\omega_X$ is surjective 
{\cred on $X'$} by 
the snake lemma and induction on $n$.   
Therefore, the image 
\[
I_n := {\rm Im}(T^f_{n} : f_*W_n\omega_Y \to W_n\omega_X) 
\]
is co-small.
Since $I_n$ is ${\cccred T_{n}}$-stable, we get $\tau(W_n\omega_X) \subseteq I_n$ by the 
 minimality of $\tau(W_n\omega_X)$. 
Therefore, 
\[
\tau_n(\omega_X) = 
(({\cccred \mathbf{R}}^{n-1})^*)^{-1}(\tau(W_n\omega_X)) \subseteq (({\cccred \mathbf{R}}^{n-1})^*)^{-1}(I_n). 
\]
Thus (1) holds. 

Let us show (2). 
{\cccred Note that, by \cite[\href{https://stacks.math.columbia.edu/tag/0AWN}{Tag 0AWN}]{stacks-project}, $\omega_Y$ is torsion-free.}
Since $f$ is birational, the vertical arrows in Diagram (\ref{eq:diagram-trace-wnomega}) are injective (as they are injective after removing the non-regular loci of $X$ and $Y$), and so the statement follows by (1) and diagram chase. Specifically, pick $a \in \tau_n(\omega_X)$. Then by (1), 
\[
({\cccred \mathbf{R}}^{n-1})^*(a) = T^f_n(b)
\]
for some $b \in f_*W_n\omega_Y$. Since ${\cccred \mathbf{V}}^*(({\cccred \mathbf{R}}^{n-1})^*(a)) = 0$ and $F_*T^f_{n-1}$ is injective, we get that ${\cccred \mathbf{V}}^*(b)=0$, and so $b = ({\cccred \mathbf{R}}^{n-1})^*(c)$ for some $c \in f_*\omega_Y$. Finally, as $({\cccred \mathbf{R}}^{n-1})^* : \omega_X \to W_n\omega_X$ is injective, we obtain $a = T^f_1(c)$. In particular, $a \in f_*\omega_Y$ as desired.

Now (3) follows immediately from (2), because, as in Definition \ref{d-q-tau-omega}(2), there exists $n>0$ such that $\tau_n(\omega_X)=\tau^q(\omega_X)=\omega_X$.
\qedhere

\end{proof}

\begin{definition}\label{d-q-F-rat2}
We work in the general setting (Notation \ref{n-non-local}).  
For an integer $n>0$, 
we say that $R$ is {\em $n$-quasi-$F$-rational} if $R$ is Cohen-Macaulay and $\omega_R = \tau_n(\omega_R)$. 
We say that $R$ is {\em quasi-$F$-rational} 
if $R$ is $n$-quasi-$F$-rational for some integer $n>0$. 
\end{definition}

\begin{remark}\label{r-q-F-rat2}
\begin{enumerate}
\item Recall that we have the inclusions $\tau_n(\omega_R) \subseteq \tau_{n+1}(\omega_R) \subseteq \omega_R$ 
(\ref{e1-q-tau-omega}). 
Therefore, if $R$ is $n$-quasi-$F$-rational, then $R$ is $(n+1)$-quasi-$F$-rational. 
\item 
{\cred 
Assume that $R$ is Cohen-Macaulay. 
Then $\Supp\,(\omega_R / \tau_n(\omega_R))$ coincides with the non-$n$-quasi-$F$-rational locus.} 
\end{enumerate}
\end{remark}

\begin{corollary}\label{cor: corr test sub 0^*}
We work in the general setting (Notation \ref{n-non-local}).  
Take an integer $n>0$. 
Then the following are equivalent.
\begin{enumerate}
    \item $R$ is $n$-quasi-$F$-rational. 
    \item $R_\fp$ is $n$-quasi-$F$-rational for every prime ideal $\fp$ of $R$. 
    \item $R_\m$ is $n$-quasi-$F$-rational for every maximal ideal $\m$ of $R$.
\end{enumerate} 
\end{corollary}

\begin{proof}
Assume (1).
Let $\fp$ be a prime ideal of $R$. 
Since $\tau_n(\omega_R)$ commutes with localisation (Theorem \ref{thm:test submodule locali}(1)), we get 
\[
\tau_n(\omega_{R_\fp}) = \tau_n(\omega_R) \otimes_R R_{\fp} =\omega_R \otimes_R R_{\fp} = \omega_{R_\fp}. 
\]
Therefore, $R_\fp$ is $n$-quasi-$F$-rational. 
Thus the implication $(1) \Rightarrow (2)$ holds. 
The implication $(2) \Rightarrow (3)$ is clear.

Let us prove the remaining implication $(3) \Rightarrow (1)$. 
Assume $(3)$. 
Then $R$ is Cohen-Macaulay. 
For every maximal ideal $\m$ of $R$, we have $\tau_n(\omega_R) \otimes_R R_\m=\tau_n(\omega_{R_\m})=\omega_{R_\m} = \omega_R \otimes_R R_\m$, which implies $\tau_n(\omega_R)=\omega_R$. 
Thus (1) holds. 
\end{proof}

\begin{corollary}\label{cor:qFrat to pseud-rat}
We work in the general setting (Notation \ref{n-non-local}).  
If $R$ is quasi-$F$-rational, then $\tau^q(\omega_R)=\omega_R$ 
and $R$ is pseudo-rational.
\end{corollary}

\begin{proof}
Recall that we have 
$\tau_n(\omega_R)=\tau^q(\omega_R)$ for $n \gg  0$ (Definition \ref{d-q-tau-omega}). 
Since $R$ is quasi-$F$-rational, 
we obtain $\tau^q(\omega_R) = \tau_n(\omega_R)=\omega_R$ for some $n>0$. 
By Proposition \ref{prop: rat submod}(3), $R$ is pseudo-rational.
\end{proof}

\subsection{Quasi-tight closure in local cohomologies (no boundary)}\label{ss q tight}
{In the above subsection, we defined quasi-$F$-rationality of Cohen-Macaulay ring{\cred s} by requiring that $\tau_n(\omega_R) = \omega_R$. 
{\cred We will give a characterisation of quasi-$F$-rationality 
by local cohomology (Corollary \ref{cor:equivalent-definitions-of-quasi-F-rationality}).} 
To this end, it is convenient to 
introduce the notion of quasi-tight closure. 



In this paragraph, {we work in the general setting (Notation \ref{n-non-local})}. 
Take $c \in R^{\circ}$ and integers $e \geq 0$ and $n>0$. 
We define the $W_nR$-modules $Q^e_{R, n}$ and $Q^{e, c}_{R, n}$,
and the $W_nR$-module homomorphisms $\Phi^e_{R, n}$ and $\Phi^{e, c}_{R, n}$ 
by the following diagram in which all the squares are pushouts: 
\begin{equation}\label{e-def-Q^{e,c}}
\begin{tikzcd}
W_nR \arrow[r, "F^e"] \arrow[d, "{\cccred \mathbf{R}}^{n-1}"'] & F_*^eW_nR \arrow[d] \arrow[r, "{ \cdot F_*^e[c]}"] & F_*^eW_nR \arrow[d]\\
R \arrow[r, "\Phi^e_{R, n}"] \arrow[rr, "\Phi^{e, c}_{R, n}"', bend right] & Q^e_{R, n} \arrow[r] & Q^{e, c}_{R, n}. 
\end{tikzcd}
\end{equation}

{We start by defining the quasi-tight closure.}

\begin{definition}\label{d-q-F-rat}
We work in the local setting (Notation \ref{n-local}). 
Take integers $n >0$ and $e \geq 0$. 
\begin{enumerate}
    \item 
    For $c \in R^{\circ}$, we set 
    \[
    \wt{K^{e, c}_{n}} := \Ker \left(
    H^d_\m(W_nR) \xrightarrow{F^e} H^d_\m(F^e_*W_nR) \xrightarrow{\cdot F^e_*[c]} H^d_\m(F^e_*W_nR)\right), 
    \]
which is a $W_nR$-submodule of $H^d_\m(W_nR)$. 
We define 
$\wt{0_n^{*}}$ by 
\[
\wt{0_n^{*}} := \bigcup_{\substack{c \in R^{\circ},\\ e_0 \in \Z_{>0}}} \bigcap_{e \geq e_0} \widetilde{K^{e, c}_{n}} \subseteq H^d_\m(W_nR). 
\]
Equivalently, given $z \in H^d_\m(W_nR)$, 
    we have $z \in \wt{0_n^{*}}$ if and only if 
    there exist $c \in R^{\circ}$ and $e_0 \in \Z_{>0}$ such that 
    $z$ is contained in $\wt{K^{e, c}_{n}}$ for every integer $e \geq e_0$. 
    Note that $\wt{0_n^{*}}$ is a $W_nR$-submodule of $H^d_\m(W_nR)$ 
    (Remark \ref{r-q-F-rat}). 
    \item 
   For $c \in R^{\circ}$, we set 
    \[
    K^{e, c}_{n} := \Ker\left( H^d_\m(R) \xrightarrow{\Phi^{e, c}_{R,n}} H^d_\m(Q^{e,c}_{R,n})\right), 
    \]
    which is an  $R$-submodule of $H^d_\m(R)$. 
    We define {the {\em $n$-quasi-tight closure}} $0_n^{*}$ by 
    \[
0_n^{*} := \bigcup_{\substack{c \in R^{\circ},\\ e_0 \in \Z_{>0}}} \bigcap_{e \geq e_0} K^{e, c}_{n} \subseteq H^d_\m(R). 
\]
Equivalently, given $z \in H^d_\m(R)$, 
we have $z \in 0_n^{*}$ if and only if 
there exist $c \in R^{\circ}$ and $e_0 \in \Z_{>0}$ 
such that $z$ is contained in $K^{e, c}_{n}$ for every integer $e \geq e_0$. 
Note that $0_n^{*}$ is an $R$-submodule of $H^d_\m({\cccred R})$ 
    (Remark \ref{r-q-F-rat}). 
    \item We define the  $R$-submodule $0^{q*}$ of $H^d_\m(R)$ by
    \[
    0^{q*}:= \bigcap_{n=1}^{\infty} 0^{*}_n, 
    \]
    which we call the {\em quasi-tight closure}. 
\end{enumerate}
\end{definition}

\begin{remark}\label{r-q-F-rat}
We use the same notation as in Definition \ref{d-q-F-rat}. 
Let us show that $\wt{0^*_n}$ is a $W_nR$-submodule of $H^d_{\m}(W_nR)$. 
Since it is clear that $\wt{0^*_n}$ is stable under scalar product, 
we only check that $\wt{0^*_n}$ is stable under sum. 
Take $z_1, z_2 \in \wt{0^*_n}$. 
Then there exist $c_1, c_2 \in R^{\circ}$ and $e_1, e_2 \in \Z_{>0}$ such that 
$z_1 \in \bigcap_{e \geq e_1} \wt{K^{e, c_1}_{n}}$ and 
$z_2 \in \bigcap_{e \geq e_1} \wt{K^{e, c_2}_{n}}$. 
By $ \wt{K^{e, c_1}_{n}} \subseteq \wt{K^{e, c}_{n}}$ and $ \wt{K^{e, c_2}_{n}} \subseteq \wt{K^{e, c}_{n}}$ for $c:=c_1c_2 \in R^{\circ}$, 
we obtain 
\[
z_1 + z_2 \in \bigcap_{e \geq e_0} \wt{K^{e, c}_{n}}
\]
for $e_0 := \max\{e_1, e_2\}$. 
Similarly, we see that $0^*_n$ is an $R$-submodule of $H^d_{\m}(R)$. 
\end{remark}

\begin{remark}\label{remark: easy remark}
We work in the local setting (Notation \ref{n-local}). 
\begin{enumerate}
    \item 
    Consider the commutative diagram
    \[
    \begin{tikzcd}
H^d_\m(F_*W_{n-1}R) \arrow[r, equal] \arrow[d, "{\cccred \mathbf{V}}"] & 
H^d_\m(F_*W_{n-1}R) \arrow[r, equal] \arrow[d, "{F^e{\cccred \mathbf{V}}}"] & H^d_\m(F_*W_{n-1}R) \arrow[d, "{(\cdot F^e_*[c]) \circ F^e{\cccred \mathbf{V}}}"] \\
H^d_\m(W_nR) \arrow[r, "F^e"] \arrow[d, twoheadrightarrow, "{\cccred \mathbf{R}}^{n-1}"] & 
H^d_\m(F^e_*W_nR) \arrow[r, "{\cdot F^e_*[c]}"] \arrow[d] & 
H^d_\m(F^e_*W_nR) \arrow[d] \\
H^d_\m(R) \arrow[r] \arrow[rr, "\Phi^{e, c}_{R, n}", bend right] & H^d_\m(Q^e_{X, n}) \arrow[r] & H^d_\m(Q^{e, c}_{X,n}), 
\end{tikzcd}
    \]
    in which all the vertical sequences are exact 
    {\cred (indeed, the leftmost sequence is the natural exact sequence, which induces the other two exact sequences, because    
    $Q^e_{R, n}$ and $Q^{e, c}_{R, n}$ are the pushouts,  
    and each of $F^e : W_nR \to F_*^eW_nR$ and 
    $\cdot F_*^e[c] : F_*^eW_nR \to F_*^eW_nR$ is injective)}.
    By diagram chase, we get 
    \[
    {\cccred \mathbf{R}}^{n-1}(\wt{K^{e, c}_{n}})=K^{e, c}_{n}.
    \]
    {Specifically, given $\alpha \in K^{e,c}_n$, we can find $\wt \alpha \in H^d_\m(W_nR)$ such that ${\cccred \mathbf{R}}^{n-1}(\wt \alpha) = \alpha$. Then the exactness of the rightmost vertical sequence allows us to find $\beta \in H^d_\m(F_*W_{n-1}R)$ such that $\alpha - {\cccred \mathbf{V}}\beta \in \wt{K^{e, c}_{n}}$. By construction, ${\cccred \mathbf{R}}^{n-1}(\alpha - {\cccred \mathbf{V}}\beta) = \wt \alpha$.}
    \item 
    Take $\wt{z} \in H^d_{\m}(W_nR)$ and  $z \in H^d_\m(R)$ satisfying ${\cccred \mathbf{R}}^{n-1}(\wt{z})=z$. 
    Then $z \in 0_n^{*}$ 
    if and only if there exist $c \in R^{\circ}$ and $e_0 \in \Z_{>0}$ 
    such that given an integer $e \geq e_0$, 
    we can find 
    $w_e \in H^d_\m(F_*W_{n-1}R)$ satisfying 
    \[
    ({F^e_*[c]}) \cdot F^e(\wt{z}-{\cccred \mathbf{V}}(w_e))=0. 
    \]
    \item We have 
    \[
    H^d_\m(R) \supseteq 0_1^* \supseteq  0_2^* \supseteq 0_3^* \supseteq  \cdots. 
    \]
    Indeed, this holds by the following factorisation: 
    \[
    \Phi^{e,c}_{R, n} : H^d_\m(R) \xrightarrow{\Phi^{e, c}_{R, n+1}} H^d_\m(Q^{e, c}_{R,n+1}) \to H^d_\m(Q^{e, c}_{R,n}).
    \]
    Since $H^d_\m(R)$ is an Artinian $R$-module, 
    we have $0^{q*}=0_n^{*}$ for $n \gg 0$.
    \item It holds that 
    \[
    0^*_1 = \wt{0^*_1} = 0^*_{H^d_{\m}(R)}, 
    \]
    where the first and second equalities follow from 
    Definition \ref{d-q-F-rat} and Proposition \ref{p wt0 vs tight}, respectively. 
\end{enumerate}
\end{remark}

{We recall the following standard lemma.}
\begin{lemma}\label{lem: ann local coh}
Let $(R, \m)$ be an equi-dimensional reduced  Noetherian local ring 
{\cred admitting a dualising complex} 
and let  $M$ be a finitely generated $R$-module. 
Then there exists $r \in R^{\circ}$ such that 
\[
r \cdot H^i_\m(M)=0
\]
for all $i < \dim R$.
\end{lemma}

\begin{proof}
Set $d := \dim R$. 
The assertion holds by applying \cite{stacks-project}*{\href{https://stacks.math.columbia.edu/tag/0EFC}{Tag 0EFC}}
for $s := d-1, T:= {\cccred \mathbf{V}}(\m) = \{\m\}$, and $T':= \bigcup_{f \in R^{\circ}} {\cccred \mathbf{V}}(f)$, 
which is applicable because  we have $\dim (R/\mathfrak p)_{\mathfrak q} =d > d-1$ for 
every $\mathfrak p \not\in T'$ and every $\mathfrak q \in T$. 
\end{proof}

\begin{proposition}\label{prop: tight cl l=n}
We work in the local setting (Notation \ref{n-local}).
Then the following hold. 
\begin{enumerate}
    \item ${\cccred \mathbf{R}}(\wt{0_{n+1}^{*}}) \subseteq \wt{0_n^{*}}$. 
    \item ${\cccred \mathbf{V}}^{-1}(\wt{0_n^{*}}) = F_*\wt{0_{n-1}^{*}}$.
\end{enumerate}
\end{proposition}
\noindent {Statement (1) says that the first commutative diagram below exists, whilst Statement (2) implies that the second diagram exists and is Cartiesian:
\[
\begin{tikzcd}
H^d_\m(W_{n+1}R) \ar{r}{{\cccred \mathbf{R}}} & H^d_\m(W_{n}R)\\
\wt{0_{n+1}^{*}} \ar[hook]{u} \ar{r}{{\cccred \mathbf{R}}} & \wt{0_{n}^{*}}, \ar[hook]{u}
\end{tikzcd}
\qquad
\begin{tikzcd}
H^d_\m(F_*W_{n-1}R) \ar{r}{{\cccred \mathbf{V}}} & H^d_\m(W_{n}R)\\
F_*\wt{0_{n-1}^{*}} \ar[hook]{u} \ar{r}{{\cccred \mathbf{V}}} & \wt{0_{n}^{*}}. \ar[hook]{u}
\end{tikzcd}
\]
}
\begin{proof}
{\cred The assertion (1) follows from 
\[
R(\wt{0_{n+1}^{*}}) = 
R\Big( \bigcup_{\substack{c \in R^{\circ},\\ e_0 \in \Z_{>0}}} \bigcap_{e \geq e_0} \widetilde{K^{e, c}_{n+1}}\Big) 
\subseteq 
\bigcup_{\substack{c \in R^{\circ},\\ e_0 \in \Z_{>0}}} \bigcap_{e \geq e_0} 
R(\widetilde{K^{e, c}_{n+1}})
\subseteq \bigcup_{\substack{c \in R^{\circ},\\ e_0 \in \Z_{>0}}} \bigcap_{e \geq e_0} 
\widetilde{K^{e, c}_{n}} = \wt{0_n^{*}}.
\]
}
Let us show (2). {\cred By a commutative diagram 
\[
\begin{tikzcd}[column sep=2.5cm]
H^d_\m(F_*W_{n-1}R) \arrow[r, "F^e"] \arrow[d, "{\cccred \mathbf{V}}"] & 
H^d_\m(F^{e+1}_*W_{n-1}R) \arrow[r, "{\cdot F^e_*[c]} = {\cdot F^{e+1}_*[c^p]} "]  \arrow[d, "{\cccred \mathbf{V}}"]& 
H^d_\m(F^{e+1}_*W_{n-1}R) \arrow[d, "{\cccred \mathbf{V}}"]\\
H^d_\m(W_nR) \arrow[r, "F^e"]  & 
H^d_\m(F^e_*W_nR) \arrow[r, "{\cdot F^e_*[c]}"]  & 
H^d_\m(F^e_*W_nR), 
\end{tikzcd}
\]
we get ${\cccred \mathbf{V}}(F_*\wt{K^{e, c^p}_{n-1}}) \subseteq \wt{K^{e, c}_{n}}$ 
for every $c \in R^{\circ}$, 
which implies 
\[
{\cccred \mathbf{V}}(F_*\wt{0_{n-1}^{*}}) = {\cccred \mathbf{V}}\Big( \bigcup_{\substack{c \in R^{\circ},\\ e_0 \in \Z_{>0}}} \bigcap_{e \geq e_0} F_*\widetilde{K^{e, c}_{n-1}}
\Big)  
\overset{(\star)}{=} 
 {\cccred \mathbf{V}}\Big(\bigcup_{\substack{c \in R^{\circ},\\ e_0 \in \Z_{>0}}} \bigcap_{e \geq e_0} F_*\widetilde{K^{e, c^p}_{n-1}}\Big) 
 \]
 \[
 \subseteq 
\bigcup_{\substack{c \in R^{\circ},\\ e_0 \in \Z_{>0}}} \bigcap_{e \geq e_0} {\cccred \mathbf{V}}(F_* \widetilde{K^{e, c^p}_{n-1}}) \subseteq 
\bigcup_{\substack{c \in R^{\circ},\\ e_0 \in \Z_{>0}}} \bigcap_{e \geq e_0} \wt{K^{e, c}_{n}} = \wt{0_n^{*}}, 
\]
where ($\star$) follows from 
$\widetilde{K^{e, c}_{n-1}} \subseteq \widetilde{K^{e, c^p}_{n-1}}$. 
Hence it holds that} $F_*\wt{0_{n-1}^{*}} \subseteq {\cccred \mathbf{V}}^{-1}(\wt{0_n^{*}})$.
In what follows, we prove the opposite inclusion 
\[
{\cccred \mathbf{V}}^{-1}(\wt{0_n^{*}}) \subseteq F_*\wt{0_{n-1}^{*}}.
\]
Take $F_*z \in {\cccred \mathbf{V}}^{-1}(\wt{0_n^{*}})$, that is, ${\cccred \mathbf{V}}(F_*z) \in \wt{0_n^{*}}$. 
Then there exist $c \in R^{\circ}$ and $e_0>0$ such that 
${\cccred \mathbf{V}}(F_*z) \in \bigcap_{e \geq e_0} \wt{K^{e, c}_{n}}$. 
Moreover, recall that we can pick $c' \in R^{\circ}$ satisfying $c' \cdot H^{d-1}_\m(R)=0$ (Lemma \ref{lem: ann local coh}), which implies that 
\[
(F_*^e c')  \cdot H^{d-1}_\m(F_*^eR)=(F_*^e c')  \cdot F_*^eH^{d-1}_\m(R)=F_*^e(c' \cdot H^{d-1}_\m(R))=0
\]
for every $e \geq 0$.  
The exact sequence 
\[
H^{d-1}_{\m}(R) \to   H^d_{\m}(F_*W_{n-1}R) \xrightarrow{{\cccred \mathbf{V}}}  H^d_{\m}(W_nR) 
\]
induces  the following commutative diagram in which each horizontal sequence is exact: 
\[
\begin{tikzcd}
H^{d-1}_{\m}(R) \arrow[r] \arrow[d, "F^e"] & 
{\cccred \mathbf{V}}^{-1}(\wt{0_n^{*}}) \arrow[r, "{\cccred \mathbf{V}}"]  \arrow[d, "F^e"]& 
\wt{0^*_n} \arrow[d, "F^e"]\\
H^{d-1}_{\m}(F_*^eR) \arrow[r] \arrow[d, "{\cdot F^e_*[c]}"] & 
F_*^e{\cccred \mathbf{V}}^{-1}(\wt{0_n^{*}}) \arrow[r, "{\cccred \mathbf{V}}"]  \arrow[d, "{\cdot F^e_*[c]}"]& 
F_*^e\wt{0^*_n} \arrow[d, "{\cdot F^e_*[c]}"]\\
H^{d-1}_{\m}(F_*^eR) \arrow[r] \arrow[d, "{\cdot F^e_*[c']=0}"] & 
F_*^e{\cccred \mathbf{V}}^{-1}(\wt{0_n^{*}}) \arrow[r, "{\cccred \mathbf{V}}"]  \arrow[d, "{\cdot F^e_*[c']}"]& 
F_*^e\wt{0^*_n} \arrow[d, "{\cdot F^e_*[c']}"]\\
H^{d-1}_{\m}(F_*^eR) \arrow[r] & 
F_*^e{\cccred \mathbf{V}}^{-1}(\wt{0_n^{*}}) \arrow[r, "{\cccred \mathbf{V}}"] & 
F_*^e\wt{0^*_n}. 
\end{tikzcd}
\]
By the diagram chase starting with $F_*z \in {\cccred \mathbf{V}}^{-1}(\wt{0_n^{*}})$ (located in the centre at the top), we get $(F_*^e[cc']) \circ F^e (F_*z)=0$ for every $e \geq e_0$, and hence 
$F_*z \in F_*\wt{0^*_{n-1}}$. 
Thus (2) holds. 
\end{proof}

{The following result is the first step towards proving that the Teichm\"uller lift of a fixed power of a test element is a {quasi-}test element.} 
\begin{proposition}\label{prop: ker V and kerF}
We work in the local setting (Notation \ref{n-local}).
Take $t \in R^{\circ}$ satisfying  $t \cdot 0^*_{H^d_\m(R)}=0$. 
Then, for every $n>0$ and every $m\geq 0$,   the following hold. 
\begin{enumerate}
\item $[t^2] \cdot \widetilde{0_n^*}=0$. 
\item $[t^2] \cdot \mathrm{Ker}(F^m \colon H^d_\m(W_nR) \to H^d_\m(F^m_*W_nR))=0$. 
\end{enumerate}
\end{proposition}

\begin{proof}
We have the following inclusion (Definition \ref{d-q-F-rat}(1)): 
\[
 \mathrm{Ker}(F^m \colon H^d_\m(W_nR) \to H^d_\m(F^m_*W_nR))\subseteq \widetilde{0_n^*}. 
\]
Hence (1) implies (2).

Thus it is enough to show (1). Specifically, 
we will show that 
$[t^2] \cdot \wt{0_n^{*}}=0$ by an increasing induction on $n$.
The base case $n=1$ holds  by 
\[
[t] \cdot \wt{0_1^{*}}=[t] \cdot 0_1^{*}=t \cdot 0^*_{H^d_\m(R)}=0
\]
(Remark \ref{remark: easy remark}(4)). 
Now assume $n \geq 2$ and consider the exact sequence
\begin{equation}\label{e2: ker V and kerF}
H^d_\m(F_*W_{n-1}R) \xrightarrow{{\cccred \mathbf{V}}} H^d_\m(W_nR) \xrightarrow{{\cccred \mathbf{R}}^{n-1}} H^d_\m(R).
\end{equation}
Fix $z \in \wt{0_n^{*}}$. 
Since ${\cccred \mathbf{R}}^{n-1}(z) \in \wt{0_1^{*}}$ (Proposition \ref{prop: tight cl l=n}(1)), 
we get that ${\cccred \mathbf{R}}^{n-1}([t]z) = [t]{\cccred \mathbf{R}}^{n-1}(z) =0$.
Thus the exact sequence (\ref{e2: ker V and kerF}) enables us to find 
\[
F_*z' \in F_*H^d_\m(W_{n-1}R)= H^d_\m(F_*W_{n-1}R)
\]
satisfying ${\cccred \mathbf{V}}(F_*z')=[t]z$, 
where $F_*z' \in F_*H^d_\m(W_{n-1}R)$ denotes the same element as 
$z' \in H^d_\m(W_{n-1}R)$ (in order to clarify the $R$-module structure see  Subsubsection \ref{sss-F_*M}). 
As $F_*z' \in {\cccred \mathbf{V}}^{-1}(\wt{0^*_n}) = F_*\wt{0^*_{n-1}}$ (Proposition \ref{prop: tight cl l=n}(2)), 
we get $z' \in \wt{0^*_{n-1}}$, and so, by the induction hypothesis, $[t^2]z' = 0$. 
We thus obtain 
\[
[t^2]z=[t]{\cccred \mathbf{V}}(F_*z')={\cccred \mathbf{V}}([t]F_*z')={\cccred \mathbf{V}}(F_*([t^p]z'))=0. 
\]
Hence (1) holds.
\end{proof}

\begin{definition}\label{definition: small F-stable}
We work in the local setting (Notation \ref{n-local}). 
Fix $n \in \Z_{>0}$. 
Let $S$ be a $W_nR$-submodule  of $H^d_\m(W_nR)$.
\begin{enumerate}
    \item We say that $S$ is {\em $F$-stable} if $F(S) \subseteq F_*S$, 
    where $F_*S \subseteq H^d_\m(F_*W_nR) = F_*(H^d_\m(W_nR))$ denotes the same subset as $S \subseteq H^d_\m(W_nR)$ via the identification 
    $H^d_\m(F_*W_nR) = F_*H^d_\m(W_nR) =H^d_{\m}(W_nR)$. 
    \item We say that $S$ is {\em small} if there exists $c \in R^{\circ}$ such that $[c] \cdot S=0$. 
\end{enumerate}
\end{definition}

\noindent 
{\cred 
Opposed to the case when $n=1$, 
there exists a proper $W_nR$-submodule of $H^d_{\m}(W_nR)$ which is not small (e.g., if  $R := \F_p$ and $n=2$, then 
$W_nR = \Z/p^2\Z$ and $H^d_{\m}(W_nR) = H^0_{\m}(\Z/p^2\Z) = \Z/p^2\Z$, and hence $S := p\Z/ p^2\Z \subsetneq \Z/p^2\Z = H^d_{\m}(W_nR)$ is 
a proper $W_2R$-submodule of $H^d_{\m}(W_nR)$ which is not small). 
}

\begin{proposition}\label{prop: tight cl l=n2}
We work in the local setting (Notation \ref{n-local}).  
Then $\wt{0_n^{*}}$ is the largest $W_nR$-submodule of $H^d_\m(W_nR)$ 
that is small and $F$-stable. 
\end{proposition}

\begin{proof}
By Definition \ref{d-q-F-rat}, $\wt{0_n^{*}}$ is $F$-stable, and  
it follows from 
Remark \ref{r-test-ideal-summary}(1) and 
Proposition \ref{prop: ker V and kerF}(1) that $\wt{0^*_n}$ is small. 
{To show that $\wt{0_n^{*}}$ is the largest such submodule,} pick any small ${\cred F}$-stable $W_nR$-submodule $M$ of $H^d_\m(W_nR)$ and take $z \in M$. 
It is enough to show that $z \in \wt{0_n^{*}}$. 
Since $M$ is small, there exists $c \in R^{\circ}$ such that $[c] \cdot M=0$.
For every $e>0$,  we get 
\[
(F_*^e[c]) \cdot F^e(z) \in (F_*^e[c]) \cdot F^e(M) 
\subseteq (F_*^e[c]) \cdot F_*^eM = F_*^e([c] \cdot M) =0,
\]
which implies  $z \in \wt{0^{*}_n}$. 
\qedhere
\end{proof}

{The following proposition is fundamental as it shows that the infinite sum and intersection in the definition of $\wt{0^{*}_n}$ stabilise for any choice of $c \in (t^4) \cap R^\circ$, where $t$ is a test element in $R^{\circ}$.}

\begin{proposition}\label{prop: test elem tilda}
We work in the local setting (Notation \ref{n-local}).  
Take $t \in R^{\circ}$ satisfying $t \cdot 0^*_{H^d_\m(R)} =0$. 
Fix $n \in \Z_{>0}$ and $c \in (t^4) \cap R^{\circ}$.
Then the following hold. 
\begin{enumerate}
    \item $H^d_\m(W_nR) \supseteq \wt{K^{0, c}_{n}} \supseteq \wt{K^{1, c}_{n}} \supseteq 
    \cdots \supseteq \wt{K^{e, c}_{n}} \supseteq \wt{K^{e+1, c}_{n}} \supseteq   \cdots$\\[-1em]
    \item There exists an integer $e_1>0$ such that 
    \[
    \wt{0^{*}_n}=\wt{K^{e, c}_{n}}
    \]
    for every integer $e \geq e_1$.\\[-1em] 
        \item $H^d_\m(R) \supseteq K^{0, c}_{n} \supseteq K^{1, c}_{n} \supseteq 
    \cdots \supseteq K^{e, c}_{n} \supseteq K^{e+1, c}_{n} \supseteq   \cdots.$
\end{enumerate}
\end{proposition}

\begin{proof}
Let us show (1). 
There exists $c' \in R^{\circ}$ such that $c=c't^4$.
Fix $z \in \wt{K^{e+1, c}_{n}}$. We then  have
\[
F((F_*^e[c't^2]) \cdot F^e(z))=
(F_*^{e+1}[c'^pt^{2p}]) \cdot (F^{e+1}(z))=
(F_*^{e+1}[c_1]) \cdot (F_*^{e+1}[c]) \cdot (F^{e+1}(z)) =0
\]
for $c_1 := c'^{p-1}t^{2p-4} \in R^{\circ}$, where the last equality is guaranteed by 
$z \in \wt{K^{e+1, c}_{n}}$. 
By Proposition \ref{prop: ker V and kerF}(2), we get 
\[
0 = (F_*^e[t^2]) \cdot (F_*^e[c't^2]) \cdot F^e(z) = (F_*^e [c't^4]) \cdot F^e(z) = 
 (F_*^e [c]) \cdot F^e(z), 
\]
which implies $z \in \wt{K^{e, c}_{n}}$. 
Thus (1) holds. {The proof {of (1)} may be visualised by the following diagram:
\[
\begin{tikzcd}[column sep = large]
 W_nR \ar{d}{=} \ar{r}{F^e} & F^e_*W_nR \ar{d}{=} \ar{d}{=} \ar[bend left = 25]{rr}{\cdot F^{e+1}_*[c_1c] \circ F} \ar{r}{\cdot F^e_*[c't^2]} & F^e_*W_nR \ar{d}{\cdot F^e_*[t^2]} \ar{r}{F} & F^{e+1}_*W_nR \\
 W_nR \ar{r}{F^e} & F^e_*W_nR \ar{r}{\cdot F_*^e[c]} & F^{e}_*W_nR.
\end{tikzcd}
\]}


Let us show (2). 
Since $H^d_\m(W_nR)$ is an Artinian $W_nR$-module, 
(1) enables us to find $e_1>0$ such that 
\[
\wt{K^{e, c}_{n}}=\wt{K^{e_1, c}_{n}}
\]
for 
every $e \geq e_1$. 
Then Definition \ref{d-q-F-rat}(1) implies  
\[
\wt{0_n^{*}} = \bigcup_{\substack{c \in R^{\circ},\\ e_0 \in \Z_{>0}}} \bigcap_{e \geq e_0} \widetilde{K^{e, c}_{n}} 
\supseteq \bigcap_{e \geq e_1} \widetilde{K^{e, c}_{n}} 
=  \widetilde{K^{e_1, c}_{n}}. 
\]
Conversely, if $z \in \wt{0^*_n}$ and $e>0$, then we get 
\begin{align*}
F_*^e[c] \cdot F^e(z) 
\in F_*^e[c] \cdot F^e(\wt{0^*_n}) &\overset{{\rm (i)}}{\subseteq} F_*^e[c] \cdot F_*^e\wt{0^{*}_n} \\
&=F_*^e(  [c] \cdot \wt{0^{*}_n} )= F_*^e( [c't^2] \cdot [t^2] \cdot \wt{0^{*}_n})
\overset{{\rm (ii)}}{=} 0,
\end{align*}
where (i) and (ii) follow from Proposition \ref{prop: tight cl l=n2} and 
Proposition \ref{prop: ker V and kerF}(1), respectively. 
Hence $z \in \bigcap_{e >0} \widetilde{K^{e, c}_{n}} = \widetilde{K^{e_1, c}_{n}}$. Thus (2) holds. 

The assertion (3) follows from (1) and ${\cccred \mathbf{R}}^{n-1}( \widetilde{K^{e, c}_{n}}) = K^{e, c}_{n}$ (Remark \ref{remark: easy remark}(1)). 
\end{proof}


\begin{theorem}\label{thm: lift of 0^*,n}
We work in the local setting (Notation \ref{n-local}).  
Fix $n \in \Z_{>0}$. 
Then  
\[
{\cccred \mathbf{R}}^{n-1}(\wt{0^{*}_n})=0^{*}_n.
\]
\end{theorem}

\begin{proof}
{\cred It holds that} 
\[
{\cred 
{\cccred \mathbf{R}}^{n-1}(\wt{0^{*}_n}) = 
{\cccred \mathbf{R}}^{n-1}\Big( 
 \bigcup_{\substack{c \in R^{\circ},\\ e_0 \in \Z_{>0}}} \bigcap_{e \geq e_0} \wt{K^{e_0, c}_{n}} 
\Big) \subseteq 
 \bigcup_{\substack{c \in R^{\circ},\\ e_0 \in \Z_{>0}}} \bigcap_{e \geq e_0} {\cccred \mathbf{R}}^{n-1}(\wt{K^{e_0, c}_{n}})  
\overset{(\star)}{=} 
\bigcup_{\substack{c \in R^{\circ},\\ e_0 \in \Z_{>0}}} \bigcap_{e \geq e_0}  K^{e_0, c}_{n} = 0^{*}_n,}
\]
where the  equality $(\star)$ follows from  Remark \ref{remark: easy remark}(1).

Then it is enough to show ${\cccred \mathbf{R}}^{n-1}(\wt{0^{*}_n}) \supseteq 0^{*}_n$. 
Take $z \in 0^*_n$ and $t \in R^{\circ}$ satisfying $t \cdot 0^*_{H^d_\m(R)}=0$.
There exist  $c \in R^{\circ}$ and $e_1>0$ 
such that $z \in \bigcap_{e \geq e_1} K^{e, c}_{n}$. 
Replacing $c$ by $ct^4$ (which is allowed by $K^{e, c}_{n} \subseteq K^{e, ct^4}_{n}$), 
we may assume $\wt{K^{e, c}_{n}}=\wt{0^*_n}$ for $e \gg 0$ (Proposition \ref{prop: test elem tilda}(2)).
For $e \gg 0$, we obtain 
\[
z \in \bigcap_{e \geq e_1} K^{e, c}_{n} \subseteq K^{e, c}_{n} 
\overset{(\star\star)}{=}{\cccred \mathbf{R}}^{n-1}(\wt{K^{e, c}_{n}})
={\cccred \mathbf{R}}^{n-1}(\wt{0^*_n}), 
\]
where 
the  equality $(\star\star)$ follows from  Remark \ref{remark: easy remark}(1). 
\qedhere 


\end{proof}

Below, (1) shows that the infinite sum and intersection in the definition of quasi-tight closure $0^*_n$ stabilise, whilst (2) shows that the triviality of $K^{e_2,t^3}_n$, as opposed to just $K^{e_2,t^4}_n$, is already enough to deduce that $0^*_n=0$. {Although not needed in 
this paper, 
this stronger result will be used in future work.} 


\begin{proposition}\label{prop: test element non-tilda}
We work in the local setting (Notation \ref{n-local}). 
Take $t \in R^{\circ}$ satisfying  $t \cdot 0^{*}_{H^d_\m(R)}=0$. 
Fix $n \in \Z_{>0}$ and $c \in (t^4) \cap R^{\circ}$.
Then the following hold. 
\begin{enumerate}
    \item There exists an integer $e_1>0$ such that 
    $0_n^{*}=K^{e, c}_{n}$ for every integer $e \geq e_1$. 
    \item If there is an integer $e_2>0$ satisfying $K^{e_2, t^3}_{n}=0$, then $0_n^{*}=0$.
\end{enumerate}
\end{proposition}

\begin{proof}
The assertion (1) follows from ${\cccred \mathbf{R}}^{n-1}(\wt{K^{e, c}_{n}})=K^{e, c}_{n}$ (Remark \ref{remark: easy remark}(1)), 
Proposition \ref{prop: test elem tilda}(2), and Theorem \ref{thm: lift of 0^*,n}. 

Let us show (2). 
Take $z \in 0^{*}_n$.
Then there exists $\wt{z} \in \wt{0^{*}_n}$ such that ${\cccred \mathbf{R}}^{n-1}(\wt{z})=z$ 
(Theorem \ref{thm: lift of 0^*,n}). 
We can find an integer $e$ such that  $e \geq e_2+2$ and $\wt{z} \in \wt{K^{e, t^4}_{n}}$,  i.e., 
$(F_*^e[t^4]) \cdot F^e(\wt{z})=0$ 
(Proposition \ref{prop: test elem tilda}(2)).
Therefore, 
\[
F^{e-e_2}(F_*^{e_2}[t] \cdot F^{e_2}(\wt{z}))=(F_*^e[t]^{p^{e-e_2}}) \cdot F^e(\wt{z})=0,
\]
which implies $(F_*^{e_2}[t^3]) \cdot F^{e_2}(\wt{z})=
(F_*^{e_2}[t^2]) \cdot (F_*^{e_2}[t])\cdot F^{e_2}(\wt{z})=
0$ (Proposition \ref{prop: ker V and kerF}(2)). 
In particular, $z \in K^{e_2, t^3}_{n}=0$, and hence $z=0$.
\end{proof}


\begin{proposition}\label{prop: chara F-rat}
We work in the local setting (Notation \ref{n-local}). 
Assume that 
\[
{\cccred \mathbf{V}} \colon H^d_\m(F_*W_mR) \to H^d_\m(W_{m+1}R)
\]
is injective for every integer $m>0$. 
Then $0^*_{H^d_\m(R)}=0$ if and only if $\wt{0^{*}_n}=0$. 
\end{proposition}

\begin{proof}
By $0^*_{H^d_\m(R)} =0^*_1 = \wt{0^{*}_1}$ (Remark \ref{remark: easy remark}(4)), 
the only-if part is clear. 
Assume that  $\wt{0^{*}_n}=0$ for some  $n>0$. 
By Proposition \ref{prop: tight cl l=n}(2), we get 
\[
{\cccred \mathbf{V}}^{n-1}(F_*^{n-1}0^*_{H^d_\m(R)})=
{\cccred \mathbf{V}}^{n-1}(F_*^{n-1}\wt{0^*_1}) 
\subseteq \wt{0^{*}_n}=0,
\]
which implies  $0^*_{H^d_\m(R)}=0$ by the injectivity 
of ${\cccred \mathbf{V}} \colon H^d_\m(F_*W_mR) \to H^d_\m(W_{m+1}R)$. 
\end{proof}

{\cred 
\noindent 
By the same argument as above, 
we can prove that $0^*_{H^d_\m(R)} =0$ 
if $\wt{0^{*}_n}=0$ and ${\cccred \mathbf{V}}^{n-1} \colon H_{\m}^d(F^{n-1}R) \to H^d_{\m}(W_nR)$ is injective.}

\subsection{
{Duality between quasi-test submodules and quasi-tight closure}}
We start by relating quasi-test submodules with quasi-tight closure.
\begin{proposition}\label{prop: corr 0^* and tau: no pair}
We work in the local setting (Notation \ref{n-local}).    
Take an integer $n>0$. 
Then {the following equality} of $W_nR$-submodules of
$(W_n\omega_R)^{\vee} = H^d_{\m}(W_nR)$ holds:
\[
 \left(\frac{W_n\omega_R}{\tau(W_n\omega_R)}\right)^\vee = \wt{0^{*}_n}.
\]
{In particular, 
\[
\tau(W_n\omega_R)^{{\cred \vee}} = 
\frac{H^d_\m(W_nR)}{\wt{0^{*}_n}}.
\]
}
\end{proposition}

\begin{proof}
Fix a test element $t \in R^{\circ}$. 
By Proposition \ref{prop: test submod}(4), 
there exists $e_0 \geq 0$ such that 
$\tau(W_n\omega_R)$ is the image of
\[
\PsiT^e_n\circ (\cdot F^e_*[t^4]) \colon  
F^e_*W_n\omega_R \to W_n\omega_R.
\]
We now apply the Matlis duality functor $\Hom_{W_nR}(-, E)$.  
Since $(-)^{\vee} := \Hom_{W_nR}(-, E)$ is an exact functor (recall that $E$ is an injective $W_nR$-module), 
the {diagrams} 
\[
T_n^{e, t^4} : F^e_*W_n\omega_R \twoheadrightarrow \tau(W_n\omega_R) \hookrightarrow W_n\omega_R 
\]
\[
0 \to \tau(W_n\omega_R) \to W_n\omega_R \to \frac{W_n\omega_R}{\tau(W_n\omega_R)} \to 0
\]
induce the following {diagrams}:  
\[
H^d_{\m}(F^e_*W_nR)\hookleftarrow  \tau(W_n\omega_R)^{\vee} \twoheadleftarrow  H^d_{\m}(W_nR) 
\]
\begin{equation} \label{eq:ses-duality-Wnomega}
0 \leftarrow \tau(W_n\omega_R)^{\vee} \leftarrow H^d_{\m}(W_nR)\leftarrow\left(\frac{W_n\omega_R}{\tau(W_n\omega_R)}\right)^\vee \leftarrow 0, 
\end{equation}
{where the latter is exact.} 
For $e \gg 0$, we get 
\begin{eqnarray*}
    \left(\frac{W_n\omega_R}{\tau(W_n\omega_R)}\right)^\vee 
&=& \Ker( H^d_{\m}(W_nR) \to \tau(W_n\omega_R)^{\vee})\\[-0.5em]
&=& \Ker((\PsiT^e_n\circ (\cdot F^e_*[t^4]))^{\vee} : H^d_{\m}(W_nR) \to  H^d_{\m}(F^e_*W_nR))\\[0.3em]
&=& 
\Ker( (\cdot F_*^e[t^4]) \circ F^e: H^d_{\m}(W_nR) \to  H^d_{\m}(F^e_*W_nR))\\[0.2em]
&\overset{{\rm (i)}}{=}& \bigcap_{e \geq 0} \wt{K^{e, t^4}_{n}}\\[-0.1em]
&\overset{{\rm (ii)}}{=}& \wt{0^*_n}, 
\end{eqnarray*}
where (i) and (ii) follow from 
Proposition \ref{prop: test elem tilda}. {The above equality and (\ref{eq:ses-duality-Wnomega}) immediately yield the in-particular part.}
\qedhere 



\end{proof}

\begin{theorem}\label{thm: corr 0^* and tau}
We work in the local setting (Notation \ref{n-local}).   
Take an integer $n>0$. 
Then {the following equality} of $R$-submodules of
$(\omega_R)^{\vee} = H^d_{\m}(R)$ holds:
\[
\left(\frac{\omega_R}{\tau_n(\omega_R)}\right)^\vee =  0^{*}_n.
\]
{In particular,
\[
\tau_n(\omega_R)^{{\cred \vee}} = \frac{H^d_\m(R)}{{0^{*}_n}}.
\]
}
\end{theorem}

\begin{proof}
Since $\tau_n(\omega_R)$ is the inverse image of $\tau(W_n\omega_R)$ 
by $({\cccred \mathbf{R}}^{n-1})^* : \omega_R \hookrightarrow W_n\omega_R$ 
(Definition \ref{d-q-tau-omega}), 
we obtain the following commutative diagram:  
\[
\begin{tikzcd}
\omega_R \arrow[r, "({\cccred \mathbf{R}}^{n-1})^*", hook] \arrow[d, twoheadrightarrow] & W_n\omega_R\arrow[d, twoheadrightarrow]\\
\frac{\omega_R}{\tau_n(\omega_R)} \arrow[r, hook] &\frac{W_n\omega_R}{\tau(W_n\omega_R)}.
\end{tikzcd}
\]
Applying the Matlis dual functor $(-)^{\vee} := \Hom_{W_nR}(-, E)$ (which is exact) to this diagram, we get 
\[
\begin{tikzcd}
H^d_{\m}(R) \arrow[r, "{\cccred \mathbf{R}}^{n-1}", twoheadleftarrow] \arrow[d, hookleftarrow] & H^d_{\m}(W_nR) \arrow[d, hookleftarrow]\\
(\frac{\omega_R}{\tau_n(\omega_R)})^{\vee} \arrow[r, twoheadleftarrow] & (\frac{W_n\omega_R}{\tau(W_n\omega_R)})^{\vee}.
\end{tikzcd}
\]
By $(\frac{W_n\omega_R}{\tau(W_n\omega_R)})^{\vee} = \wt{0_n^*}$ (Proposition \ref{prop: corr 0^* and tau: no pair}) 
and ${\cccred \mathbf{R}}^{n-1}(\wt{0_n^*}) = 0^*_n$ (Theorem \ref{thm: lift of 0^*,n}), 
we obtain $(\frac{\omega_R}{\tau_n(\omega_R)})^{\vee} = 0^*_n$. 
{The in-particular part follows by Matlis duality applied to the short exact sequence:
\[
0 \to \tau_n(\omega_R) \to \omega_R \to \frac{\omega_R}{\tau_n(\omega_R)} \to 0. \qedhere
\] }


\end{proof}

Next, we provide an explicit definition of $\tau_n(\omega_R)$ in terms of $Q^{e,c}_{R,n}$. {First, we need the following remark.}
\begin{remark} \label{remark:inequalities-for-quasi-f-module-with-test-element}
Note that we have the following inclusions:
\begin{multline*}
{\rm Im}\Big(\Hom_{W_nR}(Q^{e,c}_{R,n}, W_n\omega_R) \xrightarrow{(\Phi^{e, c}_{R,n})^*} \omega_R \Big) \\
\subseteq {\rm Im}\Big(\Hom_{W_nR}(Q^{e',c}_{R,n}, W_n\omega_R) \xrightarrow{(\Phi^{e', c}_{R,n})^*} \omega_R \Big)
\end{multline*}
for every $e' \geq e$ if 
{$t \in R^{\circ}$ is a test element} and 
$c \in (t^4) \cap R^\circ$, 
{\cccred where $(\Phi^{e, c}_{R,n})^*$ (and similarly $(\Phi^{e', c}_{R,n})^*$) denotes the map induced by the dual of $\Phi^{e, c}_{R,n}$:
\[
\Hom_{W_nR}(Q^{e,c}_{R,n}, W_n\omega_R) \xrightarrow{(\Phi^{e, c}_{R,n})^*} \Hom_{W_n R} (R, W_n \omega_{R}) \simeq \omega_R.
\]
}
Indeed, this can be checked after localisation, in which case it follows by 
Proposition \ref{prop: test elem tilda} and the fact that  
\[
{\rm Im}\Big(\Hom_{W_nR}(Q^{e,c}_{R,n}, W_n\omega_R) \xrightarrow{(\Phi^{e, c}_{R,n})^*} \omega_R \Big)
\]
is Matlis dual to 
${\rm Im}\Big( H^d_\m(R) \xrightarrow{\Phi^{e, c}_{R,n}} H^d_\m(Q^{e,c}_{R,n})\Big) = \frac{H^d_\m(R)}{{K^{e,c}_n}}$. 
{Moreover,} since $\omega_R$ is a finitely gerenated $R$-module, 
there exists an integer $e_1 >0$ such that 
\begin{multline*}
{\rm Im}\Big(\Hom_{W_nR}(Q^{e,c}_{R,n}, W_n\omega_R) \xrightarrow{(\Phi^{e, c}_{R,n})^*} \omega_R \Big) \\
= {\rm Im}\Big(\Hom_{W_nR}(Q^{e_1,c}_{R,n}, W_n\omega_R) \xrightarrow{(\Phi^{e_1, c}_{R,n})^*} \omega_R \Big)
\end{multline*}
for every integer $e \geq e_1$. 
\end{remark}


\begin{corollary} \label{cor:Qdefs-of-quasi-tight-submodule}
We work in the general setting (Notation \ref{n-non-local}). 
Then the following 
hold.
{\setlength{\leftmargini}{2em} \begin{enumerate}
    \item $\tau_n(\omega_R) = \bigcap_{c \in R^\circ} \bigcap_{e_0>0} \sum_{e \geq e_0} {\rm Im}\Big(\Hom_{W_nR}(Q^{e,c}_{R,n}, W_n\omega_R) \xrightarrow{(\Phi^{e, c}_{R,n})^*} \omega_R \Big)$.\\[-0.5em]
    \item $\tau_n(\omega_R) = \bigcap_{c \in R^\circ} \sum_{e >0 } {\rm Im}\Big(\Hom_{W_nR}(Q^{e,c}_{R,n}, W_n\omega_R) \xrightarrow{(\Phi^{e, c}_{R,n})^*} \omega_R \Big)$.\\[-0.5em]
    \item 
    {Fix a test element $t \in R^{\circ}$ and $c \in (t^4) \cap R^{\circ}$. 
    Then there exists an integer $e_0 >0$ such that 
    \[
    \tau_n(\omega_R) = {\rm Im}\Big(\Hom_{W_nR}(Q^{e,c}_{R,n}, W_n\omega_R) \xrightarrow{(\Phi^{e, c}_{R,n})^*} \omega_R \Big)
    \]
    for every integer $e \geq e_0$.} 
\end{enumerate}}
\end{corollary}


\begin{proof}
First, we observe that (3) immediately implies (1). Indeed,
\begin{align*}
\bigcap_{c \in R^\circ} \bigcap_{e_0>0} \sum_{e \geq e_0} & {\rm Im}\Big(\Hom_{W_nR}(Q^{e,c}_{R,n}, W_n\omega_R) \xrightarrow{(\Phi^{e, c}_{R,n})^*} \omega_R \Big) \\ &= \bigcap_{c \in t^4R^\circ} \bigcap_{e_0>0} \sum_{e \geq e_0} {\rm Im}\Big(\Hom_{W_nR}(Q^{e,c}_{R,n}, W_n\omega_R) \xrightarrow{(\Phi^{e, c}_{R,n})^*} \omega_R \Big)  \\
&\overset{\mathclap{\ref{remark:inequalities-for-quasi-f-module-with-test-element}}}{=} \bigcap_{c \in t^4R^\circ} \sum_{e >0 }  {\rm Im}\Big(\Hom_{W_nR}(Q^{e,c}_{R,n}, W_n\omega_R) \xrightarrow{(\Phi^{e, c}_{R,n})^*} \omega_R \Big) \\
&\overset{(3)}{=} \tau_n(\omega_R).
\end{align*}
The same argument shows that (3) implies (2).

In what follows, we fix a test element $t \in R^\circ$, take any $c \in (t^4) \cap R^\circ$, and aim for showing (3). By Remark \ref{remark:inequalities-for-quasi-f-module-with-test-element} and Noetherianity, it is enough to prove that
\[
\tau_n(\omega_R) = 
{\rm Im}\Big(\Hom_{W_nR}(Q^{e,c}_{R,n}, W_n\omega_R) \xrightarrow{(\Phi^{e, c}_{R,n})^*} \omega_R \Big)
\]
{\cred for every $e \gg 0$.} 
In order to verify this, we may assume that $R$ is a {\cred complete} 
local ring by Theorem \ref{thm:test submodule locali} 
(cf.\ 
Lemma \ref{l-W_nR-local-cplt}).
In this case,
\[
{\rm Im}\Big(\Hom_{W_nR}(Q^{e,c}_{R,n}, W_n\omega_R) \xrightarrow{(\Phi^{e, c}_{R,n})^*} \omega_R \Big)
\]
is Matlis dual to
\begin{equation} \label{eq:taun-dual-new-def}
{\rm Im}\Big( H^d_\m(R) \xrightarrow{\Phi^{e, c}_{R,n}} H^d_\m(Q^{e,c}_{R,n})\Big) = \frac{H^d_\m(R)}{{K^{e,c}_n}}.
\end{equation}

Now, Theorem \ref{thm: corr 0^* and tau} and Proposition \ref{prop: test element non-tilda}(1) imply
\begin{equation} \label{eq:hidden-equation-kecn2}
\tau_n(\omega_R) = \Big(\frac{H^d_\m(R)}{{K^{e,c}_n}}\Big)^\vee.
\end{equation}
for every $e \gg 0$ (depending on $c \in (t^4) \cap R^\circ$). Thus (3) follows by (\ref{eq:taun-dual-new-def}) {and (\ref{eq:hidden-equation-kecn2})}.
\qedhere 

\end{proof}

\begin{corollary} \label{cor:equivalent-definitions-of-quasi-F-rationality}
We work in the general setting (Notation \ref{n-non-local}). 
Fix an integer $n >0$. 
Then $R$ is $n$-quasi-$F$-rational if and only if $R$ is Cohen-Macaulay and one of the following equivalent conditions hold.
\begin{enumerate}
\item $\tau_n(\omega_R) = \omega_R$.
\item For every maximal ideal $\m$ of $R$ and for every $c \in R^\circ$, the following map
\[
 \Phi^{e, c}_{R_\m,n} \colon H^d_\m(R_\m) \to H^d_\m(Q^{e,c}_{R_\m,n})
\]
is injective $\begin{aligned}[t] &\text{{\rm (a):} for some } e > 0, \text{ or equivalently} \\ &\text{{\rm (b):} for every } e\gg 0 \text{ (depending on $c$)}. \end{aligned}$\\[-0.3em]
\item For every $c \in R^\circ$, the following map
\[
\Hom_{W_nR}(\Phi^{e, c}_{R,n}, W_n\omega_R) \colon \Hom_{W_nR}(Q^{e,c}_{R,n}, W_n\omega_R) \to \omega_R 
\]
is surjective $\begin{aligned}[t] &\text{{\rm (a):} for some } e > 0, \text{ or equivalently} \\ &\text{{\rm (b):} for every } e\gg 0 \text{ (depending on $c$)}. \end{aligned}$
\end{enumerate}
\end{corollary}

\begin{proof}
Note that (1) together with Cohen-Macauliness of $R$ is exactly our definition of $n$-quasi-$F$-rationality. By coherence of $Q^{e,c}_{R,n}$ and 
Theorem \ref{thm:test submodule locali}, we may assume that $R$ is a {\cred complete} local ring with maximal ideal $\m$. 

{Then, by Matlis duality, we get the equivalences (2a) $\Leftrightarrow$ (3a) and 
{\cred (2b)} $\Leftrightarrow$ (3b).}
The equivalence of (1) and (3a) follows from Corollary \ref{cor:Qdefs-of-quasi-tight-submodule}(2), whilst the equivalence of (1) and  (3b) follows from Corollary \ref{cor:Qdefs-of-quasi-tight-submodule}(3). 
{\cred Here Corollary \ref{cor:Qdefs-of-quasi-tight-submodule} is applicable, because  we may assume $c \in (t^4) \cap R^{\circ}$ for a test element $t \in R^{\circ}$ in order to prove these equivalences 
(cf.\ Remark \ref{remark:inequalities-for-quasi-f-module-with-test-element})}. 
\qedhere 




    
\end{proof}

{One can use the above theorem to show that quasi-$F$-rationality agree with various non-uniform versions thereof (see Remark \ref{remark:non-uniform-quasi-F-regularity} for a similar discussion in the context of quasi-$F$-regularity). }\\

{Now we explain how to reconstruct $\tau(W_n\omega_R)$ from $\tau_n(\omega_R)$ {when $R$ is Cohen-Macaulay}. In particular, the following result shows that $\tau(W_n\omega_R) = W_n\omega_R$ if and only if $R$ is $F$-rational in the usual sense.}

\begin{proposition}
We work in the general setting (Notation \ref{n-non-local}). 
{\cred Then the following hold.}
\begin{enumerate}
\item {\cred ${\cccred \mathbf{R}}^*(W_{n-1}\omega_R) \subseteq W_n\omega_R$ and 
${\cccred \mathbf{V}}^*(W_n\omega_R) \subseteq F_*W_{n-1}\omega_R$.} 
\item If  $R$ is Cohen-Macaulay, then  
there exists the following commutative diagram in which each horizontal sequence is exact and the vertical homomorphisms are natural injections:
\[
\begin{tikzcd}
0 \ar{r} & \tau_n(\omega_R) \ar[hook]{d} \ar{r} & \tau(W_n\omega_R) \ar[hook]{d}\ar{r} & F_*\tau(W_{n-1}\omega_R) \ar[hook]{d} \ar{r} & 0 \\
0 \ar[r] & \omega_{{\cred R}} \ar{r}{({\cccred \mathbf{R}}^{n-1})^*} & W_n\omega_{{\cred R}} \ar{r}{{\cccred \mathbf{V}}^*} & F_*W_{n-1}\omega_{{\cred R}}  {\cred  \ar{r}} & {\cred 0}. 
\end{tikzcd}
\]
\end{enumerate}
\end{proposition}
\noindent Specifically, the proof below shows that the graded pieces of the natural ${\cccred \mathbf{V}}$-filtration on $\tau(W_n\omega_R)$ are isomorphic to $\tau_{n}(\omega_R), F_*\tau_{n-1}(\omega_R), \ldots, F^{n-1}_*\tau(\omega_R)$.

\begin{proof}
{\cred Taking localisation and completion,}  we may assume that $(R,\m)$ is a {\cred complete} local ring.

{\cred 
Let us show (1). 
By Proposition \ref{prop: tight cl l=n}, we have the following commutative diagrams 
\[
\begin{tikzcd}
H^d_\m(W_{n}R) \ar{r}{{\cccred \mathbf{R}}} & H^d_\m(W_{n-1}R)\\
\wt{0_{n}^{*}} \ar[hook]{u} \ar{r}{{\cccred \mathbf{R}}} & \wt{0_{n-1}^{*}}, \ar[hook]{u}
\end{tikzcd}
\qquad
\begin{tikzcd}
H^d_\m(F_*W_{n-1}R) \ar{r}{{\cccred \mathbf{V}}} & H^d_\m(W_{n}R)\\
F_*\wt{0_{n-1}^{*}} \ar[hook]{u} \ar{r}{{\cccred \mathbf{V}}} & \wt{0_{n}^{*}}. \ar[hook]{u}
\end{tikzcd}
\]
Applying Matlis duality, we get 
\[
R^* : \tau(W_{n-1}\omega_R) \overset{(\star)}{=} \Big( \frac{H^d_{\m}(W_{n-1}R)}{\wt{0_{n-1}^{*}}}\Big)^{\vee} 
\to  \Big( \frac{H^d_{\m}(W_{n}R)}{\wt{0_{n}^{*}}}\Big)^{\vee} 
\overset{(\star)}{=} \tau(W_{n}\omega_R), 
\]
\[
{\cccred \mathbf{V}}^* : \tau(W_{n}\omega_R) \overset{(\star)}{=} 
 \Big( \frac{H^d_{\m}(W_{n}R)}{\wt{0_{n}^{*}}}\Big)^{\vee} 
 \to  \Big( \frac{H^d_{\m}(F_*W_{n-1}R)}{F_*\wt{0_{n-1}^{*}}}\Big)^{\vee} \overset{(\star)}{=}  F_*
\tau(W_{n-1}\omega_R), 
\]
where ($\star$) holds by Proposition \ref{prop: corr 0^* and tau: no pair}.
Thus (1) holds.}


{\cred Let us show (2). By (1),} we have a complex like above:
\[
0 \to \tau_n(\omega_R)  \to \tau(W_n\omega_R) \to F_*\tau(W_{n-1}\omega_R) \to 0
\]
and it remains to show that it is exact. 
Consider the exact sequence
\[
H^d_\m(F_*W_{n-1}R) \xrightarrow{{\cccred \mathbf{V}}} H^d_\m(W_{n}R) \xrightarrow{{\cccred \mathbf{R}}^{n-1}} H^d_\m(R) \to 0. 
\]
By Theorem \ref{thm: lift of 0^*,n} and Proposition \ref{prop: tight cl l=n}, we get the following commutative diagram
\[
\begin{tikzcd}
 H^d_\m(F_*W_{n-1}R)  \ar{r}{{\cccred \mathbf{V}}}  & H^d_\m(W_{n}R) \ar{r}{{\cccred \mathbf{R}}^{n-1}} & H^d_\m(R) \ar{r} & 0 \\
 F_*\wt{0^*_{n-1}} \arrow[ur, phantom, "\urcorner", very near start] \arrow[hook]{u} \ar{r}{{\cccred \mathbf{V}}} & \wt{0^*_n} \arrow[hook]{u} \ar{r}{{\cccred \mathbf{R}}^{n-1}} & 0^*_n \ar{r} \arrow[hook]{u} & 0,
\end{tikzcd}
\]
where the left square is Cartesian. 
{Since} $R$ is Cohen-Macaulay, 
{\cred we get $H^{d-1}_\m(R)=0$, and hence} 
the homomorphism ${\cccred \mathbf{V}} \colon H^d_\m(F_*W_{n-1}R)  \to H^d_\m(W_{n}R)$  is injective. Thus the required diagram exists by applying Snake lemma, Matlis duality,  Theorem \ref{prop: corr 0^* and tau: no pair}, and Theorem \ref{thm: corr 0^* and tau}. 
{\cred Thus (2) holds.}
\end{proof}


\begin{proposition}\label{prop: another chara of qFrat}
We work in the local setting (Notation \ref{n-local}).  
Then the following hold. 
\begin{enumerate}
    \item 
If 
$\tau_n(\omega_R)=\omega_R$ for some integer $n>0$, then the induced $W_{m+n-1}R$-module homomorphism 
\[
{\cccred \mathbf{R}}^{n-1} \colon \wt{0_{m+n-1}^{*}} \to \wt{0_m^{*}}
\]
is zero for every integer  $m >0$.
\item The following are equivalent.
\begin{enumerate}
\renewcommand{\labelenumii}{(\roman{enumii})}
    \item $\tau^q(\omega_R)=\omega_R$.  
    \item $\varprojlim_{n} \wt{0_n^{*}}=0$. 
    \item $\varinjlim_{n} \tau(W_n\omega_R)=\varinjlim_{n} W_n\omega_R$, where each inductive system is induced by ${\cccred \mathbf{R}}^*$.
\end{enumerate}
\end{enumerate}
\end{proposition}
{\noindent We emphasise that $R^* : W_n\omega_R \to W_{n+1}\omega_R$ is injective for every $n>0$.}

\begin{proof}
Let us show (1). 
We prove that 
\[
{\cccred \mathbf{R}}^{n-1} \colon \wt{0_{m+n-1}^{*}} \to \wt{0_m^{*}}
\]
is zero by induction on $m \geq 1$. 
The base case $m=1$ of this induction follows from 
${\cccred \mathbf{R}}^{n-1}(\wt{0_n^{*}})=0_n^{*} = (\omega_R/\tau_n(\omega_R))^{\vee}=0$ (Theorem \ref{thm: lift of 0^*,n}, Theorem \ref{thm: corr 0^* and tau}). 
Assume  $m \geq 2$. Consider the following commutative diagram in which each horizontal sequence is exact: 
\begin{equation*}
\begin{tikzcd}
H^d_\m(F_*W_{m+n-2}R) \arrow{r}{{\cccred \mathbf{V}}} \arrow{d}{{\cccred \mathbf{R}}^{n-1}} & H^d_\m(W_{m+n-1}R) \arrow{r}{{\cccred \mathbf{R}}^{m+n-2}} \arrow{d}{{\cccred \mathbf{R}}^{n-1}} & H^d_\m(R) \arrow[d,equal] \\
H^d_\m(F_*W_{m-1}R) \arrow{r}{{\cccred \mathbf{V}}} & H^d_\m(W_mR) \arrow{r}{{\cccred \mathbf{R}}^{m-1}} & H^d_\m(R).
\end{tikzcd}
\end{equation*}
We have ${\cccred \mathbf{R}}^{m+n-2}(\wt{0_{m+n-1}^{*}})=0_{m+n-1}^{*} \subseteq 0^*_n =0$ 
as $m+n-1 \geq n$ (Remark \ref{remark: easy remark}(3), Theorem \ref{thm: lift of 0^*,n}). 
{\cred This, together with the upper exact sequence in the above diagram, 
implies $\wt{0_{m+n-1}^{*}} \subseteq {\cccred \mathbf{V}}({\cccred \mathbf{V}}^{-1}(\wt{0_{m+n-1}^{*}}))$.} 
By the following equality (Proposition \ref{prop: tight cl l=n}): 
\[
{\cccred \mathbf{V}}^{-1}(\wt{0_{m+n-1}^{*}}) = F_*\wt{0_{m+n-2}^{*}},
\]
it suffices to show ${\cccred \mathbf{R}}^{n-1}(\wt{0_{m+n-2}^{*}})=0$, 
which follows from the induction hypothesis. Thus (1) holds. 

Let us show (2). 
By (1), we get the implication (i) $\Rightarrow$ (ii).
To prove the implication (ii) $\Rightarrow$ (iii), 
{it is enough to prove that \[
\varinjlim_n\, 
 (W_n\omega_R/\tau(W_n\omega_R)) = 0
\]
as $\varinjlim_n$ is exact, which in turn follows from} 
$\wt{0^*_n} 
= (W_n\omega_R /\tau(W_n\omega_R))^{\vee}$ 
(Proposition \ref{prop: corr 0^* and tau: no pair}) {given that the Matlis duality functor $(-)^{\vee}$ turns colimits into limits}.    
Let us show (iii) $\Rightarrow$ (i). 
Pick $a \in \omega_R$. 
By (iii),  there exists $n>0$ such that $({\cccred \mathbf{R}}^*)^{n-1}(a) \in \tau(W_n\omega_R)$. 
Therefore, we have $a \in \tau_n(\omega_R)$, and thus $\omega_R=\tau^q(\omega_R)$.
\end{proof}

\subsection{Quasi-$+$-rationality}

In this subsection, we briefly discuss the concept of quasi-$+$-rationality and constructions of quasi-$+$-test submodules. These constructions could be done analogously to the above subsections by ways of quasi-$+$-closure. However, we will not pursue such an approach for two reasons. First, we do not understand quasi-$+$-closure as well as quasi-tight closure. Second, we wish to present a different, more geometric, approach to the problem, which is useful in other contexts (see, for example, {Subsection \ref{ss QF^eS}}). 
The main goal of this subsection is to show that quasi-$+$-rationality implies pseudo-rationality, which, in turn, provides an alternative proof of the same implication for quasi-$F$-rationality.

In this paragraph, we work in the general setting (Notation \ref{n-non-local}). Take an integer $n>0$ 
and a finite extension  $f \colon R \hookrightarrow S$ of reduced rings 
{(here a {\em finite extension} $f \colon R \hookrightarrow S$ is an injective ring homomorphism such that
 $S$ is a finitely generated $R$-module)}.  
We define a $W_nR$-module $Q^f_{R, n}$ and 
a $W_nR$-module homomorphism 
$\Phi^f_{R, n}$ by the following pushout diagram: 
\[
\begin{tikzcd}
W_nR \arrow[r, "W_nf"] \arrow[d, "{\cccred \mathbf{R}}^{n-1}"] & W_nS \arrow[d]\\
R  \arrow[r, "\Phi^f_{R, n}"] & Q^f_{R,n}. 
\end{tikzcd}
\]

\begin{definition}
We work in the general setting (Notation \ref{n-non-local}). 
We define the \emph{{\cred $n$-}quasi-$+$-test submodule} of $\omega_R$  by the following formula:
\begin{align*}
\tau_{+,n}(\omega_R) &:= \bigcap_{f : R \hookrightarrow S} \tau_{f,n}(\omega_R)\qquad  \text{ for } \\
\tau_{f,n}(\omega_R) &:= {\rm Im}\Big((\Phi^f_{R, n})^*  \colon \Hom_{W_nR}(Q^f_{R, n}, W_n\omega_R) \to \omega_R\Big),
\end{align*}
where $(\Phi^f_{R, n})^* := \Hom_{W_nR}(\Phi^f_{R, n}, W_n\omega_R)$ and the intersection is taken over all finite extensions $f \colon R \hookrightarrow S$ of reduced rings.
\end{definition}

\begin{definition} We work in the general setting (Notation \ref{n-non-local}). 
We say that $R$ is \emph{{\cred $n$-}quasi-$+$-rational} if $\tau_{+,n}(\omega_R) = \omega_R$. This is the same as saying that for every finite extension $f \colon R \hookrightarrow S$ of reduced rings, the following map is surjective:
\[
(\Phi^f_{R, n})^* \colon \Hom_{W_nR}(Q^f_{R, n}, W_n\omega_R) \to \omega_R.
\]
{\cred We say that $R$ is \emph{quasi-$+$-rational} if $R$ is $n$-quasi-$+$-rational for some integer $n>0$.} 
\end{definition}

\begin{proposition}
We work in the general setting (Notation \ref{n-non-local}). 
{\cred Fix an integer $n>0$.} 
Assume that $R$ is {\cred $n$-}quasi-$F$-rational. Then $R$ is 
{\cred $n$-}quasi-$+$-rational.
\end{proposition}

\begin{proof}
The proof is analogous to \cite[Proposition 4.9]{TWY}.
\end{proof}

To demonstrate that quasi-$+$-rationality implies pseudo-rationality, we introduce the quasi-$+$-test submodule of $W_n\omega_R$.
\begin{definition}\label{d tau+ T_nwomega def}
We work in the general setting (Notation \ref{n-non-local}). We define the \emph{quasi-$+$-test $W_nR$-submodule} of $W_n\omega_R$  by the following formula:
\begin{align*}
\tau_+(W_n\omega_R) &:= \bigcap_{f : R \hookrightarrow S} \tau_f(W_n\omega_R)\qquad  \text{ for }\\ 
\tau_f(W_n\omega_R) &:= {\rm Im}\Big(T^f_{n}  \colon W_n\omega_S \to W_n\omega_R\Big),
\end{align*}
where $T^f_{n} := \Hom_{W_nR}({W_nf}, W_n\omega_R)$ and the intersection is taken over all finite extensions $f \colon R \hookrightarrow S$ of reduced rings.
\end{definition}

\begin{proposition} \label{prop:witt-+-from-full-witt}
We work in the general setting (Notation \ref{n-non-local}). Then \[
\tau_{+,n}(\omega_R) = (({\cccred \mathbf{R}}^{n-1})^*)^{-1}(\tau_+(W_n\omega_R)).
\]
In other words, there exists the following pullback diagram in which the vertical arrows are natural inclusions
\[
\begin{tikzcd}
 \tau_{+,n}(\omega_R) \ar[hook]{d} \ar[hook]{r} & \tau_+(W_n\omega_R) \ar[hook]{d} \\
 \omega_R \ar[hook]{r}{({\cccred \mathbf{R}}^{n-1})^*} & W_n\omega_R.
\end{tikzcd}
\]    

\end{proposition}
\begin{proof}
Let $f \colon R \to S$ be a finite extension of reduced rings. Consider the following diagram
\[
\begin{tikzcd}
0 \ar{r} & F_*W_{n-1}R \ar{r} & W_nS \arrow[r] & Q^f_{R,n} \ar{r} & 0 \\
0 \ar{r} & F_*W_{n-1}R \ar{u}{=} \ar{r}{{\cccred \mathbf{V}}} & W_nR \arrow[u, "W_nf"] \arrow[r, "{\cccred \mathbf{R}}^{n-1}"] & R \arrow[u, "\Phi^f_{R, n}"]  \ar{r} & 0.
\end{tikzcd}
\]    
By applying $\Hom_{W_nR}(-, W_n\omega_R)$, we get:
\[
\begin{tikzcd}
 F_*W_{n-1}\omega_R \ar{d}{=}  & \ar{l} W_n\omega_S \ar{d}{T^f_{n}} \arrow[l] & \ar{l} \Hom_{W_nR}(Q^f_{R,n}, W_n\omega_R) \ar{d}{(\Phi^f_{R, n})^*}    & \ar{l} 0 \\
 F_*W_{n-1}\omega_R   & \ar{l}{{\cccred \mathbf{V}}^*} W_n\omega_R  & \omega_R \arrow{l}{({\cccred \mathbf{R}}^{n-1})^*} & \ar{l} 0.
 \arrow[phantom, from=1-2, to=2-3, "(\star)" description]
\end{tikzcd}
\]    
Then $(\star)$ is a pullback square by diagram chase in which it is essential that the left vertical arrow is an equality. This immediately implies that:
\[
\tau_{+,n}(\omega_R) =  (({\cccred \mathbf{R}}^{n-1})^*)^{-1}(\tau_+(W_n\omega_R)). \qedhere
\] 
\end{proof}

\begin{remark} There are a couple of results we leave unproven in this subsection. For starters, we do not argue that $\tau(W_n\omega_R)$ and $\tau_+(\omega_R)$ are calculated by a single finite extension, and, in particular, that they are stable under localisation and completion. This should follow by the same argument as in 
{\cred \cite[Proposition 4.13]{TWY}.}

Second, with a little bit of work, one should be able to establish that there exists the following commutative diagram in which each row is exact and vertical maps are natural injections:
\[
\begin{tikzcd}
0 \ar{r} & \tau_{+,n}(\omega_R) \ar[hook]{d} \ar{r} & \tau_+(W_n\omega_R) \ar[hook]{d}\ar{r} & F_*\tau_+(W_{n-1}\omega_R) \ar[hook]{d} \ar{r} & 0 \\
0 \ar[r] & \omega_R \ar{r}{({\cccred \mathbf{R}}^{n-1})^*} & W_n\omega_R \ar{r}{{\cccred \mathbf{V}}^*} & F_*W_{n-1}\omega_R. & \hphantom{a}
\end{tikzcd}
\]    
Specifically the graded pieces of the natural ${\cccred \mathbf{V}}$-filtration on $\tau_+(W_n\omega_R)$ should be isomorphic to $\tau_{+,n}(\omega_R), F_*\tau_{+,n-1}(\omega_R), \ldots, F^{n-1}_*\tau_+(\omega_R)$.
\end{remark}

We will also need the generalisation of Hochster-Huneke's vanishing to the Witt setting.
\begin{lemma} \label{lem:HHWitt} We work in the general setting (Notation \ref{n-non-local}). Fix an integer $n>0$. Let $\pi \colon X \to \Spec R$ be a proper morphism, {\cred where $X$ is a reduced scheme}. There there exists a finite {surjective} morphism $f \colon Y \to X$ over $R$ such that {\cred $Y$ is a reduced scheme and} the pullback maps
\[
H^i(X, W_n\cO_X) \xrightarrow{f^*} H^i(Y, W_n\cO_Y)
\]
are zero for every $i>0$.
\end{lemma}
\begin{proof}
The case of $n=1$ follows from \cite[Introduction, Theorem 1.5]{Bhatt_2012}. In general, we argue by an increasing induction on $n$. Let $f \colon Y \to X$ be a finite morphism over $R$ such that
\[
H^i(X, \cO_X) \xrightarrow{f^*} H^i(Y, \cO_Y)
\]
is zero for every $i>0$. Assume by induction that there exists a finite  {surjective} morphism $g \colon Z \to Y$ over $R$ such that 
\[
H^i(Y, W_{n-1}\cO_Y) \xrightarrow{g^*} H^i(Z, W_{n-1}\cO_Z)
\]
is also zero for every $i>0$. Then the composition
\[
H^i(X, W_n\cO_X) \xrightarrow{(f \circ g)^*} H^i(Z, W_n\cO_Z)
\]
is zero for every $i>0$ by tracing through the following commutative diagram 
\[
\begin{tikzcd}
0 \ar{r} & {\ccred f_*g_*}F_*W_{n-1}\cO_Z \ar{r}{{\cccred \mathbf{V}}} 
& {\ccred f_*g_*}W_n\cO_Z \ar{r}{{\cccred \mathbf{R}}^{n-1}} 
& {\ccred f_*g_*}\cO_Z \ar{r} & 0   \\
0 \ar{r} & {\ccred f_*}F_*W_{n-1}\cO_Y \ar{u}{g^*} \ar{r}{{\cccred \mathbf{V}}} 
& {\ccred f_*}W_n\cO_Y \ar{u}{g^*} \ar{r}{{\cccred \mathbf{R}}^{n-1}} 
& {\ccred f_*}\cO_Y \ar{u}{g^*} \ar{r} & 0   \\
0 \ar{r} & F_*W_{n-1}\cO_X \ar{u}{f^*} \ar{r}{{\cccred \mathbf{V}}} & W_n\cO_X \ar{u}{f^*} \ar{r}{{\cccred \mathbf{R}}^{n-1}} & \cO_X \ar{u}{f^*} \ar{r} & 0   
\end{tikzcd}
\]
with $H^i$ applied to it.
\end{proof}

\begin{remark} \label{remark:witt-connected-fibres}
We work in the general setting (Notation \ref{n-non-local}). Fix an integer $n>0$ and let $\pi \colon Y \to X$ be a projective 
morphism over $R$ {\cred such that 
$Y$ and $X$ are reduced schemes and $\pi_*\MO_Y = \MO_X$, i.e., 
the induced homomorphism $\MO_X \to \pi_*\MO_Y$ is an isomorphism}. Then $\pi_*W_n\cO_Y = W_n\cO_X$. Indeed, we have the natural commutative diagram
\[
\begin{tikzcd}
0 \ar{r} & \pi_*F_*W_{n-1}\cO_Y \ar{r}{\pi_*{\cccred \mathbf{V}}} & \pi_*W_n\cO_Y \ar{r}{\pi_* {\cccred \mathbf{R}}^{n-1}} & \pi_*\cO_Y & \hphantom{a} \\
0 \ar{r} & F_*W_{n-1}\cO_X \ar{u}{=} \ar{r}{{\cccred \mathbf{V}}} & W_n\cO_X \ar{u}{\pi^*} \ar{r}{{\cccred \mathbf{R}}^{n-1}} & \ar{u}{=} \cO_X \ar{r} & 0, 
\end{tikzcd}
\]
in which the left vertical arrow is an equality by ascending induction on $n$. Then $\pi_*{\cccred \mathbf{R}}^{n-1}$ is surjective by diagram chase, and so $\pi^* \colon W_n \cO_X \to \pi_*W_n\cO_Y$ is an equality by Five Lemma.
\end{remark}

By combining Lemma \ref{lem:HHWitt} and Remark \ref{remark:witt-connected-fibres}, we 
see that the ind-systems $\{R\pi_*\cO_X\}_{\pi \colon X \to \Spec R}$ and $\{\pi_*\cO_X\}_{\pi \colon X \to \Spec R}$ taken over all alterations and all finite surjective morphisms, respectively, are equivalent. By Matlis duality, this yields the following result.

\begin{proposition} \label{prop:full-witt-pushforward} We work in the general setting (Notation \ref{n-non-local}). Let $\pi \colon X \to \Spec R$ be a projective birational morphism  
{\cred of integral schemes}. Then $\tau_{+}(W_n\omega_R) \subseteq \pi_*W_n\omega_X$.
\end{proposition}

\begin{proof}
{\cred Taking a localisation, we may assume that $(R, \m)$ is a local ring.} 
{\cred Recall that $d = \dim R = \dim X$.} 
By Lemma \ref{lem:HHWitt} {applied $d$ times\footnote{we warn the reader that it is \emph{not} true in general that a map in the derived category is zero if and only if it is zero on each cohomology} (see also \cite[Lemma 2.2]{Bhatt_2012}, cf.\ the proof of \cite[Theorem 0.4]{Bhatt_2012}),} we get that there exists a finite surjective morphism  $f \colon Y \to X$ over $R$ such that the pullback 
{\[
f^* \colon R^{>0}\pi_*W_n\cO_X \to R^{>0}\rho_*W_n\cO_Y
\]
is {\cred the} zero map,} 
where $\rho := \pi \circ f \colon Y \xrightarrow{f} X \xrightarrow{\pi} \Spec R$. {Let $D$ denote the derived category of $W_nR$-modules.} In particular, by applying  $\Hom_{D}(R\pi_*W_n\cO_X, -)$ 
to the exact triangle (cf.\ the proof of \cite[Theorem 0.3]{Bhatt_2012}):
\[
\rho_*W_n\cO_Y \to R\rho_*W_n\cO_Y \to R^{>0}\rho_*W_n\cO_Y \xrightarrow{+1},
\]
we get an exact sequence of abelian groups:  
\begin{align*}
\Hom_D(R\pi_*W_n\MO_X, \rho_*W_n\cO_Y) &\to 
\Hom_D(R\pi_*W_n\MO_X, R\rho_*W_n\cO_Y) \\
&\to 
\Hom_D(R\pi_*W_n\MO_X, R^{>0}\rho_*W_n\cO_Y). 
\end{align*}
By using this sequence, we get that
\[
f^* \colon R\pi_*W_n\cO_X \to R\rho_*W_n\cO_Y
\]
factors through $\rho_*W_n\cO_Y$. {Choose one such factorisation}.
By Remark \ref{remark:witt-connected-fibres}, we have that $\rho_*W_n\cO_Y = \theta_*W_n\cO_W$ where $\theta \colon W \to \Spec R$ is the finite part of the Stein factorisation of $\rho \colon Y \to \Spec R$:
\[
\begin{tikzcd}
Y \ar{r}{f} \ar{d} \ar{rd}{\rho} & X \ar{d}{\pi} \\
W \ar{r}{\theta} & \Spec R.
\end{tikzcd}
\]
To conclude, we have a sequence of maps:
\[
\theta^* \colon {W_n}R \xrightarrow{\pi^*} R\pi_*W_n\cO_X \xrightarrow{f^*} \rho_*W_n\cO_Y = \theta_*W_n\cO_W.
\]
In particular, 
we get an induced surjection:
\[
{\rm Im}\big(\pi^* \colon H^d_\m({W_n}R) \to H^d_\m(R\pi_*W_n\cO_X)\big) \twoheadrightarrow {\rm Im}\big(\theta^* \colon H^d_\m({W_n}R) \to H^d_\m(\theta_*W_n\cO_W)\big). 
\]
{By applying Matlis duality $(-)^{\vee} = \Hom_{W_nR}(-,E)$, 
\[
{\cred \tau_\theta(W_n\omega_R) \otimes_{W_nR} W_n\widehat{R}  \hookrightarrow \pi_*W_n\omega_X  \otimes_{W_nR} W_n\widehat{R},} 
\]
{\cred where $\widehat{R}$ denotes the $\m$-adic completion of $R$.}
Here the Matlis dual $H^d_\m(R\pi_*W_n\cO_X)^{\vee}$ can be computed as follows:}
\[
H^d_\m(R\pi_*W_n\cO_X)^{\vee} \simeq 
\Big( H^0R\Hom_{{W_nR}}(R\pi_*W_n\cO_X, {W_n}\omega^\bullet_{{R}})
\Big) \otimes_{W_nR} W_n\widehat{R}  \simeq \pi_*W_n\omega_{{X}} \otimes_{W_nR} W_n\widehat{R} , 
\]
where the first isomorphism is \cite[Proposition 2.2]{KTTWYY1} and the second one is Grothendieck duality.
Therefore, we get inclusions 
 $\tau_+(W_n\omega_R) \subseteq  
 {\cred \tau_\theta(W_n\omega_R)}  \subseteq \pi_*W_n\omega_X$, 
{\cred where the first inclusion follows from Definition \ref{d tau+ T_nwomega def} and the second one holds because $W_nR \to W_n\widehat{R}$ is faithfully flat}.
\end{proof}

\begin{theorem} \label{thm:quasi+rational-are-pseudorational} We work in the general setting (Notation \ref{n-non-local}). Let $\pi \colon X \to \Spec R$ be a projective birational morphism {\cred of integral schemes}. Then $\tau_{+,n}(\omega_R) \subseteq \pi_*\omega_{{\cred X}}$. In particular, if $R$ is quasi-$+$-rational, then

it is pseudo-rational.
\end{theorem}
\begin{proof}
Since $W_n\omega_X$ is $S_1$ 
{\cred \cite[Tag 0AWK]{stacks-project}}, 
the vertical arrows in the following diagram are injective
\[
\begin{tikzcd}
0 \ar{r} & \pi_*\omega_X \ar[hook]{d} \ar{r}{({\cccred \mathbf{R}}^{n-1})^*} & \pi_* W_n\omega_X \ar[hook]{d} \ar{r}{{\cccred \mathbf{V}}^*} & F_*\pi_*W_{n-1}\omega_X \ar[hook]{d} \\
0 \ar{r} & \omega_R \ar{r}{({\cccred \mathbf{R}}^{n-1})^*} & W_n\omega_R \ar{r}{{\cccred \mathbf{V}}^*} & F_*W_{n-1}\omega_R.
\end{tikzcd}
\]
By diagram chase, we can thus see that
\[
\pi_*\omega_X = (({\cccred \mathbf{R}}^{n-1})^*)^{-1}(\pi_*W_n\omega_X).
\]
Therefore, by Proposition \ref{prop:witt-+-from-full-witt} and Proposition \ref{prop:full-witt-pushforward}, we get that
\[
\tau_{+,n}(\omega_R) = (({\cccred \mathbf{R}}^{n-1})^*)^{-1}(\tau_+(W_n\omega_R)) \subseteq (({\cccred \mathbf{R}}^{n-1})^*)^{-1}(\pi_*W_n\omega_X) =  \pi_*\omega_X, 
\]
concluding the proof of the theorem.
\end{proof}

\section{Quasi-test ideals {and quasi-test submodules for pairs}}

This section is the central part of this paper, in which we shall introduce 
the quasi-test ideal $\tau^q(R, \Delta)$, 
where $R$ is  an $F$-finite normal integral domain of characteristic $p>0$ 
and $\Delta$ is  a $\Q$-divisor on $\Spec R$ such that $K_R +\Delta$ is $\Q$-Cartier. 
Set $D :=K_R+\Delta$. 
In Subsection \ref{ss log test submod} and  Subsection \ref{ss log tight closure}, we shall introduce the {\em $n$-quasi-test submodule} $\tau_n(\omega_R, D) \subseteq \omega_R(-D) := \cHom_R(R(D), \omega_R)$  and
the {\em $n$-quasi-tight closure} $0^*_{D, n} \subseteq H^d_{\m}(R(D))$, respectively.  {We warn the reader again that $R(D) = R(\rdown{D})$ and $\omega_R(-D) = \omega_R(\rup{-D})$ 
 in our notation {(cf.\ Remark \ref{r-W_n omega D}(1))}.}

These subsections are the log versions of Subsection \ref{ss q test submod} and
Subsection \ref{ss q tight}. 
Then the {\em $n$-quasi-test ideal} $\tau_n(R, \Delta)$ is defined by 
\[
\tau_n(R, \Delta) := \tau_n(\omega_R, K_R+\Delta). 
\]
The {\em quasi-test ideal} $\tau^q(R, \Delta)$ is given  by 
\[
\tau^q(R, \Delta) := \bigcup_{n\geq 1} \tau_n(R, \Delta) =\tau_n(R, \Delta) 
\]
for $n \gg 0$, where the latter equality holds by $\tau_1(R, \Delta) \subseteq \tau_2(R, \Delta) \subseteq \cdots$. 
We then establish some fundamental properties of $\tau_n(R,\Delta)$ (Subsection \ref{ss q-test ideal}): 
\begin{enumerate}
\item 
$\tau_n(R, \Delta)$ is compatible with localisations and completions (Remark \ref{r q test ideal clear}). 
\item 
If $E$ is an effective Cartier divisor on $\Spec R$, then
there exists a rational number $\epsilon >0$ 
such that $\tau_n(R, \Delta) = \tau_n(R, \Delta + \epsilon E)$  (Theorem \ref{thm: suff small qFr}). 
\item $(R, \Delta)$ is $n$-quasi-$F$-regular if and only if $\tau_n(R, \Delta) = R$  
(Proposition \ref{prop: test ideal pair non-tilda}). 
\end{enumerate}
In order to show (2) and (3), we prove the following  duality 
for the case 
when $R$ is a local ring (Proposition \ref{prop: test ideal pair non-tilda}): 
    \[
\Hom_R(R(\rup{-\Delta}) / \tau_n(R, \Delta), E)  = 0^{*}_{K_R+\Delta, n}.
\] 





\subsection{Quasi-test submodules (pairs)}\label{ss log test submod}



\begin{notation}\label{n-normal-non-local}
Let $R$ be an $F$-finite Noetherian normal integral domain of characteristic $p>0$. 
Set $X:=\Spec R$. 
\end{notation}

\begin{definition}
We use Notation \ref{n-normal-non-local}. 
Let $D$ be a $\Q$-divisor on $X$. 
Take $g \in R(-(p^e-1)D)$, which induces a $W_nR$-module homomorphism: 
\[
W_nR \xrightarrow{\cdot [g]} W_nR(- (p^e-1)D). 
\]
By applying $\Hom_{W_nR}(-, W_n\omega_R)$ to  the  $W_nR$-module homomorphisms
\[
W_nR \xrightarrow{F^e} F_*^eW_nR \xrightarrow{\cdot F_*^e[g]} F_*^eW_nR(- (p^e-1)D), 
\]
we obtain $W_nR$-module homomorphisms 
\[
T^{e, g}_{n} : 
F_*^eW_n\omega_R( (p^e-1)D)
\xrightarrow{ \cdot F_*^e[g]}
F_*^eW_n\omega_R 
\xrightarrow{T^e_n} 
W_n\omega_R, 
\]
where we used 
\begin{itemize}
\item 
$\Hom_{W_nR}(F_*^eW_nR, W_n\omega_R) \simeq F_*^e W_n\omega_R$, and 
\item 
$\Hom_{W_nR}(F_*^eW_nR(- (p^e-1)D), W_n\omega_R) \simeq 
F_*^eW_n\omega_R((p^e-1)D)$. 
\end{itemize}
For an integer $e \geq 0$, we set 
\[
C^e_n(D):=\big\{\PsiT^{e, g}_{n}  \,\big|\, g \in R(-(p^e-1)D)\big\} \subseteq 
\Hom_{W_nR}(F^e_*W_n\omega_R((p^e-1)D),W_n\omega_R). 
\]  
\end{definition}

\begin{definition}\label{d-log-cosmall-stable}
We use Notation \ref{n-normal-non-local}. 
Take an integer $n>0$. 
Let $D$ be an effective $\Q$-divisor on $X$. 
Let $M$ be a $W_nR$-submodule of $W_n\omega_R$.
\begin{enumerate}
\item We say that $M$ is {\em co-small} if there exists $c \in R^{\circ}$ such that 
$[c] \cdot W_n\omega_R \subseteq M$. 
\item 
We say that $M$ is {\em $C_n(D)$-stable} 
if the inclusion 
\[
\varphi(F^e_*M) \subseteq M
\]
holds for every integer $e \geq 0$ and every $\varphi \in C^e_n(D)$. 
Note that $F_*^eM \subseteq F_*^e W_n\omega_R \subseteq F_*^e W_n\omega_R( (p^e-1)D)$ 
(Remark \ref{r-W_nomega-j^*}). 
\end{enumerate}
\end{definition}

{Note that if $(p^e-1)D = \div(g)$ for some $g \in R$ (for example, this is always true when $R$ is local and $(p^e-1)D$ is Cartier), then in Definition \ref{d-log-cosmall-stable}(2), it would be enough for us to only consider $\varphi = T^{e,g}_n$ for this particular choice of $g$.}


\begin{remark}\label{r-T^{e,g}-composition}
Take $T^{e, g}_{n} \in C^e_n(D)$ and $T^{e', g'}_{n} \in C^{e'}_n(D)$. 
Let us construct the natural composition $T^{e, g}_{n} \circ (F^e_*T^{e', g'}_{n})$.
Recall that $g \in R(-(p^e-1)D)$ and $g' \in R(-(p^{e'}-1)D)$. 
Consider the composition
\begin{alignat*}{3}
W_nR \xrightarrow{F^e}\ &F_*^eW_nR \xrightarrow{\cdot F_*^e[g]} F_*^eW_nR(-(p^e-1)D) && \\
\xrightarrow{F^{e'}}\ &F_*^{e+e'}W_nR(-p^{e'}(p^e-1)D) \xrightarrow{\cdot F_*^{e+e'}[g']} 
&&F_*^{e+e'}W_nR(-p^{e'}(p^e-1)D - (p^{e'}-1)D) \\[0.5em]
&  &&\mathllap{=\ }F_*^{e+e'}W_nR(-(p^{e+e'}-1)D). 
\end{alignat*}
Applying $\Hom_{W_nR}(-, W_n\omega_R)$, we obtain $T^{e, g}_{n} \circ (F^e_*T^{e'}_{n, g'})$: 
\[
F_*^{e+e'}W_nR((p^{e+e'}-1)D) \xrightarrow{F^e_*T^{e', g'}_{n}} 
F_*^{e}W_nR((p^{e}-1)D) \xrightarrow{T^{e, g}_{n}} W_n\omega_R. 
\]
Moreover, we have $T^{e, g}_{n} \circ (F^e_*T^{e', g'}_{n}) = T^{e+e',g^{p^{e'}}g'}_{n}$ 
by 
\[
((\cdot F^{e+e'}_*[g']) \circ F^{e'}) \circ ((\cdot F^e_*[g]) \circ F^e) = (\cdot F^{e+e'}_*[g^{p^e}g']) \circ F^{e+e'}. 
\]
\end{remark}
\begin{proposition}\label{prop: first prop of C_n(D)-sub}
We use Notation \ref{n-normal-non-local}. 
Let $D$ be an effective $\Q$-divisor on $X$. 
Fix a test element $t \in R^{\circ}$.
Take $d \in R^{\circ}$ such that $D \leq \mathrm{div}(d)=:D'$.
Then the following hold. 
\begin{enumerate}
    \item $W_n\omega_R(-D)$ is a co-small $C_n(D)$-stable $W_nR$-submodule of $W_n\omega_R$.  
    \item For every $C_n(D')$-stable co-small  $W_nR$-submodule $M$ of $W_n\omega_R(-D')$, 
    $[d^{-1}] \cdot M$ is a co-small ${\cccred T_n}$-stable $W_nR$-submodule  of $W_n\omega_R$. 
    \item For every $C_n(D)$-stable co-small $W_nR$-submodule $M$ of $W_n\omega_R$, 
    it holds that 
    \[
    [t^2d]\cdot W_n\omega_R \subseteq M.
    \]
\end{enumerate}
\end{proposition}

\begin{proof}
Let us show (1). 
Take an integer $e \geq 0$ and $g \in R(-(p^e-1)D)$.
Then $\PsiT^{e, g}_{n} \in C_n(D)$ induces the following $W_nR$-module homomorphisms:  
\[
F^e_*W_n\omega_R(-D) \xrightarrow{\cdot F^e_*[g]} F^e_*W_n\omega_R(-p^eD) \xrightarrow{\PsiT_{n}^e}  W_n\omega_R(-D).  
\]
Therefore, $W_n\omega_R(-D)$ is $C_n(D)$-stable.
Since $[d] \cdot W_n\omega_R = W_n\omega_R(-\div(d)) \subseteq W_n\omega_R(-D)$,  we get that
$W_n\omega_R(-D)$ is co-small. 
Thus (1) holds.

Let us show (2). 
We have 
\[
[d^{-1}] \cdot M \subseteq [d^{-1}] \cdot W_n\omega_R(-D') 
= W_n\omega_R(\div(d)-D')  =W_n\omega_R.
\]
Moreover, it holds that 
\begin{align*}
    &[d] \cdot \PsiT_{n}^e(F^e_*([d^{-1}]\cdot M)) = \PsiT_{n}^e(F^e_*([d^{p^e-1}]\cdot M)) \subseteq M, 
\end{align*}
because  $M$ is $C_n(D')$-stable and $d^{p^e-1} \in R(-(p^e-1)D')$.
Therefore, $[d^{-1}]\cdot M$ is ${\cccred T_n}$-stable and co-small. 
Thus (2) holds.

Let us show (3). 
We now reduce the problem to the case when $D=D'$. 
By $D \leq D'$, we get the natural {\ccred homomorphism} $C^e_n(D') \to C^e_n(D)$ given by $R(-(p^e-1)D') \subseteq R(-(p^e-1)D)$. 
Specifically, 
if $T_n^{e, h} \in C^e_n(D')$ with $h \in R(-(p^e-1)D') \subseteq R(-(p^e-1)D)$, 
then we get  
\[
\begin{tikzcd}
 F_*^eW_n\omega_R((p^e-1)D') \arrow[rr, bend left, "{T_n^{e, h}}"'] \arrow[r, "{\cdot F_*^e[h]}"]  & F_*^eW_n\omega_R \arrow[r, "{T^e_n}"] & W_n\omega_R\\
 F_*^eW_n\omega_R((p^e-1)D) \arrow[u, hook] \arrow[r, "{\cdot F_*^e[h]}"]  & F_*^eW_n\omega_R. \arrow[u, equal] 
\end{tikzcd}
\]
Since $M$ is $C_n(D)$-stable, $M$ is also $C_n(D')$-stable.
Therefore, we may assume $D=D'$. 

Since $W_n\omega_R(-D)=[d]\cdot W_n\omega_R$ is $C_n(D)$-stable by $(1)$, 
also $M':=M \cap [d] \cdot W_n\omega_R$ is  $C_n(D)$-stable and co-small. 
Replacing $M$ by $M'$, 
we may assume that $M \subseteq [d] \cdot W_n\omega_R$.
By $(2)$, $[d^{-1}]\cdot M$ is ${\cccred T_n}$-stable and co-small.
By Proposition \ref{prop: first prop of F-stable}(2), we have
\[
[t^2d] \cdot W_n\omega_R = [d] \cdot ( [t^2] \cdot W_n\omega_R) \subseteq [d] \cdot ([d^{-1}] \cdot M)=M.
\]
Thus (3) holds. 
\end{proof}



\begin{proposition}\label{prop: test submodule gen pair}
We use Notation \ref{n-normal-non-local}. 
Let $D$ be an effective $\Q$-divisor on $X$ and take $d \in R^{\circ}$ such that $D \leq \mathrm{div}(d)$. 
Fix a test element $t \in R^{\circ}$  and $c \in (dt^2) \cap R^{\circ}$. 
Then the following hold. 
\begin{enumerate}
    \item There exists a smallest co-small $C_n(D)$-stable 
    $W_nR$-submodule $\tau(W_n\omega_R,D)$ of $W_n\omega_R$. 
    Moreover, $\tau(W_n\omega_R,D) \subseteq W_n\omega_R(-D)$. 
    \item For $h \in R^{\circ}$, it holds that 
    \[
    \tau(W_n\omega_R,D+\mathrm{div}(h))=[h] \cdot \tau(W_n\omega_R,D).
    \]
    \item     For every co-small $W_nR$-submodule $M$ of $W_n\omega_R$, it holds that 
    \[
    \tau(W_n\omega_R,D)=\sum_{e \geq 0} \sum_{\varphi \in C^e_n(D)} \varphi(F^e_*([c] \cdot M))
    \]
    \item If $f \in R^{\circ}$ and $D=a \mathrm{div}(f)$ for some $a \in \Q_{>0}$, then 
    \[
    \tau(W_n\omega_R,D)=\sum_{e \geq 0} \PsiT_n^e(F^e_*([cf^{\rup{a(p^e-1)}}] \cdot M)) 
    \]
    for every co-small $W_nR$-submodule $M$ of $W_n\omega_R$.
    \item 
    Assume 
$D=a \mathrm{div}(f)$ and $c \in (f^{\rup{a}}) \cap (dt^2) \cap R^{\circ}$  
for some $a \in \Q_{>0}$ and $f \in R^{\circ}$. 
\begin{enumerate}
\item For every integer $e \geq 0$, the inclusion 
\begin{align*}
T_n^{e}(&F^{e}_*([cf^{\rup{a(p^{e}-1)}}] \cdot \tau(W_n\omega_R, D))) \\
&\quad \subseteq 
T_n^{e+1}(F^{e+1}_*([cf^{\rup{a(p^{e+1}-1)}}] \cdot \tau(W_n\omega_R, D))) 
\end{align*}
holds. 
\item  There exists an integer $e_0 \geq 0$ such that 
\[
\tau(W_n\omega_R,D)= T_n^{e}(F^{e}_*([cf^{\rup{a(p^{e}-1)}}] \cdot M)) 
\]
holds for every $e \geq e_0$ and every co-small $C_n(D)$-stable $W_nR$-submodule $M$ of $W_n\omega_R$. 
\item 
If $L$ is a co-small $W_nR$-submodule of $W_n\omega_R$, 
then there exists an integer $e_L \geq 0$ such that 
\[
\tau(W_n\omega_R,D)= T_n^{e}(F^{e}_*([cf^{\rup{a(p^{e}-1)}}] \cdot L)) 
\]
holds for every $e \geq e_L$. 
\end{enumerate}
\end{enumerate}
\end{proposition}

We call $\tau(W_n\omega_R, D)$ the {\em quasi-test $(W_nR$-$)$submodule} of $(W_n\omega_R, D)$. 

\begin{proof}
Let us show (1). 
Set 
\[
\tau(W_n\omega_R,D) := \bigcap_M M \subseteq W_n\omega_R,  
\]
where $M$ runs over all the $C_n(D)$-stable co-small $W_nR$-submodules $M$ of $W_n\omega_R$. 
By Proposition \ref{prop: first prop of C_n(D)-sub}(3), 
$\tau(W_n\omega_R,D)$  is co-small. 
Proposition \ref{prop: first prop of C_n(D)-sub}(1) implies 
the inclusion 
$\tau(W_n\omega_R,D) \subseteq W_n\omega_R(-D)$. 
Thus (1) holds. 


Let us show (2). 
Since $\tau(W_n\omega_R,D+\mathrm{div}(h))$ is $C_n(D+\mathrm{div}(h))$-stable, 
it follows from a similar argument to the one of Proposition \ref{prop: first prop of C_n(D)-sub} 
that $[h^{-1}]\cdot \tau(W_n\omega_R,D+\mathrm{div}(h))$ is $C_n(D)$-stable. 
By (1), we get $\tau(W_n\omega_R,D+\mathrm{div}(h))  \supseteq [h] \cdot \tau(W_n\omega_R,D)$. 
The opposite inclusion holds by (1) and the fact that  $[h] \cdot \tau(W_n\omega_R,D)$ is 
$C_n(D+\mathrm{div}(h))$-stable.
 Thus (2) holds. 


Let us show (3). 
We denote by $N$ the right hand side of (3). 
By ${\rm id}_{W_n\omega_R} \in C^0_n(D)$ (where ${\rm id}_{W_n\omega_R} : W_n\omega_R \to W_n\omega_R$ denotes the identity map), we obtain 
\[
N \supseteq {\rm id}_{W_n\omega_R}(F_*^0[c] \cdot M) = [c] \cdot M. 
\]
Since $[c] \cdot M$ is co-small, so is $N$. 
For every $\PsiT^{e', h}_n \in C^{e'}_n(D)$, we have  
\begin{align*}
     \PsiT^{e', h}_n(F^{e'}_*N) 
     &=\sum_{e\geq 0} \sum_{\varphi \in C^e_{n}(D)} \PsiT^{e', h}_n(F^{e'}_* \varphi(F^e_*[c] \cdot M)) \\
    &=\sum_{e\geq 0} \sum_{g \in R(-(p^e-1)D)} \PsiT^{e', h}_n(F^{e'}_* T^{e, g}_n(F^e_*[c] \cdot M)) \\
    &\overset{(\star)}{=} \sum_{e\geq 0} \sum_{g \in R(-(p^e-1)D)} 
    T^{e+e', h^{p^{e'}}\!\!\!g}_n(F^{e+e'}_*[c] \cdot M)) 
    \subseteq N, 
\end{align*}
where $(\star)$ follows from Remark \ref{r-T^{e,g}-composition}. 
Therefore, $N$ is $C_n(D)$-stable. Then (1) implies 
\[
\tau(W_n\omega_R,D) \subseteq N. 
\]
The opposite inclusion $N \subseteq \tau(W_n\omega_R,D)$ 
holds by 
\begin{align*}
    N
    &= \sum_{e \geq 0} \sum_{\varphi \in C^e_n(D)} \varphi(F^e_*[c] \cdot M) \\
    & \subseteq \sum_{e \geq 0} \sum_{\varphi \in C^e_n(D)} \varphi(F^e_*[c] \cdot W_n\omega_R) \\
    & \overset{{\rm (i)}}{\subseteq} \sum_{e \geq 0} \sum_{\varphi \in C^e_n(D)} \varphi(F^e_*\tau(W_n\omega_R,D)) \\
    & \overset{{\rm (ii)}}{\subseteq} \tau(W_n\omega_R,D), 
\end{align*}
where (i) follows from Proposition \ref{prop: first prop of C_n(D)-sub}(3) and 
(ii) holds by the fact that $\tau(W_n\omega_R,D)$ is $C_n(D)$-stable. 
Thus (3) holds. 

Let us show (4). 
Since $D = a \div(f)$, we have $\rup{a(p^e-1)}\div(f) -(p^e-1)D \geq 0$, 
Hence we get  $f^{\rup{a(p^e-1)}} \in R(-(p^e-1)D)$, which shows that
\[
\PsiT^e_n(F^e_*[f^{\rup{a(p^e-1)}}] \cdot -) \in C^e_n(D). 
\]
Then 
the assertion (3) implies
\[
\sum_{e \geq 0} \PsiT^e_n(F^e_*[cf^{\rup{a(p^e-1)}}] \cdot M) \subseteq \tau(W_n\omega_R,D). 
\]
Let us prove the opposite inclusion. 
For 
every $e \geq 0$ and every $g \in R(- (p^e-1)D)$, 
we get 
\[
\mathrm{div}(g)+ \mathrm{div}(f) \geq  (p^e-1)D+\mathrm{div}(f) 
= (a (p^e-1) +1) \div(f) \geq  \rup{a(p^e-1)}\mathrm{div}(f), 
\]
which implies 
\[
\PsiT^{e, g}_{n}(F^e_*[cf]M)= 
\PsiT^e_{n}(F^e_*[cfg]M)\subseteq \PsiT_n^e(F^e_*[cf^{\rup{a(p^e-1)}}]M).
\]
As $cf \in (dt^2) \cap R^{\circ}$, 
we obtain 
\begin{align*}
    \tau(W_n\omega_R,D) 
    &\overset{{\rm (3)}}{=} \sum_{e \geq 0} \sum_{\PsiT^{e, g}_{n} \in C^e_n(D)} \PsiT^{e, g}_{n}(F^e_*[cf] \cdot M) \\
    &\subseteq \sum_{e \geq 0} \PsiT_n^e(F^e_*[cf^{\rup{a(p^e-1)}}]M).
\end{align*}
Thus (4) holds. 


Let us show (5). 
By a similar argument to the proofs of Proposition \ref{prop: test submod}(4)(5), 
it is easy to see that (a) implies (b) and (c). 
Hence it is enough to prove (a). 
By $c \in (f^{\rup{a}})$, 
it holds that 
\begin{eqnarray*}
\mathrm{div}(c^{p-1})+\mathrm{div}(f^{p\rup{a(p^e-1)}}) 
&\geq& (p-1)\rup{a} \div(f) + p\rup{a(p^e-1)}\div(f)\\
&\geq& (\rup{a(p-1)}+p\rup{a(p^e-1)})\div(f) \\
&\geq& \rup{a(p^{e+1}-1)} \div(f). 
\end{eqnarray*}
We then obtain
\begin{align*}
T^{e+1}_n(F^{e+1}_*([cf^{\rup{a(p^{e+1}-1)}}]\cdot\tau(W_n\omega_R, D)))
&\supseteq  T^{e+1}_n(F^{e+1}_*([c^pf^{p\rup{a(p^e-1)}}]\cdot \tau(W_n\omega_R, D))) \\
&=T^e_n(F^e_*([cf^{\rup{a(p^e-1)}}]\cdot T_n(F_*\tau(W_n\omega_R, D)))) \\
&\supseteq T^e_n(F^e_*([cf^{\rup{a(p^e-1)}}]\cdot \tau(W_n\omega_R, D))), 
\end{align*}
where the last inclusion follows from (1) and ($\sharp$) below. 
\begin{enumerate}
\item[($\sharp$)] $T_n(F_*\tau(W_n\omega_R, D))$ is co-small and $C_n(D)$-stable. 
\end{enumerate}

Let us show $(\sharp)$. 
Since it is easy to see that $T_n(F_*\tau(W_n\omega_R, D))$ is co-small, 
we only prove that $T_n(F_*\tau(W_n\omega_R, D))$ is $C_n(D)$-stable. 
Pick an integer $e \geq 0$ and $g \in R(-(p^e-1)D)$. 
By 
\[
g^p \in R(-p(p^e-1)D) \subseteq R(-(p^e-1)D), 
\]
we get 
\begin{eqnarray*}
T_n^e( F_*^e[g] \cdot F_*^eT_n(F_*\tau(W_n\omega_R, D)))
&=& T_n^{e+1}(F_*^{e+1} ([g^p] \tau(W_n\omega_R, D)))\\
&=& T_n(F_*( T_n^e (F_*^e([g^p] \tau(W_n\omega_R, D)) ))\\
&\subseteq & T_n(F_*(\tau(W_n\omega_R, D))). 
\end{eqnarray*}
Thus $(\sharp)$ holds.

Then (a) holds, which completes the proof of (5). 
\end{proof}

\begin{definition}\label{d-tau-W_nomega-non-eff}
We use Notation \ref{n-normal-non-local}. 
Let $D$ be a $\Q$-divisor on $X$.
Take $h \in R^{\circ}$ such that $D+\div(h)$ 
is effective. 
We define the $W_nR$-submodule $\tau(W_n\omega_R,D)$ of 
$W_n\omega_R \otimes_{W_nR} W_n(K(R))$ by
\[
\tau(W_n\omega_R,D):=[h^{-1}] \tau(W_n\omega_R,D+\div(h)) \subseteq 
W_n\omega_R \otimes_{W_nR} W_n(K(R)). 
\]
This 
is independent of the choice of $h \in R^{\circ}$ by Proposition \ref{prop: test submodule gen pair}(2). 
{\ccred Note that the natural map 
$W_n\omega_R \to W_n\omega_R \otimes_{W_nR} W_n(K(R))$ is injective, 
because $W_n\omega_R$ is $S_2$ and $W_n\omega_R \otimes_{W_nR} W_n(K(R))$ coincides with the stalk of at the generic point of $\Spec (W_nR)$.}
\end{definition}


\begin{proposition}\label{prop: test submod pair compat}
We use Notation \ref{n-normal-non-local}. 
Let $D$ be a  $\Q$-divisor on $X$.
\begin{enumerate}
    \item 
    If $S$ is a multiplicatively closed subset $S$ of $R$, then  it holds that  
    \[
    \tau(W_n\omega_R,D) \otimes_{W_nR} W_n(S^{-1}R) 
    \simeq \tau(W_n\omega_{S^{-1}R},j^*D), 
    \]
    where $j : \Spec S^{-1}R \to \Spec R$ denotes the induced morphism. 
    \item If $(R,\m)$ is a local ring, then it holds that  
    \[
    \tau(W_n\omega_R,D) \otimes_{W_nR} W_n\widehat{R} 
    \simeq \tau(W_n\omega_{\widehat{R}},\pi^*D),  
    \]
    where $\widehat R$ denotes the $\m$-adic completion of $R$ and $\pi : \Spec \widehat{R} \to \Spec R$ is the induced morphism. 
\end{enumerate}
\end{proposition}

\begin{proof}
Taking $h\in R^{\circ}$ such that $D+\div[h]$ is effective, 
we may assume that $D$ is effective (Definition \ref{d-tau-W_nomega-non-eff}). 

Let us show (1). 
Pick $d \in R^{\circ}$ satisfying $D \leq \mathrm{div}(d)$. 
Fix a test element $t \in R^{\circ}$ and 
$c \in (dt^2) \cap R^{\circ}$.
Note that the conditions in Notation \ref{n-normal-non-local} are preserved under taking localisation. 
We have
\begin{align*}
    \tau(W_n\omega_R,D) \otimes_{W_nR} W_n(S^{-1}R)
    &= \left(\sum_{e \geq 0} \sum_{\varphi \in C^e_n(D)} \varphi(F^e_*[c] \cdot W_n\omega_R)\right) \otimes_{W_nR} W_n(S^{-1}R) \\
    &= \sum_{e \geq 0} \sum_{\varphi \in C^e_n(D)} 
    (\varphi \otimes_{W_nR} W_n(S^{-1}R))(F^e_*[c] \cdot W_n\omega_{S^{-1}R}) \\
    &\subseteq \sum_{e \geq 0} \sum_{\psi \in C^e_n(j^*D)} \psi(F^e_*[c] \cdot W_n\omega_{S^{-1}R}) \\
    &=\tau(W_n\omega_{S^{-1}R},j^*D), 
\end{align*}
where the above inclusion holds because 
we have $\varphi \otimes_{W_nR} W_n(S^{-1}R) \in C^e_n(j^*D)$ 
for $\varphi \in  C^e_n(D)$. 
It is enough to show the opposite inclusion 
$\tau(W_n\omega_R,D) \otimes_{W_nR} W_n(S^{-1}R) \supseteq 
\tau(W_n\omega_{S^{-1}R},j^*D)$. 
Take $e \geq 0$, $z \in W_n\omega_{S^{-1}R}$, and $\psi \in C^e_n(D|_{U})$. 
It suffices to show $\psi(F_*^e[c] \cdot z) \in \tau(W_n\omega_R,D) \otimes_{W_nR} W_n(S^{-1}R)$. 
We have $s_1 z \in W_n\omega_R$ for some $s_1 \in S$. 
We can write $\psi=\PsiT^{e, g}_{n}$ for some 
\[
g \in \MO_{\Spec S^{-1}R}(-(p^e-1)j^*D) = j^*R(-(p^e-1)D) = 
R(-(p^e-1)D) \otimes_R S^{-1}R. 
\]
Then there exists $s_2 \in S$ such that 
$s_2g \in R(-(p^e-1)D)$.  
In particular, $T^{e, c_2g}_{n} \in C^e_n(D)$. 
Therefore, we have 
\[
[s_1s_2] \cdot \psi(F^e_*([c]\cdot z)) 
=T^{e, g}_{n}(F^e_*([c s_1^{p^e}s_2^{p^e}] \cdot z)) 
= T^{e, s_2g}_{n}(F^e_*([c]\cdot  [s_1^{p^e-1}s_2^{p^e-1}] \cdot (s_1z))) 
\]
\[
\in T^{e, s_2g}_{n}(F_*([c]\cdot  W_n\omega_R)) \subseteq 
\tau(W_n\omega_R,D). 
\]
Hence $\psi(F^e_*([c]\cdot z)) \in \tau(W_n\omega_R,D) \otimes_{W_nR} W_n(S^{-1}R)$. 
Thus (1) holds. 


Let us show (2).
The inclusion $\pi^*\tau(W_n\omega_R,D) \subseteq \tau(W_n\omega_{\widehat{R}},\pi^*D)$ follows from a similar argument to the one of (1).
Then 
we have the induced injective  $W_n\widehat{R}$-module homomorphism 
\[
\iota \colon \pi^*\tau(W_n\omega_R,D) \hookrightarrow \tau(W_n\omega_{\widehat{R}},\pi^*D)
\]
of finitely generated $W_n\widehat{R}$. 
By the Artin-Rees lemma and completeness, it is enough to show that $\iota$ is surjective after taking modulo $W_n(\m^{l})W_n\omega_{\widehat{R}} \cap \tau(W_n\omega_{\widehat{R}},\pi^*D)$ for every positive integer $l$.
If $\psi \in C^e_n(\pi^*D)$, then there exists $g \in \sO_{X'}(-(p^e-1)\pi^*D)$ such that $\psi=T^e_g$.
By $\sO_{X'}( -(p^e-1)\pi^*D) \simeq \pi^*R(-(p^e-1)D)$, there exists $g' \in R(-(p^e-1)D)$ such that $g-g' \in (\m^l)^{[p^e]}$. 
Then we have $\alpha:=[g]-[g'] \in W_n((\m^l)^{[p^e]}\widehat{R})$.
In particular, we have
\[
T^{e, \alpha}(F^e_*[c]W_n\omega_{\widehat{R}}) \subseteq W_n(\m^l)W_n\omega_{\widehat{R}}.
\]
Therefore, we have
\begin{align*}
\psi(F^e_*[c]\cdot W_n\omega_{\widehat{R}}) 
&= T^{e, g}_{n}(F^e_*[c]\cdot W_n\omega_{\widehat{R}}) \\
&\equiv T^{e, g'}_{n}(F^e_*[c]\cdot W_n\omega_{\widehat{R}}) \hspace{1cm} \mod W_n(\m^l)W_n\omega_{\widehat{R}} \cap \tau(W_n\omega_{\widehat{R}},\pi^*D) \\
&\subseteq \pi^*\tau(W_n\omega_R,D).
\end{align*}
Thus, we obtain the surjectivity of $\iota$ after taking modulo $W_n(\m^{l})W_n\omega_{\widehat{R}} \cap \tau(W_n\omega_{\widehat{R}},\pi^*D)$.
\qedhere


\end{proof}

\begin{proposition}\label{prop: fund maps}
We use Notation \ref{n-normal-non-local}. 
Let $D$ be a $\Q$-divisor on $X$. 
Then, for every integer $n \geq 2$, 
the $W_nR$-module homomorphism ${\cccred \mathbf{R}}^* \colon W_{n-1}\omega_R \to W_n\omega_R$  
induces the following $W_nR$-module homomorphism: 
\[
{\cccred \mathbf{R}}^* \colon \tau(W_{n-1}\omega_R,D) \to \tau(W_{n}\omega_R,D).
\]
\end{proposition}

\begin{proof}
We may assume that $D$ is effective.
It is enough to show that 
$({\cccred \mathbf{R}}^*)^{-1}(\tau(W_n\omega_R,D))$ is co-small and $C_{n-1}(D)$-stable (as these  imply $\tau(W_{n-1}\omega_R,D) \subseteq ({\cccred \mathbf{R}}^*)^{-1}(\tau(W_n\omega_R,D))$ 
by Proposition \ref{prop: test submodule gen pair}(1)). 
Since $\tau(W_n\omega_R,D)$ is co-small, 
we get $[c] \cdot W_n\omega_R \subseteq \tau(W_n\omega_R,D)$ for some $c \in R^{\circ}$, 
which implies 
\[
{\cccred \mathbf{R}}^*([c] \cdot W_{n-1}\omega_R)  = [c] \cdot {\cccred \mathbf{R}}^*(W_{n-1}\omega_R) \subseteq \tau(W_n\omega_R,D), 
\]
i.e., $[c] \cdot W_{n-1}\omega_R \subseteq ({\cccred \mathbf{R}}^*)^{-1}(\tau(W_n\omega_R,D))$. 
Thus $({\cccred \mathbf{R}}^*)^{-1}(\tau(W_n\omega_R,D))$ is co-small.

It suffices to show that $({\cccred \mathbf{R}}^*)^{-1}(\tau(W_n\omega_R,D))$ is $C_{n-1}(D)$-stable. 
Take $e \geq 0$ and $T^{e, g}_{n-1} \in C^e_{n-1}(D)$, where $g \in R(-(p^e-1)D)$. 
We have 
\[
{\cccred \mathbf{R}}^* \circ T^{e, g}_{n-1} = T^{e, g}_{n} \circ {\cccred \mathbf{R}}^{*}
\]
for $T^{e, g}_{n} \in C^e_n(D)$, 
because we have $T^{e, g}_{m}(-) = \PsiT^e_{m}(F^e_*[g]\cdot-)$ for $m \in \{n-1, n\}$. 
Therefore, we get 
\begin{align*}
    {\cccred \mathbf{R}}^* \circ T^{e, g}_{n-1}( F_*^e({\cccred \mathbf{R}}^*)^{-1}(\tau(W_n\omega_R,D))) 
    &=
    T^{e, g}_{n} \circ {\cccred \mathbf{R}}^*  ( ({\cccred \mathbf{R}}^*)^{-1}(F_*^e\tau(W_n\omega_R,D))) \\
    &\subseteq  T^{e, g}_{n}(F^e_*\tau(W_{n}\omega_R,D)) \subseteq \tau(W_{n}\omega_R,D), 
\end{align*}
where the last inclusion follows from the fact that 
$\tau(W_{n}\omega_R,D)$ is $C_n(D)$-stable 
(Proposition \ref{prop: test submodule gen pair}(1)). 
Hence $({\cccred \mathbf{R}}^*)^{-1}(\tau(W_n\omega_R,D))$ is $C_{n-1}(D)$-stable. 
\qedhere 



\end{proof}

\begin{definition}\label{d-taun-omega-log}
We use Notation \ref{n-normal-non-local}. 
Let $D$ be a $\Q$-divisor on $X$.
\begin{enumerate}
    \item 
    We set 
    \[
    \tau_n(\omega_R,D) := ( ({\cccred \mathbf{R}}^{n-1})^*)^{-1}(\tau(W_n\omega_R,D)), 
    \]
    which is the  $R$-submodule of $\omega_R( -D)$ 
    obtained as the inverse image of $\tau(W_n\omega_R,D)$     by 
    \[
    ({\cccred \mathbf{R}}^{n-1})^* : \omega_R( -D) \to W_n\omega_R(-D). 
    \]
    We call $\tau_n(\omega_R, D)$ the {\em $n$-quasi-test $R$-submodule  of $\omega_R$ along $D$}. 
We have the following inclusions (Remark \ref{r-taun-omega-log}): 
\begin{equation}\label{e1-taun-omega-log}
\tau_1(\omega_{{\ccred R}}, D) \subseteq 
\tau_2(\omega_{{\ccred R}}, D) \subseteq \cdots \subseteq 
\tau_n(\omega_{{\ccred R}}, D) \subseteq 
\tau_{n+1}(\omega_{{\ccred R}}, D) \subseteq \cdots \subseteq \omega_R(- D).  
\end{equation}
    \item 
    We set $\tau^q(\omega_R,D) :=\sum_{n=1}^{\infty} \tau_n(\omega_R,D)$. 
    Since (\ref{e1-taun-omega-log}) is an ascending chain of $R$-submodules of     $\omega_R(- D)$ , there exists $n_0>0$ such that 
    \[
    \tau^q(\omega_R,D)  =\sum_{n=1}^{\infty} \tau_n(\omega_R,D) 
    = {\ccred \bigcup_{n=1}^{\infty}\tau_n(\omega_R,D)} = \tau_n(\omega_R, D) \subseteq  \omega_R(- D)
    \]
    for every integer $n \geq n_0$. 
    In particular, $\tau^q(\omega_R,D)$ an $R$-submodule of $\omega_R$. 
    We call $\tau^q(\omega_R, D)$ the {\em quasi-test $R$-submodule of $\omega_R$ along $D$}. 
\end{enumerate}
\end{definition}

\begin{remark}\label{r-taun-omega-log}
We use the same notation as in Definition \ref{d-taun-omega-log}. 
We have the following commutative diagram in which each horizontal sequence is exact: 
\[
\begin{tikzcd}
0 \arrow[r] & 
\omega_R(-D) \arrow[r, "({\cccred \mathbf{R}}^{n-1})^*"] \arrow[d, equal]&
W_n\omega_R(-D) \arrow[r, "{\cccred \mathbf{V}}^*"] \arrow[d, "{\cccred \mathbf{R}}^*"]&
F_*W_{n-1}\omega_R(-pD) \arrow[d, "{\cccred \mathbf{R}}^*"]\\
0 \arrow[r] & 
\omega_R(-D) \arrow[r, "({\cccred \mathbf{R}}^{n})^*"] &
W_{n+1}\omega_R(-D) \arrow[r, "{\cccred \mathbf{V}}^*"] &
F_*W_{n}\omega_R(-pD). 
\end{tikzcd}
\]
This diagram, together with Proposition \ref{prop: fund maps}, induces the following commutative diagram 
in which each horizontal sequence is exact: 
    \[
    \xymatrix{
    0 \ar[r] & \tau_n(\omega_R,D) \ar[r]^-{({\cccred \mathbf{R}}^{n-1})^*} &
    \tau(W_n\omega_R,D) \ar[r]^-{{\cccred \mathbf{V}}^*} \ar[d]^-{{\cccred \mathbf{R}}^*} & 
    F_*W_{n-1}\omega_R(-pD) 
    \ar[d]^-{{\cccred \mathbf{R}}^*} \\
    0 \ar[r] & \tau_{n+1}(\omega_R,D) \ar[r]^-{({\cccred \mathbf{R}}^{n})^*} & \tau(W_{n+1}\omega_R,D) \ar[r]^-{{\cccred \mathbf{V}}^*} & 
        F_*W_{n}\omega_R(-pD).  
    }
    \]
Therefore, we obtain an inclusion $\tau_n(\omega_R,D) \subseteq \tau_{n+1} (\omega_R,D)$.
\end{remark}

\begin{proposition}\label{prop: stability of tau omega}
We use Notation \ref{n-normal-non-local}. 
Take a $\Q$-Cartier $\Q$-divisor $D$ and 
an effective $\Q$-divisor $E$ on $X$.
Then the following hold. 
\begin{enumerate}
\item 
For every integer $n \geq 1$, 
there exists a rational number $\epsilon_n>0$ such that 
\[
\tau(W_n\omega_R,D)=\tau(W_n\omega_R,D+\epsilon_n E)
\]
\item 
There exists a rational number $\epsilon>0$ such that 
\[
\tau_n(\omega_R,D)=\tau_n(\omega_R,D+\epsilon E) \qquad \text{and}\qquad \tau^q(\omega_R,D)=\tau^q(\omega_R,D+\epsilon E)
\]
for every integer $n>0$. 
\end{enumerate}
\end{proposition}

\begin{proof}
We first reduce the problem to the case when   
\begin{enumerate}
\renewcommand{\labelenumi}{(\roman{enumi})}
\item $D=a\div(f)$ for some $a \in \Q_{>0}$ and $f \in R^\circ$, 
\item $E=\mathrm{div}(g)$ for some $g \in R^{\circ}$, and  
\item $a<1$. 
\end{enumerate}
Enlarging $E$, we may assume that $E$ is an effective Cartier divisor. 
Taking a suitable affine open cover of $X$, we may assume that (i) and (ii) hold. 
By $D = a\div(f) = \frac{a}{r}\div(f^r)$ for every $r \in \Z_{>0}$, 
we may assume (iii). 

Let us show (1). Fix $n \in \Z_{>0}$. 
For $\epsilon \in \Q_{>0}$, the inclusion 
\[
\tau(W_n\omega_R,D+\epsilon E) \subseteq \tau(W_n\omega_R,D)
\]
is clear by definition. 
Pick a test element $t \in R^{\circ}$ and $c \in R^\circ \cap (fgt^2) \cap (f^{\rup{a}})$. 
By {applying} Proposition \ref{prop: test submodule gen pair}(5)(c) {for $L = W_n\omega_R$}, 
there exists an integer $e >0$ (depending on $n$) such that 
\[
\tau(W_n\omega_R,D)= 
\PsiT^e_n(F^e_*[cg][f^{\rup{a(p^e-1)}}]W_n\omega_R).
\]
Set 
$\epsilon_n := 
 \frac{1}{p^{e}-1}$. 
We then get 
\begin{align*}
\tau(W_n\omega_R,D) &= \PsiT^e_n(F^e_*[cg][f^{\rup{a(p^e-1)}}]W_n\omega_R) \\[0.3em]
&= \PsiT^e_n(F^e_*[c][f^{\rup{a(p^e-1)}}g^{\rup{\epsilon_n(p^e-1)}}]W_n\omega_R)\\[-0.2em]
&\overset{{\small(\star)}}{\subseteq}  \tau(W_n\omega_R,D+\epsilon_n E),
\end{align*}
where $(\star)$ can be proven  as follows. 
By Proposition \ref{prop: test submodule gen pair}(3), 
it is enough to show 
\[
f^{\rup{a(p^e-1)}}g^{\rup{\epsilon_n(p^e-1)}} \in R( -(p^e-1)(D+\epsilon_n E)), 
\]
that is, $\div(f^{\rup{a(p^e-1)}}g^{\rup{\epsilon_n(p^e-1)}}) \geq (p^e-1)(D+\epsilon_n E)$, 
which is assured by 
\[
\rup{a(p^e-1)} \div(f) \geq (p^e-1)D\quad \text{and}\quad 
\rup{\epsilon_n(p^e-1)} \div(g)  = \div(g) = \epsilon_n  (p^e-1) E. 
\]
Thus (1) holds.


Let us show (2).
For every $\epsilon>0$, the inclusion 
$\tau^q(\omega_R,D+\epsilon E) \subseteq \tau^q(\omega_R,D)$ is clear. 
Fix an integer $n_0>0$ satisfying  
$\tau^q(\omega_R,D)=\tau_{n_0}(\omega_R,D)$ (Definition \ref{d-taun-omega-log}). 
By (1) {and Definition \ref{d-taun-omega-log}},  there exists $\epsilon_0 >0$ (depending on $n_0$) such that $\tau_{n_0}(\omega_R,D+\epsilon_0E)=\tau_{n_0}(\omega_R, D)$.
Then we get 
\[
\tau^q(\omega_R,D+\epsilon_0 E) \supseteq \tau_{n_0}(\omega_R,D+\epsilon_0 E)=\tau_{n_0}(\omega_R,D)=\tau^q(\omega_R, D). 
\]
This completes the proof of $\tau^q(\omega_R, D) =\tau^q(\omega_R, D +\epsilon_0 E)$. 

Again by Definition \ref{d-taun-omega-log}, we can find $n_1>0$ such that 
\[
\tau_n(\omega_R, D)=
\tau^q(\omega_R, D) =\tau^q(\omega_R, D +\epsilon_0 E) =\tau_n(\omega_R, D +\epsilon_0 E) 
\]
for every $n > n_1$. 
Set 
\[
\epsilon := \min \{ \epsilon_0, \epsilon_1, ..., \epsilon_{n_1}\}, 
\]
where each of $\epsilon_1, ..., \epsilon_{n_1}$ is as in (1). 
Then, for every $1 \leq n \leq n_1$, 
we get $\tau_n(\omega_R, D)=\tau_n(\omega_R, D +\epsilon E)$  by (1). 
Thus (2) holds. 
\end{proof}


{At this point, we do not know yet if the condition that $\tau_n(\omega_R, {K_R+\Delta}
) = R$ is equivalent to $(R,\Delta)$ being {$n$-}quasi-$F$-regular, when $K_R+\Delta$ is $\bQ$-Cartier. In order to verify this, we need to develop the theory of quasi-tight closure {\cred for pairs}.}

\subsection{Quasi-tight closure in local cohomologies (pairs)}\label{ss log tight closure}


\begin{notation}\label{n-normal-local}
Let $(R, \m)$ be an $F$-finite Noetherian normal local ring of characteristic $p>0$. 
Set $d:= \dim R$ and $X:=\Spec R$.
\end{notation}

We use Notation \ref{n-normal-local}. 
Take $c \in R^{\circ}$, a $\Q$-divisor $D$ on $X$, and 
integers $n >0$ and $e \geq 0$. 
We define the $W_nR$-modules $Q^e_{R, D, n}$ and $Q^{e, c}_{R, D, n}$
and the $W_nR$-module homomorphisms $\Phi^e_{R, D, n}$ and $\Phi^{e, c}_{R, 
D,n}$ 
by the following diagram in which all the squares are pushouts: 
\begin{equation}\label{e-def-Q^{e,c}_D}
\begin{tikzcd}
W_nR(D) \arrow[r, "F^e"] \arrow[d, "{\cccred \mathbf{R}}^{n-1}"'] & F_*^e(W_nR(p^eD)) \arrow[d] \arrow[r, "{ \cdot F_*^e[c]}"] & F_*^e(W_nR(p^eD)) \arrow[d]\\
R(D) \arrow[r, "\Phi^e_{R, D, n}"] \arrow[rr, "\Phi^{e, c}_{R, D, n}"', bend right] & Q^e_{R, D, n} \arrow[r] & Q^{e, c}_{R, D, n}. 
\end{tikzcd}
\end{equation}

\begin{definition}\label{d-q-F-rat-log}
We use Notation \ref{n-normal-local}. 
Take integers $n >0$ and $e \geq 0$. 
Let $D$ be a $\Q$-divisor on $X$. 
\begin{enumerate}
    \item 
    For $c \in R^{\circ}$, we set 
    \[
    \hspace{3.3em} \wt{K^{e, c}_{D, n}} := \Ker \left(
    H^d_\m(W_nR(D)) \xrightarrow{F^e} H^d_\m(F^e_*W_nR(p^eD)) \xrightarrow{\cdot F^e_*[c]} H^d_\m(F^e_*W_nR(p^eD))\right), 
    \]
which is a $W_nR$-submodule of $H^d_\m(W_nR(D))$. 
Then we define 
$\wt{0_{D, n}^{*}}$ by 
\[
\wt{0_{D, n}^{*}} := \bigcup_{\substack{c \in R^{\circ},\\ e_0 \in \Z_{>0}}} \bigcap_{e \geq e_0} \widetilde{K^{e, c}_{D, n}} 
\subseteq H^d_\m(W_nR(D)). 
\]
Equivalently, given $z \in H^d_\m(W_nR(D))$, 
    we have $z \in \wt{0_{D, n}^{*}}$ if and only if 
    there exist $c \in R^{\circ}$ and $e_0 \in \Z_{>0}$ such that 
    $z$ is contained in $\wt{K^{e, c}_{D, n}}$ for every integer $e \geq e_0$. 
    Note that $\wt{0_{D, n}^{*}}$ is a $W_nR$-submodule of $H^d_\m(W_nR(D))$ 
    (cf.\ Remark \ref{r-q-F-rat}). 
    \item 
   For $c \in R^{\circ}$, we set 
    \[
    K^{e, c}_{D, n} := \Ker\left( H^d_\m(R(D)) \xrightarrow{\Phi^{e, c}_{R,D,n}} H^d_\m(Q^{e, c}_{R,D,n})\right),
    \]
    which is an  $R$-submodule of $H^d_\m(R(D))$. 
    We define the {\em $n$-quasi-tight closure} $0_{D,n}^{*}$ by 
    \[
0_{D, n}^{*} := \bigcup_{\substack{c \in R^{\circ},\\ e_0 \in \Z_{>0}}} \bigcap_{e \geq e_0} K^{e, c}_{D, n} \subseteq H^d_\m(R(D)). 
\]
Equivalently, given $z \in H^d_\m(R)$, 
we have $z \in 0_{D, n}^{*}$ if and only if 
there exist $c \in R^{\circ}$ and $e_0 \in \Z_{>0}$ 
such that $z$ is contained in $K^{e, c}_{D, n}$ for every integer $e \geq e_0$. 
Note that $0_{D, n}^{*}$ is an $R$-submodule of $H^d_\m(R(D))$ 
    (cf.\ Remark \ref{r-q-F-rat}). 
\end{enumerate}
\end{definition}


\begin{remark}[{cf.\ Remark \ref{remark: easy remark}}]\label{remark: easy remark-log}
We use Notation \ref{n-normal-local}. 
Take a $\Q$-divisor $D$ on $X$. 
\begin{enumerate}
    \item 
    Consider the commutative diagram
    \[
    \begin{tikzcd}
H^d_\m(F_*W_{n-1}R) \arrow[r, equal] \arrow[d, "{\cccred \mathbf{V}}"] & 
H^d_\m(F_*W_{n-1}R) \arrow[r, equal] \arrow[d, "{F^e{\cccred \mathbf{V}}}"] & H^d_\m(F_*W_{n-1}R) \arrow[d, "{(\cdot F^e_*[c]) \circ F^e{\cccred \mathbf{V}}}"] \\
H^d_\m(W_nR(D)) \arrow[r, "F^e"] \arrow[d, twoheadrightarrow, "{\cccred \mathbf{R}}^{n-1}"] & 
H^d_\m(F^e_*W_nR(p^eD)) \arrow[r, "{\cdot F^e_*[c]}"] \arrow[d] & 
H^d_\m(F^e_*W_nR(p^eD)) \arrow[d] \\
H^d_\m(R(D)) \arrow[r] & H^d_\m(Q^e_{R,D, n}) \arrow[r] & H^d_\m(Q^{e,c}_{R,D,n})
\end{tikzcd}
    \]
    where the vertical sequences are {\cccred exact, as in Remark~\ref{remark: easy remark}}.
    {\cccred By diagram chase similar to that in Remark \ref{remark: easy remark},} we get 
    \[
    {\cccred \mathbf{R}}^{n-1}(\wt{K^{e, c}_{D, n}})=K^{e, c}_{D, n}.
    \]
    
    \item We have 
    \[
    H^d_\m(R(D)) \supseteq 0_{D, 1}^* \supseteq  0_{D,2}^* \supseteq 0_{D, 3}^* \supseteq  \cdots. 
    \]
    Indeed, this holds by the following factorisation: 
    \[
    H^d_\m(R(D)) \to H^d_\m(Q^{e, c}_{R,D, n+1}) \to H^d_\m(Q^{e, c}_{R,D, n}).
    \]
\item 
By Proposition \ref{p wt0 vs tight} and Definition \ref{d-q-F-rat-log}, it holds that 
\[
\wt{0^*_{D,1}}=0^*_{D, 1} = 0^{*\{D\}}_{H^d_{\m}(R(D))}. 
\]
\end{enumerate}
\end{remark}

\begin{proposition}\label{prop: tight cl l=n-log}
We use Notation \ref{n-normal-local}. 
Take a  $\Q$-Cartier $\Q$-divisor $D$ on $X$. 
Then the following hold. 
\begin{enumerate}
    \item ${\cccred \mathbf{R}}(\wt{0_{D, n+1}^{*}}) \subseteq \wt{0_{D, n}^{*}}$. 
    \item ${\cccred \mathbf{V}}^{-1}(\wt{0_{D, n}^{*}}) = F_*\wt{0_{pD, n-1}^{*}}$.
\end{enumerate}
\end{proposition}
\noindent {Statement (1) says that the first commutative diagram below exists, whilst Statement (2) implies that the second diagram exists and is Cartesian:
\[
\begin{tikzcd}
H^d_\m(W_{n+1}R(D)) \ar{r}{{\cccred \mathbf{R}}} & H^d_\m(W_{n}R(D))\\
\wt{0_{D,n+1}^{*}} \ar[hook]{u} \ar{r}{{\cccred \mathbf{R}}} & \wt{0_{D,n}^{*}}, \ar[hook]{u}
\end{tikzcd}
\qquad
\begin{tikzcd}
H^d_\m(F_*W_{n-1}R(pD)) \ar{r}{{\cccred \mathbf{V}}} & H^d_\m(W_{n}R(D))\\
F_*\wt{0_{pD,n-1}^{*}} \ar[hook]{u} \ar{r}{{\cccred \mathbf{V}}} & \wt{0_{D,n}^{*}}. \ar[hook]{u}
\end{tikzcd}
\]
}
\begin{proof}
By Definition \ref{d-q-F-rat-log}, (1) holds. 
Let us show (2). 
By Definition \ref{d-q-F-rat-log}, we get 
\[
\wt{0_{D, n}^{*}} \supseteq {\cccred \mathbf{V}}(F_*\wt{0_{pD, n-1}^{*}}),
\]
that is,  ${\cccred \mathbf{V}}^{-1}(\wt{0_{D, n}^{*}}) \supseteq F_*\wt{0_{pD, n-1}^{*}}$. 
In what follows, we prove the opposite inclusion 
\[
{\cccred \mathbf{V}}^{-1}(\wt{0_{D, n}^{*}}) \subseteq F_*\wt{0_{D, n-1}^{*}}.
\]
Take $F_*z \in {\cccred \mathbf{V}}^{-1}(\wt{0_{D, n}^{*}})$, that is, 
${\cccred \mathbf{V}}(F_*z) \in \wt{0_{D, n}^{*}}$. 
Then there exist $c \in R^{\circ}$ and $e_0>0$ such that 
${\cccred \mathbf{V}}(F_*z) \in \bigcap_{e \geq e_0} \wt{K^{e, c}_{D, n}}$. 
Since $R$ is a local ring and $D$ is a $\Q$-Cartier $\Q$-divisor on $\Spec R$, 
there are finitely many possibilities for $R(p^eD)$, {up to isomorphisms,} when $e \geq 0$. 
Hence we can find $c' \in R^{\circ}$ such that $c' \cdot H^{d-1}_{\m}(R(p^eD))=0$ for every $e \geq 0$ (Lemma \ref{lem: ann local coh}), which implies that
\[
(F_*^e c')  \cdot H^{d-1}_\m(F_*^eR(p^eD))=(F_*^e c')  \cdot F_*^eH^{d-1}_\m(R(p^eD))=F_*^e(c' \cdot H^{d-1}_\m(R(p^eD)))=0
\]
for every $e \geq 0$.  
Given a $\Q$-divisor $E$, we have the exact sequence 
\[
H^{d-1}_{\m}(R(E)) \to   H^d_{\m}(F_*W_{n-1}R(pE)) \xrightarrow{{\cccred \mathbf{V}}}  H^d_{\m}(W_nR(E)).  
\]
Hence we get  the following commutative diagram in which each horizontal sequence is exact: 
\[
\begin{tikzcd}
H^{d-1}_{\m}(R(D)) \arrow[r] \arrow[d, "F^e"] & 
{\cccred \mathbf{V}}^{-1}(\wt{0_{D, n}^{*}}) \arrow[r, "{\cccred \mathbf{V}}"]  \arrow[d, "F^e"]& 
\wt{0^*_{D, n}} \arrow[d, "F^e"]\\
H^{d-1}_{\m}(F_*^eR(p^eD)) \arrow[r] \arrow[d, "{\cdot F^e_*[c]}"] & 
F_*^e{\cccred \mathbf{V}}^{-1}(\wt{0_{p^eD, n}^{*}}) \arrow[r, "{\cccred \mathbf{V}}"]  \arrow[d, "{\cdot F^e_*[c]}"]& 
F_*^e\wt{0^*_{p^eD, n}} \arrow[d, "{\cdot F^e_*[c]}"]\\
H^{d-1}_{\m}(F_*^eR(p^eD)) \arrow[r] \arrow[d, "{\cdot F^e_*[c']=0}"] & 
F_*^e{\cccred \mathbf{V}}^{-1}(\wt{0_{p^eD, n}^{*}}) \arrow[r, "{\cccred \mathbf{V}}"]  \arrow[d, "{\cdot F^e_*[c']}"]& 
F_*^e\wt{0^*_{p^eD, n}} \arrow[d, "{\cdot F^e_*[c']}"]\\
H^{d-1}_{\m}(F_*^eR(p^eD)) \arrow[r] & 
F_*^e{\cccred \mathbf{V}}^{-1}(\wt{0_{p^eD, n}^{*}}) \arrow[r, "{\cccred \mathbf{V}}"] & 
F_*^e\wt{0^*_{p^eD, n}}. 
\end{tikzcd}
\]
By the diagram chase starting with $F_*z \in {\cccred \mathbf{V}}^{-1}(\wt{0_{D, n}^{*}})$ (located in the centre at the top), we get $(F_*^e[cc']) \circ F^e (F_*z)=0$ for every $e \geq e_0$, and hence 
$F_*z \in F_*\wt{0^*_{pD, n-1}}$. 
Thus (2) holds. 
\end{proof}

The choice of $t \in \tau(R,\{p^eD\})$ is {motivated} by Proposition \ref{p wt0 vs tight} and Definition \ref{d-test-ideal}.

\begin{proposition}\label{prop: ker V and ker F: D}
We use {Notation \ref{n-normal-local}}. 
Take  a $\Q$-Cartier $\Q$-divisor $D$ on $X$ 
and  $t \in R^\circ$ such that $t$ is contained in 
$\tau(R,\{p^eD\})$ for every $e \geq 0$.
Fix integers $n>0$ and $m \geq 0$. 
Then the following hold. 
\begin{enumerate}
\item $[t^{2}]\cdot \wt{0_{D, n}^{*}}=0$. 
\item 
$[t^2]\cdot \mathrm{Ker}(F^m \colon H^d_\m(W_nR(D)) \to H^d_\m({\ccred F^m_*}W_nR(p^{m}D)))=0.$
\end{enumerate}
\end{proposition}

\begin{proof}
We have the following inclusion (Definition \ref{d-q-F-rat-log}(1)): 
\[
 \mathrm{Ker}(F^m \colon H^d_\m(W_nR(D)) \to H^d_\m(F^m_*W_nR(p^mD)))\subseteq \widetilde{0_{D, n}^*}. 
\]
Hence (1) implies (2). 

Thus, it is enough to prove (1). 
Specifically, we will show that 
$[t^2] \cdot \wt{0_{D, n}^{*}}=0$ by an increasing induction on $n$ {(for $D$ replaced by any $p^k$-th multiple of $D$)}. 
The base case $n=1$ holds  by 
\[
[t] \cdot \wt{0_{D, 1}^{*}}=[t] \cdot 0_{D, 1}^{*}=t \cdot 0^{*\{D\}}_{H^d_\m(R(D))}=0
\]
(Remark \ref{remark: easy remark-log}(3)). 
Now assume $n \geq 2$ and consider the exact sequence
\begin{equation}\label{e2: ker V and kerF :D}
H^d_\m(F_*W_{n-1}R(pD)) \xrightarrow{{\cccred \mathbf{V}}} H^d_\m(W_nR(D)) \xrightarrow{{\cccred \mathbf{R}}^{n-1}} H^d_\m(R(D)).
\end{equation}
Fix $z \in \wt{0_{D, n}^{*}}$. 
Since ${\cccred \mathbf{R}}^{n-1}(z) \in \wt{0_{D, 1}^{*}}$ (Proposition \ref{prop: tight cl l=n-log}(1)), 
we get that ${\cccred \mathbf{R}}^{n-1}([t]z) = [t]{\cccred \mathbf{R}}^{n-1}(z) =0$.
Then the exact sequence (\ref{e2: ker V and kerF :D}) enables us to find 
\[
F_*z' \in F_*H^d_\m(W_{n-1}R(pD))= H^d_\m(F_*W_{n-1}R(pD))
\]
satisfying ${\cccred \mathbf{V}}(F_*z')=[t]z$. 
As $F_*z' \in {\cccred \mathbf{V}}^{-1}(\wt{0^*_{D, n}}) = F_*\wt{0^*_{pD, n-1}}$ (Proposition \ref{prop: tight cl l=n-log}(2)), 
we get $z' \in \wt{0^*_{pD, n-1}}$, and so, by the induction hypothesis, $[t^2]z' = 0$. 
We thus obtain 
\[
[t^2]z=[t]{\cccred \mathbf{V}}(F_*z')={\cccred \mathbf{V}}([t]F_*z')={\cccred \mathbf{V}}(F_*([t^p]z'))=0. 
\]
Hence (1) holds.
\end{proof}

\begin{proposition}\label{prop: ker V and ker F: pair}
We use Notation \ref{n-normal-local}. 
Fix $t \in \tau(R) \cap R^{\circ}$ and an integer $n>0$. 
Let $D$ be a  $\Q$-Cartier $\Q$-divisor on $X$. 
Assume that there exist 
$g \in R^{\circ}$ and a $\Q$-Cartier Weil divisor $D'$ on $X$ 
such that 
$D \leq D' \leq D+\mathrm{div}(g)$.
Then the following hold. 
\begin{enumerate}
\item $[t^{2}g]\cdot \wt{0_{D, n}^{*}}=0$. 
\item For every integer $m \geq 0$, 
\[
[t^2g]\cdot \mathrm{Ker}(F^m \colon H^d_\m(W_nR(D)) \to H^d_\m(W_nR(p^{m}D)))=0.  
\]
\end{enumerate}
\end{proposition}

\begin{proof}
We have the following commutative diagram 
\[
\begin{tikzcd}
    H^d_\m(W_nR(D)) \arrow[r, "\iota"] \arrow[d, "F^e"]  & H^d_\m(W_nR(D')) \arrow[d, "F^e"] \\
H^d_\m(F^e_*W_nR(p^eD)) \arrow[r, "\iota'"] & H^d_\m(F^e_*W_nR(p^eD')),
\end{tikzcd}
\]
where $\iota$ and $\iota'$ are the $W_nR$-module homomorphisms 
induced by the natural inclusions 
$W_nR(D) \hookrightarrow W_nR(D')$ and 
$W_nR(p^eD) \hookrightarrow W_nR(p^eD')$, respectively.
Thus $\iota(\wt{0^{*}_{D, n}}) \subseteq \wt{0^{*}_{D', n}}$.
By Proposition \ref{prop: ker V and ker F: D}(1), we have 
$[t^2] \cdot \wt{0^{*}_{D', n}} =0$, 
which implies $[t^2] \cdot \wt{0^{*}_{D, n}} \subseteq \mathrm{Ker}(\iota)$. 
By the factorisation: 
\[
\cdot [g]: W_nR(D) \hookrightarrow W_nR(D') \hookrightarrow W_nR(D+ \div(g))
\xrightarrow{\simeq, \cdot [g] }W_nR(D), 
\]
we get $[gt^2] \cdot \wt{0^{*}_{D, n}}=0$. 
Thus (1) holds. 
The assertion (2) follows from (1) and 
$\mathrm{Ker}(F^m) \subseteq \wt{0^{*}_{D, n}}$ (Definition \ref{d-q-F-rat-log}). 
\end{proof}

\begin{proposition}\label{prop: test elem pair tilda}
We use Notation \ref{n-normal-local}. 
Let $D$ be a $\Q$-Cartier $\Q$-divisor on $X$ 
such that $D=a \cdot \mathrm{div}(g)$ for some  $a \in \Q$ and $g \in R^{\circ}$. 
Take a test element $t \in R^{\circ}$  and 
$t_D \in R^\circ$ satisfying  $t_D \in \tau(R,\{p^eD\})$ for every integer $e\geq 0$. 
Fix $c \in R^{\circ} \cap ((g^2t^4) \cup (t_D^4))$. 
Then 
the following hold. 
\begin{enumerate}
    \item 
    $H^d_{\m}(W_nR(D)) \supseteq \wt{K^{0,c}_{D,n}}
    \supseteq \wt{K^{1,c}_{D,n}} \supseteq \cdots \supseteq \wt{K^{e,c}_{D,n}}
    \supseteq \wt{K^{e+1,c}_{D,n}} \supseteq \cdots.$ 
    \item 
    There exists an integer $e_0 \geq 0$ such that 
    \[
\wt{0_{D, n}^{*}}=\wt{K^{e,c}_{D, n}}
\]
for every integer $e \geq e_0$. 

\end{enumerate}
\end{proposition}

\begin{proof}
We first treat the case when  $c \in R^\circ \cap (g^2t^4)$. 
In this case, we have $c=c'g^2t^4$ for some $c' \in R^{\circ}$.

Let us show (1). 
Set $D_e:=\rup{ap^e}  \mathrm{div}(g)$ and take $z \in \wt{K^{e+1,c}_{D, n}}$. 
Since we have $cc_1 = c'^pg^pt^{2p}$ for $c_1 := c'^{p-1}g^{p-2}t^{2p-4} \in R^\circ$, 
we get 
\[
0 = F^{e+1}_*[cc_1] \cdot F^{e+1}(z)=
F_*^{e+1}[c'^pg^pt^{2p}] \cdot F^{e+1}(z)=
F(F_*^e[c'gt^2] \cdot F^{e}(z)). 
\]
As $z \in H^d_{\m}(R(D))$, we have  
\[
F_*^e[c'gt^2] \cdot F^{e}(z) \in 
F_*^e[c'gt^2] \cdot H^d_{\m}(F_*^eR(p^eD)) \subseteq 
H^d_{\m}(F_*^eR(p^eD)). 
\]
By $p^eD \leq D_e \leq p^eD+\mathrm{div}(g)$,  
we may apply Proposition \ref{prop: ker V and ker F: pair}(1), 
and hence  $[t^2g] \cdot \wt{0^*_{{p^e}D, n}}=0$. 
Then 
\[
0 = F_*^e[t^2g] \cdot F_*^e[c'gt^2] \cdot F^{e}(z) =F_*^e[c'g^2t^4] \cdot F^{e}(z) =
F_*^e[c] \cdot F^e(z). 
\]
Therefore, $z \in \wt{K^{e,c}_{D, n}}$. Thus (1) holds.

Let us show (2). 
Since $H^d_\m(W_nR(D))$ is an Artinian $W_nR$-module, 
(1) enables us to find an integer $e_0 >0$ such that $\wt{K^{e, c}_{D, n}} = \wt{K^{e_0, c}_{D, n}}$ 
for every $e \geq e_0$. 
Definition \ref{d-q-F-rat-log}(1) implies  
\[
\wt{0^{*}_{D, n}} \supseteq \bigcap_{e \geq 0}\wt{K^{e, c}_{D, n}} 
= \wt{K^{e_0, c}_{D, n}}. 
\]
Conversely, if $z \in \wt{0^{*}_{D, n}}$ and $e \geq 0$, then we get 
\begin{align*}
F_*^e[c] \cdot F^e(z) 
&\in F_*^e[c] \cdot F^e(\wt{0^{*}_{D, n}}) \subseteq 
F_*^e[c] \cdot F_*^e\wt{0^{*}_{p^eD, n}} \\
& =F_*^e(  [c'g^2t^4] \cdot \wt{0^{*}_{p^eD, n}} )
= F_*^e( [c't^2] \cdot [gt^2] \cdot \wt{0^{*}_{p^eD, n}})=0,
\end{align*}
where the last equality follows from Proposition \ref{prop: ker V and ker F: pair}(1), 
which is applicable by $p^eD \leq D_e \leq p^eD + \div(g)$. 
Therefore, $z \in \bigcap_{e \geq 0} \widetilde{K^{e, c}_{D, n}} = \widetilde{K^{e_0, c}_{D, n}}$. Thus (2) holds. 


As for the case when 
$c \in R^\circ \cap (t_D^4)$, 
we may apply the same argument as above 
by using 
 Proposition \ref{prop: ker V and ker F: D} 
instead of Proposition \ref{prop: ker V and ker F: pair}. 
\end{proof}

\begin{theorem}\label{thm: lifting thm pair}
We use Notation \ref{n-normal-local}. 
Let $D$ be a $\Q$-Cartier $\Q$-divisor on $X$.
Then 
\[
{\cccred \mathbf{R}}^{n-1}(\wt{0^{*}_{D, n}})=0^{*}_{D, n}.
\]
\end{theorem}

\begin{proof}
First, we prove that ${\cccred \mathbf{R}}^{n-1}(\wt{0^{*}_{D, n}}) \subseteq 0^{*}_{D, n}$.
Take $z \in \wt{0^{*}_{D, n}}$. 
Then there exist $e_0 \geq 0$ and $c \in R^{\circ}$ such that 
$z \in \bigcap_{e \geq e_0} \wt{K^{e,c}_{D,n}}$. 
We get 
\[
{\cccred \mathbf{R}}^{n-1}(z) \in {\cccred \mathbf{R}}^{n-1}\left( \bigcap_{e \geq e_0} \wt{K^{e,c}_{D,n}}\right) 
\subseteq \bigcap_{e \geq e_0}{\cccred \mathbf{R}}^{n-1}\left(\wt{K^{e,c}_{D,n}}\right)  = 
\bigcap_{e \geq e_0}  K^{e,c}_{D,n}, 
\]
which implies ${\cccred \mathbf{R}}^{n-1}(z) \in 0^{*}_{D, n}$. {Here the last equality follows from Remark \ref{remark: easy remark-log}(1).} 
This completes the proof of the inclusion 
${\cccred \mathbf{R}}^{n-1}(\wt{0^{*}_{D, n}}) \subseteq 0^{*}_{D, n}$.

Now, it is enough to show the opposite inclusion 
${\cccred \mathbf{R}}^{n-1}(\wt{0^{*}_{D, n}}) \supseteq 0^{*}_{D, n}$. 
We have $D = a \cdot \div(g)$ for some $a \in \Q$ and $g \in R^{\circ}$. 
Take $z \in 0^*_{D, n}$ and $t \in \tau(R) \cap R^{\circ}$. 
There exist $c \in R^{\circ}$ and $e_1>0$ 
such that $z \in \bigcap_{e \geq e_1} K^{e, c}_{D, n}$. 
Replacing $c$ by $cg^2t^4$ (this is allowed by $K^{e, c}_{D, n} \subseteq K^{e, cg^2t^4}_{D, n}$), 
we may assume $\wt{K^{e, c}_{D, n}}=\wt{0^{*}_{D, n}}$ for $e \gg 0$ 
(Proposition \ref{prop: test elem pair tilda}(2)). 
For $e \gg 0$, we obtain 
\[
z \in \bigcap_{e \geq e_1} K^{e, c}_{D, n} \subseteq K^{e, c}_{D, n} 
\overset{(\star)}{=}{\cccred \mathbf{R}}^{n-1}(\wt{K^{e, c}_{D, n}})
={\cccred \mathbf{R}}^{n-1}(\wt{0^{*}_{D, n}}), 
\]
where 
the  equality $(\star)$ follows from  (Remark \ref{remark: easy remark-log}(1)). 
Thus ${\cccred \mathbf{R}}^{n-1}(\wt{0^{*}_{D, n}}) \supseteq 0^{*}_{D, n}$. 
\end{proof}


\begin{proposition}\label{prop: test element pair non-tilda}
We use Notation \ref{n-normal-local}. 
Let $D$ be a $\Q$-Cartier $\Q$-divisor on $X$ 
such that $D=a \cdot \mathrm{div}(g)$ for some  $a \in \Q$ and $g \in R^{\circ}$.
Take   $t \in \tau(R) \cap R^{\circ}$ and $t_D \in R^\circ$ 
satisfying $t_D \in \tau(R,\{p^eD\})$ for every integer $e\geq 0$. 
Fix $c \in R^{\circ} \cap ((g^2t^4) \cup (t_D^4))$.
Then  the following hold. 
\begin{enumerate}
    \item $
    H^d_{\m}(R(D)) \supseteq {K^{0,c}_{D,n}}
    \supseteq {K^{1,c}_{D,n}} \supseteq \cdots \supseteq {K^{e,c}_{D,n}}
    \supseteq {K^{e+1,c}_{D,n}} \supseteq \cdots.$
    \item 
    There exists an integer $e_1\geq 0$ such that 
    \[
    0^{*}_{D, n}=K^{e,c}_{D, n}
    \]
    for every integer $e \geq e_1$. 
    \item If 
    there exists an integer $e_2 \geq 0$  
    satisfying  $K^{e_2,c}_{D, n}=0$, then $0^{*}_{D, n}=0$.
\end{enumerate}
\end{proposition}

\begin{proof}
The assertions $(1)$ and $(2)$ follow from Proposition \ref{prop: test elem pair tilda}, 
Theorem \ref{thm: lifting thm pair}, and ${\cccred \mathbf{R}}^{n-1}(\wt{K^{e,c}_{D,n}})=K^{e,c}_{D,n}$ (Remark \ref{remark: easy remark-log}(1)). 
\qedhere 

\end{proof}

{The following {lemma} 
will be essential later in comparing the construction of $\tau(R,\Delta)$ with the definition of quasi-$F$-regularity.}

\begin{lemma} \label{lem:forgotten-definition-log}
{We use Notation \ref{n-normal-non-local} and let $D$ be an effective $\bQ$-divisor. Take $t_D \in R^\circ$ 
satisfying $t_D \in \tau(R,\{p^eD\})$ for every integer $e\geq 0$. Further take $f \in R^{\circ}$ and set $c := ft_D^2$. 
Fix integers $e_0 \geq 0$ and $n>0$. Then
\[
K^{e,f^{p^{e-e_0}}}_{D,n} \subseteq K^{e_0,c}_{D,n}.
\]
for every $e \geq e_0$.}
\end{lemma}

\begin{proof}Recall first that $K_{D, n}^{e, f^{p^{e-e_0}}}$ is the kernel of   
\[
\Phi^{e, f^{p^{e-e_0}}}_{R, D, n}  : 
H^d_\m(R(D)) \to  H^d_{\m}(Q_{R, D, n}^{e, f^{p^{e-e_0}}}).
\]
Take  $z \in K_{D, n}^{e, f^{p^{e-e_0}}}$ and 
a lift $\wt{z} \in \wt{K_{D, n}^{e, f^{p^{e-e_0}}}}$ of $z$ 
by the surjection 
\[
{\cccred \mathbf{R}}^{n-1} : \wt{K_{D, n}^{e, f^{p^{e-e_0}}}} \twoheadrightarrow K_{D, n}^{e, f^{p^{e-e_0}}}
\]
(Remark \ref{remark: easy remark-log}(1)). 
Now, for $e \geq e_0$, we have 
\[
0 = F_*^e[f]^{p^{e-e_0}} \cdot F^e(\wt{z}) 
= 
F^{e-e_0}(F_*^{e_0}[f] \cdot F^{e_0}(\wt{z})). 
\]
In view of \[
F_*^{e_0}[f] \cdot F^{e_0}(\wt{z}) \in F_*^{e_0}\Ker(F^{e-e_0} : H^d_\m(R(p^{e_0}D)) \to H^d_{\m}(R(p^eD))),
\]
Proposition \ref{prop: ker V and ker F: D}(2) 
implies that 
\[
F_*^{e_0}[c]  \cdot F^{e_0}(\wt{z})=F_*^{e_0}[ft_D^2] \cdot F^{e_0}(\wt{z})= 0,
\]
which in turn yields $\wt{z} \in \wt{K^{e_0,c}_{D, n}}$. 
Therefore, $z ={\cccred \mathbf{R}}^{n-1}(\widetilde z) \in 
{\cccred \mathbf{R}}^{n-1}(\wt{K^{e_0,c}_{D, n}}) =K^{e_0,c}_{D, n}$.    
    \end{proof}

\subsection{Duality and quasi-test ideals}\label{ss q-test ideal}

{We are ready to define and study quasi-test ideals. We warn the reader again about our notation in which $\omega_R(-(K_R+\Delta)) = R(\rup{-\Delta})$ {(cf.\ Remark \ref{r-W_n omega D}(1))}.}

\begin{definition}\label{d q test ideal}
We use Notation \ref{n-normal-non-local}. 
Let $\Delta$ be a $\Q$-divisor on $X$ such that $K_R+\Delta$ is $\Q$-Cartier. 
Take an integer $n>0$. 
Set 
\begin{eqnarray*}
\tau(W_nR,\Delta)&:=&\tau(W_n\omega_R,K_R+\Delta)\subseteq W_n\omega_R(-(K_R+\Delta)),\\
\tau_n(R, \Delta)&:=& \tau_n(\omega_R,K_R+\Delta) \subseteq 
R(\lceil-\Delta\rceil) \simeq \omega_R(-(K_R+\Delta)),\\
\tau^q(R, \Delta)&:=& \tau^q(\omega_R,K_R+\Delta) \subseteq 
R(\lceil-\Delta \rceil) \simeq \omega_R(-(K_R+\Delta)). 
\end{eqnarray*}
In particular, 
$\tau_n(R, \Delta)$ and $\tau^q(R, \Delta)$ are $R$-submodules of 
$R(\lceil-\Delta \rceil)$, 
which we call the {\em $n$-quasi-test ideal} and 
the {\em quasi-test ideal} of $(R, \Delta)$, 
respectively. 
\end{definition}

\noindent
{\ccred We warn the reader that $\tau(W_nR,\Delta)$ is not an ideal of $W_nR$ even if $R$ is regular and $\Delta =0$.}

\begin{remark}\label{r q test ideal clear}
We use Notation \ref{n-normal-non-local}. 
Let $\Delta$ be a $\Q$-divisor on $X$ such that $K_R+\Delta$ is $\Q$-Cartier. 
We here summarise some properties which can be immediately obtained from the earlier parts. 
\begin{enumerate}
\item 
By Definition \ref{d-taun-omega-log}, it holds that $\tau_n(R, \Delta) = (({\cccred \mathbf{R}}^{n-1})^*)^{-1}(\tau(W_nR, \Delta))$ 
for 
\[
({\cccred \mathbf{R}}^{n-1})^* : \omega_R(-(K_R+\Delta)) \to W_n\omega_R(-(K_R+\Delta)). 
\]
\item 
We have 
\[
\tau_1(R, \Delta) \subseteq \tau_2(R, \Delta) \subseteq \cdots 
\subseteq 
\tau_n(R, \Delta) \subseteq \tau_{n+1}(R, \Delta) \subseteq \cdots \subseteq R(-\Delta)
\]
and $\tau^q(R, \Delta) = \bigcup_{n>0} \tau_n(R, \Delta) = \tau_n(R, \Delta)$ for every $n \gg 0$ (Definition \ref{d-taun-omega-log}). 
\item 
Quasi-test ideals commute with localisations and completions 
(Proposition \ref{prop: test submod pair compat}). 
Specifically, the following hold. 
\begin{enumerate}
    \item 
    If $S$ is a multiplicatively closed subset of $R$, then 
\begin{equation}\label{e1-test submod compat-log2}
\tau(W_n\omega_R, \Delta) \otimes_{W_nR} W_n(S^{-1}R) \simeq \tau(W_n\omega_{S^{-1}R}, \Delta')
\end{equation}
    \begin{equation}\label{e2-test submod compat-log2}
\tau_n(R, \Delta) \otimes_{R} S^{-1}R \simeq \tau_n(S^{-1}R, \Delta')
\end{equation}
    \begin{equation}\label{e3-test submod compat-log2}
\tau^q(R, \Delta) \otimes_{R} S^{-1}R \simeq \tau^q(S^{-1}R, \Delta')
\end{equation}
    where $\Delta'$ denotes the pullback of $\Delta$ to $\Spec (S^{-1}R)$. 
    \item If $(R,\m)$ is a local ring, then
\begin{equation}\label{e4-test submod compat-log2}
\tau(W_n\omega_R, \Delta) \otimes_{W_nR} W_n(\widehat{R}) \simeq \tau(W_n\omega_{\widehat{R}}, \widehat{\Delta})
\end{equation}
    \begin{equation}\label{e5-test submod compat-log2}
\tau_n(R, \Delta) \otimes_{R} \widehat{R} \simeq \tau_n(\widehat{R}, \widehat{\Delta})
\end{equation}
    \begin{equation}\label{e6-test submod compat-log2}
\tau^q(R, \Delta) \otimes_{R} \widehat{R} \simeq \tau^q(\widehat{R},\widehat{\Delta})
\end{equation}
    where $\widehat{R}$ denotes the $\m$-adic completion of $R$ and 
    $\widehat{\Delta}$ is the pullback of $\Delta$ to $\Spec \widehat{R}$. 
\end{enumerate}
\end{enumerate}
\end{remark}

{The goal of the next proposition is two-fold: we show that the ideal $\tau(W_nR,\Delta)$ can be computed using quasi-tight closure and, to this end, we simplify its definition (cf.\ Proposition \ref{prop: test submodule gen pair}(3)) by showing that it is equal to $\tau^c(W_nR,\Delta)$ below. This equality is intuitively clear, as $\tau(W_nR,\Delta)$ and $\tau^c(W_nR,\Delta)$ differ by a ``fixed constant'' which can be subsumed by perturbations inherent in the theory of quasi-test ideals.}


\begin{proposition}\label{p-tau-WnX}
We use Notation \ref{n-normal-non-local}. 
Fix an integer $n>0$. 
Let $\Delta$ be a $\Q$-divisor on $X$ such that $K_R+\Delta$ is $\Q$-Cartier. 
Take  $a \in \Q$ and $g \in R^{\circ}$ 
satisfying $K_R+\Delta=a \cdot \mathrm{div}(g)$. 
Pick $t \in R^{\circ} \cap \tau(R)$  and $c \in R^{\circ} \cap (g^2t^4)$. 
Set 
\[
\tau^c(W_nR, \Delta) := \sum_{e \geq 0} \PsiT_{n}^{e, c}(F^e_*W_n\omega_R(-p^e(K_R+\Delta))), 
\]
where 
\[
T_{n}^{e, c} : F^e_*W_n\omega_R(-p^e(K_R+\Delta)) \xrightarrow{\cdot F_*^e[c]} 
F^e_*W_n\omega_R(-p^e(K_R+\Delta))  \xrightarrow{T^e_n} W_n\omega_R(-(K_R+\Delta))
\]
is the $W_nR$-module homomorphism obtained by applying $\mathcal Hom_{W_nR}(-, W_n\omega_R)$ to 
\[
W_nR(K_R+\Delta) \xrightarrow{F^e} F^e_*W_nR(p^e(K_R+\Delta)) \xrightarrow{\cdot F^e_*[c]} F^e_*W_nR(p^e(K_R+\Delta)). 
\]
Then the following hold. 
\begin{enumerate}
\item {$I^{0, c}_{K_X+\Delta} \subseteq I^{1, c}_{K_X+\Delta} \subseteq I^{2, c}_{K_X+\Delta} \subseteq\cdots $ 
for  $I^{e, c}_{K_X+\Delta} :=  \PsiT_{n}^{e, c}(F^e_*W_n\omega_R(-p^e(K_R+\Delta)))$.} 
\item 
$\tau(W_nR,\Delta) = \tau^c(W_nR, \Delta)$. 
\item If $(R, \m)$ is a local ring, then the equality 
    \[
    \left(\frac{W_n\omega_R(-(K_R+\Delta))}{\tau(W_nR,\Delta)}\right)^\vee = \wt{0^{*}_{K_R+\Delta, n}}
    \]
    of $W_nR$-submodules of $(W_n\omega_R(-(K_R+\Delta)))^{\vee}= H^d_{\m}(W_nR(K_R+\Delta))$ holds. 
\end{enumerate}
\end{proposition}
\begin{proof} 
We start by proving the following claim.
\begin{claim*}
If $(R, \m)$ is a local ring, then the equality 
    \[
    \left(\frac{W_n\omega_R(-(K_R+\Delta))}{\tau^c(W_nR,\Delta)}\right)^\vee = \wt{0^{*}_{K_R+\Delta, n}}
    \]
    of $W_nR$-submodules of $(W_n\omega_R(-(K_R+\Delta)))^{\vee}= H^d_{\m}(W_nR(K_R+\Delta))$ holds. 
\end{claim*}
\begin{proof}[Proof of Claim]
Set $D:=K_R+\Delta$. 
We define $I^{e, c}_{D, n}$ and $C^{e, c}_{D, n}$ as the image and the cokernel of 
\[
T^{e, c}_n : 
 F^e_*W_n\omega_R(-p^eD) \xrightarrow{\cdot F_*^e[c]} 
F^e_*W_n\omega_R(-p^eD)  \xrightarrow{T^e_n} W_n\omega_R(-D). 
\]
By definition, we have 
\[
\tau^c(W_nR,\Delta) = \sum_{e \geq 0} I^{e, c}_{D, n}
\]
and $T^{e, c}_n$ is obtained by applying the Matlis duality functor 
$\Hom_{W_nR}(-, E)$ to 
\begin{eqnarray*}
 F^{e, c} : H^d_\m(W_nR(K_R+\Delta)) &\to&  H^d_\m(W_nR(p^e(K_R+\Delta)))\\
z &\mapsto& F_*^e[c] \cdot F^e(z). 
\end{eqnarray*}
By  Proposition \ref{prop: test elem pair tilda}, we have 
\[
H^d_{\m}(W_nR(D)) \supseteq 
\wt{K^{0, c}_{D, n}} \supseteq 
\wt{K^{1, c}_{D, n}} \supseteq \cdots \supseteq 
\wt{K^{e_0, c}_{D, n}}
= 
\wt{K^{e_0+1, c}_{D, n}}=
\cdots = \wt{0^{*}_{D, n}}. 
\]
We now apply the Matlis duality functor $\Hom_{W_nR}(-, E)$. 
Then the kernel $\wt{K^{e, c}_{D, n}} = \Ker F^{e, c}$ turns into the cokernel $C^{e, c}_{D, n} = \Coker T_n^{e, c}$. 
Hence we get 
\[
W_n\omega_R(-D) \twoheadrightarrow
C^{0, c}_{D, n}  \twoheadrightarrow
C^{1, c}_{D, n}  \twoheadrightarrow \cdots  \twoheadrightarrow
C^{e_0, c}_{D, n}  \overset{\simeq}{\twoheadrightarrow} 
C^{e_0+1, c}_{D, n}  \overset{\simeq}{\twoheadrightarrow} \cdots 
= (\wt{0^{*}_{D, n}})^{\vee}. 
\]
Since $I^{e, c}_{D, n} = \Ker(W_n\omega_R(-D) \twoheadrightarrow C^{e, c}_{D, n})$, it holds that 
\[
I^{0, c}_{D, n} \subseteq 
I^{1, c}_{D, n}\subseteq \cdots \subseteq 
I^{e_0, c}_{D, n} = 
I^{e_0+1, c}_{D, n} = \cdots. 
\]
In particular, 
\[
\tau^c(W_nR,\Delta) = \sum_{e \geq 0} I^{e, c}_{D, n} = I^{e, c}_{D, n}
\]
for $e \gg 0$, which implies 
\begin{equation}\label{e1 tau-WnX}
\left(\frac{W_n\omega_R(-D)}{\tau^c(W_nR,\Delta)}\right)^\vee
=
\left(\frac{W_n\omega_R(-D)}{I^{e, c}_{D, n}}\right)^\vee 
=
\left( C^{e, c}_{D, n} \right)^\vee  = \wt{K^{e, c}_{D, n}} =  \wt{0^{*}_{D, n}}. 
\end{equation}
Thus Claim holds. 
\end{proof}

Now {(1) follows from the above argument and} (3) is immediate from (2) and Claim, and hence it is enough to prove (2). To this end, by taking a localisation at a prime ideal, we may assume that $R$ is a local ring 
(here 
$\tau^c$ localises as it is defined as a colimit, and 
$\tau$ localises by Proposition \ref{prop: test submod compat} and  Proposition \ref{prop: test submod pair compat}).  
The problem is thus reduced to the case when $K_R$ and $K_R+\Delta$ are effective. 
Then we have $a \geq 0$ and 
\begin{eqnarray*}
\tau^c({W_n}R, \Delta) &=& \sum_{e \geq 0} \PsiT_{n}^{e, c}(F^e_*W_n\omega_R(-p^e(K_R+\Delta))) \subseteq W_n\omega_R(-(K_R+\Delta)),\\
\tau(W_nR, \Delta) &=& 
\sum_{e \geq 0} \PsiT_{n}^{e, cg^{\rup{a(p^e-1)}}}(F^e_*W_n\omega_R(-(K_R+\Delta))) \subseteq W_n\omega_R(-(K_R+\Delta)). 
\end{eqnarray*}
where the latter equality holds by 
Proposition \ref{prop: test submodule gen pair}$(4)$. 

We first prove the inclusion 
$\tau(W_nR,\Delta) \supseteq \tau^c({W_n}R, \Delta)$.
By $\div(g) \geq 0$ {and $a(p^e-1) +1 \geq \rup{a(p^e-1)}$}, 
we get 
\begin{eqnarray*}
    -p^e(K_R+\Delta)-\mathrm{div}(g) 
    &=& -(K_R+\Delta) -a(p^e-1) \cdot \div(g)- \div (g)\\
    &\leq& -(K_R+\Delta) -\rup{a(p^e-1)} \cdot \div(g)\\
    &=& -(K_R+\Delta)-\mathrm{div}(g^{\rup{a(p^e-1)}}). 
\end{eqnarray*}
We have 
$\tau^c({W_n}R, \Delta) = \tau^{cg}({W_n}R, \Delta)$ by {Claim}. 
Then it holds that
\begin{align*}
\tau^c({W_n}R, \Delta) &= \tau^{cg}({W_n}R, \Delta)\\
    &= \PsiT^{e,cg}_{n}(F^e_*W_n\omega_R(-p^e(K_R+\Delta))) \\
    &= \PsiT^{e, c}_{n}(F^e_*W_n\omega_R(-p^e(K_R+\Delta)-\mathrm{div}(g))) \\
    &\subseteq \PsiT^{e,c}_{n}(F^e_*W_n\omega_R(-(K_R+\Delta)-\mathrm{div}(g^{\rup{a(p^e-1)}}))) \\
    &= \PsiT^{e,cg^{\rup{a(p^e-1)}}}_{n}(F^e_*W_n\omega_R(-(K_R+\Delta))) \subseteq \tau(W_n\omega_R,K_R+\Delta),
\end{align*}
as required. 

The opposite inclusion 
$\tau(W_nR,\Delta)  \subseteq  \tau^c(W_nR, \Delta)$ follows from 
\begin{align*}
&\PsiT^{e,cg^{\rup{a(p^e-1)}}}_{n}(F^e_*W_n\omega_R(-(K_R+\Delta)))\\
&\quad = \PsiT^{e,c}_{n}(F^e_*W_n\omega_R(-(K_R+\Delta)-\mathrm{div}(g^{\rup{a(p^e-1)}})))\\
&\quad = \PsiT^{e,c}_{n}(F^e_*W_n\omega_R(-(K_R+\Delta)-\rup{a(p^e-1)}\mathrm{div}(g)))\\
&\quad \subseteq  \PsiT^{e,c}_{n}(F^e_*W_n\omega_R(-(K_R+\Delta)-a(p^e-1)\mathrm{div}(g)))\\
&\quad =  \PsiT^{e,c}_{n}(F^e_*W_n\omega_R(-p^e(K_R+\Delta))). 
\end{align*}
This concludes the proof of (2) and the entire proposition.
\qedhere

\end{proof}

\begin{proposition}\label{prop: test submod constant mult}
We use Notation \ref{n-normal-non-local}. 
Let $\Delta$ be a $\Q$-divisor on $X$ such that $K_R+\Delta$ is $\Q$-Cartier. 
Take $h \in R^{\circ}$.
Then the following hold. 
\begin{enumerate}
\item $\tau(W_nR, \Delta+\div(h)) = [h] \cdot \tau(W_nR, \Delta)$. 
\item $\tau_n(R, \Delta+\div(h)) = h \cdot \tau_n(R, \Delta)$. 
\item $\tau^q(R, \Delta+\div(h)) = h \cdot \tau^q(R, \Delta)$. 
\end{enumerate}
\end{proposition}

\begin{proof}
We may assume that there exist $a \in \Q$ and $g \in R^{\circ}$ 
satisfying $K_R+\Delta=a \cdot \mathrm{div}(g)$. 
Take a test element $t \in R^{\circ}$  and $c \in R^{\circ} \cap (h^2g^2t^4)$.
By Proposition \ref{p-tau-WnX}(3), 
the assertion (1) follows from 
\begin{align*}
\tau(W_nR,\Delta+\mathrm{div}(h))
&=\sum_{e \geq 0} T^{e,c}_{n}(F^e_*W_n\omega_R(-p^e(K_R+\Delta+\mathrm{div}(h))) \\
&= \sum_{e \geq 0} [h] \cdot T^{e,c}_{n}(F^e_*W_n\omega_R(-p^e(K_R+\Delta))) \\
&=[h] \cdot \tau(W_nR,\Delta). 
\end{align*}
By Remark \ref{r q test ideal clear}, (1) implies (2) and (2) implies (3). 
\end{proof}

{\begin{theorem}\label{thm: corr 0^* and tauPAIR}
In addition to Notation \ref{n-normal-non-local}, assume that $(R,\m)$ is a local ring. 
Let $\Delta$ be a  $\Q$-divisor on $X$ such that $K_R+\Delta$ is $\Q$-Cartier. Take an integer $n>0$. Then the equality 
    \[
\left(\frac{R(\rup{-\Delta})}{\tau_n(R, \Delta)}\right)^\vee = 0^{*}_{K_R+\Delta, n}.
\] 
of $R$-submodules of 
$R(\rup{-\Delta})^{\vee} =\omega_R(-(K_R+\Delta))^{\vee}= H^d_{\m}(R(K_R+\Delta))$ holds. 

In particular,
\[
\tau_n(R,\Delta)^{{\ccred \vee}} = 
\frac{H^d_\m(R(K_R+\Delta))}{0^{*}_{K_R+\Delta,n}}.
\]
\end{theorem}}
\begin{proof}
    Since $\tau_n(R, \Delta)$ is the inverse image of $\tau(W_nR, \Delta)$ 
by $({\cccred \mathbf{R}}^{n-1})^* : \omega_R(-(K_R+\Delta)) \hookrightarrow W_n\omega_R(-(K_R+\Delta))$ (Remark \ref{r q test ideal clear}), 
we obtain the following commutative diagram:  
\[
\begin{tikzcd}
\omega_R(-(K_R+\Delta)) \arrow[r, "({\cccred \mathbf{R}}^{n-1})^*", hook] \arrow[d, twoheadrightarrow] & W_n\omega_R(-(K_R+\Delta))\arrow[d, twoheadrightarrow]\\
\frac{\omega_R(-(K_R+\Delta))}{\tau_n(R, \Delta)} \arrow[r, hook] &\frac{W_n\omega_R(-(K_R+\Delta))}{\tau(W_nR, \Delta)}.
\end{tikzcd}
\]
By applying the Matlis dual functor $(-)^{\vee} := \Hom_{W_nR}(-, E)$ (which is exact) to this diagram, we get 
\[
\begin{tikzcd}
H^d_{\m}(R(K_R+\Delta)) \arrow[r, "{\cccred \mathbf{R}}^{n-1}", twoheadleftarrow] \arrow[d, hookleftarrow] & H^d_{\m}(W_nR(K_R+\Delta)) \arrow[d, hookleftarrow]\\
(\frac{\omega_R(-(K_R+\Delta))}{\tau_n(R, \Delta)})^{\vee} \arrow[r, twoheadleftarrow] & (\frac{W_n\omega_R(-(K_R+\Delta))}{\tau(W_nR, \Delta)})^{\vee}.
\end{tikzcd}
\]
Since 
\[
\left(\frac{W_n\omega_R(-(K_R+\Delta))}{\tau(W_nR, \Delta)}\right)^{\vee} = \wt{0^{*}_{K_R+\Delta, n}} \qquad \text{ and } \qquad {\cccred \mathbf{R}}^{n-1}\left(\wt{0^{*}_{K_R+\Delta, n}}\right) =
0^{*}_{K_R+\Delta, n}
\]
by Proposition \ref{p-tau-WnX}
and Theorem \ref{thm: lifting thm pair}, respectively, 
we obtain the following:
\begin{align*}
\left(\frac{R(\rup{-\Delta})}{\tau_n(R, \Delta)}\right)^{\vee} &=
\left(\frac{\omega_R(-(K_R+\Delta))}{\tau_n(R, \Delta)}\right)^{\vee} = 
{\cccred \mathbf{R}}^{n-1}\left(\left(\frac{W_n\omega_R(-(K_R+\Delta))}{\tau(W_nR, \Delta)}\right)^{\vee}\right)\\
&={\cccred \mathbf{R}}^{n-1}(\wt{0^{*}_{K_R+\Delta, n}}) = 0^*_{K_R+\Delta, n}.
\end{align*}
This concludes the proof.
\end{proof}

{Next, we provide an explicit definition of $\tau_n(R,\Delta)$ 
(Corollary \ref{cor:Qdefs-of-quasi-tight-submodule-PAIR}) 
in terms of $Q^{e,c}_{R,K_R+\Delta,n}$ (recall the notation from (\ref{e-def-Q^{e,c}_D})). {First, we need the following {lemma}.}

\begin{lemma} \label{lemma:inequalities-for-quasi-f-module-with-test-elementPAIR}
We use Notation \ref{n-normal-non-local} and let $D$ be an effective $\bQ$-divisor. Write $D = a \cdot {\rm div}(g)$ for $a \in \bQ$ and $g \in R^\circ$. Take $t \in \tau(R) \cap R^{\circ}$ and $t_D \in R^\circ$ 
satisfying $t_D \in \tau(R,\{p^eD\})$ for every integer $e\geq 0$. 
Fix $c \in R^{\circ} \cap ((g^2t^4) \cup (t_D^4))$. 
Then the following assertions are valid. 
\begin{enumerate}
\item 
The inclusion 
\begin{multline*}
{\rm Im}\Big(\Hom_{W_nR}(Q^{e,c}_{R,D,n}, W_n\omega_R) \xrightarrow{(\Phi^{e, c}_{R,D,n})^*} \omega_R(-D) \Big) \\
\subseteq {\rm Im}\Big(\Hom_{W_nR}(Q^{e',c}_{R,D,n}, W_n\omega_R) \xrightarrow{(\Phi^{e', c}_{R,D,n})^*} \omega_R(-D) \Big)
\end{multline*}
holds for every $e' \geq e$. 
\item 
There exists an integer $e_1>0$ such that the equality 
\begin{multline*}
{\rm Im}\Big(\Hom_{W_nR}(Q^{e,c}_{R,D,n}, W_n\omega_R) \xrightarrow{(\Phi^{e, c}_{R,D,n})^*} \omega_R(-D) \Big) \\
= {\rm Im}\Big(\Hom_{W_nR}(Q^{e_1,c}_{R,D,n}, W_n\omega_R) \xrightarrow{(\Phi^{e_1, c}_{R,D,n})^*} \omega_R(-D) \Big)
\end{multline*}
holds for every integer $e \geq e_1$. 
\end{enumerate}
\end{lemma}

\begin{proof}
It is enough to prove (1), because (2) holds by (1) and the fact that $\omega_R$ is a finitely generated $R$-module. 
In order to show (1), we may assume that $R$ is a local ring. 
Let $\m$ be its maximal ideal. 
Then the assertion follows from 
Proposition \ref{prop: test element pair non-tilda}{(1)} and the fact that  
\[
{\rm Im}\Big(\Hom_{W_nR}(Q^{e,c}_{R,D,n}, W_n\omega_R) \xrightarrow{(\Phi^{e, c}_{R,D,n})^*} \omega_R(-D) \Big)
\]
is Matlis dual to 
${\rm Im}\Big( H^d_\m(R(D)) \xrightarrow{\Phi^{e, c}_{R,D,n}} H^d_\m(Q^{e,c}_{R,D,n})\Big) = \frac{H^d_\m(R)}{{K^{e,c}_{D,n}}}$. 
\end{proof}

\begin{corollary} \label{cor:Qdefs-of-quasi-tight-submodule-PAIR}
We use Notation \ref{n-normal-non-local}. Let $\Delta$ be a $\bQ$-divisor on $X$ such that $K_R+\Delta$ is $\bQ$-Cartier and set $D = K_R + \Delta$. Take $a \in \bQ$ and $g \in R^\circ$ with $D = a \cdot {\rm div}(g)$.
Then the following 
hold.
{\setlength{\leftmargini}{2em} \begin{enumerate}
    \item $\tau_n(R,\Delta) = \bigcap_{c \in R^\circ} \bigcap_{e_0>0} \sum_{e \geq e_0} {\rm Im}\Big(\Hom_{W_nR}(Q^{e,c}_{R,D,n}, W_n\omega_R) \xrightarrow{(\Phi^{e, c}_{R,D,n})^*} \omega_R(-D) \Big)$,\\[-0.5em]
    \item $\tau_n(R,\Delta) = \bigcap_{c \in R^\circ} \sum_{e >0 } {\rm Im}\Big(\Hom_{W_nR}(Q^{e,c}_{R,D,n}, W_n\omega_R) \xrightarrow{(\Phi^{e, c}_{R,D,n})^*} \omega_R(-D) \Big)$,\\[-0.5em]
    \item  $\tau_n(R,\Delta) = \bigcap_{c \in R^\circ} \sum_{e_0>0} \bigcap_{e \geq e_0} {\rm Im}\Big(\Hom_{W_nR}(Q^{e,c^{p^{e-e_0}}}_{R,D,n}, W_n\omega_R) \xrightarrow{(\Phi^{e, c^{p^{e-e_0}}}_{R,D,n})^*} \omega_R(-D) \Big)$,\\[-0.5em]
    \item 
    {Take $t \in \tau(R) \cap R^{\circ}$ and $t_D \in R^\circ$ 
satisfying $t_D \in \tau(R,\{p^eD\})$ for every integer $e\geq 0$. Fix $c \in R^{\circ} \cap ((g^2t^4) \cup (t_D^4))$. 
    Then there exists an integer $e_0 >0$ such that 
    \[
    \tau_n(R,\Delta) = {\rm Im}\Big(\Hom_{W_nR}(Q^{e,c}_{R,D,n}, W_n\omega_R) \xrightarrow{(\Phi^{e, c}_{R,D,n})^*} \omega_R(-D) \Big)
    \]
    for every integer $e \geq e_0$.} 
\end{enumerate}}
\end{corollary}


\begin{proof}
We start by proving (4). 
By Lemma \ref{lemma:inequalities-for-quasi-f-module-with-test-elementPAIR},   
{\ccred there exists an integer $e_0>0$ such that 
\begin{multline*}
\bigcup_{e>0} {\rm Im}\Big(\Hom_{W_nR}(Q^{e,c}_{R,D,n}, W_n\omega_R) \xrightarrow{(\Phi^{e, c}_{R,D,n})^*} \omega_R(-D) \Big) \\
= {\rm Im}\Big(\Hom_{W_nR}(Q^{e,c}_{R,D,n}, W_n\omega_R) \xrightarrow{(\Phi^{e, c}_{R,D,n})^*} \omega_R(-D) \Big)
\end{multline*}
for every $e \geq e_0$. 
Fix $e \geq e_0$.} 
It is enough to prove that
\[
\tau_n({\ccred R, \Delta}) = {\rm Im}\Big(\Hom_{W_nR}(Q^{e,c}_{R,D,n}, W_n\omega_R) \xrightarrow{(\Phi^{e, c}_{R,D,n})^*} \omega_R(-D) \Big).
\]
In order to verify this, we may assume that $R$ is a {\ccred complete} local ring by (\ref{e2-test submod compat-log2}) and 
{\ccred (\ref{e5-test submod compat-log2})}. 
In this case,
\begin{equation} \label{eq:taun-dual-new-defPAIR00}
{\rm Im}\Big(\Hom_{W_nR}(Q^{e,c}_{R,D,n}, W_n\omega_R) \xrightarrow{(\Phi^{e, c}_{R,D,n})^*} \omega_R(-D) \Big)
\end{equation}
is Matlis dual to
\begin{equation} \label{eq:taun-dual-new-defPAIR}
{\rm Im}\Big( H^d_\m(R(D)) \xrightarrow{\Phi^{e, c}_{R,D,n}} H^d_\m(Q^{e,c}_{R,D,n})\Big) = \frac{H^d_\m(R(D))}{{K^{e,c}_{D,n}}}.
\end{equation}

Now, Theorem \ref{thm: corr 0^* and tauPAIR} and Proposition \ref{prop: test element pair non-tilda}(2) imply
\begin{equation} \label{eq:hidden-equation-kecn2PAIR}
\tau_n(R,\Delta) = \Big(\frac{H^d_\m(R(D))}{{K^{e,c}_{D,n}}}\Big)^\vee.
\end{equation}
for every $e \gg 0$ (depending on $c$). Thus (4) follows by (\ref{eq:taun-dual-new-defPAIR}) {and (\ref{eq:hidden-equation-kecn2PAIR})}.

Second, 
{assertion (1) holds, because}  
\begin{align*}
\bigcap_{c \in R^\circ} \bigcap_{e_0>0} \sum_{e \geq e_0}\ & {\rm Im}\Big(\Hom_{W_nR}(Q^{e,c}_{R,D,n}, W_n\omega_R) \xrightarrow{(\Phi^{e, c}_{R,D,n})^*} \omega_R(-D) \Big) \\ &= \bigcap_{c \in g^2t^4R^\circ} \bigcap_{e_0>0} \sum_{e \geq e_0} {\rm Im}\Big(\Hom_{W_nR}(Q^{e,c}_{R,D,n}, W_n\omega_R) \xrightarrow{(\Phi^{e, c}_{R,D,n})^*} \omega_R(-D) \Big)  \\
&\overset{\mathclap{\ref{lemma:inequalities-for-quasi-f-module-with-test-elementPAIR}}}{=} \bigcap_{c \in g^2t^4R^\circ} \sum_{e >0 }  {\rm Im}\Big(\Hom_{W_nR}(Q^{e,c}_{R,D,n}, W_n\omega_R) \xrightarrow{(\Phi^{e, c}_{R,D,n})^*} \omega_R(-D) \Big) \\
&\overset{{(4)}
}{=} \tau_n(R, \Delta).
\end{align*}
Then (2) holds by the same argument. 

Thus we are left with showing (3). 
We denote its right hand side by $\tau'_n(R,\Delta)$. By (2), we immediately have that $\tau'_n(R,\Delta) \subseteq \tau_n(R,\Delta)$. 
It suffices to show the opposite inclusion $\tau'_n(R,\Delta) \supseteq \tau_n(R,\Delta)$.  
Fix $t_D \in R^\circ$ satisfying $t_D \in \tau(R,\{p^eD\})$ for every integer $e\geq 0$. 
{For every $e \geq e_0$ and $c \in R^{\circ}$}, 
we have the 
the following inclusion (Lemma \ref{lem:forgotten-definition-log}): 
\[
K^{e,c^{p^{e-e_0}}}_{D,n} \subseteq K^{e_0,ct^2_D}_{D,n}. 
\]
By applying  Matlis duality to this inclusion, 
we get 
{\begin{multline*}
{\rm Im}\Big(\Hom_{W_nR}(Q^{e,c^{p^{e-e_0}}}_{R,D,n}, W_n\omega_R) \xrightarrow{(\Phi^{e, c^{p^{e-e_0}}}_{R,D,n})^*} \omega_R(-D) \Big) \\
\supseteq {\rm Im}\Big(\Hom_{W_nR}(Q^{e_0,ct^2_D}_{R,D,n}, W_n\omega_R) \xrightarrow{(\Phi^{e_0, ct^2_D}_{R,D,n})^*} \omega_R(-D) \Big), 
\end{multline*}
{because (\ref{eq:taun-dual-new-defPAIR00}) and (\ref{eq:taun-dual-new-defPAIR}) are Matlis dual to each other}.} 
Therefore, 
\[
\tau'_n(R,\Delta) \supseteq  \bigcap_{c \in R^\circ} \sum_{e_0 >0 } {\rm Im}\Big(\Hom_{W_nR}(Q^{e_0, c{t^2_D}}_{R,D,n}, W_n\omega_R) \xrightarrow{(\Phi^{e_0, c{t^2_D}}_{R,D,n})^*} \omega_R(-D) \Big) 
\overset{{\rm (2)}}{=} \tau_n(R,\Delta).
\]
Thus (3) holds. 
\qedhere 
\end{proof}


\begin{proposition}\label{prop: test ideal pair non-tilda}
We use Notation \ref{n-normal-non-local}. 
Let $\Delta$ be a  $\Q$-divisor on $X$ such that $K_R+\Delta$ is $\Q$-Cartier.
Then the following hold. 
\begin{enumerate}
\item 
$(R, \Delta)$ is $n$-quasi-$F$-regular if and only if $\tau_n(R, \Delta)=R$.
\item 
$(R, \Delta)$ is quasi-$F$-regular if and only if $\tau^q(R, \Delta)=R$.
\end{enumerate}

\end{proposition}
\noindent {In particular, 
{$n$-quasi-$F$-regularity and quasi-$F$-regularity commute with localisation}.} 


\begin{proof}
Since (2) holds by (1), let us  show (1). 
It follows from 
Corollary \ref{cor:Qdefs-of-quasi-tight-submodule-PAIR}(3) 
that 
\[
\tau_n(R,\Delta) = \bigcap_{c \in R^\circ} \sum_{e_0>0} \bigcap_{e \geq e_0} {\rm Im}\Big(\Hom_{W_nR}(Q^{e}_{R,D_c,n}, W_n\omega_R) \xrightarrow{(\Phi^e_{R,D_c,n})^*} \omega_R(-D) \Big),
\]
where $D_c := D+ (1/p^{e_0}){\rm div}(c)$. 
In particular, by Definition \ref{d-IQFS}, we get that 
{$\tau_n(R,\Delta)=R$}  
if and only if for every effective Cartier divisor $E$, there exists $e_0>0$ such that $(R,\Delta + \frac{1}{p^{e_0}}E)$ is $n$-quasi-$F^e$-{\cred split} 
for every $e \geq e_0$ (equivalently, for every $e \geq 0$). But this is clearly equivalent to the original definition of $n$-quasi-$F$-regularity (see Definition \ref{d-QFR}).
\qedhere


\end{proof}

\begin{theorem}\label{thm:chara qFr pair}
{We use Notation \ref{n-normal-non-local}.} 
Let  $\Delta$ an effective $\Q$-divisor on $X$ such that $\rdown{\Delta}=0$ 
and $D:=K_R+\Delta$ is $\Q$-Cartier. 
Take $a \in \Q$ and $g \in R^{\circ}$ with $D=a \cdot \mathrm{div}(g)$.
Pick $t \in \tau(R) \cap R^{\circ}$ and $t_D \in R^\circ$ satisfying $t_D \in \tau(R,\{p^e\Delta\})$ for every integer $e \geq 0$.
Fix $c \in R^{\circ} \cap ((g^2t^4) \cup (t_D^4))$ and an integer $n>0$. 
Then the following are equivalent.
\begin{enumerate}
    \item $(R, \Delta)$ is $n$-quasi-$F$-regular.\\[-0.9em] 
    \item There exists $e \in \bZ_{> 0}$ such that 
    $(R,\Delta+(1/p^{e})\mathrm{div}(c))$ is $n$-quasi-$F^{e}$-split. \\[-0.9em]
    \item {
    There exists an integer $e_0>0$ such that 
    \begin{equation}\label{e1 chara qFr pair}
    (\Phi^{e, c}_{R,D,n})^* \colon \Hom_{W_nR}(Q^{e,c}_{R,D,n}, W_n\omega_R) \xrightarrow{} \omega_R(-D)
    \end{equation}
    is surjective for every $e \geq e_0$.\\[-0.9em]}
    \item $0^{*}_{D_{\m}, n}=0$ for every maximal ideal $\m$ of $R$, 
    where $D_{\m}$ denotes the pullback of $D$ to $\Spec R_{\m}$. 
\end{enumerate}
{In particular, $(R,\Delta)$ is quasi-$F$-regular if and only if there exist integers $n>0$ and $e >0$ such that $(R,\Delta + (1/p^e){\rm div}(c))$ is $n$-quasi-$F^e$-split.}
\end{theorem}
\begin{proof}
Note that $0^*_{D_{\m},n} = K^{e,c}_{D_{\m},n}$ for every $e \gg 0$ by Proposition \ref{prop: test element pair non-tilda}. Thus the equivalence  (3) $\Leftrightarrow$ (4) is immediate by Noetherian induction and Matlis duality 
(cf.\ Matlis duality between (\ref{eq:taun-dual-new-defPAIR00}) and (\ref{eq:taun-dual-new-defPAIR})). 
{The equivalence (2) $\Leftrightarrow$  (3) 
follows from Lemma \ref{lemma:inequalities-for-quasi-f-module-with-test-elementPAIR}(1), 
because  $(R,\Delta+(1/p^{e})\mathrm{div}(c))$ is $n$-quasi-$F^{e}$-split 
if and only if (\ref{e1 chara qFr pair}) is surjective.} 
Finally, the equivalence  (1) $\Leftrightarrow$  (3) immediately holds by Corollary \ref{cor:Qdefs-of-quasi-tight-submodule-PAIR}(4) and Proposition \ref{prop: test ideal pair non-tilda}(1).

{The ``in-particular'' part is immediate from the equivalence of (1) $\Leftrightarrow$ (2).}
\end{proof}

\begin{theorem}\label{thm:chara qFR pair non-Q-Gor}
We use {Notation \ref{n-normal-non-local}.} 
Let $\Delta$ be an effective $\Q$-divisor on $X$ such that $K_R+\Delta$ is $\Q$-Cartier.
Then the following hold. 
\begin{enumerate}
\item Fix an integer $n>0$. 
Then $(R, \Delta)$ is $n$-quasi-$F$-regular if and only if for every $f \in R^{\circ}$, there exists an integer $e> 0$ such that $(R,\Delta+(1/p^e)\mathrm{div}(f))$ is $n$-quasi-$F^e$-split. 
\item  $(R, \Delta)$ is 
quasi-$F$-regular if and only if for every $f \in R^{\circ}$, there exists an integer $e>0$ such that  $(R,\Delta+(1/p^e)\mathrm{div}(f))$ is quasi-$F^e$-split.
\end{enumerate}

\end{theorem}
{\noindent The reader should compare this result with Definitions \ref{d-QFR} and \ref{d-weak-QFR}. We emphasise that in (2), for every $f \in R^\circ$, we just need to find integers $n>0$ and $e>0$ such that $(R,\Delta + (1/p^e){\rm div}(f))$ is $n$-quasi-$F^e$-split. In particular, we can take $n \gg e$, which is not possible in the original definition.} 
\begin{proof}
Assertion (1) is immediate from Corollary \ref{cor:Qdefs-of-quasi-tight-submodule-PAIR}(2) and Proposition \ref{prop: test ideal pair non-tilda}(1).  Now we prove (2).

{
First, suppose that $(R,\Delta)$ is $n$-quasi-$F$-regular for some integer $n>0$. Then, up to multiplying $f$ by an element of $R^\circ$, 
we may assume that $f \in R^\circ \cap ((g^2t^4) \cup (t^4_D))$ with notation as in Theorem \ref{thm:chara qFr pair}. Then the same theorem implies that $(R,\Delta+(1/p^e){\rm div}(f))$ is $n$-quasi-$F^e$-split for some $e>0$. 

As for the other implication, pick $f \in R^\circ \cap ((g^2t^4) \cup (t^4_D))$. Then there exist $e>0$ and $n>0$ such that $(R,\Delta+(1/p^e){\rm div}(f))$ is $n$-quasi-$F^e$-split, 
which implies that $(R,\Delta)$ is $n$-quasi-$F$-regular by Theorem \ref{thm:chara qFr pair} again.} 
\end{proof}
 We can finally give the proof of a theorem from the introduction.  

\begin{theorem}[= Theorem \ref{intro:thm:best-def-quasi-F-regular}]\label{thm:best-def-quasi-F-regular-Section4}
{\cred Let $R$ be an $F$-finite normal Noetherian $\bQ$-Gorenstein domain over $\F_p$. Let $D$ be an {effective} {Cartier} divisor on $\Spec R$ such that 
the inclusion $I_D \subseteq \tau(R)$ holds for the defining ideal $I_D$ of $D$. 
Then $R$ is quasi-$F$-regular 
if and only if there exist integers $e>0$ and $n>0$ 
such that $( R, (4/p^e)D)$ 
is $n$-quasi-$F^e$-split.} 
\end{theorem}



\begin{proof}
Taking a suitable finite affine open cover of $\Spec R$, 
we may assume that $D = {\rm div}(g)$. 
{\cred 
It follows from $gR = I_D \subseteq \tau(R)$ that} 
$c := g^{4}$ satisfies the assumptions of Theorem \ref{thm:chara qFr pair}. 
By 
{\cred $(4/p^e)D = (4/p^e) {\rm div}(g) = (1/p^4) {\rm div}(c)$}, we can conclude the proof by the ``in-particular'' part of Theorem \ref{thm:chara qFr pair}.
\end{proof}

\begin{theorem}\label{thm:chara qFr pair2}
We use Notation \ref{n-normal-non-local}. 
Let  $\Delta$ be an effective $\Q$-divisor on $X$ such that 
$K_R+\Delta$ is $\Q$-Cartier. 
Then $(R, \Delta)$ is quasi-$F$-regular if and only if $(R, \Delta)$ is feebly quasi-$F$-regular. 
\end{theorem}
\begin{proof}
{This is immediate from Theorem \ref{thm:chara qFR pair non-Q-Gor}. For the covenience of the reader, we spell the argument out. First, quasi-$F$-regularity implies feeble quasi-$F$-regularity. Now, assume that $(R,\Delta)$ is feebly quasi-$F$-regular. 
By Definition \ref{d-weak-QFR} and Definition \ref{d QFR prelim}, 
for every effective $\Q$-divisor $E$ on $X$, 
there exist 
$n \in \Z_{>0}$ and $\epsilon \in \Q_{>0}$  
such that $(X, \Delta + \epsilon E)$ is $n$-quasi-$F^e$-split 
for all $e \in \Z_{>0}$. In particular, for every $f \in R^{\circ}$, there exists an integer $e>0$ such that  $(R,\Delta+(1/p^e)\mathrm{div}(f))$ is $n$-quasi-$F^e$-split for some $n>0$. Thus $(R,\Delta)$ is quasi-$F$-regular by Theorem \ref{thm:chara qFR pair non-Q-Gor}(2).}
\qedhere

\end{proof}

{\begin{remark} \label{remark:non-uniform-quasi-F-regularity}
When working with quasi-$F^\infty$-splittings, one needs to consider its two variants: the usual one and the uniform one, depending on the order of the quantifiers between $n$ and $e$. This difference is fundamental. For example, cones over {supersingular} elliptic curves are always quasi-$F^\infty$-split, but they are never uniformly quasi-$F^\infty$-split. However, this difference disappears when working with quasi-$F$-rationality or quasi-$F$-regularity. This idea is essential for many parts of our paper.

In particular, one can show much stronger statements than that quasi-$F$-regularity and feeble quasi-$F$-regularity agree. For example, it is now easy to see from Theorem \ref{thm:chara qFr pair} that $(R,\Delta)$ is quasi-$F$-regular if and only if for every effective Weil divisor $E$, there exists a rational number $\epsilon > 0$ such that for every $e \gg 0$ there exists $n >0$ such that $(R,\Delta+\epsilon E)$ is $n$-quasi-$F^e$-split. 
\end{remark}}

\begin{theorem}\label{thm: suff small qFr}
We use Notation \ref{n-normal-non-local}. 
Let $\Delta$ be a $\Q$-divisor on $X$ such that $K_R+\Delta$ is $\Q$-Cartier.
Let $E$ be an effective $\Q$-divisor on $X$. 
Then the following hold. 
\begin{enumerate}
    \item 
    If $E$ is $\Q$-Cartier, 
    then there exists a rational number $\epsilon >0$ such that 
\[
\tau_n(R, \Delta)=\tau_n(R, \Delta+\epsilon E)\qquad\text{and}\qquad 
\tau^q(R, \Delta)=\tau^q(R, \Delta+\epsilon E)
\]
for every integer $n>0$. 
\item 
Take $n>0$. 
If $(R, \Delta)$ is $n$-quasi-$F$-regular, then 
there exists a rational number $\epsilon >0$ such that  $(R,\Delta+\epsilon E)$ is 
$n$-quasi-$F$-regular.     
\item 
If $(R, \Delta)$ is quasi-$F$-regular, then 
there exists a rational number $\epsilon >0$ such that  $(R,\Delta+\epsilon E)$ is quasi-$F$-regular.     
\end{enumerate}
\end{theorem}

\begin{proof}
The assertion (1) follows from Proposition \ref{prop: stability of tau omega}(2). In order to show (2) and (3), we may assume, by enlarging $E$, 
that $E$ is Cartier. 
 By Proposition \ref{prop: test ideal pair non-tilda}, (1) implies (2) and (3). 
 \qedhere


\end{proof}


The following result provides a reformulation of the splitting definition of a quasi-$F^e$-splitting from \cite[Proposition 3.19]{TWY} {and, together with Proposition \ref{prop:  qFreg 0-map} below, will be used in the proof that three-dimensional quasi-$F$-regular singularities are Cohen-Macaulay.}

\begin{proposition}\label{prop:qFs-divisor-split}
We use Notation \ref{n-normal-local}.
Let $D$ a $\Q$-divisor on $\Spec R$. 
Take $c \in R^\circ$ and integers $n >0$ and $e >0$. 
Assume that $(R,\{D\}+(1/p^e)\mathrm{div}(c)))$ is $n$-quasi-$F^e$-split. 
Then there exists a $W_nR$-module homomorphism $\alpha \colon F^e_*W_nR(p^eD) \to W_n\omega_R(-K_R+\rdown{D})$ such that  the following diagram is commutative:
\begin{equation}\label{e1:qFs-divisor-split}
\begin{tikzcd}
    W_nR(D) \arrow[r,"{(\cdot F^e_*[c])}\circ F^e"] \arrow[d,"{\cccred \mathbf{R}}^{n-1}"'] & F^e_*W_nR(p^eD) \arrow[ldd,"\alpha"] \\
    \mathllap{R(D)=\,\, }\omega_R(-K_R+\rdown{D}) 
    \arrow[d, "({\cccred \mathbf{R}}^{n-1})^*"'] \\
    W_n\omega_R(-K_R+\rdown{D})
\end{tikzcd}
\end{equation}
\end{proposition}

\begin{proof}
For every $\Q$-divisor $D'$ on $\Spec R$, we set
\[
\varphi_{D'}:=(\cdot F^e_*[c]) \circ F^e \colon W_nR(D') \to F^e_*W_nR(p^eD')
\]
Then $\varphi_{D'}$ coincides with the composition 
\[
W_nR(D') \xrightarrow{F^e} F^e_*W_nR(p^eD') \hookrightarrow F^e_*W_nR(p^eD'+\mathrm{div}(c)) \xrightarrow{\cdot F^e_*[c]} F^e_*W_nR(p^eD'),
\]
where the second homomorphism is the natural injection.
Since $(R,\{D\}+(1/p^e)\mathrm{div(c)})$ is $n$-quasi-$F^e$-split, there exists a $W_nR$-module homomorphism $\beta \colon F^e_*W_nR(p^e\{D\}+\mathrm{div}(c)) \to W_n\omega_R(-K_R)$  such that  the following diagram is commutative \cite[Proposition 3.19]{TWY}:
\[
\begin{tikzcd}
    W_nR(\{D\} +(1/p^e)\mathrm{div(c)}) \arrow[r,"F^e"] \arrow[d,"{\cccred \mathbf{R}}^{n-1}"'] & F^e_*W_nR(p^e\{D\}+\mathrm{div}(c)) \arrow[ldd,"\beta"] \\
    R \arrow[d,"({\cccred \mathbf{R}}^{n-1})^*"'] \\
    W_n\omega_R(-K_R)
\end{tikzcd}
\]
Therefore, if we define a $W_nR$-module homomorphism $\alpha'$ as
\[
\alpha' \colon F^e_*W_nR(p^e\{D\}) \xrightarrow{\times F^e_*[c]^{-1}} F^e_*W_nR(p^e\{D\}+\mathrm{div}(c)) \xrightarrow{\beta} W_n\omega_R(-K_R),
\]
then we obtain the following commutative diagram:
\begin{equation}\label{eq:qFs-new-def}
\begin{tikzcd}
    W_nR(\{D\}) \arrow[r,"\varphi_{\{D\}}"] \arrow[d,"{\cccred \mathbf{R}}^{n-1}"] & F^e_*W_nR(p^e\{D\}) \arrow[ldd,"\alpha'"] \\
    R \arrow[d] \\
    W_n\omega_R(-K_R)
\end{tikzcd}
\end{equation}
By tensoring (\ref{eq:qFs-new-def}) with $W_nR(\rdown{D})$ and taking 
the $S_2$-fications, we obtain the required commutative diagram.
\end{proof}


\begin{proposition}\label{prop:  qFreg 0-map}
We use Notation \ref{n-normal-local}. 
{Take integers $n>0$ and $m>0$.} 
Let $\Delta$ be an effective $\Q$-divisor on $X$ such that $(R, \Delta)$ 
is $n$-quasi-$F$-regular. 
{Assume that 
the inequality 
{$ \{p^rD\} \leq \Delta$} holds 
for every integer $r \geq 0$.}  
Then the induced $W_{m+n-1}R$-module homomorphism
\[
{\cccred \mathbf{R}}^{n-1} \colon \wt{0^*_{D,m+n-1}} \to \wt{0^*_{D,m}}
\]
is zero. 
\end{proposition}

\begin{proof}
Let us show the assertion by induction on $m$. {We emphasise that the assumptions of the proposition are valid for $D$ replaced by $p^kD$ for any integer $k \geq 0$.}

First, we treat the base case $m=1$ of this induction. 
Fix $t_D \in R^{\circ}$ such that $t_D \in \tau(R, \{p^e  D\})$ 
for every integer $e \geq 0$. 
Then there exist $c \in R^\circ \cap (t_D^{4})$ and an integer $e_0 \geq 0$ such that $\wt{0^*_{D,n}}=\wt{K^{e,c}_{D,n}}$ for every integer $e \geq e_0$ 
(Proposition {\ref{prop: test elem pair tilda}}).
Since $(R, \{D\})$ is quasi-$F$-regular, 
we can find an integer $e \geq e_0$ 
such that $(R,\{D\}+(1/p^{e})\mathrm{div}(c))$  is quasi-$F^{e}$-split (Theorem \ref{thm:chara qFr pair}). 
By Proposition \ref{prop:qFs-divisor-split}, there exists a $W_nR$-module homomorphism $\alpha \colon F^e_*W_nR(p^eD) \to W_n\omega_R(\rdown{D}-K_R)$ 
{fitting inside the commutative diagram as in} (\ref{e1:qFs-divisor-split}). 
By taking $H^d_{\m}$, we obtain the following commutative diagram:
\[
\begin{tikzcd}
    H^d_{\m}(W_nR(D)) \arrow[r,"{(\cdot F^e_*[c])}\circ F^e"] \arrow[d,"{\cccred \mathbf{R}}^{n-1}"] & H^d_{\m}(F^e_*W_nR(p^{e}D)) \arrow[ldd,"\alpha'"] \\
    H^d_{\m}(R(D)) \arrow[d,"\beta"] \\
    H^d_{\m}(W_n\omega_R(-K_R+\rdown{D})).
\end{tikzcd}
\]
By ${(\cdot F^e_*[c])}\circ F^e (\wt{0^*_{D,n}})=0$, 
we get $\beta \circ {\cccred \mathbf{R}}^{n-1}(\wt{0^*_{D,n}})=0$, that is, 
\[
{\cccred \mathbf{R}}^{n-1}(\wt{0^*_{D,n}}) \subseteq \mathrm{Ker}(\beta).
\]
It suffices to show that $\beta$ is injective. 
Applying the Matlis duality $\Hom_{W_nR}(-, E)$ to $\beta$, 
we get 
\[
{\cccred \mathbf{R}}^{n-1} : W_nR(K_R-\rdown{D}) \to R(K_R-\rdown{D}),
\]
which is surjective. Hence $\beta$ is injective. 
This completes the proof for the case when $m=1$. 

Next, we assume $m \geq 2$. 
Consider the following commutative diagram in which each horizontal sequence is exact:
\[
\begin{tikzcd}
    H^d_{\m}(F_*W_{m+n-2}R(pD)) \arrow[r,"{\cccred \mathbf{V}}"] \arrow[d,"{\cccred \mathbf{R}}^{n-1}"] & H^d_{\m}(W_{m+n-1}R(D)) \arrow[r,"{\cccred \mathbf{R}}^{m+n-2}"] \arrow[d,"{\cccred \mathbf{R}}^{n-1}"] & H^d_{\m}(R(D)) \arrow[d,equal] \\
    H^d_{\m}(F_*W_{m-1}R(pD)) \arrow[r, "{\cccred \mathbf{V}}"] & H^d_{\m}(W_mR(D)) \arrow[r, "{\cccred \mathbf{R}}^{m-1}"] & H^d_{\m}(R(D)).
\end{tikzcd}
\]
We have 
\[
{\cccred \mathbf{R}}^{m+n-2}(\wt{0^*_{D,m+n-1}}) 
\overset{{\rm (i)}}{=} 0^*_{D,m+n-1}
\overset{{\rm (ii)}}{\subseteq} 0^*_{D,n}
\overset{{\rm (iii)}}{=} {\cccred \mathbf{R}}^{n-1}(\wt{0^*_{D, n}}) 
\overset{{\rm (iv)}}{=} 0, 
\]
where (i) and (iii) hold by Theorem \ref{thm: lifting thm pair}, 
(ii) follows from $m+n-1 \geq n$, and 
(iv) is valid by the case when $m=1$, 
which has been settled already. 
By the exactness of the upper horizontal sequence and 
${\cccred \mathbf{V}}^{-1}(\wt{0^*_{D,m+n-1}})=F_*\wt{0^*_{pD,m+n-2}}$ (Proposition \ref{prop: tight cl l=n-log}(2)), 
it is enough to show that ${\cccred \mathbf{R}}^{n-1}(\wt{0^*_{pD,m+n-2}})=0$. 
This holds by the induction hypothesis. 
\qedhere
\end{proof}


\section{Quasi-test ideals via alterations}

In this section, we investigate the relationship between quasi-test ideals and alterations. Our goal is to prove the implication: quasi-$F$-regular 
$\Rightarrow$ klt (Theorem \ref{thm:qFR to klt pair}). This  is established as a consequence of the inclusion $\tau_n(X, \Delta) \subseteq \mathcal J(X, \Delta)$ (Corollary \ref{cor:trace under bir}).

{Moreover}, {in the log $\bQ$-Gorenstein case, we will show that} for a sufficiently high finite cover $f: Y \to X$, the $n$-quasi-test ideal $\tau_n(X, \Delta)$ coincides with the image of the trace map 
of 
$\Phi^f_{X, K_X+\Delta, n} : \mathcal{O}_X(K_X+\Delta) \to Q^f_{X, K_X+\Delta, n}$ (Theorem \ref{thm: chara via alteration, quasi}). This observation leads to the conclusion that $(X, \Delta)$ is $n$-quasi-$F$-regular if and only if $(R, \Delta)$ is $n$-quasi-$+$-regular (Theorem \ref{thm: rel q+R, qFrat and qFR}).


As it is not natural to assume in this section that our schemes are affine, 
we start by globalising the local concepts introduced earlier 
(Subsection \ref{ss global q test}).




\subsection{Quasi-test ideals in the global case}\label{ss global q test}

In this subsection, we globalise some notions obtained in earlier parts. 

\begin{definition}\label{d-global-test-submod}
We use Notation \ref{n-global}. 
Take an integer $n>0$ and a $\Q$-divisor $D$ on $X$. 
For the generic point $\xi$ of $X$, we denote the stalk of $W_n\omega_X$ at $\xi$ by
$W_n\omega_{X, \xi}$, and we 
let  $\underline{W_n\omega_{X, \xi}}$ be the corresponding constant sheaf on $X$. Set $\underline{\omega_{X, \xi}} := \underline{W_1\omega_{X, \xi}}$. 
\begin{enumerate}
\item 
$\tau(W_n\omega_X, D)$ is defined as the coherent $W_n\MO_X$-submodule 
of $\underline{W_n\omega_{X, \xi}}$ such that 
for every affine open subset $U$ of $X$ and $R := \Gamma(U, \MO_X)$, 
the following equality 
\[
\Gamma(U, \tau(W_n\omega_X, D)) = \tau(W_n\omega_{R}, D|_U), 
\]
of  $R$-submodules of the $R$-module 
$W_n\omega_{X, \xi} = W_n\omega_{R} \otimes_{W_n{R}} W_n(K(R))$ holds. 
The existence of $\tau(W_n\omega_X, D)$ is assured by Proposition \ref{prop: test submod pair compat}. 
\item 
$\tau_n(\omega_X, D)$ is defined as the coherent $\MO_X$-submodule of 
$\underline{\omega_{X, \xi}}$ such that 
for every affine open subset $U$ of $X$ and $R := \Gamma(U, \MO_X)$, 
the following equality 
\[
\Gamma(U, \tau_n(\omega_X, D)) = \tau_n(\omega_R, D|_U), 
\]
of $R$-submodules of the $R$-module 
$\omega_{X, \xi} = \omega_{R} \otimes_{R} K(R)$ holds. 
Equivalently, $\tau_n(\omega_X, D)$ is defined as 
\[
\tau_n(\omega_X, D) = (({\cccred \mathbf{R}}^{n-1})^*)^{-1} (\tau(W_n\omega_X, D)), 
\]
where 
$({\cccred \mathbf{R}}^{n-1})^* : \underline{\omega_{X, \xi}} \to \underline{W_n\omega_{X, \xi}}$. 
\end{enumerate} 
\end{definition}

\begin{definition}\label{d-global-test-ideal}
We use Notation \ref{n-global}. 
Let $\Delta$ be a $\Q$-divisor on $X$ such that $K_X+\Delta$ is $\Q$-Cartier. 
Take an integer $n>0$. 
Set 
\begin{align*}
\tau(W_nX,\Delta)&:=\tau(W_n\omega_X,K_X+\Delta) \subseteq W_n\omega_X(-(K_X+\Delta))\\
\tau_n(X, \Delta)&:= \tau_n(\omega_X,K_X+\Delta) \subseteq 
\MO_X({\rup{-\Delta}}) \simeq \omega_X(-(K_X+\Delta))\\
\tau^q(X, \Delta)&:= \tau^q(\omega_X,K_X+\Delta) \subseteq 
{{\ccred \MO_X}}({\rup{-\Delta}}) \simeq \omega_{{\ccred X}}(-(K_X+\Delta)). 
\end{align*} 
In particular, 
$\tau_n(X, \Delta)$ and $\tau^q(X, \Delta)$ are coherent $\MO_X$-submodules of $\MO_X({\rup{-\Delta}})$, 
which we call the {\em $n$-quasi-test ideal} and 
the {\em quasi-test ideal} of $(X, \Delta)$, 
respectively. 
\end{definition}

 For an affine open subset $U$ of $X$ and $R := \Gamma(U, \MO_X)$, it is easy to see that 
\begin{align*}
\Gamma(U, \tau(W_nX, \Delta)) &= \tau(W_nR, \Delta|_{U}),\\ 
\Gamma(U, \tau_n(X, \Delta)) &= \tau_n(R, \Delta|_{U}),\\  
\Gamma(U, \tau^q(X, \Delta)) &= \tau^q(R, \Delta|_{U}).
\end{align*}

\begin{proposition}\label{p-test-tensorH}
We use Notation \ref{n-global}. 
Take a $\Q$-divisor $\Delta$ such that $D:=K_X+\Delta$ is $\Q$-Cartier. 
Let $H$ be a Cartier divisor on $X$. 
Then the following hold. 
\begin{enumerate}
\item 
$\tau(W_n\omega_X, D+H) = \tau(W_n\omega_X, D) \otimes_{W_n\MO_X} W_n\MO_X(-H)$. 
\item 
$\tau_n(\omega_X, D+H) = \tau_n(\omega_X, D) \otimes_{\MO_X} \MO_X(-H)$. 
\item 
$\tau(W_nX, \Delta + H) = \tau(W_nX, \Delta) \otimes_{W_n\MO_X} W_n\MO_X(-H)$. 
\item 
$\tau_n(X, \Delta + H) = \tau_n(X, \Delta) \otimes_{\MO_X} \MO_X(-H)$. 
\end{enumerate}
\end{proposition}

\begin{proof}
Let us show (1). 
Recall that 
$\tau(W_n\omega_X, D+H)\subseteq \underline{W_n\omega_{X, \xi}}$ 
and $\tau(W_n\omega_X, D)\subseteq \underline{W_n\omega_{X, \xi}}$ 
(Definition \ref{d-global-test-submod}(1)), 
which implies 
\[
\tau(W_n\omega_X, D) \otimes_{W_n\MO_X} W_n\MO_X(-H) \subseteq 
\underline{W_n\omega_{X, \xi}}\otimes_{W_n\MO_X} W_n\MO_X(-H)
= \underline{W_n\omega_{X, \xi}}. 
\]
Therefore, both sides of (1) are coherent $W_n\MO_X$-submodules 
of the constant sheaf $\underline{W_n\omega_{X, \xi}}$. 
Taking an affine open cover of $X$ on which $H$ is principal, 
we may assume that $X = \Spec R$ and 
$H = \div(h)$ for some $h \in K(R)^{\times}$. 
There exist $h_1, h_2 \in R^{\circ}$ satisfying $h=h_1/h_2$. 
Hence the problem is reduced to the case when $h \in R^{\circ}$, in which it holds that 
\[
{ \tau(W_n\omega_X, D) \otimes_{W_n\MO_X} W_n\MO_X(-H) = [h] \cdot \tau(W_n\omega_X, D).} 
\]
Then the assertion follows from Proposition \ref{prop: test submodule gen pair}(2), and so (1) holds.

The assertion (2) follows from (1) and 
$\tau_n(\omega_X, D) = (({\cccred \mathbf{R}}^{n-1})^*)^{-1} (\tau(W_n\omega_X, D))$ (Definition \ref{d-global-test-submod}(2)). 
Then (1) and (2) imply (3) and (4), respectively (Definition \ref{d-global-test-ideal}). 
\end{proof}

\begin{corollary}\label{cor:openness-qFr-locus}
We use Notation \ref{n-global}.
Take an effective $\Q$-divisor $\Delta$ such that $K_X+\Delta$ is $\Q$-Cartier and $\rdown{\Delta}=0$.
Fix an integer $n \geq 1$ and a point $x \in X$. 
{Then} $(\Spec \cO_{X,x},\Delta_x)$ is $n$-quasi-$F$-regular 
(resp.\ quasi-$F$-regular) if and only if 
$\tau_n(X,\Delta)_x=\cO_{X,x}$ (resp.\ $\tau^q(X,\Delta)_x=\cO_{X,x}$). 
In particular, the $n$-quasi-$F$-regular locus and quasi-$F$-regular locus of $X$ are open subsets of $X$.
\end{corollary}

\begin{proof}
We take 
an affine open subset $U = \Spec R$ of $X$ such that $x \in U$.
By construction, we have $\tau_n(X,\Delta)_x=\tau_n(U,(\Delta|_{U}))_x = \tau_n(R, (\Delta|_U))_x$. 
Thus, by Proposition \ref{prop: test ideal pair non-tilda}, $(\Spec \cO_{X,x},\Delta_x)$ is $n$-quasi-$F$-regular if and only if $\tau_n(X,\Delta)_x=\cO_{X,x}$.
Therefore, the $n$-quasi-$F$-regular locus of $X$ coincides with $\{x \in X \ \mid\ \tau_n(X,\Delta)_x=\cO_{X,x}\}$, thus it is open, as desired.
The remaining assertion  follows from a similar argument as above.
\end{proof}

{\cred

\begin{proposition}\label{etale-cov}
We use Notation \ref{n-global}. 
Let $f \colon Y \to X$ be an \'etale morphism, where $Y$ is an integral normal Noetherian scheme. 
Fix an integer $n>0$ and let $D$ be a $\Q$-Cartier $\Q$-divisor on $X$.
Then the following hold. 
\begin{enumerate}
\item 
$(W_nf)^*(\tau(W_n\omega_X,D)) 
=  \tau\!\bigl(W_n\omega_Y,f^*D\bigr)$, where $W_nf : (Y, W_n\MO_Y) \to (X, W_n\MO_X)$ denotes the induced \'etale morphism of Noetherian schemes. 
\item 
$f^*\tau_n(\omega_X,D) 
=\tau_n\!\bigl(\omega_Y,f^*D\bigr)$. 
\item 
$f^*\tau_n(X, \Delta) 
= \tau_n(Y, f^*\Delta)$ 
for an effective $\Q$-divisor $\Delta$ on $X$ such that $(K_X+\Delta)$ is $\Q$-Cartier. 
\end{enumerate}
\end{proposition}


\begin{proof}
We may assume that $X$ and $Y$ are affine: $X = \Spec R$ and $Y = \Spec S$.

Let us show (1). 
By Definition \ref{d-tau-W_nomega-non-eff}, we may assume that $D$ is effective.
Choose $a \in \Q_{>0}$ and $g \in R^\circ$ such that $D=a\,\mathrm{div}(g)$.
Take $d \in R^\circ$ with $D \leq \mathrm{div}(d)$.
Pick $t' \in R^\circ$ such that $R_{t'}$ is regular; then $S_{t'}$ is also regular.
By \cite{HH94}*{Theorem~5.10}, there exists a positive integer $l$ such that $t:=(t')^l$ is a test element for both $R$ and $S$.
Set $c:=dt^2$.
By Proposition \ref{prop: test submodule gen pair}(4), we obtain
\[
\tau(W_n\omega_R,D) \otimes_{W_n(R)} W_n(S)
=\sum_{e \geq 0} T^e_n\!\bigl(F^e_*([c\,g^{\rup{a(p^e-1)}}] \cdot W_n\omega_R)\bigr) \otimes_{W_n(R)} W_n(S).
\]
We have 
\[
(F^e \colon W_n(R) \to F^e_*W_n(R)) \otimes_{W_n(R)} W_n(S)
= (F^e \colon W_n(S) \to F^e_*W_n(S))
\]
by \cite{KTY22}*{Lemma~2.4}. 
By Grothendieck duality, it follows that
\[
(T^e_n \colon F^e_*W_n\omega_R \to W_n\omega_R) \otimes_{W_n(R)} W_n(S)
= (T^e_n \colon F^e_*W_n\omega_S \to W_n\omega_S)
\]
for every integer $e\geq 0$ (here we used the fact that an \'etale morphism is decomposed into an open immersion and a finite morphism).
Therefore, 
\[
\sum_{e \geq 0} T^e_n(F^e_*([c\,g^{\rup{a(p^e-1)}}] \cdot W_n\omega_R)) \otimes_{W_n(R)} W_n(S)
\simeq
\sum_{e \geq 0} T^e_n(F^e_*([c \, g^{\rup{a(p^e-1)}}] \cdot W_n\omega_S)),
\]
which is isomorphic to $\tau(W_n\omega_S,f^*D)$ by Proposition \ref{prop: test submodule gen pair}(4).
Therefore $\tau(W_n\omega_R,D)\otimes_{W_n(R)} W_n(S)=\tau(W_n\omega_S,f^*D)$. 
Thus (1) holds. 

Let us show (2). 
By the definition of $\tau_n(\omega_R,D)$, the natural map $\omega_R \to W_n\omega_R$ induces an injection
\[
({\cccred \mathbf{R}}^{n-1})^*\colon  \omega_R/\tau_n(\omega_R,D) \hookrightarrow W_n\omega_R/\tau(W_n\omega_R,D).
\]
Tensoring with $W_n(S)$ over $W_n(R)$ and using the first assertion, we obtain an injection
\[
({\cccred \mathbf{R}}^{n-1})^*\colon\omega_S/(\tau_n(\omega_R,D) \otimes_R S) \hookrightarrow W_n\omega_S/\tau(W_n\omega_S,f^*D).
\]
Hence
\[
 \tau_n(\omega_S,f^*D)
=(({\cccred \mathbf{R}}^{n-1})^*)^{-1}\!\bigl(\tau(W_n\omega_S,f^*D)\bigr)
=\tau_n(\omega_R,D) \otimes_R S. 
\]
Thus (2) holds. The assertion (3) immediately follows from (2). 
\end{proof}
}

\subsection{Behaviour of quasi-test ideals under alterations}

\begin{proposition}\label{prop: under fin pair}
We use Notation \ref{n-global}.
Let $f \colon Y \to X$ be a finite surjective morphism from an integral normal 
 scheme $Y$,
and let $D$ be a $\Q$-Cartier $\Q$-divisor on $X$.
Then 
\[
T^f_n(\tau(W_n\omega_Y,f^*D))=\tau(W_n\omega_X,D), 
\]
where $T^f_n : f_* W_n\omega_Y \to W_n\omega_X$ denotes the trace map, obtained from {the pullback homomorphism} $W_n\MO_X \to f_*W_n\MO_Y$ by 
applying $\cHom_{W_n\MO_X}(-, W_n\omega_X)$.  
\end{proposition}

\begin{proof}
We may assume that $X=\Spec R$, $Y=\Spec S$, and 
$D$ is effective.
Set $J_Y:=\mathrm{Im}(T^f_n\colon f_*W_n\omega_Y(-f^*D) \to W_n\omega_X(-D))$. 
For the time being, we finish the proof by assuming Claim below. 

\begin{claim*}
$J_Y$ is a co-small $W_n\MO_X$-submodule of $W_n\omega_X(-D)$.
\end{claim*}

\noindent 
After taking a suitable affine cover of $X$, 
we may assume that $D=a\mathrm{div}(g)$ for some $a \in \Q_{\geq 0}$ and $g \in R^{\circ}$.
By Proposition \ref{prop: test submodule gen pair}(4), there exists $c \in R^{\circ}$ such that
\begin{align*}
    \tau(W_n\omega_X,D)&=\sum_{e \geq 0} \PsiT^e_{X,n}(F^e_*[cg^{\rup{a(p^e-1)}}] \cdot J_Y), \\
    \tau(W_n\omega_Y,f^*D)&=\sum_{e \geq 0} \PsiT^e_{Y,n}(F^e_*[cg^{\rup{a(p^e-1)}}] \cdot W_n\omega_Y(-f^*D)),
\end{align*}
where $T_{X, n}^e$ and $T^e_{Y, n}$ denote the trace maps of $X$ and $Y$, respectively: 
\[
T^e_{X, n}: F^e_*W_n\omega_X \to W_n\omega_X, \qquad 
T^e_{Y, n}: F^e_*W_n\omega_Y \to W_n\omega_Y.
\]
Therefore, we have
\begin{align*}
    T^f_n(f_*\tau(W_n\omega_Y,f^*D))
    &= \sum_{e \geq 0} T^f_n(\PsiT^e_{Y,n}(F^e_*[cg^{\rup{a(p^e-1)}}] \cdot W_n\omega_Y(-f^*D))) \\
    &= \sum_{e \geq 0}\PsiT^e_{X,n} (T^f_n(F^e_*[cg^{\rup{a(p^e-1)}}] \cdot W_n\omega_Y(-f^*D))) \\
    &= \sum_{e \geq 0}\PsiT^e_{X,n}(F^e_*[cg^{\rup{a(p^e-1)}}] \cdot T^f_n( W_n\omega_Y(-f^*D))) \\
    &= \sum_{e \geq 0}\PsiT^e_{X,n}(F^e_*[cg^{\rup{a(p^e-1)}}] \cdot J_Y)=\tau(W_n\omega_X,D).
\end{align*}
Thus  it is enough to prove Claim. 
By shrinking $X$, we may assume that 
\begin{itemize}
    \item $X$ and $Y$ are regular, 
    \item $D=0$, and 
    \item $\sO_X \to f_*\sO_Y$ splits as an $\MO_X$-module homomorphism (note that the induced field extension $K(X) \hookrightarrow K(Y)$ splits as a $K(X)$-linear map).
\end{itemize}
Then we obtain the following commutative diagram in which each horizontal sequence is  exact: 
\[
\xymatrix{
0 \ar[r] & f_*\omega_Y \ar[r] \ar[d]^-{T_{1}^f} &  f_*W_n\omega_Y \ar[r] \ar[d]^-{T^f_n} & f_*F_*W_{n-1}\omega_Y \ar[r] \ar[d]^-{F_*T_{n-1}^f} & 0 \\
0 \ar[r] & \omega_X \ar[r] & W_n\omega_X \ar[r] & F_*W_{n-1}\omega_{{\cccred X}} \ar[r] & 0.
}
\]
Since $\sO_X \to f_*\sO_Y$ splits, $T_{1}^f$ is surjective.
By the snake lemma and the induction on $n$, $T^f_n$ is surjective for every integer $n>0$. 
This completes the proof of Claim. 
\end{proof}

\begin{theorem}\label{thm: chara via alteration}
We use Notation \ref{n-global}. 
Let $D$ be a $\Q$-Cartier $\Q$-divisor on $X$.
Then the following hold. 
\begin{enumerate}
    \item For every integer $n>0$ and every alteration $f \colon Y \to X$ from an integral normal scheme $Y$, 
    it holds that 
    \begin{equation}\label{e1-chara via alteration}
    \tau(W_n\omega_X,D) \subseteq \mathrm{Im}(f_*W_n\omega_Y(-f^*D) \to W_n\omega_X(-D)).
    \end{equation}
    \item For every integer $n>0$, 
    there exists a finite separable morphism $f \colon Y \to X$ 
    from an integral  normal  scheme $Y$ such that 
    \[
    \tau(W_n\omega_X,D) =\mathrm{Im}(f_*W_n\omega_Y(-f^*D) \to W_n\omega_X(-D)).
    \]
\end{enumerate}
\end{theorem}

\begin{proof}
We may assume that $X=\Spec R$. 

Let us show (1). 
We may assume that $K_X$ {and $D$ are}
effective. 
We consider the diagram
    \begin{equation}\label{e2-chara via alteration}
\xymatrix{
0 \ar[r] & f_*\omega_Y(-f^*D) \ar[r]^-{({\cccred \mathbf{R}}^{n-1})^*} \ar[d]^{T^{f}_1} & f_*W_n\omega_Y(-f^*D) \ar[r]^-{{\cccred \mathbf{V}}^*} \ar[d]^{T^{f}_n} & 
f_*F_*W_{n-1}\omega_Y(-f^*D) \ar[d]^{T^{f}_{n-1}} \\
0 \ar[r] & \omega_X(-D) \ar[r]^-{({\cccred \mathbf{R}}^{n-1})^*} & W_n\omega_X(-D) \ar[r]^-{{\cccred \mathbf{V}}^*} & F_*W_{n-1}\omega_X(-D).
}
\end{equation}
First, we prove that the image $J_n$ of $T^{f}_n$ is co-small and $C_n(D)$-stable.
Since $f$ is finite over a suitable open subset on $X$, 
$J_n$ is co-small by the proof of Proposition \ref{prop: under fin pair}. 
{To prove that it is $C_n(D)$-stable,} take $\varphi \in C^e_n(D)$. 
Then there exists $g \in \sO_X(-(p^e-1)D)$ such that 
$\varphi=\PsiT^e_{X,n} \circ (\cdot F^e_*[g])$: 
\[
\varphi=\PsiT^e_{X,n} \circ (\cdot F^e_*[g]) : F_*^eW_n\omega_X((p^e-1)D) \xrightarrow{\cdot F_*^e[g]} F_*^eW_n\omega_X \xrightarrow{T^e_{X,n}} W_n\omega_X. 
\]
For
\[
\varphi_Y:=\PsiT^e_{Y,n} \circ (\cdot F^e_*[g]) \in \Hom_{W_n\MO_Y}(F^e_*W_n\omega_Y((p^e-1)f^*D),W_n\omega_Y),
\]
we have
\begin{align*}
    \varphi(F^e_*J_n) 
    &= \varphi \circ T^{f}_n(F^e_*W_n\omega_Y(-f^*D)) \\
    &=T^{f}_n(\varphi_Y(F^e_*W_n\omega_Y(-f^*D))) \\
    &\subseteq T^{f}_n(W_n\omega_Y(-f^*D))=J_n,
\end{align*}
where the inclusion follows from 
Proposition \ref{prop: first prop of C_n(D)-sub}$(1)$. 
Therefore, we have
\[
\tau(W_n\omega_X,D) \subseteq J_n
\]
by minimality (Proposition \ref{prop: test submodule gen pair}(1)).  
Thus (1) holds.


Let us show (2). 
We first treat the case when $D = 0$. 
By (1) and Noetherian induction, we may assume that $R$ is a local ring. 
However, being local will not be stable under the following argument, and hence we only assume that $R$ is a semi-local ring, that is, $R$ has only finitely many maximal ideals $\m_1, ...,  \m_r$. 
Set $\m$ to be its Jacobson radical ($\m := \m_1 \cap \cdots \cap \m_r$). 
For $X_i := \Spec \MO_{X, \m_i}$, 
\begin{itemize}
\item[(A)] it holds that 
$H^d_{\m}(W_n\MO_X) = \bigoplus_{i=1}^rH^d_{\m_1}(W_n\MO_{X_1}) \oplus \cdots \oplus H^d_{\m_r}(W_n\MO_{X_r})$, and 
 \item[(B)] we set 
 \[
 \wt{0^*_{X, n}} := \bigoplus_{i=1}^r\wt{0^*_{X_i, n}}\subseteq 
\bigoplus_{i=1}^r H^d_{\m_i}(W_n\MO_{X_i}) = H^d_{\m}(W_n\MO_X), 
 \]
\end{itemize}
where (A) follows from $H^d_{\m}(W_nR) \simeq H^d_{\widehat{\m}}(W_n\widehat{R})$ \cite[Proposition 3.5.4(d)]{BH93} and 
$\widehat{R} \simeq \widehat{R_{\m_1}} \times \cdots \times  \widehat{R_{\m_r}}$ 
\cite[Theorem 17.7]{Nag62}. 
Consider the following statement: 
\begin{enumerate}
\item[$(\star_n)$] 
there exists a separable finite cover $f \colon Y \to X$ such that $f^*(\wt{0^{*}_{X, n}})=0$ for  the pullback map
\[
f^* : H^d_\m(W_n\sO_X) \to H^d_\m(f_*W_n\sO_Y). 
\]
\end{enumerate}

We now prove $(\star_n)$ by induction on $n$.
The base case $(\star_{1})$ of this induction 
follows from \cite[Theorem 3.2(b)]{bst}. 
Assume $n \geq 2$. 
By $(\star_{1})$, 
we can find 
a separable finite cover $f \colon Y \to X$ such that 
$f^*(\wt{0^{*}_{X, 1}})=0$. Recall that we have 
\[
H^d_{\m}(f_*\MO_Y) \overset{{\rm (i)}}{=} H^d_{\m \MO_Y}(\MO_Y) 
\overset{{\rm (ii)}}{=} H^d_{\m_Y}(\MO_Y)
\]
where $\m_Y$ denotes the Jacobson radical of $Y$, 
(i) follows from 
\cite{stacks-project}*{\href{https://stacks.math.columbia.edu/tag/0952}{Tag 0952}}
, and (ii) holds by $\sqrt{\m \MO_Y} = \m_Y$.  
By applying $(\star_{n-1})$ to $Y$, 
there exists a separable finite cover $g \colon Z \to Y$ such that 
$g^*(\wt{0^{*}_{Y, n-1}})=0$ for 
\[
g^* : H^d_{\m_Y}(W_{n-1}\sO_{Y}) \to H^d_{\m_Y}(g_*W_{n-1}\sO_{Z}). 
\]

Pick $x \in \wt{0^{*}_{X, n}} \subseteq H^d_\m(W_n\sO_X)$.
By the choice of $f$, we have ${\cccred \mathbf{R}}^{n-1}(f^*x) = f^*{\cccred \mathbf{R}}^{n-1}(x)=0$.
Then there exists $y \in H^d_\m(f_*F_*W_{n-1}\sO_Y)$ such that ${\cccred \mathbf{V}}y=f^*x$.
By (A) and Proposition \ref{prop: tight cl l=n}, 
we obtain $y \in F_*\wt{0^{*}_{Y, n-1}} \subseteq  H^d_\m(f_*\sO_Y)$.
By the choice of $g$, we have $g^*y=0$. 
Therefore, $\wt{0^{*}_{X,n}}$ maps to zero via the composite pullback map
\[
g^*f^* : H^d_\m(W_n\sO_X) \to H^d_\m(f_*g_*W_n\sO_Z).
\]
This completes the proof of $(\star_n)$. 

Let us settle the case when $D=0$ by using $(\star_n)$. 
Note that  $(\star_n)$ induces the map 
\[
f^* : 
H^d_\m(W_n\sO_X) 
\twoheadrightarrow
\frac{H^d_\m(W_n\sO_X)}{\wt{0^{*}_{n}}} \to H^d_\m(f_*\sO_Y).
\]
By applying the Matlis duality functor $\Hom_{W_nR}(-, E)$, 
we get the following 
\[
W_n\omega_X \hookleftarrow \tau(W_n\omega_X) \leftarrow f_*W_n\omega_Y,
\]
where the middle term can be computed by Proposition \ref{prop: corr 0^* and tau: no pair}. 
This, together with (1), completes the proof for the case when $D=0$. 


Let us reduce the general case to this case. 
There exists a separable finite cover $f \colon Y \to X$ from an integral normal scheme $Y$ such that $f^*D$ is Cartier \cite[Lemma 4.5]{bst}. 
Since the case when $D=0$ has been settled already, 
we can find  a finite separable cover $g \colon Z \to Y$ such that
\[
\tau(W_n\omega_Y)=\mathrm{Im}(g_*W_n\omega_Z \to W_n\omega_Y).
\]
By the following isomorphism (Proposition \ref{p-test-tensorH}):
\begin{align*}
    \tau(W_n\omega_Y,f^*D) &\simeq \tau(W_n\omega_Y) \otimes_{W_n\MO_Y} W_n\sO_Y(-f^*D),
\end{align*}
we get 
\[
\tau(W_n\omega_Y,f^*D)=\mathrm{Im}(g_*W_n\omega_Z(-g^*f^*D) \to W_n\omega_Y(-f^*D)).
\]
It follows from  Proposition \ref{prop: under fin pair} that 
\[
\tau(W_n\omega_X,D)=\mathrm{Im}(f_*g_*W_n\omega_Z(-g^*f^*D) \to W_n\omega_X(-D)).
\]
\end{proof}

\begin{corollary}\label{cor:trace under bir}
We use Notation \ref{n-global}.  
Let $f \colon Y \to X$ be an alteration 
from an integral normal scheme $Y$. 
Let $\Delta$ be a  $\Q$-divisor on $X$ such that $K_X+\Delta$ is $\Q$-Cartier. Then the following hold. 
\begin{enumerate}
\item 
$\tau_n(X, \Delta) \subseteq 
(({\cccred \mathbf{R}}^{n-1})^*)^{-1}(\mathrm{Im}(f_*W_n\omega_Y(-f^*(K_X+\Delta)) \to W_n\omega_X(-(K_X+\Delta))))$ 
\item 
$ \tau_n(X,\Delta) \subseteq \mathcal{J}(X,\Delta)$.  
\end{enumerate}
\end{corollary}

\begin{proof}
By Definition \ref{d-global-test-ideal}, we have 
\[
\tau_n(X, \Delta) = 
\tau_n(\omega_X, K_X+\Delta). 
\]
The assertion (1) follows immediately from Theorem \ref{thm: chara via alteration}(1) by taking the inverse images by $(({\cccred \mathbf{R}}^{n-1})^*)^{-1}$ (cf. (\ref{e2-chara via alteration})). 

Let us show (2). Set $D:=K_X+\Delta$. 
Assume that $f$ is birational. 
Then all the vertical arrows in (\ref{e2-chara via alteration}) are injective. 
Hence the left square in the  diagram (\ref{e2-chara via alteration}) consists of injections. 
In what follows, we consider each $W_n\MO_X$-module in this square as a $W_n\MO_X$-submodule of $W_n\omega_X(-D)$.  
We then   get 
\[
\tau_n(\omega_X, D) \overset{{\rm (i)}}{\subseteq} \omega_X(-D) \cap f_*W_n\omega_Y(-f^*D)
\overset{{\rm (ii)}}{\subseteq} f_*\omega_Y(-f^*D), 
\]
where (i) holds by (1) and (ii) follows by chasing  the diagram (\ref{e2-chara via alteration}) 
(specifically, $\zeta \in \omega_X(-D) \cap f_*W_n\omega_Y(-f^*D) \Rightarrow T^f_{n-1} \circ {\cccred \mathbf{V}}^*(\zeta)=0 \Rightarrow {\cccred \mathbf{V}}^*(\zeta)=0 \Rightarrow \zeta \in f_*\omega_Y(-f^*D)$). 
Therefore, we get $\tau_n(\omega_X,K_X+\Delta) \subseteq \mathcal{J}(X,\Delta)$.
Thus (2) holds. 
\qedhere 



\end{proof}

\begin{theorem}\label{thm:qFR to klt pair}
We use 
Notation \ref{n-global}.
Let $\Delta$ be an effective $\Q$-divisor on $X$ such that $K_X+\Delta$ is $\Q$-Cartier. 
If $(X, \Delta)$ is quasi-$F$-regular, 
then $(X, \Delta)$ is klt. 
\end{theorem}

\begin{proof}
It follows from Corollary \ref{cor:trace under bir}(2) that 
\[
\tau_n(X, \Delta) 
\subseteq 
\mathcal J(X, \Delta) \subseteq \MO_X.  
\]
Assume that $(X, \Delta)$ is quasi-$F$-regular. 
Then $\tau_n(X, \Delta) = \MO_X$ for some $n>0$ (Proposition \ref{prop: test ideal pair non-tilda}(1)), 
which implies $\mathcal J(X, \Delta) = \MO_X$, that is, $(X, \Delta)$ is klt. 
\end{proof}

\subsection{Quasi-$F$-regularity and quasi-$+$-regularity}

{In what follows, we recall the definition of quasi-$+$-regularity from \cite{TWY}.}
\begin{dfn}\label{d-Q+R}
We use Notation \ref{n-global}. 
Let  $\Delta$ be a $\Q$-divisor on $X$. Take an integer $n>0$ 
and a finite surjective morphism $f: Y \to X$ from an integral normal  scheme $Y$.  
We define a $W_n\MO_X$-module $Q^f_{X, \Delta, n}$ and 
a $W_n\MO_X$-module homomorphism 
$\Phi^f_{X, \Delta, n}$ by the following pushout diagram: 
\[
\begin{tikzcd}
W_n\MO_X(\Delta) \arrow[r, "f^*"] \arrow[d, "{\cccred \mathbf{R}}^{n-1}"] & f_*W_n\MO_Y(f^*\Delta) \arrow[d]\\
\MO_X(\Delta)  \arrow[r, "\Phi^f_{X, \Delta, n}"] & {Q^f_{X, \Delta, n}}
\end{tikzcd}
\]
\begin{enumerate}
\item 
We say that $(X, \Delta)$ is {\em globally $n$-quasi-$f$-regular} if 
$\rdown{\Delta} =0$ and the induced map 
\begin{align}\label{e1-Q+R}
    \Hom_{W_n\MO_X}(\Phi^f_{X, \Delta, n}, W_n\omega_X(-K_X)) :
&\Hom_{W_n\MO_X}(Q^f_{X, \Delta, n}, W_n\omega_X(-K_X)) \\
&\longrightarrow
\Hom_{W_n\MO_X}(\MO_X, W_n\omega_X(-K_X))  \nonumber
\end{align}
is surjective, where 
the map (\ref{e1-Q+R}) is obtained by 
applying 
\[
\Hom_{W_n\MO_X}(-, W_n\omega_X(-K_X))
\]
to $\Phi^f_{X, \Delta,n} : \MO_X =\MO_X(\Delta) \to Q^f_{X, \Delta, n}$.  
\item 
We say that $(X, \Delta)$ is {\em globally $n$-quasi-+-regular} if 
$\rdown{\Delta} =0$ and it is $n$-quasi-$f$-regular for every finite surjective morphism $f: Y \to X$ from an integral normal scheme $Y$. 
We say that $(X, \Delta)$ is {\em globally quasi-+-regular} if 
$(X, \Delta)$ is globally $n$-quasi-+-regular for some $n \in \Z_{>0}$. 
\end{enumerate}
When $X$ is affine, we simply say that $(X, \Delta)$ is 
{\em $n$-quasi-$+$-regular} 
(resp.\ {\em quasi-$+$-regular}) if 
it is globally $n$-quasi-$+$-regular (resp.\ globally quasi-$+$-regular). 
\end{dfn}

\begin{theorem}\label{thm: chara via alteration, quasi}
We use Notation \ref{n-global}. 
Let $\Delta$ be a $\Q$-divisor on $X$ such that $K_X+\Delta$ is $\Q$-Cartier.
For an alteration $f \colon Y \to X$ from an integral normal scheme $Y$, 
we define the $\MO_X$-submodule $I^f_n(\Delta)$ of $\sO_X(\rup{-\Delta})$ 
as the image of the $W_n\MO_X$-module homomorphism 
\[
(\Phi^f_{X, K_X+\Delta, n})^* : 
\mathcal{H}om(Q^f_{X,K_X+\Delta,n},W_n\omega_X)  \to \mathcal{H}om(\sO_X(K_X+\Delta),W_n\omega_X)=\sO_X(\rup{-\Delta})
\]
induced by applying $\mathcal Hom_{W_n\MO_X}(-, W_n\omega_X)$  to 
$\Phi^f_{X, K_X+\Delta, n} : \sO_X(K_X+\Delta) \to Q^f_{X,K_X+\Delta,n}$.
Then the following hold. 
\begin{enumerate}
    \item For every integer $n>0$ and 
    every alteration $f \colon Y \to X$ from an integral normal scheme $Y$, 
    it holds that 
    \[
    \tau_n(X,\Delta) \subseteq I^f_n(\Delta). 
    \]
    \item For every integer $n >0$, 
    there exists a finite separable morphism $f \colon Y \to X$ from an integral normal scheme $Y$ such that
    \[
    \tau_n(X,\Delta) =I^f_n(\Delta).
    \]
\end{enumerate}
\end{theorem}

\begin{proof}
By Definition \ref{d-global-test-submod} and Definition \ref{d-global-test-ideal}, 
we  have 
$\tau_n(X,\Delta)=\tau_n(\omega_X,K_X+\Delta)$ 
and 
\[
(({\cccred \mathbf{R}}^{n-1})^*)^{-1}(\tau(W_n\omega_X,K_X+\Delta)) =  \tau_n(\omega_X, K_X+\Delta), 
\]
where 
\[
({\cccred \mathbf{R}}^{n-1})^* : \omega_X(-(K_X+\Delta)) \to W_n\omega_X(-(K_X+\Delta)). 
\]
Fix an alteration $f \colon Y \to X$ from an integral normal scheme $Y$ such that $f^*(K_X+\Delta)$ is Cartier. 
By Theorem \ref{thm: chara via alteration}, 
it is enough to show that
\[
(({\cccred \mathbf{R}}^{n-1})^*)^{-1}(\mathrm{Im}(f_*W_n\omega_Y(-f^*(K_X+\Delta)) \to 
W_n\omega_X(-(K_X+\Delta)))) 
=I^f_n(\Delta). 
\]
Consider the following commutative diagram in which each horizontal sequence is exact: \begin{equation*}
    \begin{tikzcd}
        0 \arrow{r} & F_*W_{n-1}\sO_X(p(K_X+\Delta)) \arrow{r}{{\cccred \mathbf{V}}} \arrow[d, equal] & W_n\sO_X(K_X+\Delta) \arrow[r, "{\cccred \mathbf{R}}^{n-1}"] 
        \arrow[d] & \sO_X(K_X+\Delta)\arrow[d, "{\Phi^f_{X, K_X+\Delta, n}}"] \arrow[r] & 0\\
        0 \arrow[r] & F_*W_{n-1}\cO_X(p(K_X+\Delta)) \arrow[r] & f_*W_n\sO_Y(f^*(K_X+\Delta)) \arrow[r, "\psi"] & Q^f_{X,K_X+\Delta,n} \arrow[r] & 0.
    \end{tikzcd}
\end{equation*}
Taking $(-)^* :=\mathcal{H}om_{W_n\MO_X}(-,W_n\omega_X)$, we get 
\begin{equation*}
    \begin{tikzcd}
        0 \arrow{r}  &\mathcal{H}om(Q^f_{X,K_X+\Delta,n},W_n\omega_X) \arrow{r}{\psi^*} \arrow[d, "(\Phi^f_{X, K_X+\Delta, n})^*"] & f_*W_n\omega_Y(-f^*(K_X+\Delta))  \arrow[d, "T_n^f"] \\
        0 \arrow[r] & \omega_X(-(K_X+\Delta))
        \arrow[r, "({\cccred \mathbf{R}}^{n-1})^*"] & W_n\omega_X(-(K_X+\Delta)). 
    \end{tikzcd}
\end{equation*}
By construction {(specifically, because the leftmost vertical arrow in the first diagram of this proof is an equality)}, we get that $\Coker (\psi^*) \to \Coker (({\cccred \mathbf{R}}^{n-1})^*)$ 
is injective. 
Therefore, we obtain 
\[
(({\cccred \mathbf{R}}^{n-1})^*)^{-1}(\mathrm{Im}(T_n^f)) 
= 
{\rm Im}((\Phi^f_{X, K_X+\Delta, n})^*) = I^f_n(\Delta),
\]
as required. 
\qedhere 
\end{proof}


\begin{theorem}\label{thm: rel q+R, qFrat and qFR}
We use Notation \ref{n-normal-non-local}.  
Let $\Delta$ be an effective $\Q$-divisor such that 
$\rdown{\Delta}=0$ and 
$K_R+\Delta$ is $\Q$-Cartier. 
Take an integer $n>0$. 
Then $(R,\Delta)$ is $n$-quasi-$F$-regular if and only if $(R,\Delta)$ is $n$-quasi-+-regular.
\end{theorem}

\begin{proof}
Set $X := \Spec R$. 
Note that the map (\ref{e1-Q+R}) is the same as 
\begin{align}\label{e1-rel q+R, qFrat and qFR}
    \Hom_{W_n\MO_X}(\Phi^f_{X, K_X+\Delta, n}, W_n\omega_X):
&\Hom_{W_n\MO_X}(Q^f_{X, K_X+\Delta, n}, W_n\omega_X) \\ 
&\longrightarrow 
\Hom_{W_n\MO_X}(\MO_X(K_X+\Delta), W_n\omega_X).  \nonumber
\end{align}
Let $I^f_n(\Delta)$ be the image of this $W_n\MO_X$-module homomorphism. 
It is enough to prove that the following statements are equivalent. 
\begin{enumerate}
\renewcommand{\labelenumi}{(\roman{enumi})}
\item $(X,\Delta)$ is $n$-quasi-$F$-regular. 
\item $\tau_n(X, \Delta) = \MO_X$. 
\item $I^f_n(\Delta) = \MO_X$ for every finite surjective  morphism $f \colon Y \to X$ from an integral normal scheme $Y$. 
\item  $(X,\Delta)$ is $n$-quasi-+-regular
\end{enumerate}
The equivalences ${\rm (i)} \Leftrightarrow {\rm (ii)}$ 
and ${\rm (iii)} \Leftrightarrow {\rm (iv)}$ 
follow from Proposition \ref{prop: test ideal pair non-tilda} and Definition \ref{d-Q+R}, respectively. 
We get ${\rm (ii)} \Leftrightarrow {\rm (iii)}$  by Theorem \ref{thm: chara via alteration, quasi} 
{\ccred (indeed,  
the implication ${\rm (ii)} \Rightarrow {\rm (iii)}$ 
(resp.\ ${\rm (ii)} \Leftarrow {\rm (iii)}$)
follows from Theorem \ref{thm: chara via alteration, quasi}(1) (resp.\ Theorem \ref{thm: chara via alteration, quasi}(2))).}
\qedhere

\end{proof}

 \begin{proposition}\label{prop:reduction}
Let $S$ be an integral  
scheme which is flat and of finite type over $\Z$. 
Let $\pi \colon X \to S$ be a 
morphism of finite type from an integral  normal scheme $X$.
Let $\Delta$ be a 
$\Q$-divisor on $X$ such that $K_X+\Delta$ is $\Q$-Cartier.
Then there exists an open dense subset $U$ of $S$ such that 
\[
\mathcal{J}(X_{\mu},\Delta_{\mu})=\tau^q(X_{\mu},\Delta_{\mu})=\tau(X_{\mu},\Delta_{\mu})
\]
for every closed point $\mu \in U$, 
where $X_\mu := \pi^{-1}(\mu)$ and 
$\Delta_{\mu}=\Delta|_{X_{\mu}}$. 
\end{proposition}

\begin{proof}
By \cite{Tak04}*{Theorem 3.2} (cf.\ \cite[Theorem 11]{FDTT15}), there exists an open dense subset $U$ of $S$ such that
\begin{equation}\label{eq:1}
    \mathcal{J}(X_{\mu},\Delta_{\mu})=\tau(X_{\mu},\Delta_{\mu})
\end{equation}
for every closed point $\mu \in U$.
By the definition of $\tau^q(X_{\mu}, \Delta_{\mu})$, we have
\begin{equation}\label{eq:2}
    \tau(X_{\mu},\Delta_{\mu})=\tau_1(X_{\mu},\Delta_{\mu}) \subseteq
\tau^q(X_{\mu},\Delta_{\mu}).
\end{equation}
On the other hand, by Corollary \ref{cor:trace under bir}, we have
\begin{equation}\label{eq:3}
    \tau^q(X_{\mu},\Delta_{\mu}) \subseteq \mathcal{J}(X_{\mu},\Delta_{\mu}).
\end{equation}
Combining (\ref{eq:1}, \ref{eq:2}, \ref{eq:3}), we obtain
\[
\mathcal{J}(X_{\mu},\Delta_{\mu})=\tau^q(X_{\mu},\Delta_{\mu})=\tau(X_{\mu},\Delta_{\mu}).
\]
\end{proof}

\begin{remark}
{Last, we point out that one can also prove that $\bQ$-Gorenstein quasi-$+$-regular rings are klt by a similar argument to that of Theorem \ref{thm:quasi+rational-are-pseudorational}. In particular, this gives a somewhat different proof that $\bQ$-Gorenstein quasi-$F$-regular rings are klt.}
\end{remark}

\section{Cone correspondence for quasi-
$F$-regularity}

In this section, we establish 
the cone correspondences for $n$-quasi-$F^e$-splitting (Theorem \ref{thm:height equality}) 
and quasi-$F$-regularity (Theorem \ref{thm: corr quasi-F-reg}). 
In most parts,  the arguments are identical to 
the one of the cone correspondence for $n$-quasi-$F$-splitting established in \cite{KTTWYY2}.

\begin{notation}\label{n-cone-QFS}
\begin{enumerate}
\item 
Let $k$ be an $F$-finite field 
of characteristic $p>0$. 
Let $X$ be a projective normal variety over $k$
with $\dim X \geq 1$  and $H^0(X, \MO_X)=k$. 
Let $D$ be an ample $\Q$-Cartier $\Q$-divisor. 
Set 
\[
R := R(X, D) := \bigoplus_{d \in \Z_{\geq 0}} H^0(X, \MO_X(dD))t^d \subseteq {K(X)}[t],
\]
which is a $\Z_{\geq 0}$-graded subring of the standard $\Z_{\geq 0}$-graded polynomial ring ${K(X)}[t]$. 
Note that $R$ is a finitely generated $\Z_{\geq 0}$-graded  $k$-algebra 
and $\Spec R$ is an affine normal variety \cite{KTTWYY2}*{Theorem A.4}. 
Let $D = \sum_{i=1}^r \frac{\ell_i}{d_i}D_i$ be the irreducible decomposition, 
where $\ell_i$ and $d_i$ are coprime integers satisfying $d_i>0$ for each $1 \leq i \leq r$. 
Set $D':=\sum_{i=1}^r \frac{d_i-1}{d_i}D_i$.  
\item 
For the graded maximal ideal $\m := \bigoplus_{d >0} H^0(X, \MO_X(dD))t^d \subseteq R$, 
we set $U := \Spec R \setminus \{\m\}$. 
For $n \in \Z_{>0}$, let 
\[
W_nR = \bigoplus_{e \in p^{-(n-1)}\Z_{\geq 0}} (W_nR)_e
\]
be  the $p^{-(n-1)}\Z_{\geq 0}$-graded ring structure induced by 
Proposition \cite{KTTWYY2}*{Proposition 7.1}: 
\[
(W_nR)_e := \{ (r_0, r_1, ..., r_{n-1}) \in W_nR \,|\, 
r_0 \in R_e, r_1 \in R_{pe}, ..., r_{n-1} \in R_{p^{n-1}e}\}, 
\]
where $R_{c} :=0$ for $c \in p^{-(n-1)}\Z_{\geq 0} \setminus \Z_{\geq 0}$. 
Set $(W_nR)_{>0} := \bigoplus_{e \in p^{-(n-1)}\Z_{>0}} W_nR_{{e}}$, 
which is a graded primary ideal such that $\sqrt{(W_nR)_{>0}}$ is a  maximal ideal. 
\item 
For the induced isomorphism $U \simeq \Spec_X\left(\bigoplus_{d \in \Z} \MO_X(dD)t^d\right)$ \cite{KTTWYY2}*{Theorem A.4}, 
let $\rho \colon U \to X$ be the natural morphism 
(cf.\ \cite{KTTWYY2}*{Notation 7.7}). 
We note that $\rho$ is flat over the regular locus of $X$ by \cite{KTTWYY2}*{Proposition A.7(1)}.
Let $E$ be a $\Q$-divisor on $X$ such that $\Supp(\{E\})$ and $\Supp(D')$ have no common irreducible component. 
Then we define the $\Q$-divisor $\wt{E}$ on $\Spec R$ such that
the restriction of $\wt{E}$ to $U':=\rho^{-1}(X_{\reg})$ coincides with $(\rho|_{U'})^*(E|_{X_{\reg}})$.
Therefore, if $E=\sum_i a_i E_i$ is the irreducible decomposition, 
then 
$\wt{E}=\sum_i a_i \wt{E_i}$ is also the irreducible decomposition,
and in particular, 
$\rdown{\wt{E}}=\wt{\rdown{E}}$.
Furthermore, for $q \in \Z$, we have
\begin{align*}
    R(\wt{E}+qK_R) 
    &= H^0(
     X_{\reg}
    ,\rho'_*\sO_{U'}(\wt{E}|_{U'}+qK_{U'}))\\ &
    =
    H^0(
     X_{\reg}
    ,\rho'_*\sO_{U'}(\rdown{\wt{E}}|_{U'}+qK_{U'})) \\
    &=H^0(X_{\reg},\sO_{X_{\reg}}( \lfloor E \rfloor |_{X_{\reg}}) \otimes \rho'_*\sO_{U'}(qK_{U'})) \\
    &= \bigoplus_{d \in \Z} H^0(X,\sO_X(E+q(K_X+D')+dD))\cdot t^d \\ &\subseteq K(X)[t,t^{-1}],
\end{align*}
where $\rho'=\rho|_{U'}$ (cf.~\cite{KTTWYY2}*{7.9}).
Here, the third equality follows from \cite[Proposition 7.12 (2)]{KTTWYY2} and the equality
\[
\rdown{E+q(K_X+D')+dD}=\rdown{E}+\rdown{q(K_X+D')+dD},
\]
 which holds by the assumption on supports.
In particular,  $R(\wt{E}+qK_R)$ has the natural $\Z$-graded $R$-module structure. 
We further define the $p^{-(n-1)}\Z_{\geq 0}$-graded structure on $W_nR(\wt{E}+qK_R)$ by a similar way to the one of \cite{KTTWYY2}*{Proposition 7.5}. 
In particular, $W_nR(\wt{E}+qK_R)$ 
is a graded $W_nR$-module.
\end{enumerate}
\end{notation}

\begin{proposition}\label{prop:coh of Q}
We use Notation \ref{n-cone-QFS}. 
Fix $d, n \in \Z_{>0}$. 
Then, for every $e \in \Z_{>0}$, the following commutative diagram 
consists of graded $W_nR$-module homomorphisms of $p^{-(n+e-1)}\Z$-graded $W_nR$-modules
\begin{equation}\label{e1:coh of Q}
\begin{tikzcd}[column sep=1.5in]
H^d(U, \MO_U(K_U+\wt{E}|_U)) \arrow[r, "\Phi^e_{U, K_U+\wt{E}|_U, n}"] \arrow[d, "\simeq"] & H^d(Q^e_{U,K_U+\wt{E}|_U,n}) \arrow[d, "\simeq"] \\
H^{d+1}_\m(R(K_R+\wt{E})) \arrow[r, "\Phi^e_{R, K_R+\wt{E}, n}"] & H^{d+1}_\m(Q^e_{R,K_R+\wt{E},n}), 
\end{tikzcd}
\end{equation}
where all the maps are the natural ones, each vertical map is an isomorphism, 
and  the graded structures are given in Remark \ref{r:coh of Q}. 
\end{proposition}

\begin{proof}
The assertion follows from a similar argument to the proof of \cite{KTTWYY2}*{Proposition 7.14}.
\end{proof}

\begin{remark}\label{r:coh of Q}
We use the same notation as in Proposition \ref{prop:coh of Q}. 
The graded structure on the $W_nR$-modules in (\ref{e1:coh of Q}) are defined as follows. 
\begin{enumerate}
\item Since $R(K_R+\wt{E})$ is a $\Z$-graded $R$-module (Notation \ref{n-cone-QFS}(3)), 
the $\Z$-graded $R$-module structure on $H^{d+1}_\m(R(K_R+\wt{E}))$ is introduced by 
\cite[7.6]{KTTWYY2}. 
\item 
We define the ${p^{-(n+e-1)}}\Z$-graded $W_nR$-module structure on $Q^e_{R,K_R+\wt{E},n}$ 
by the same way as in  
\cite[Remark 7.13]{KTTWYY2}. 
The $p^{-(n+e-1)}\Z$-graded $W_nR$-module structure on $H^{d+1}_\m(Q^e_{R,K_R+\wt{E},n})$ is introduced by \cite[7.6]{KTTWYY2}. 
\item 
We have the $\Z$-graded $R$-module structure 
$\MO_U(K_U+\wt{E}|_U)$ induced by the restriction of that of $\wt{R(K_R+\wt{E})}$, 
where the $\Z$-graded structure on $\wt{R(K_R+\wt{E})}$ is given in (1). 
Then the $\Z$-graded $R$-module structure on $H^d(U, \MO_U(K_U+\wt{E}|_U))$  is introduced by 
\cite[7.6]{KTTWYY2}. 
We get the $p^{-(n+e-1)}\Z$-graded $W_nR$-module structure on $H^d(Q^e_{U,K_U+\wt{E}|_U,n})$ 
{in the same way as in (2).} 
\end{enumerate}
\end{remark}

\begin{proposition}\label{prop:commutativity}
We use Notation \ref{n-cone-QFS}. 
Set $d := \dim X$.  
Then, for every $e \in \Z_{>0}$, we have the commutative diagram
\[
\begin{tikzcd}[column sep=1.5in]
H^d(U, \MO_U(K_U +\wt{E}|_{U}))_0 \arrow[r, "{H^d(U, \Phi^e_{U, K_U+\wt{E}|_{U}, n})_0}"]  & H^d(U, Q^e_{U,K_U+\wt{E}|_{U},n})_0  \\
H^d(X,\cO_X(K_X+D'+E)) \arrow[r, "{H^d(X, \Phi^e_{X, K_X+D'+E, n})}"] \arrow[u, "\simeq"] & H^d(X,Q^e_{X,K_X+D'+E,n}).
\arrow[u, "\simeq"]
\end{tikzcd}
\]
\end{proposition}

\begin{proof}
The assertion follows from a similar argument to the proof of \cite{KTTWYY2}*{Proposition 7.15}.
\end{proof}

\begin{theorem}\label{thm:height equality}
We use Notation \ref{n-cone-QFS}.
Let $\Delta$ be an effective $\Q$-divisor on $X$ such that $\rdown{\Delta}=0$ and $\Supp(\Delta)$ and $\Supp(D')$ have no common irreducible component.
Take integers $n>0$ and $e>0$. 
Then 
 $(\Spec R_{\m},\wt{\Delta})$ is $n$-quasi-$F^e$-split if and only if $(X, D'+\Delta)$ is $n$-quasi-$F^e$-split. 
In particular, 
$(\Spec R_{\m},\wt{\Delta})$ is quasi-$F^e$-split if and only if $(X, D'+\Delta)$ is quasi-$F^e$-split. 
\end{theorem}

\begin{proof}
The assertion follows from a similar argument to  the proof of \cite{KTTWYY2}*{Theorem 7.16} by using Proposition \ref{prop:coh of Q} and 
Proposition \ref{prop:commutativity}.
\end{proof}

\begin{corollary}\label{cor: ample cartier case}
Let $k$ be an $F$-finite field of characteristic $p>0$ 
and let $X$ be a projective normal variety over $k$. 
Let $\Delta$ be an effective $\Q$-divisor on $X$ such that $-(K_X+\Delta)$ is $\Q$-Cartier and ample.
Assume that there exists 
an effective Weil divisor $H$ on $X$ such that (1)-(3) hold. 
\begin{enumerate}
\item  $H \sim -\ell (K_X+\Delta)$ for some integer $\ell >0$.  
\item  $(X \setminus \Supp(H), \{p^e\Delta\}|_{X \setminus \Supp(H)})$ is strongly $F$-regular for every integer $e \geq 0$.
\item For every integer $m>0$, 
there exists an integer $e>0$ such that 
\[
(X,\Delta+ (m/p^e) H)
\]
is quasi-$F^e$-split. 
\end{enumerate}
Then $(X,\Delta)$ is globally quasi-$F$-regular.
\end{corollary}

\begin{proof}
We may assume $H^0(X,\sO_X)=k$. 
Replacing $H$ by $rH$ for some $r \in \Z_{>0}$, 
the problem is reduced to the case when $H$ is Cartier. 
Set $D:=H$ and we use Notation \ref{n-cone-QFS}. 
Since $D=H$ is Cartier, $D'=0$ and every fibre of 
$\rho : U = \Spec_X \mathcal R \to X$ is smooth, where 
\[
\mathcal{R}:=\bigoplus_{d \in \Z} \cO_X(dH)t^d. 
\]
As $H$ is effective, 
we get $t \in R$ 
by $R = \bigoplus_{d \geq 0} H^0(X, \MO_X(dD))t^d \subseteq K(X)[t]$ (Notation \ref{n-cone-QFS}). 
Let $\wt{H}$ be the Weil divisor on $\Spec{R}$ corresponding to $H$.

Let us show that $\wt{H}=\mathrm{div}(t)$.
We take an affine open  cover $\{V_i\}_{i \in I}$ of $X$ 
such that, for every $i \in I$, 
we have $H|_{V_i}=\mathrm{div}_{V_i}(f_i)$ for some $f_i  \in \cO_X(V_i)$.
Recall that $\rho$ is induced by the natural injection $\cO_X \to \mathcal{R}$. 
For $U_i:=\rho^{-1}(V_i)$, we get $\wt{H}|_{U_i}=\mathrm{div}_{U_i}(f_i)$.
By $f_i \in \cO_{V_i}(-H)$ and $f_i^{-1} \in \cO_{V_i}(H)$, 
both $t^{-1}f_i$ and $tf_i^{-1}$ are contained in $\mathcal{R}(U_i)=\cO_{U}(U_i)$.
Therefore, we have $\mathrm{div}_{U_i}(t)=\mathrm{div}_{U_i}(f_i)=\wt{H}|_{U_i}$. 
Since $\mathrm{div}(t)$ coincides with $\wt{H}$ on $U$, so does on $\Spec(R)$.

Therefore, $\rho \colon U \to X$ induces a smooth morphism
\[
\Spec{R_t} = U \backslash \Supp(\wt{H}) \to X \backslash \Supp(H).
\]
%
By (2), $(\Spec R_t,\{p^e\wt{\Delta}\})$ is strongly $F$-regular for every $e \geq 0$. 
Then $\tau(R, \{p^e\wt{\Delta}\})_t = \tau(R_t, \{p^e\wt{\Delta}\}) =R_t$. 
Since $\wt\Delta$ is a $\Q$-divisor, there exists $m>0$ such that 
$t^m \in \tau(R, \{p^e\wt{\Delta}\})$ for every $e \geq 0$. 
By (3), there exists $e_0>0$ such that
$(X,\Delta+\frac{3m}{p^{e_0}}H)$ 
is  quasi-$F^{e_0}$-split. 
Then $(\Spec R_\m, \wt{\Delta}+(1/p^{e_0})\mathrm{div}(t^{3m}))$ is quasi-$F^{e_0}$-split by Theorem \ref{thm:height equality}.
By Theorem \ref{thm:chara qFr pair}, $(\Spec R_\m,\wt{\Delta})$ is quasi-$F$-regular. 
Pick an effective Weil divisor $E$ on $X$.
Since $(\Spec R_\m,\wt{\Delta})$ is $n$-quasi-$F$-regular for some $n>0$, 
there exists a rational number $\epsilon >0$ 
such that $(\Spec R_\m,\wt{\Delta}+\epsilon \wt{E})$ is $n$-quasi-$F^{e}$-split for every $e\geq 0$ (Definition \ref{d-QFR}). 
Again by Theorem \ref{thm:height equality}, 
$(X, \Delta + \epsilon E)$ is $n$-quasi-$F^e$-split for every $e\geq 0$, where 
$n$ is fixed and $e$ is independent of $n$. 
Therefore, $(X, \Delta)$ is globally quasi-$F$-regular (Definition \ref{d-QFR}). 
\end{proof}

\begin{theorem}\label{thm: corr quasi-F-reg}
Let $k$ be an $F$-finite field of characteristic $p>0$. 
Let $X$ be a projective normal variety over $k$
with $\dim X \geq 1$  and $H^0(X, \MO_X)=k$. 
Take a $\Q$-divisor $B$ on $X$ such that 
all the coefficients of $B$ are contained in $\{ 1- \frac{1}{m} \,|\, m \in \Z_{>0}\}$ 
 and $-(K_X+B)$ is $\Q$-Cartier and ample. 
Set 
\[
R := R(X, -(K_X+B)) := \bigoplus_{d \in \Z_{\geq 0}} H^0(X, \MO_X(-d(K_X+B)))t^d \subseteq {K(X)}[t],
\]
which is a $\Z_{\geq 0}$-graded subring of the standard $\Z_{\geq 0}$-graded polynomial ring ${K(X)}[t]$. 
Set $\m := \bigoplus_{d >0} H^0(X, \MO_X(-d(K_X+B)))t^d \subseteq R$. 
Then $R_{\m}$ is quasi-$F$-regular if and only if $(X, B)$ is globally quasi-$F$-regular.
\end{theorem}

\begin{proof}
Set $D := -(K_X+B)$. 
In what follows, we use Notation \ref{n-cone-QFS}. Then $B  =D'$. 
We take an ample effective Cartier divisor $H$ on $X$ such that 
\begin{itemize}
    \item $\Supp(H)$ and $\Supp(B)$ have no common irreducible components, 
    \item $(X \backslash \Supp(H),B|_{X \backslash \Supp(H)})$ is simple normal crossing, and 
    \item $H \sim -\ell (K_X+B)  =\ell D$ for some integer $\ell >0$.
\end{itemize}

First, we assume that $R_\m$  is quasi-$F$-regular.
Then we can find a rational number $\epsilon>0$ such that 
$-(K_X+B+\epsilon H)$ is still ample and $(\Spec R_\m,\epsilon\wt{H})$ is quasi-$F$-regular by Theorem \ref{thm: suff small qFr}(3).
In particular, there exists $n \in \Z_{>0}$ such that $(\Spec R_\m,\epsilon\wt{H})$ is $n$-quasi-$F^e$-split for every integer $e>0$. 
Applying Theorem \ref{thm:height equality} by setting 
$\wt{\Delta}:=\epsilon H$,
$(X,B+\epsilon H)$ is $n$-quasi-$F^e$-split for every $e>0$. 
For every integer $m>0$, there exists $e \in \Z_{>0}$ such that 
$m/p^e \leq \epsilon$, 
and hence $(X, B+ (m/p^e)H)$ is quasi-$F^e$-split. 
Thus $(X,B)$ is globally quasi-$F$-regular by Corollary \ref{cor: ample cartier case}.

Next, we assume that $(X,B)$ is globally quasi-$F$-regular.
By the same argument as in the second paragraph of the proof of Corollary \ref{cor: ample cartier case}, we get $\widetilde{D}  = \div(t)$.
Since we have $H \sim \ell D$, 
it holds that $ \div(f) = \widetilde{H} \sim 
\ell \wt{D}=  \div(t^{\ell})$ 
for some homogeneous element $f \in K(X)[t]$ of degree $\ell$. 
We have the induced morphism: 
\[
\Spec R_f = U \setminus \Supp(\wt{H}) \to X \setminus \Supp(H). 
\]
Since $(X \backslash \Supp(H),B|_{X \backslash \Supp(H)})$ is simple normal crossing, 
$U \setminus \Supp(\wt{H})$ has \'etale locally toric singularities \cite[Proposition B.6]{KTTWYY2}, and hence $R_f$ is strongly $F$-regular, i.e., $\tau(R)_f = \tau(R_f) = R_f$. 
Therefore, there exists an integer $m > 0$ such that 
$f^m \in \tau(R)$. Then $f^m$ is a test element of $R_{\m}$ (Remark \ref{r-test-ideal-summary}). 
Since $(X,B)$ is globally quasi-$F$-regular, there exist $n \in \Z_{>0}$ and $e\in\Z_{>0}$ 
such that $(X,B+(3m/p^e)H)$ 
is $n$-quasi-$F^e$-split. 
By $\widetilde{H} = \div(f)$ and Theorem \ref{thm:height equality}, $(\Spec R_\m, (1/p^e)\mathrm{div}(f^{3m}))$ is $n$-quasi-$F^e$-split. 
Then it follows from  Theorem \ref{thm:chara qFr pair} that $R_\m$ is quasi-$F$-regular.
\qedhere


\end{proof}

\section{Quasi-$F$
-regularity for Fano varieties and klt singularities}

In this section, we apply our results to Fano varieties and klt singularities.

\subsection{Fano varieties}

It is natural to ask which Fano varieties are globally quasi-$F$-regular. 
As a fundamental result, we prove that 
quasi-$F$-regularity is equivalent to quasi-$F$-splitting 
for $\Q$-factorial strongly $F$-regular Fano varieties (Corollary \ref{qFs to qFr for am Fano}). 
This result is meaningful, because it is much easier to check whether a given variety is quasi-$F$-split.  
As a consequence, klt del Pezzo surfaces are globally quasi-$F$-regular 
when $p>5$ (Corollary \ref{c klt dP QFR}). 

\begin{lemma}\label{l-QFR-descent}
Let $k \subseteq k'$ be a field extension of $F$-finite fields of characteristic $p>0$. 
Let $X$ be a projective normal variety over $k$ 
such that $H^0(X, \MO_X)=k$ and $X \times_k k'$ is normal. 
Take an effective $\Q$-divisor $\Delta$ on $X$. 
Then the following hold. 
\begin{enumerate}
    \item
    If $(X \times_k k', \Delta \times_k k')$ is globally $n$-quasi-$F^e$-split, 
    then so is $(X, \Delta)$. 
    \item 
    If $(X \times_k k', \Delta \times_k k')$ is globally $n$-quasi-$F$-regular, 
    then so is $(X, \Delta)$. 
\end{enumerate}     
\end{lemma}

\begin{proof}
Let us show (1). 
Set $X' := X \times_k k'$ and $\Delta' := \Delta \times_k k'$. 
For the induced morphism $\pi : X' = X \times_k k' \to X$, 
{we have that $\pi^*K_X = K_{X'}$, so} 
we get the following commutative diagram: 
\[
\begin{tikzcd}[column sep = large]
\MO_X(K_X) \arrow[r, "\Phi^e_{X, \Delta, n}"]\arrow[d, "\alpha"] & Q^e_{X, K_X+\Delta, n} \arrow[d]\\
\pi_*\MO_{X'}(K_{X'}) \arrow[r, "\Phi^e_{X, \Delta, n}"] & \pi_*Q^e_{X', K_X+\Delta', n}
\end{tikzcd}
\]
Since $\pi$ is an affine morphism, 
we obtain the following commutative diagram by taking $H^d(X, -)$: 
\[
\begin{tikzcd}[column sep = huge]
H^d(X, \MO_X(K_X)) \arrow[r, "H^d(\Phi^e_{X, \Delta, n})"]\arrow[d, "\beta"] & 
H^d(X, Q^e_{X, K_X+\Delta, n}) \arrow[d]\\
H^d(X', \MO_{X'}(K_{X'})) \arrow[r, "H^d(\Phi^e_{X', \Delta', n})"] & 
H^d(X', Q^e_{X', K_X+\Delta', n}).
\end{tikzcd}
\]
Since $\alpha$ is a split injection, so is $\beta$. 
Therefore, if $H^d(\Phi^e_{X', \Delta', n})$ is injective, then so is $H^d(\Phi^e_{X, \Delta, n})$. 
Thus (1) holds. 
The assertion (2) follows from (1). 
\end{proof}

\begin{lemma}\label{l QFS times A1}
Let $k$ be an $F$-finite field of characteristic $p>0$. 
Let $X$ be a normal variety over $k$ and take an effective $\Q$-divisor $\Delta$ on $X$. 
Fix an integer $n>0$. 
Assume that $(X, \Delta)$ is 
$n$-quasi-$F$-split. Then $(X \times_k \mathbb A^1_k, \Delta \times_k \mathbb A^1_k)$ is 
$n$-quasi-$F$-split. 
\end{lemma}

\begin{proof}
    Since $(X, \Delta)$ is $n$-quasi-$F$-split, 
there exists 
a $W_n\MO_X$-module homomorphism 
\[
\sigma \colon F_*W_n\cO_X(p\Delta) \to \cO_X
\]
such that $\sigma(F_*1)=1$.
Set $Y:=\Spec k[t]$.
Then $X \times_k Y$ is the relative affine spectrum of the $\cO_X$-algebra $\cO_X[t]$ over $X$. 
We equip $\MO_X[t]$ with the natural graded ring structure given by 
\[
\MO_X[t] = \bigoplus_{d \in \Z_{\geq 0}} \MO_X t^d. 
\]
Then $W_n(\MO_X[t])$ has a $\frac{1}{p^{n-1}}\Z_{\geq 0}$-graded ring structure (\cite[Proposition 7.1]{KTTWYY2}): 
\[
W_n(\MO_X[t]) =\bigoplus_{e \in \frac{1}{p^{n-1}}\Z_{\geq 0}} W_n(\MO_X[t])_e, 
\]
 where
\[
W_n(\MO_X[t])_e := 
\{ (f_0, \ldots, f_{n-1}) \in W_n(\MO_X[t]) \,|\, f_0 \in \MO_Xt^e, 
f_1 \in \MO_X t^{pe}, \ldots , f_{n-1} \in {\cred \MO_X}t^{p^{n-1}e}\}.
\]
Here we set $t^c :=0$ for $c \in \Q \setminus \Z$. 
Define $s := [t] = (t, 0, ..., 0) \in W_n(\MO_X[t])$ to be the Teichm\"uller lift of $t$.

\begin{claim*}
The following equality 
\[
(W_n\MO_X)[s] = \bigoplus_{d \in \Z_{\geq 0}} W_n(\MO_X[t])_d
\]
of subrings of $W_n(\MO_X[t])$ holds.  
\end{claim*}

\begin{proof}[Proof of Claim]
The inclusion 
$(W_n\MO_X)[s] \subseteq \bigoplus_{d \in \Z_{\geq 0}} W_n(\MO_X[t])_d$ 
holds by 
$W_n\MO_X = W_n(\MO_X[t])_0$ and $s \in W_n(\MO_X[t])_1$. 
Let us show the opposite inclusion: 
$(W_n\MO_X)[s] \supseteq \bigoplus_{d \in \Z_{\geq 0}} W_n(\MO_X[t])_d$. 
Take $d \in \Z_{\geq 0}$ and 
$f \in W_n(\MO_X[t])_d$. 
We can write $f = \sum_{i=0}^{n-1}{\cccred \mathbf{V}}^i([a_it^{p^id}]) = $ 
for some $a_i \in \MO_X$. 
Fix $i \in \{0, 1, ..., n-1\}$. 
It suffices to show ${\cccred \mathbf{V}}^i([a_it^{p^id}]) \in (W_n\MO_X)[s]$, which follows from 
\[
{\cccred \mathbf{V}}^i([a_it^{p^id}]) ={\cccred \mathbf{V}}^i( [a_i] \cdot F^i[t^d]) = [t^d] \cdot {\cccred \mathbf{V}}^i([a_i]) = s^d \cdot {\cccred \mathbf{V}}^i([a_i]) \in (W_n\MO_X)[s]. 
\]
This completes the proof of Claim. 
\end{proof}

We have the following decomposition into the integral and non-integral parts: 
\begin{eqnarray*}
W_n(\MO_X[t]) &=& \bigoplus_{e \in \frac{1}{p^{n-1}}\Z_{\geq 0}} W_n(\MO_X[t])_e\\
&=& \left(\bigoplus_{d \in  \Z_{\geq 0}} W_n(\MO_X[t])_d\right)
\oplus 
\left(\bigoplus_{e \in \frac{1}{p^{n-1}}\Z_{\geq 0} \setminus \Z} W_n(\MO_X[t])_e\right)\\
&=& (W_n\MO_X)[s]
\oplus 
\left(\bigoplus_{e \in \frac{1}{p^{n-1}}\Z_{\geq 0} \setminus \Z} W_n(\MO_X[t])_e\right), 
\end{eqnarray*}
where the last equality holds by Claim. 
Since each direct summand is a $(W_n\MO_X)[s]$-submodule, 
the projection 
\[
\pi : W_n(\MO_X[t]) \to (W_n\MO_X)[s]
\]
is a $(W_n\MO_X)[s]$-module homomorphism. 
Consider the additive homomorphism  
\begin{eqnarray*}
\tau : 
F_*((W_n\cO_X)[s]) &\to&  \cO_X[t]\\
F_*(as^i) &\mapsto& 
\begin{cases}
    \sigma(F_*a)t^{{\frac{i}{p}}} \qquad \text{if}\quad 
    i \in p \Z \\
    0 \hspace{30mm}\text{otherwise},
\end{cases}
    \end{eqnarray*}
where $a \in W_n\cO_X$. 
Then $\tau$ is a $(W_n\MO_X)[s]$-module homomorphism, 
because the following holds for $b \in W_n\MO_X$: 
\[
\tau(bs^j F_*(as^i)) = \tau (F_*( aF(b)s^{i+pj})) = \sigma(F_*(aF(b))) t^{{\frac{i+pj}{p}}} 
\]
\[
= (bt^j) \cdot  \sigma(F_*(a)) t^{{\frac{i}{p}}} = (bs^j)\cdot \tau(F_*(as^i)),
\]
where we set $t^c :=0$ for $c \in \Q \setminus \Z$. 

Take the composite $(W_n\MO_X)[s]$-module homomorphism: 
\[
\theta : F_*(W_n(\cO_X[t])) \xrightarrow{F_*\pi} F_*((W_n\cO_X)[s]) \xrightarrow{\tau} \cO_X[t]. 
\]
We have $\theta (F_*1) = \tau \circ (F_*\pi)(F_*1) = \tau(F_*1) = 1$. 
After replacing $\cO_X$ by the function field $K(X)$ of $X$, 
we also have $(W_nK(X))[s]$-module homomorphisms: 
\[
\theta: F_*W_n(K(X)[t])  \xrightarrow{F_*\pi} F_*W_n(K(X))[s]  \xrightarrow{\tau} K(X)[t]. 
\]
It is enough to check the inclusion $\theta(F_*W_n(\cO_X(p\Delta)[t])) \subseteq \MO_X[t]$. 
This follows from $\pi(W_n(\cO_X(p\Delta)[t])) \subseteq W_n(\cO_X(p\Delta))[s]$ 
and $\sigma(F_*W_n\cO(p\Delta)) \subseteq \MO_X$.  
\qedhere




\end{proof}

\begin{lemma}\label{l-QFS-pure-tr-ext}
Let $k$ be an $F$-finite field of characteristic $p>0$ and take a purely transcendental extension $k \subseteq k':=k(t_1, ..., t_m)$ of finite degree. 
Let $X$ be a projective normal variety over $k$ with $H^0(X, \MO_X)=k$. 
Take an integer $n >0$ and an effective $\Q$-divisor $\Delta$ on $X$. 
Then $(X, \Delta)$ is  $n$-quasi-$F$-split if and only if 
    $(X \times_k k', \Delta \times_k k')$ is  $n$-quasi-$F$-split. 
\end{lemma}

\begin{proof}
By induction on $m$, we may assume that $k'  = k(t_1)$. 
The \lq\lq if" part holds by Lemma \ref{l-QFR-descent}. 
Conversely, assume that $(X, \Delta)$ is  $n$-quasi-$F$-split. 
Then $(X \times_k \mathbb A^1_k, \Delta \times_k \mathbb A^1_k)$  is $n$-quasi-$F$-split (Lemma \ref{l QFS times A1}), 
and hence so is $(X \times_k k', \Delta \times_k k')$ 
by taking the localisation $(-) \times_{\mathbb A^1_k} \Spec k'$. 
\qedhere

\end{proof}

\begin{lemma}\label{l divCM bc}
Let $k$ be a field and let $X$ be a normal variety over $k$. 
Assume that $X$ is divisorially Cohen-Macaulay. 
Then $X \times_k \mathbb A^1_k$ and $X \times_k k(t)$ are divisorially Cohen-Macaulay, where $k(t) := K(\mathbb A^1_k)$. 
\end{lemma}

\begin{proof}
It is enough to show that $X \times_k \mathbb A^1_k$ is divisorially Cohen-Macaulay, 
which follows from the fact that the pullback map 
\begin{eqnarray*}
\pi^*: {\rm Cl}(X) &\to& {\rm Cl} (X \times_k \mathbb A^1_k)\\
D &\mapsto& \pi^*D
\end{eqnarray*}
 is an isomorphism 
 for the projection $\pi : X \times_k \mathbb A^1_k \to X$ 
 \cite[Ch.\ II, Proposition 6.6]{hartshorne77}. 
\end{proof}

{The following is the key theorem of this section. The reader should note that the first assumption is satisfied, for example, when $\Delta$ has standard coefficients.}

\begin{theorem}\label{thm:qFs to qFr}
Let $k$ be an $F$-finite field of characteristic $p>0$ and 
let $X$ be a divisorially Cohen-Macaulay projective normal variety 
over $k$. 
Let $\Delta$ be an effective $\Q$-divisor such that 
\begin{enumerate}
\item $\{ p^r\Delta\} \leq \Delta$ for every 
integer $r \geq 0$,  
\item $-(K_X+\Delta)$ is $\Q$-Cartier and ample, 
\item $(X, \Delta)$ is strongly $F$-regular, and 
\item $(X, \Delta)$ is quasi-$F$-split. 
\end{enumerate}
Then 
$(X,\Delta)$ is globally quasi-$F$-regular.
\end{theorem}

\begin{proof}
Set $d:=\dim X$. 
We have $\rdown{\Delta}=0$. 
Replacing $k$ by $H^0(X, \MO_X)$, the problem is reduced to the case when $H^0(X, \MO_X)=k$. 
By Lemma \ref{l-QFR-descent}, Lemma \ref{l-QFS-pure-tr-ext}, 
and Lemma \ref{l divCM bc}, 
we may replace $(X, \Delta)$ by $(X \times_k k', \Delta \times_k k')$ 
for a purely transcendental field extension $k \subseteq k'$ of finite degree. 
By the Bertini theorem for the generic member \cite[Theorem 4.9 and Theorem 4.19]{TanakaBertini}, 
we can find  an integer $\ell >0$ and 
a normal prime Cartier divisor $H_0$ 
such that 
\begin{enumerate}
\renewcommand{\labelenumi}{(\roman{enumi})}
\item $H_0 \sim  -\ell(K_X+\Delta)$, 
\item $H_0 \not\subseteq \Supp \Delta$, and 
\item $(X, \Delta + aH_0)$ is strongly $F$-regular for every rational number $0 \leq a <1$. 
\end{enumerate}

Fix an integer $s_0 > 0$. 
For  $H:=s_0H_0$, we set 
\[
D_e:= K_X+\Delta+(1/p^e)H\qquad\text{and} \qquad D_{\infty} := K_X+\Delta. 
\]
By Corollary \ref{cor: ample cartier case}, 
it is enough to show that 
$(X,\Delta+(1/p^e)H)$ is quasi-$F^e$-split for some $e>0$. 
Then the problem is reduced to showing $(\star)$ below (see \cite[Lemma 3.10]{TWY}).
\begin{enumerate}
\item[$(\star)$] There exist integers  $n>0$ and $e>0$ such that 
\[
H^{d-1}(X, B^e_{X, D_e, n})=0.
\]
\end{enumerate}
Recall that,  given integers $e>0$, $n>0$, and a $\Q$-divisor $D$, 
a coherent $W_n\MO_X$-module $B^e_{X, D, n}$ is defined by the following exact sequence (cf.\ \cite[(3.0.2), (3.1.1)]{TWY}): 
\[
0 \to W_n\MO_X{(D)} \xrightarrow{F^e} F_*^e W_n\MO_X{(p^eD)} \to B^e_{X, D, n} \to 0. 
\]
Set $Q_{X, \Delta, n} := Q^1_{X, \Delta, n}$ and 
$B_{X, \Delta, n} := B^1_{X, \Delta, n}$. 
By (4), there is an integer $m_0>0$ such that $(X, \Delta)$ is $m_0$-quasi-$F$-split. 


\setcounter{step}{0}

\begin{step}\label{s1:qFs to qFr}
$B_{X, p^rD_e, n}$ {and $B_{X, p^rD_{\infty}, n}$ are} Cohen-Macaulay for every triple $(e, n, r) \in \Z^3$ satisfying  $e \geq 0$, $n>0$, and $r \geq 0$. 
\end{step}

\begin{proof}[Proof of Step \ref{s1:qFs to qFr}]
Fix integers $e \geq 0$ and $r \geq 0$. 
By \cite[Proposition 3.9]{KTTWYY1}, 
it is enough to show that 
$(X, \{p^rD_e\})$ is strongly $F$-regular, 
because this implies that $(X, \{p^rD_e\})$  is naively keenly $F$-pure. 
It holds that 
\[
\{ p^rD_e\} = \{ p^r(K_X+\Delta + (1/p^e)H)\} 
\overset{{\rm (ii)}}{=} \{ p^r\Delta\} + \{ (p^r/p^e)H\} 
\overset{{\rm (1)}}{\leq}  \Delta + \{ s_0p^r/p^e\} H_0. 
\]
Thus $(X, \{ p^rD_e\})$ is strongly $F$-regular by (iii). 
{Similarly,  $B_{X, p^rD_{\infty}, n}$  is Cohen-Macaulay, because 
$(X, \{p^r \Delta\})$ is strongly $F$-regular by (1) and (3).} 
This completes the proof of Step \ref{s1:qFs to qFr}. 
\end{proof}

\begin{step}\label{s2:qFs to qFr}
There exists an integer $n_0>0$ such that 
\[
H^j(X, \MO_X(p^rD_{\infty})) 
=
H^j(X, B_{X, p^rD_{\infty}, n})=0.
\]
for every triple $(j, n, r) \in \Z^3$ satisfying 
$j<d$, $n\geq n_0$, and $r \geq 0$. 
\end{step}

\begin{proof}[Proof of Step \ref{s2:qFs to qFr}]
By the Fujita vanishing theorem, there exists an integer $r_1>0$ such that 
\[
H^j(X, \MO_X(p^rD_{\infty}))=H^j(X, B_{X, p^rD_{\infty}, 1})=0
\]
for every $j<d$ and every $r \geq r_1$, {where the latter equality holds 
by Serre duality and the fact that $B_{X, p^rD_{\infty}, 1}$ is Cohen-Macaulay (Step \ref{s1:qFs to qFr}).} 
We have the following exact sequence \cite[Lemma 3.8]{KTTWYY1}:  
\[
0 \to F_*B_{X, p^{r+1}D_{\infty}, n} \to B_{X, p^rD_{\infty}, n+1} \to B_{X, p^rD_{\infty}, 1} \to 0. 
\]
Therefore, it holds that 
\begin{equation}\label{e1-s1:qFs to qFr}
H^j(X, \MO_X(p^rD_{\infty}))=H^j(X, B_{X, p^rD_{\infty}, n})=0
\qquad \text{if} \,\,j<d,\,\, n>0,\,\, r \geq r_1. 
\end{equation}

Fix $j <d$. 
Set $n_0 := r_1 + m_0$ and take an integer $n$ 
satisfying $n \geq n_0$. 
It is enough to prove that
\begin{equation}\label{e2-s1:qFs to qFr}
H^j(X, \MO_X(p^rD_{\infty})) = H^j(X, B_{X, p^rD_{\infty}, 
n-r})=0
\end{equation}
for every $r \geq 0$ by descending induction on $r$. 
The base case $r \geq r_1$ of this induction 
has been settled by (\ref{e1-s1:qFs to qFr}). 
Fix an integer $r$ satisfying $0 \leq r < r_1$. 
It is enough to show the implication 
$(\ref{e3-s1:qFs to qFr}) \Rightarrow (\ref{e2-s1:qFs to qFr})$, where 
\begin{equation}\label{e3-s1:qFs to qFr}
H^j(X, \MO_X(p^{r+1}D_{\infty})) = H^j(X, B_{X, p^{r+1}D_{\infty}, 
n-r-1})=0.  
\end{equation}
We have the following  exact sequences  \cite[Lemma 3.8]{KTTWYY1}: 
\[
0 \to F_*B_{X, p^{r+1}D_{\infty}, m-1} \to Q_{X, p^{r}D_{\infty}, m} \to F_*\MO_X(p^{r+1}D_{\infty}) \to 0, 
\]
\[
0 \to \MO_X(p^rD_{\infty}) \to Q_{X, p^rD_{\infty}, m} \to B_{X, p^rD_{\infty}, m} \to 0. 
\]
Since $(X, \Delta)$ is $m_0$-quasi-$F$-split and $D_{\infty} = K_X+\Delta$,  the latter sequence  splits if  $m \geq m_0$. 
Therefore, we get the following implications for $m \geq m_0$: 
\begin{align*}
&H^j(X, \MO_X(p^{r+1}D_{\infty})= H^j(X, B_{X, p^{r+1}D_{\infty}, m-1}) =0\\
\implies\ &H^j(X, Q_{X, p^{r}D_{\infty}, m})= 0\\
\implies\ &H^j(X, \MO_X(p^rD_{\infty})) =  H^j(X, B_{X, p^rD_{\infty}, m})=0. 
\end{align*}
By setting $m:=n-r$, we get the required implication $(\ref{e3-s1:qFs to qFr}) \Rightarrow (\ref{e2-s1:qFs to qFr})$. 
Here the condition $m \geq m_0$  holds by 
\[
m = n-r \geq n_0 -r_1 =m_0. 
\]
This completes the proof of Step \ref{s2:qFs to qFr}. 
\end{proof}

\begin{step}\label{s3:qFs to qFr}
There exist integers $e_1>0$  and $r_0>0$ such that 
\[
H^j(X, B_{X, p^rD_e, n})=0 
\]
for every quadruple $(j, e, n, r) \in \Z^4$
satisfying $j<d, e \geq e_1, n>0$, and $r \geq r_0$. 
\end{step}

\begin{proof}[Proof of Step \ref{s3:qFs to qFr}]
By the following exact sequence \cite[Lemma 3.8]{KTTWYY1}:  
\[
0 \to F_*B_{X, p^{r+1}D_e, n} \to B_{X, p^rD_e, n+1} \to B_{X, p^rD_e, 1} \to 0, 
\]
we may assume $n=1$, i.e., the required vanishing is reduced to 
\begin{equation}\label{goal:qFs to qFr}
H^j(X, B_{X, p^rD_e, 1})=0. 
\end{equation}
In what follows, 
we shall provide some details, 
although  the assertion follows from the Fujita vanishing theorem. 

\medskip

It holds that  
\[
p^rD_e = p^r(K_X+\Delta + 
(1/p^e)H) = p^r(K_X+\Delta) + 
(s_0p^r/p^e)H_0. 
\]
Fix an integer $e^{\dagger} > 0$ such that $0 <s_0/p^{e^{\dagger}} <1$. 
Depending on $e$, we shall treat the following three cases separately. 
\begin{enumerate}
\item[(I)] $e \leq r$ ($e$ is small). 
\item[(II)] $r \leq e \leq r+1+e^{\dagger}$ ($e$ is in the middle)
\item[(III)] $e \geq r+1+e^{\dagger}$ ($e$ is large). 
\end{enumerate}
It is enough to prove Step \ref{s3:qFs to qFr} under each of (I)-(III). 
In other words, we shall find pairs $(e_{{\rm I}}, r_{{\rm I}}), (e_{{\rm II}}, r_{{\rm II}}), (e_{{\rm III}}, r_{{\rm III}})$, 
corresponding to the above cases, which assure the required vanishings, 
and then it suffices to set $e_1 := \max \{  e_{{\rm I}}, e_{{\rm II}}, e_{{\rm III}}\}$ and 
$r_0 := \max \{r_{{\rm I}}, r_{{\rm II}}, r_{{\rm III}}\}$. 


Fix an integer $\nu_0>0$ such that $\nu_0(K_X+\Delta)$ is Cartier. 
For every $r\geq 0$, we can write 
\[
p^r = a_r\nu_0 + b_r
\]
where $a_r$ and $b_r$ are the non-negative integers satisfying $0 \leq b_r <\nu_0$. It is easy to see that 
\[
\frac{p^r}{\nu_0} -1 \leq a_r \leq  \frac{p^r}{\nu_0}. 
\]
We have 
\[
p^rD_e= p^r(K_X+\Delta) + (s_0p^r/p^e)H_0 = 
b_r(K_X+\Delta) + a_r\nu_0(K_X+\Delta) + s_0p^{r-e}H_0. 
\]
\medskip

(I) 
Assume $e \leq r$. 
Then $a_r\nu_0(K_X+\Delta) + s_0p^{r-e}H_0$ is Cartier, and hence  we get  
\begin{eqnarray*}
B_{X, p^rD_e, 1} &\simeq& B_{X,  b_r(K_X+\Delta), 1}  \otimes \MO_X(a_r\nu_0(K_X+\Delta) + s_0p^{r-e}H_0). 
\end{eqnarray*}
By $0 \leq b_r <\nu_0$, there are only finitely many possibilities 
for 
$B_{X,  b_r(K_X+\Delta), n}$. 
Then it is enough to show that 
$-(a_r\nu_0(K_X+\Delta) + s_0p^{r-e}H_0)$ is sufficiently ample. 
We have 
$ \frac{a_r\nu_0}{p^r}  \geq  1- \frac{\nu_0}{p^r}  > \frac{2}{3}$ for $r \gg 0$. 
There exists 
an integer $e_{{\rm I}}>0$ such that  
\[
-\left( 
 \frac{1}{3}(K_X+\Delta) + \frac{s_0}{p^e}H_0 
 \right)
\]
is ample for every $e\geq e_{{\rm I}}$. It holds that 
\begin{eqnarray*}
a_r\nu_0(K_X+\Delta) + s_0p^{r-e}H_0 
&=& p^r \left( 
\frac{a_r\nu_0}{p^r}(K_X+\Delta) + \frac{s_0}{p^e}H_0 \right)\\
&=& p^r \left( 
\left(\frac{a_r\nu_0}{p^r}-\frac{1}{3} \right)(K_X+\Delta) + 
\left(\frac{1}{3} (K_X+\Delta) + \frac{s_0}{p^e}H_0\right) \right)
\end{eqnarray*}
Note that we may apply Serre duality by Step \ref{s1:qFs to qFr}. 
Then, by the Fujita vanishing theorem for $\Q$-divisors (cf.\ \cite[Lemma 3.5]{tanaka16_excellent}), 
there exists $r_{{\rm I}}>0$ such that the required vanishing (\ref{goal:qFs to qFr}) 
holds for $r \geq r_{{\rm I}}$ and $e \geq e_{{\rm I}}$ 
under the assumption that $e \leq r$. 

\medskip

(II) 
As in (I), we get 
\begin{eqnarray*}
B_{X, p^rD_e, 1} &\simeq& B_{X,  b_r(K_X+\Delta)+s_0p^{-(e-r)}H_0, 1}  \otimes \MO_X(a_r\nu_0(K_X+\Delta)). 
\end{eqnarray*}
If $r \leq e \leq r+1+e^{\dagger}$, then $0 \leq e-r \leq 1+e^{\dagger}$. 
Hence there are only finitely many possibilities for $ b_r(K_X+\Delta)+s_0p^{-(e-r)}H_0$. 
Then we may apply a similar argument to the one of (I). 

\medskip

(III) Assume $e \geq r+1 + e^{\dagger}$. 
Recall that $\MO_X(D) = \MO_X(\rdown{D})$. 
We have 
\[
\rdown{p^{r+1}D_e} 
= \rdown{p^{r+1}(K_X+\Delta)} + \rdown{\frac{s_0p^{r+1}}{p^e}}H_0. 
\]
By $0<s_0/p^{e^{\dagger}}<1$ and $e \geq r+1 + e^{\dagger}$, we obtain 
\[
0 < \frac{s_0\cdot p^{r+1}}{p^e}< 
\frac{p^{ e^{\dagger}} \cdot p^{r+1}}{p^{r+1+e^{\dagger}}}=1. 
\]
Therefore, $\MO_X(p^{r+1}D_e) = \MO_X(p^{r+1}(K_X+\Delta))$. 
Similarly, $\MO_X(p^rD_e) = \MO_X(p^r(K_X+\Delta))$. 
By the exact sequence  
\[
0 \to \MO_X(p^rD_e) \to F_*\MO_X(p^{r+1}D_e) \to B_{X,p^rD_e, 1} \to 0, 
\]
we get 
\[
B_{X,p^rD_e, 1} = B_{X,p^r(K_X+\Delta), 1} \quad \text{if}\quad 
e \geq r+1+e^{\dagger}.
\]
As in (II), it holds that
\begin{eqnarray*}
B_{X, p^rD_e, 1} &\simeq& B_{X,  b_r(K_X+\Delta), 1}  \otimes \MO_X(a_r\nu_0(K_X+\Delta)). 
\end{eqnarray*}
There are only finitely many possibilities for $b_r(K_X+\Delta)$. 
Then we may apply a similar argument to the one of (I). 

This completes the proof of Step \ref{s3:qFs to qFr}. 
\qedhere




\end{proof}



\begin{step}\label{s4:qFs to qFr}
There exist integers $e_2>0$ and $n_0>0$ such that 
\begin{equation}\label{e1-s4:qFs to qFr}
H^j(X, B_{X, p^rD_e, n_0})=0 
\end{equation}
for every triple $(j, e, r) \in \Z^3$ 
satisfying $j<d, e \geq e_2,$ and $ r \geq 0$. 
\end{step}

\begin{proof}[Proof of Step \ref{s4:qFs to qFr}]
Let $n_0$ be a positive integer as in Step \ref{s2:qFs to qFr}. 
Take positive integers $e_1$ and $r_0$ as in Step \ref{s3:qFs to qFr}. 
Fix $j<d$. 
If $r \geq r_0$ and $e \geq e_1$, 
then (\ref{e1-s4:qFs to qFr}) follows from Step \ref{s3:qFs to qFr}. 
If $0 \leq r \leq r_0$, then we get 
$B_{X, p^rD_e, n_0} = B_{X, p^rD_{\infty}, n_0}$ for $e \gg 0$, 
and hence  (\ref{e1-s4:qFs to qFr}) follows from Step \ref{s2:qFs to qFr}. 
This completes the proof of Step \ref{s4:qFs to qFr}. 
\qedhere
\end{proof}

\begin{step}\label{s5:qFs to qFr}
There exist integers $e_2>0$  and $n_0>0$ such that 
\[
H^j(X, B^c_{X, p^rD_e, n_0})=0 
\]
for every quadruple $(j, c, e, r) \in \Z^4$ satisfying 
$j<d, c>0, e \geq e_2$, and $r \geq 0$. 
\end{step}

\begin{proof}[Proof of Step \ref{s5:qFs to qFr}]
Take positive integers $e_2$ and $n_0$ as in Step \ref{s4:qFs to qFr}. 
We have the following exact sequence \cite[Lemma 3.14]{TWY}: 
\[
0 \to B^c_{X, p^rD_e, n_0} \to B^{c+1}_{X, p^rD_e, n_0} \to F_*^cB^1_{X, p^{c+r}D_e, n_0} \to 0. 
\]
By induction on $c\geq 1$, the assertion follows from Step \ref{s4:qFs to qFr}. 
This completes the proof of Step \ref{s5:qFs to qFr}. 
\end{proof}

Applying Step \ref{s5:qFs to qFr} by setting $j:=d-1, c:=e_2, e:=e_2$ and $r:=0$,   
we get $H^{d-1}(X, B^{e_2}_{X, D_{e_2}, n_0})=0$. 
Therefore, ($\star$) holds. 
This completes the proof of Theorem \ref{thm:qFs to qFr}. 
\end{proof}

\begin{corollary}\label{qFs to qFr for am Fano}
Let $k$ be an $F$-finite field of characteristic $p>0$ and 
let $X$ be a projective normal $\Q$-factorial variety over $k$. 
Assume that $X$ is strongly $F$-regular and $-K_X$ is ample.  
Then $X$ is quasi-$F$-split if and only if $X$ is globally quasi-$F$-regular.
\end{corollary}

\begin{proof}
Since $X$ is strongly $F$-regular and $\Q$-factorial, 
$X$ is divisorially Cohen-Macaulay \cite[Corollary 3.3]{patakfalvi-schwede14}. 
Hence the assertion follows from Theorem \ref{thm:qFs to qFr}.
\end{proof}

\begin{cor}\label{c klt dP QFR}
Let $k$ be a perfect field of characteristic $p>5$ and 
let $X$ be a klt del Pezzo surface over $k$. 
Then $X$ is globally quasi-$F$-regular. 
\end{cor}

\begin{proof}
Since $X$ is a normal surface, $X$ is divisorially Cohen-Macaulay. 
By $p>5$, $X$ is strongly $F$-regular. 
Hence the assertion follows from Theorem \ref{thm:qFs to qFr} and 
\cite[Theorem B]{KTTWYY1}. 
\end{proof}

\begin{corollary}\label{cor:QFR for log dP}
Let $k$ be a perfect field of characteristic $p>41$ and 
let $(X, \Delta)$ be a log del Pezzo surface with standard coefficients over $k$. 
Then $(X, \Delta)$ is globally quasi-$F$-regular.
\end{corollary}

\begin{proof}
Note that $X$ is divisorially Cohen-Macaulay and 
$(X, \Delta)$ is strongly $F$-regular. 
By \cite{KTTWYY2}*{Theorem C}, $(X,\Delta)$ is quasi-$F$-split. 
Then it follows from  Theorem \ref{thm:qFs to qFr} 
that $(X,\Delta)$ is globally quasi-$F$-regular.
\end{proof}

\subsection{Klt singularities}

\begin{theorem}\label{thm:QFR for klt}
Let $(X, \Delta)$ be a two-dimensional affine klt pair, 
where $\Gamma(X, \MO_X)$ is a ring essentially of finite type over a perfect field of characteristic $p>0$. 
Then $(X, \Delta)$ is quasi-$F$-regular.
\end{theorem}

\begin{proof}
Note that $X$ is $\Q$-factorial \cite[Corollary 4.11]{tanaka16_excellent}. 
By \cite[Theorem 7.12]{TWY}, $(X, \Delta)$ is feebly quasi-$F$-regular. 
Therefore, $(X, \Delta)$ is quasi-$F$-regular (Theorem \ref{thm:chara qFr pair2}). 
\qedhere

\end{proof}

\begin{theorem}\label{t-3-dim-klt}
Let $k$ be a perfect field of characteristic $p>42$ and 
let $(X, \Delta)$ be a three-dimensional $\Q$-factorial affine klt pair of finite type over $k$, where $\Delta$ has standard coefficients. 
Then $(X, \Delta)$ is quasi-$F$-regular. 
\end{theorem}

\begin{proof}
Fix a closed point $x \in X$. After replacing $X$ by an affine open neighbourhood of $x \in X$, 
there exists a projective birational morphism $f \colon Y \to X$ 
such that $x \in f(\Exc(f))$, $E := \Exc(f)$ is a prime divisor, $-(K_Y+E+f^{-1}_*\Delta)$ is ample, and $(Y, E+ f_*^{-1}\Delta)$ 
is a $\Q$-factorial plt pair (see \cite[the proof of Theorem 6.19]{KTTWYY1}).

Since $E$ is normal (see, for example, \cite[Theorem 2.11]{GNT16}), 
we can write $K_E + \Delta_E = (K_Y+E+f^{-1}_*\Delta)|_E$. 
If $f(E)$ is a point, then  $(E,\Delta_E)$ is globally quasi-$F$-regular 
(Corollary \ref{cor:QFR for log dP}), and hence $(E, \Delta_E)$ 
is globally quasi-$+$-regular \cite[Proposition 4.9]{TWY}. 
If $f(E)$ is a curve, then 
$(E,\Delta_E)$ is globally quasi-$+$-regular by \cite[Theorem 7.14]{TWY}. 
In any case, $(E, \Delta_E)$ is globally quasi-$+$-regular.

By the same argument as in \cite[the proof of Theorem 6.19]{KTTWYY1}, 
we see that $Y$ is strongly $F$-regular and divisorially Cohen-Macaulay. 
Then we can apply \cite[Corollary 6.8]{TWY}, so that 
$(Y, E+f^{-1}_*\Delta)$ is globally quasi-$+$-regular. 
Therefore, its {\cred push-forward} $(X, f_*(E+f^{-1}_*\Delta)) = (X, \Delta)$ is quasi-$+$-regular. 
Then $(X, \Delta)$ is quasi-$F$-regular by Theorem \ref{thm: rel q+R, qFrat and qFR}. 
\end{proof}

{\cccred

Contrary to Theorem \ref{t-3-dim-klt}, the analogous statement for $F$-splitting 
does not hold  as follows. 

\begin{prop}
\label{prop:non-F-pure canonical Gorenstein}
Let $k$ be an algebraically closed field of characteristic $p>5$. 
Then there exists a $\Q$-factorial canonical Gorenstein affine threefold $X$ over $k$ which is not $F$-pure. 
\end{prop}


\begin{proof}
Take a homogeneous polynomial $f(x, y, z) \in k[x, y, z]$ of degree $3$ 
such that 
$\Proj\,k[x, y, z]/(f(x, y, z))$ is a supersingular elliptic curve. 
For every integer $m>0$, we set 
\[
X_m := \Spec k[x, y, z, w]/( f(x, y, z) +w^m). 
\]
Fix an integer $n$ satisfying $n \geq p$ and $n \equiv 1 \mod 3$. 
Then $X_n$ is canonical Gorenstein affine threefold  which is not $F$-pure. 
\cite[Lemma 5.1 and Lemma 5.2]{CTW15a}. 

It is enough to show that $X_n$ is $\Q$-factorial. 
By \cite[the proof of Lemma 5.1]{CTW15a}, 
there exists a sequence 
\[
X_n =:Y_n \xleftarrow{f_{n-3}} Y_{n-3} \xleftarrow{f_{n-6}} Y_{n-6} \xleftarrow{f_{n-6}}\cdots \xleftarrow{f_1} Y_1
\]
such that (1)-(3) hold for every $m \in \{n, n-3, n-6, ..., 1\}$: 
\begin{enumerate}
\item There exists  an open cover $Y_m = X_m \cup U_m$ for some smooth threefold $U_m$. 
In particular, $Y_1$ is smooth and $Y_m$ has a unique singular point when $m >1$. 
\item $f_m : Y_{m} \to Y_{m+3}$ is the blowup at the unique singular point of $Y_{m+3}$ 
(i.e., the origin of $X_{m+3}$).   
\item $E_m := \Ex(f_m) \simeq \Proj\,k[x, y, z, w]/(f(x, y, z))$, which is the projective cone over the supersingular elliptic curve $\Proj\,k[x, y, z]/(f(x, y, z))$. 
In particular, $E_m$ is a projective normal rationally chain connected surface with $\rho(E_m)=1$. 
\end{enumerate}
Assume that $Y_m$ is $\Q$-factorial. 
By induction on $m$, it is enough to show that $Y_{m+3}$ is $\Q$-factorial. 
Since $f_m : Y_m \to Y_{m+3}$ is a projective birational morphism between quasi-projective normal threefolds with $\rho(E_m)=1$, 
we get $\rho(Y_m/Y_{m+3})=1$. 
By standard argument as in \cite[the first paragraph of the proof of Corollary 3.18]{KM98}, 
it is enough to show that every Cartier divisor $D$ on $Y_m$ satisfying $D \equiv_{f_m} 0$ is $f_m$-semi-ample, which follows from \cite[Theorem 2.9]{GNT16} 
(this is applicable, because $-K_{Y_m}$ is $f_m$-nef and $f_m$-big 
by $K_{Y_m} \sim f_m^*K_{Y_{m+3}}$). 
\end{proof}
}

\section{Miscellaneous results}

\subsection{More on quasi-$F^e$-splittings}
The goal of this subsection is to prove that a 
{Gorenstein} 
ring is quasi-$F$-split if and only if it is quasi-$F^\infty$-split. We define a $W_nR$-module $Q^e_{R, n}$ and 
a $W_nR$-module homomorphism 
$\Phi^e_{R, n}$ by the following pushout diagram: 
\[
\begin{tikzcd}
W_nR \arrow[r, "F^e"] \arrow[d, "{\cccred \mathbf{R}}^{n-1}"] & F^e_*W_nR \arrow[d]\\
R  \arrow[r, "\Phi^e_{R, n}"] & Q^e_{R,n}. 
\end{tikzcd}
\]
By applying $\Hom_{W_nR}(-, W_n\omega_R)$ we get the following diagram
\[
\begin{tikzcd}
W_n\omega_R  & \arrow{l}[swap]{T^e_n} F^e_*W_n\omega_R \arrow[d]\\
{\ccred \omega_R}   \arrow{u}{({\cccred \mathbf{R}}^{n-1})^*}  & \arrow{l}[swap]{(\Phi^e_{R, n})^*} \Hom_{W_nR}(Q^e_{R,n}, W_n\omega_R). 
\end{tikzcd}
\]
Recall that $R$ is $n$-quasi-$F^e$-split if and only if $(\Phi^e_{R, n})^*$ is surjective.

In what follows, we will denote the usual projection of Witt rings $W_nR \to W_mR$ for $n \geq m$ by ${\cccred \mathbf{R}}_{n,m}$. Moreover, we define a $W_nR$-module $Q^e_{R, n,m}$ and 
a $W_nR$-module homomorphism 
$\Phi^e_{R, n,m}$ for $n \geq m$ by the following pushout diagram: 
\[
\begin{tikzcd}
W_nR \arrow[r, "F^e"] \arrow{d}[swap]{{\cccred \mathbf{R}}_{n,m}} & F^e_*W_nR \arrow[d]\\
W_mR  \arrow[r, "\Phi^e_{R, n,m}"] & Q^e_{R,n,m}. 
\end{tikzcd}
\]
In other words, 
\[
Q^e_{R,n,m} = {\rm Coker}(F^m_*W_{n-m}R \xrightarrow{{\cccred \mathbf{V}}^m} W_nR \xrightarrow{F^e} F^e_*W_nR).
\]
By definition of a pushout, 
we have the induced $W_nR$-module homomorphism 
$Q^e_{R,n,m} \to Q^e_{R,n,l}$ for every $1 \leq l \leq m$. In particular, in the case of $l=1$, we get a factorisation:
\begin{equation} \label{eq:phie-rnm}
\begin{tikzcd}
 W_mR \ar{r}{\Phi^e_{R, n,m}} &  Q^e_{R,n,m} \ar{r} &  Q^e_{R,n}.
\end{tikzcd}
\end{equation}

\begin{remark}
For future use, we also construct the short exact sequence:
\begin{equation} \label{eq:phie-rnm-diagram}
0 \to F_*Q^e_{R, n-1,m-1} \xrightarrow{\psi} Q^e_{R,n,m} \to  F^e_*R \to 0
\end{equation}
by applying Snake Lemma to:
\[
\begin{tikzcd}
 & F^m_*W_{n-m}R \ar{r}{=} \ar{d}{F^e \circ {\cccred \mathbf{V}}^{m-1}} & F^m_*W_{n-m}R \ar{d}{F^e \circ {\cccred \mathbf{V}}^m}   &  &  \\
0 \ar{r} & F^{e+1}_*W_{n-1}R \ar{r}{{\cccred \mathbf{V}}} & F^e_*W_nR \ar{r}{{\cccred \mathbf{R}}^{n-1}} & F^e_*R \ar{r} & 0.
\end{tikzcd}
\]
\end{remark}

\begin{proposition} \label{prop:yet-new-def-of-quasiF^esplit}
We work in the general setting (Notation \ref{n-non-local}). Assume that $R$ is Cohen-Macaulay and $k$-quasi-$F^e$-split. Then
\[
(\Phi^e_{R, n,m})^* := \Hom_{W_nR}(\Phi^e_{R, n,m}, W_n\omega_R) \colon \Hom_{W_{n}R}(Q^e_{R, n,m}, W_n\omega_R) \to W_m\omega_R,
\]
is surjective for every $m \geq 1$ and $n = k + m -1$,
\end{proposition}
\begin{proof}
We may assume that $R$ is a local ring. Then, by Matlis duality, 
it suffices to verify that
\[
H^d_\m(W_mR) \xrightarrow{H^d_\m(\Phi^e_{R, n,m})} H^d_\m(Q^e_{R,n,m}) 
\]
is injective for every $m \geq 1$ and $n = k + m - 1$. We argue by ascending induction on $m$, specifically, we may assume by inductive assumption that 
\begin{equation} \label{eq:injection-qernm} 
H^d_\m(W_{m-1}R) \xrightarrow{H^d_\m(\Phi^e_{R, n-1,m-1})} H^d_\m(Q^e_{R,n-1,m-1}) \quad \text{ is injective.}
\end{equation}
Note that the base case $m=1$ of the induction 
{follows from \cite[Lemma 3.10 and Proposition 3.20]{TWY}}.

Pick $\alpha \in H^d_\m(W_mR)$ such that 
\[
H^d_\m(\Phi^e_{R,n,m})(\alpha) = 0 \in H^d_\m(Q^e_{R,n,m}).
\] 
Then consider the following diagram (see (\ref{eq:phie-rnm})):
\[
\begin{tikzcd}
Q^e_{R,n,m} \ar{r} & Q^e_{R,n} \\
W_mR  \ar{u}{\Phi^e_{R,n,m}} \ar{r}{{\cccred \mathbf{R}}^{m-1}} & R.  \ar{u}{\Phi^e_{R,n}}
\end{tikzcd}
\]
Now $H^d_\m(\Phi^e_{R,n})$ is injective, 
because $n \geq k$ and $R$ is $k$-quasi-$F^e$-split. Therefore 
\[
H^d_\m({\cccred \mathbf{R}}^{m-1})(\alpha) = 0 \in H^d_\m(R).
\]

Now, consider the following diagram (see (\ref{eq:phie-rnm-diagram}) and construction thereof):
\[
\begin{tikzcd}
0 \ar{r} & F_*Q^e_{R, n-1,m-1} \ar{r}{\psi} & Q^e_{R,n,m} \ar{r} & F^e_*R \ar{r} & 0\\
0 \ar{r}  & F_*W_{m-1}R \ar{u}{F_*\Phi^e_{R,n-1,m-1}} \ar{r}{{\cccred \mathbf{V}}} & W_mR  \ar{u}{\Phi^e_{R,n,m}} \ar{r}{{\cccred \mathbf{R}}^{m-1}} & R \ar{u}{F^e} \ar{r} & 0.
\end{tikzcd}
\]
By the above paragraph, 
we have $\alpha = H^d_{\m}({\cccred \mathbf{V}})(\beta)$ 
for some $\beta \in H^d_\m({\cccred F}_*W_{m-1}R)$. 
Moreover, $H^d_\m(F_*\Phi^e_{R,n-1,m-1})$ is injective by (\ref{eq:injection-qernm}) and $H^d_\m(\psi)$ is injective because $H^{d-1}_\m(F^e_*R)=0$ in view of Cohen-Macauliness of $R$. Since 
\[
0 = H^d_\m(\Phi^e_{R,n,m})(\alpha) = H^d_\m(\Phi^e_{R,n,m} \circ {\cccred \mathbf{V}})(\beta),
\]
it holds that $\beta =0$. 
Hence $\alpha=0$ as well, concluding the proof of the injectivity of $H^d_\m(\Phi^e_{R,n,m})$.
\qedhere 

    
\end{proof}

\begin{remark} \label{remark:qernm-diagram}
Take integers $n \geq m$ and consider the following diagram
\[
\begin{tikzcd}
0 \ar{r} & F^{m}_*W_{n-m}R \ar{r}{F^e \circ {\cccred \mathbf{V}}^m} & F^e_*W_nR \arrow[r] & Q^e_{R,n,m} \ar{r} & 0 \\
0 \ar{r} & F^{m}_*W_{n-m}R \ar{u}{=} \ar{r}{{\cccred \mathbf{V}}^m} & W_nR \arrow[u, "F^e"] \arrow[r, "{\cccred \mathbf{R}}_{n,m}"] & W_mR \arrow[u, "\Phi^e_{R, n,m}"]  \ar{r} & 0.
\end{tikzcd}
\]    
By applying $\Hom_{W_nR}(-, W_n\omega_R)$ we get:
\[
\begin{tikzcd}
 F^m_*W_{n-m}\omega_R \ar{d}{=}  & \ar{l} F^e_*W_n\omega_R \ar{d}{T^e_{n}} \arrow[l] & \ar{l} \Hom_{W_nR}(Q^e_{R,n,m}, W_n\omega_R) \ar{d}{(\Phi^e_{R, n,m})^*}    & \ar{l} 0 \\
 F^m_*W_{n-m}\omega_R   & \ar{l}{{\cccred \mathbf{V}}^*} W_n\omega_R  & W_m\omega_R \arrow{l}{{\cccred \mathbf{R}}^*_{n,m}} & \ar{l} 0.
\arrow[phantom, from=1-2, to=2-3, "(\star)" description]
\end{tikzcd}
\]    
By diagram chase, $(\star)$ is a pullback diagram. 
\end{remark}

\begin{corollary} \label{cor:yet-new-def-of-quasiF^esplit}
We work in the general setting (Notation \ref{n-non-local}). 
Suppose that $R$ is Cohen-Macaulay and $k$-quasi-$F^e$-split. Then with notation of Remark \ref{remark:qernm-diagram} for every $m \geq 1$ and $n=k+m-1$ we have that:
\[
{\cccred \mathbf{R}}^*_{n,m}(W_m\omega_R) \subseteq  {\rm Im}\big(T^e_{n} \colon F^e_*W_n\omega_R \to W_n\omega_R\big).
\]
\end{corollary}
\begin{proof}
This is immediate from Remark \ref{remark:qernm-diagram} and Proposition \ref{prop:yet-new-def-of-quasiF^esplit}.
\end{proof}
Note that the inclusion $({\cccred \mathbf{R}}^{n-1})^*(\omega_R) \subseteq  {\rm Im}\big(T^e_{n} \colon F^e_*W_n\omega_R \to W_n\omega_R\big)$ is equivalent to $R$ being $n$-quasi-$F^e$-split {when $R$ is Gorenstein}.
\begin{theorem}
We work in the general setting (Notation \ref{n-non-local}). Suppose that $R$ is 
{Gorenstein}. 
Further assume that $R$ is $m$-quasi-$F^e$-split and $k$-quasi-$F^{e'}$-split. Then $R$ is $(k+m-1)$-quasi-$F^{e+e'}$-split.
\end{theorem}
\begin{proof}
Set $n= k+m-1$ and consider the following diagram.
\[
\begin{tikzcd}[column sep = large]
&&  F^{e+e'}_*W_{n}\omega_R \ar{d}{T^{e'}_{n}} \\
& F^e_*W_{m}\omega_R \ar{r}{{\cccred \mathbf{R}}^*_{n,m}} \ar{d}{T^e_m} & F^{e}_*W_{n}\omega_R \ar{d}{T^e_{n}} \\
\omega_R \ar{r}{({\cccred \mathbf{R}}^{m-1})^*} \ar[bend right = 30]{rr}{({\cccred \mathbf{R}}^{n-1})^*} & W_{m}\omega_R \ar{r}{{\cccred \mathbf{R}}^*_{n,m}} & W_{n} \omega_R. 
\end{tikzcd}
\]
Since $R$ is $k$-quasi-$F^{e'}$-split and $m$-quasi-$F^e$-split, Corollary \ref{cor:yet-new-def-of-quasiF^esplit} yields
\begin{align}
{\cccred \mathbf{R}}^*_{n,m}(F^e_*W_{m}\omega_R) &\subseteq {\rm Im}\Big(T^{e'}_{n} \colon F^{e+e'}_*W_{n}\omega_R \to F^{e}_*W_{n}\omega_R \Big) \tag{$\dagger$}\\
({\cccred \mathbf{R}}^{m-1})^*(\omega_R) &\subseteq {\rm Im}\Big(T^{e}_{m} \colon F^{e}_*W_{m}\omega_R \to W_{m}\omega_R \Big). \tag{$\dagger\dagger$}
\end{align}
Therefore
\begin{align*}
\left({\cccred \mathbf{R}}^{n-1}\right)^*\left(\omega_R\right) &= {\cccred \mathbf{R}}^*_{n,m}\left(\left({\cccred \mathbf{R}}^{m-1}\right)^*\left(\omega_R\right)\right) \\
&\overset{(\dagger\dagger)}{\subseteq} {\cccred \mathbf{R}}^*_{n,m}\left(T^e_m\left(F^e_*W_m\omega_R\right)\right) \\[0.3em]
&= T^{e}_n\left({\cccred \mathbf{R}}^*_{n,m}\left(F^e_*W_m\omega_R\right)\right) \\
&\overset{(\dagger)}{\subseteq} T^{e}_n\left(T^{e'}_n\left(F^{e+e'}_*W_n\omega_R\right)\right) \\
&= T^{e+e'}_n\left(F^{e+e'}_*W_n\omega_R\right).
\end{align*}
{Since $R$ is Gorenstein}, $R$ is $n$-quasi-$F^{e+e'}$-split.
\end{proof}

In particular, a {Gorenstein} 
ring is quasi-$F$-split if and only if it is quasi-$F^\infty$-split.
\begin{corollary} \label{cor:quasiF=quasiFinfty}
We work in the general setting (Notation \ref{n-non-local}). Suppose that $R$ is 
{Gorenstein} 
and $n$-quasi-$F$-split. Then $R$ is $(ne-e+1)$-quasi-$F^e$-split for every integer $e>0$.
\end{corollary}
Note that this bound should be sharp by taking cones over Calabi-Yau varieties in view of \cite[Theorem 7.1]{TWY}.

\subsection{Three-dimensional quasi-$F$-regular rings are Cohen-Macaulay}\label{ss 3-dim CM}

\begin{lem}\label{lem:vanishing-loc-coh}
We use Notation \ref{n-normal-local}.
Let $D$ be a $\Q$-Cartier Weil divisor on $\Spec R$.
Then there exists $c \in R^\circ$ such that 
\[
[c] \cdot H^{d-1}_{\m}(W_n\cO_X(p^eD))=0
\]
for every $e \in \Z_{>0}$ and every $n \in \Z_{>0}$. %
\end{lem}
\begin{proof}
Since $R$ is local and $D$ is $\Q$-Cartier, there are only finitely many possibilities for $R(p^eD)$, with $e \geq 0$, up to isomorphism. 
Thus there exists $f \in R^\circ$ such that $f \cdot H^{d-1}_{\m}(R(p^eD))=0$ for every integer $e >0 $ by Lemma \ref{lem: ann local coh}.
Set $c:=f^2$. 
Let us show that $c$ satisfies the required property by induction on $n$.
The base case $n=1$ of this induction  follows from the choice of $f$.

Assume $n \geq 2$. 
Consider the exact sequence
\[
0 \to F_*W_{n-1}R(p^{e+1}D) \xrightarrow{{\cccred \mathbf{V}}} W_nR (p^eD) \xrightarrow{{\cccred \mathbf{R}}^{n-1}} R(p^eD) \to 0.
\]
Taking local cohomologies, we get the exact sequence
\[
F_*H^{d-1}_{\m}(W_{n-1}R(p^{e+1}D)) \xrightarrow{{\cccred \mathbf{V}}} H^{d-1}_{\m}(W_nR(p^eD)) \xrightarrow{{\cccred \mathbf{R}}^{n-1}} H^{d-1}_{\m}(R(p^eD)).
\]
Pick $\alpha \in H^{d-1}_{\m}(W_nR(p^eD))$.
By the choice of $f$, we have ${\cccred \mathbf{R}}^{n-1}([f]\alpha)=0$.
Therefore, there exists $\beta \in H^{d-1}_{\m}(W_{n-1}R(p^{e+1}D))$ such that 
${\cccred \mathbf{V}}(F_*\beta)=[f]\alpha$.
By the induction hypothesis, we have $[c]\beta =0$, which implies 
\[
[c]\alpha=[f]{\cccred \mathbf{V}}(F_*\beta)={\cccred \mathbf{V}}(F_*([f^p]\beta))={\cccred \mathbf{V}}(F_*([cf^{p-2}]\beta))=0,
\]
as required.
\end{proof}

\begin{lemma}\label{lem:CM lem surj}
We use Notation \ref{n-normal-local}. 
Let $\Delta$ be an effective  $\Q$-divisor on $X$ such that $(X, \Delta)$ is quasi-$F$-regular.  Take a $\Q$-Cartier $\Q$-divisor $D$ on $\Spec R$ 
such that $\{p^rD\} \leq \Delta$ 
for every integer {$r \geq 1$}. 
Fix an integer $m\geq 2$. 
Then the following hold. 
\begin{enumerate}
    \item The $W_mR$-module homomorphism 
    \[
    {\cccred \mathbf{R}}^{m-1} : H^{d-1}_{\m}(W_mR(D)) \to H^{d-1}_{\m}(R(D))
    \]
    is surjective. 
    \item The sequence 
    \[
    0 \to \omega_R(-D) \xrightarrow{({\cccred \mathbf{R}}^{m-1})^*} 
    W_m\omega_R(-D) \xrightarrow{{\cccred \mathbf{V}}^*} F_*W_{m-1}\omega_R(-pD) \to 0.
    \]
    is exact, where this is obtained by applying $\Hom_{W_mR}(-, W_m\omega_R)$ 
    to the exact sequence 
    \[
0 \to F_*W_{m-1}R(pD) \xrightarrow{{\cccred \mathbf{V}}} W_mR(D) \xrightarrow{{\cccred \mathbf{R}}^{m-1}} R(D) \to 0. 
    \]
\end{enumerate}
\end{lemma}

\begin{proof}
Fix an integer $n \geq 1$ such that $R$ is $n$-quasi-$F$-regular.

Let us show (1).
We consider the following commutative diagram 
in which each horizontal sequence is exact:
\[
\begin{tikzcd}
    0 \arrow[r]  & F_*W_{n+m-1}R(pD) \arrow[r,"{\cccred \mathbf{V}}"] \arrow[d,"{\cccred \mathbf{R}}^{n-1}"] & W_{n+m}R(D) \arrow[d,"{\cccred \mathbf{R}}^{n-1}"] \arrow[r,"{\cccred \mathbf{R}}^{n+m-1}"] & R(D) \arrow[d,equal] \arrow[r] & 0 \\
    0 \arrow[r] & F_*W_{m-1}R(pD) \arrow[r,"{\cccred \mathbf{V}}"] & W_mR(D) \arrow[r,"{\cccred \mathbf{R}}^{m-1}"] & R(D) \arrow[r] & 0.
\end{tikzcd}
\]
Taking local cohomologies, we obtain
\[
\begin{tikzcd}
    H^{d-1}_{\m}(W_{n+m}R(D)) \arrow[r] \arrow[d,"{\cccred \mathbf{R}}^{\textcolor{purple}{n}-1}"] & H^{d-1}_{\m}(R(D)) \arrow[r,"\sigma_{n+m}"] \arrow[d,equal] & H^{d}_{\m}(F_*W_{n+m-1}R(pD))  \arrow[d,"{\cccred \mathbf{R}}^{n-1}"]  \\
    H^{d-1}_{\m}(W_mR(D)) \arrow[r] & H^{d-1}_{\m}(R(D)) \arrow[r,"\sigma_{m}"] & H^d_{\m}(F_*W_{m-1}R(pD)).
\end{tikzcd}
\]
By $\mathrm{Im}(\sigma_{n+m}) ={\cccred \mathbf{V}}^{-1}(0) \subseteq {\cccred \mathbf{V}}^{-1}(\wt{0^*_{n+m}}) = F_*\wt{0^{*}_{pD,n+m-1}}$ 
(Proposition \ref{prop: tight cl l=n-log}(2)), 
we have 
\[
\mathrm{Im}(\sigma_{m})={\cccred \mathbf{R}}^{n-1}(\mathrm{Im}(\sigma_{n+m})) \subseteq {\cccred \mathbf{R}}^{n-1}(F_*\wt{0^*_{pD,n+m-1}})=0,
\]
where the last equality follows from Proposition \ref{prop:  qFreg 0-map}.
Therefore, 
$\sigma_{m}$ is zero, that is, 
the map 
\[
H^{d-1}_{\m}(W_mR(D)) \to H^{d-1}_{\m}(R(D))
\]
is surjective. 
Thus (1) holds.

Let  us show (2). 
By applying $\text{RHom}_{W_mR}(-,W_m\omega_R^{\mydot})$ to the exact sequence 
    \[
0 \to F_*W_{m-1}R(pD) \xrightarrow{{\cccred \mathbf{V}}} W_mR(D) \xrightarrow{{\cccred \mathbf{R}}^{m-1}} R(D) \to 0, 
    \]
it is enough to show that 
\[
{({\cccred \mathbf{R}}^{m-1})^*}: H^{-d +1}\text{RHom}_{W_mR}(R(D),W_m\omega_R^{\mydot}) 
\to 
H^{-d+1}\text{RHom}_{W_mR}(W_mR(D),W_m\omega_R^{\mydot}) 
\]
is injective. 
This homomorphism is obtained by applying 
the Matlis {\cred duality} functor $\Hom_{W_mR}(-, E)$ 
to  the surjection in (1). 
Thus (2) holds. 
\qedhere 
    
\end{proof}

\begin{theorem}\label{thm: CM general}
Let $R$ be an $F$-finite normal integral domain of characteristic $p>0$. 
Assume that $\dim R =3$, $R$ is quasi-$F$-regular, and $K_R$ is $\Q$-Cartier. 
Then 
$R(D)$ is Cohen-Macaulay for every $\Q$-Cartier Weil divisor $D$ on $\Spec R$.
\end{theorem}

\begin{proof}
We may assume that $R$ is a local ring. 
Let $\m$ be the maximal ideal of $R$. 
Fix a $\Q$-Cartier Weil divisor $D$ and 
an integer $n>0$ such that $R$ is $n$-quasi-$F$-regular. 
 By definition,
\[
W_{m}\omega_R(p(D-K_R))=\cHom(W_{m}R(p(K_R-D)),W_{m}\omega_R)
\]
for every integer $m \geq 1$. 
Since $W_{m}\omega_R$ is a dualising $W_mR$-module, $W_{m}\omega_R$ 
is $(S_2)$ 
\cite[\href{https://stacks.math.columbia.edu/tag/0AWK}{Tag 0AWK}]{stacks-project}. 
Therefore,  $W_{m}\omega_R(p(D-K_R))$ is also 
$(S_2)$, 
which implies 
\[
H^1_{\m}(W_{m}\omega_R(p(D-K_R)))=0
\]
for every integer $m \geq 1$.
Since $K_R-D$ is a $\Q$-Cartier Weil divisor, 
we get the following exact sequence (Lemma \ref{lem:CM lem surj}(2)): 
\[
0 \to \omega_R(D-K_R) \xrightarrow{({\cccred \mathbf{R}}^{n-1})^*} W_n\omega_R(D-K_R) 
\xrightarrow{{\cccred \mathbf{V}}^*} F_*W_{n-1}\omega_R(p(D-K_R)) \to 0.
\]
By $\omega_R(D-K_R) \simeq R(D)$, we obtain an injection
\begin{equation}\label{eq: inj}
({\cccred \mathbf{R}}^{n-1})^* :  H^2_\m(R(D)) \hookrightarrow H^2_\m(W_n\omega_R(D-K_R)).   
\end{equation}

Take $c_1 \in R^\circ$ such that 
$c_1 \cdot H^{2}_{\m}(W_nR(p^eD))=0$ for every integer $e \geq 1$, 
whose existence is guaranteed by Lemma \ref{lem:vanishing-loc-coh}. 
For $t_D \in R^{\circ} \cap \tau(R)$, we set $c_2 := {t_D^4}$ and $c :=c_1c_2$. 
Since $R$ is $n$-quasi-$F$-regular,  there exists an integer $e \geq 1$ such that $(R, (1/p^e)\mathrm{div}(c))$ is 
$n$-quasi$F^e$-split (Theorem \ref{thm:chara qFr pair}). 
By Proposition \ref{prop:qFs-divisor-split}, we have the commutative diagram 
\[
\begin{tikzcd}
    W_nR(D) \arrow[r,"{(\cdot F^e_*[c])}\circ F^e"] 
    \arrow[d,"{\cccred \mathbf{R}}^{n-1}"'] & F^e_*W_nR(p^eD) \arrow[ldd,"\alpha"] \\
    R(D) \arrow[d, "({\cccred \mathbf{R}}^{n-1})^*"'] \\
    W_n\omega_R(-K_R+D).
\end{tikzcd}
\]
Taking local cohomologies, we obtain
\[
\begin{tikzcd}
    H^2_{\m}(W_nR(D)) \arrow[r,"{(\cdot F^e_*[c])}\circ F^e"] \arrow[d,"{\cccred \mathbf{R}}^{n-1}"'] & H^2_{\m}(F^e_*W_nR(p^eD)) \arrow[ldd,"\alpha"] \\
    H^2_{\m}(R(D)) \arrow[d, "({\cccred \mathbf{R}}^{n-1})^*"'] \\
    H^2_{\m}(W_n\omega_R(-K_R+D)).
\end{tikzcd}
\]
By using (\ref{eq: inj}) and $[c] \cdot H^2_{\m}(W_nR(p^eD))=0$, we get that 
\[
{\cccred \mathbf{R}}^{n-1} : H^2_\m(W_nR(D)) \to H^2_\m(R(D))
\]
is zero.
By  Lemma \ref{lem:CM lem surj}(1), it is also surjective.
Therefore $H^2_\m(R(D))=0$.
\end{proof}

\section{Appendix: dualising complexes on Witt rings}



Given an $\F_p$-algebra $R$, we set $\Omega^i_R := \Omega^i_{R/\F_p}$.  
The purpose of  this appendix is to prove the following theorem. 

\begin{theorem}\label{thm:X-exist dualizing complex}
Let $X$ be a separated scheme of finite type over an $F$-finite Noetherian $\mathbb{F}_p$-algebra. 
For $e \geq 0$ and $n  \geq m \geq 1$, 
we define the finite morphism $\iota^e_{X, m, n} \colon W_mX \to W_nX$ as the composition 
\[
\iota^e_{X, m, n} \colon W_mX \xrightarrow{F^e} W_mX \hookrightarrow W_nX, 
\]
where $W_nX := (X, W_n\MO_X), W_mX := (X, W_m\MO_X)$, 
and 
$W_mX \hookrightarrow W_nX$ denotes the closed immersion corresponding to ${\cccred \mathbf{R}}^{n-m} \colon W_n\MO_X \to W_m\MO_X$. 
Then there exist objects 
\[
W_1\omega_X^{\mydot},\, W_2\omega_X^{\mydot},\, ...
\]
of the derived category of $W_n\cO_X$-modules and a set of isomorphisms
\[
\{\rho^e_{X,m,n} \colon W_m\omega_X^{\mydot} \xrightarrow{\simeq} (\iota^e_{X,m,n})^{!}W_n\omega_X^{\mydot} 
\,|\, e \geq 0,  n  \geq m \geq 1\}
\] 
in the derived category of $W_n\cO_X$-modules which satisfy the following.
\begin{enumerate}
\item $W_n\omega_X^{\mydot}$ is a dualising complex on 
$W_nX = (X, W_n\cO_X)$ for every integer $n>0$. 
\item 
The equality $(\iota^{e'}_{X,l,m})^{!}\rho^{e}_{X,m,n} \circ \rho^{e'}_{X,l,m}=\rho^{e+e'}_{X,l,n}$ holds for all $e,e',l,m,n \in \Z$ satisfying $e,e' \geq 0$ and $n \geq m \geq l \geq 1$.
\item 
For all $e,m,n \in \Z$ satisfying $e \geq 0$ and $n \geq m \geq 1$, 
we define the morphism $T^e_{X,m,n}$ by
\[
T^e_{X,m,n} \colon (\iota^e_{X,m,n})_*W_m\omega_X^{\mydot} \xrightarrow{\rho^e_{X,m,n}} (\iota^e_{X,m,n})_*(\iota^e_{X,m,n})^{!}W_n\omega_X^{\mydot} \xrightarrow{{\rm adj}} W_n\omega_X^{\mydot},
\]
where ${\rm adj}$ is the natural morphism induced by the fact that $(\iota^e_{X,m,n})^{!}$ is the right adjoint of $(\iota^e_{X,m,n})_*$ 
\cite[Tag 0A9Y]{stacks-project}. 
Then the equality $T^{e}_{X,m,n} \circ (\iota^{e}_{X,m,n})_*T^{e'}_{X,l,m}=T^{e+e'}_{X,l,n}$ holds for all $e,e',l,m,n \in \Z$ satisfying $e,e' \geq 0$ and $n \geq m \geq l \geq 1$.
\end{enumerate}
\end{theorem}
 
\noindent 
Note that we set $(\iota^e_{X, m, n})_* := R(\iota_{X, m, n})_*$, 
as  $\iota^e_{X, m, n}$ is a finite morphism.

The proof of Theorem \ref{thm:X-exist dualizing complex} will be reduced to the case when $X=\Spec R$ and $R$ is regular. 
In this case,  
the de Rham-Witt complex $\{W_n\Omega_R^{i}\}_{i \geq 0}$ will play a crucial role. 
We now summarise some of its properties.

\begin{nothing}[de Rham-Witt complexes] 
Let $R$ be a regular $F$-finite $\F_p$-algebra. 
In \cite{illusie_de_rham_witt}, Illusie introduced the de Rham-Witt complex 
\[
0 \to W_n\Omega_R^0 \xrightarrow{d} W_n\Omega_R^1 \xrightarrow{d} 
\cdots 
\xrightarrow{d} W_n\Omega_R^i \xrightarrow{d} W_n\Omega_R^{i+1} \xrightarrow{d} \cdots, 
\]
where each $W_n\Omega_R^i$ is a $W_nR$-module and $d$ is an additive homomorphism satisfying $d^2=0$. 
For all integers $n \geq 1$ and $i \geq 0$, we have $W_{n+1}R$-module homomorphisms 
\[
F: W_{n+1}\Omega_R^i \to  F_*W_n\Omega_R^i \qquad {\cccred \mathbf{V}} :  W_{n}\Omega_R^i \to  F_*W_{{\cred n+1}}\Omega_R^i, 
\qquad 
{\cccred \mathbf{R}}: W_{n+1}\Omega_R^i \to W_n\Omega_R^i, 
\]
which satisfies the following properties \cite[Ch.\ I, Proposition 2.18]{illusie_de_rham_witt}:
\begin{enumerate}
\item $p=F{\cccred \mathbf{V}}={\cccred \mathbf{V}}F  : W_n\Omega_R^i \to W_n\Omega_R^i$. 
\item $dF = pFd : W_{n+1}\Omega_R^i \to W_n\Omega_R^{i+1}$. 
\item ${\cccred \mathbf{V}}d = pd{\cccred \mathbf{V}} : W_{n}\Omega_R^i \to W_{n+1}\Omega_R^{i+1}$. 
\item $Fd{\cccred \mathbf{V}}=d : W_n\Omega_R^i \to W_n\Omega_R^{i+1}$. 
\end{enumerate}
Since  $R$ is regular, 
Illusie's de Rham-Witt complex $W_n\Omega_R^i$ coincides 
with the saturated de Rham-Witt complex $\mathcal W_n\Omega_R^i$ introduced by Bhatt-Lurie-Mathew 
\cite[the proof of Theorem 4.4.12]{BLM}. 
\end{nothing}

\begin{lemma}\label{lem:covered be smooth}
Let $(R,\m)$ be a regular $F$-finite local $\F_p$-algebra.
Then there exists a smooth $\F_p$-algebra $S$ and a flat ring homomorphism $\varphi \colon S \to R$ such that 
\begin{enumerate}
\item 
$\mathfrak p R= \m$ for 
$\mathfrak p := \varphi^{-1}(\m)$, and 
\item the induced $S$-module homomorphism
\[
F^e_*\Omega_S^i \otimes_S R \to F^e_*\Omega_R^i
\]
is an isomorphism for every integer $e \geq 0$ and every integer $i \geq 0$.  
\end{enumerate}  
\end{lemma}
\begin{proof}
By \cite[{the proof of} Lemma 3.1]{KTY22}, there exists 
a $p$-basis $x_1,\ldots,x_r$ of $R$ such that $x_1,\ldots,x_s$ form a regular system of parameter (i.e., $\m=(x_1, ..., x_s)$ and $s=\dim R$) for some $0 \leq s \leq r$. 
In particular, 
for every integer $e > 0$, 
the set 
\[
\{F^e_*x_1^{i_1}\cdots x_r^{i_r} \mid 0 \leq i_1,\ldots,i_r \leq p^e-1 \}
\]
is a free $R$-linear basis of the $R$-module $F^e_*R$. 
In other words, 
$\{ x_1^{i_1}\cdots x_r^{i_r} \mid 0 \leq i_1,\ldots,i_r \leq p^e-1 \}$ 
is a free $R^{p^e}$-linear basis of the $R^{p^e}$-module $R$. 
Here $R^{p^e} := \{ r^{p^e} \in R \,|\, r \in R\}$ is a subring of $R$, 
which coincides with the {image of the} $e$-th iterated Frobenius ring homomorphism 
$F^e : R \to R, r \mapsto r^{p^e}$. 

\begin{claim}\label{claim:p-basis forms a basis}
$dx_1,\ldots,dx_r$ form a free $R$-linear basis of the $R$-module $\Omega^1_R$.
\end{claim}

\begin{proof}[Proof of Claim \ref{claim:p-basis forms a basis}]
We have the following exact sequence \cite[Theorem 25.1]{Matsumura}: 
\[
\Omega^1_{R^p} \otimes_R R \to \Omega^1_{R} \to \Omega^1_{R/R^p} \to 0.
\]
The first map is zero (because $d(r^p) = pd(r^{p-1})=0$ for every $r \in R$), 
and hence we get $\Omega^1_R \xrightarrow{\simeq} \Omega^1_{R/R^p}$.
Furthermore, we consider the $R^p$-algebra homomorphism
\begin{eqnarray*}
\theta: R^p[T_1,\ldots,T_r]/(T_1^p-x_1^p,\ldots,T_r^p-x_r^p) &\to& R\\
T_i &\mapsto& x_i.
\end{eqnarray*}
Since $x_1,\ldots,x_r$ form a $p$-basis of $R$, 
$\theta$ is an isomorphism.
Thus we get $R$-module isomorphisms 
\[
\Omega^1_{R/R^p} \overset{(\star)}{\simeq} \Omega^1_{R^p[T_1,\ldots,T_r]/R^p}/(d(T_1^p-x_1^p),\ldots,d(T_r^p-x_r^p)) \simeq \Omega^1_{R^p[T_1,\ldots,T_r]/R^p}, 
\]
where $(\star)$ follows from \cite[Theorem 25.2]{Matsumura}. 
Note that the right term $\Omega^1_{R^p[T_1,\ldots,T_r]/R^p}$ 
has an $R$-linear basis  $dT_1,\ldots,dT_r$, 
and hence $dx_1,\ldots,dx_r$ form an $R$-linear basis of 
$\Omega^1_{R/R^p} (\simeq  \Omega^1_{R})$. 
This completes the proof of Claim \ref{claim:p-basis forms a basis}. 
\end{proof}

Set $S:=\F_p[X_1,\ldots,X_r]$ and consider the following $\F_p$-algebra homomorphism: 
\begin{eqnarray*}
\varphi \colon S=\F_p[X_1,\ldots,X_r] &\to& R\\
X_i &\mapsto& x_i.
\end{eqnarray*}
Then $X_1,\ldots,X_r$ form a $p$-basis of $S$ and $dX_1,\ldots,dX_r$ form a free $S$-linear basis of the $S$-module $\Omega^1_S$.
Therefore, the image of the $S$-linear basis 
\[
\{ F_*^eX_1^{i_1} \cdots X_r^{i_r} dX_j\, |\, 0 \leq i_1, ..., i_r \leq p^e-1, 1 \leq j \leq r\}
\]
of $F^e_*\Omega_S$ by the natural map $F^e_*\Omega_S \to F^e_*\Omega_R$ is the basis 
\[
\{ F_*^ex_1^{i_1} \cdots x_r^{i_r} dx_j\, |\, 0 \leq i_1, ..., i_r \leq p^e-1, 1 \leq j \leq r\}
\]
of $F^e_*\Omega_R$.
Thus we get the natural $R$-module isomorphism 
\[
F^e_*\Omega_S \otimes_S R \xrightarrow{\simeq} F^e_*\Omega_R. 
\]

Finally, it is enough  to prove that $\varphi : S \to R$ is flat.
Set $\mathfrak p':=(X_1,\ldots,X_s)$. 
Then $\mathfrak p' R=(x_1,\ldots,x_s)=\m$. 

\begin{claim}\label{claim:compare ideal}
We have $(\mathfrak p:=)\varphi^{-1}(\m) =\mathfrak p'$.
\end{claim}

\begin{proof}[Proof of Claim]
The inclusion 
$\varphi^{-1}(\m) \supseteq \mathfrak p'$ is clear.
Let us show the opposite 
inclusion $\varphi^{-1}(\m) \subseteq \mathfrak p'$. 
Take $f \in \varphi^{-1}(\m)$, i.e., $f \in S = \F_p[X_1, ..., X_r]$ and $\varphi(f) \in \m$. 
We can write 
\[
f=\sum_{v \in \Z_{\geq 0}^r} n_v X^v \in \F_p[X_1, ..., X_r],
\]
where $n_v \in \F_p$ and we set $X^v := X_1^{v_1} \cdots X_r^{v_r}$ for $v = (v_1, ..., v_r)\in \Z_{\geq 0}^r$. 
Then we have $\varphi(f)=\sum_{v \in \Z_{\geq 0}^r} n_vx^v$, where $x^v := x_1^{v_1} \cdots x_r^{v_r}$. 
Since this sum $\varphi(f)=\sum_{v \in \Z_{\geq 0}^r} n_vx^v$ is a finite sum, 
we can find an integer $e>0$ such that 
every nonzero coefficient $n_v = n_{(v_1, ..., v_r)}$ satisfies 
$0 \leq v \leq p^e-1$, i.e., $0 \leq v_1 \leq p^{e}-1, \ldots, 0 \leq v_r \leq p^e-1$. 
Note that $R = \bigoplus_{0 \leq v \leq p^e-1} R^{p^e} x^v$, i.e., 
$\{x^v \,|\, 0 \leq v \leq p^e-1\}$ is a free $R^{p^e}$-linear basis of the $R^{p^e}$-module $R$. t
By $\varphi(f) \in \m = Rx_1 + \ldots +Rx_s$, we get 
\[
\varphi(f) = \sum_{0 \leq v \leq p^e-1} m_vx^v \qquad 
\text{for}\qquad m_v \in R^{p^e},  
\]
where every nonzero coefficient $m_v$ satisfies (a) or (b) below. 
\begin{enumerate}
\item[(a)] 
$v_i >0$ for some $1 \leq i \leq s$.  
\item[(b)] $m_v \in Rx_1+ \cdots +Rx_s$. 
\end{enumerate}
By 
\[
\sum_{0 \leq v \in p^e-1} n_vx^v = 
\varphi(f) = \sum_{0 \leq v \leq p^e-1} m_vx^v,  
\]
we obtain $\F_p \ni n_v = m_v \in R^{p^e}$. 
If the case (b) occurs, then 
we would get 
$n_v = m_v \in \F_p \cap (Rx_1+ \cdots Rx_s)=0$, 
which contradicts $m_v \neq 0$. 
Therefore, $n_v (= m_v) \neq 0$ implies (a). 
Hence  we get $f \in SX_1 + \cdots + SX_s =\mathfrak p'$. 
This completes the proof of  Claim \ref{claim:compare ideal}. 
\end{proof}
By Claim \ref{claim:compare ideal}, $\varphi$ induces a local ring homomorphism $\varphi' \colon S_\mathfrak p \to R$.
As the localisation $S \to S_\mathfrak p$ is flat, 
it is enough to show that $\varphi'$ is flat.
Explicitly, we have $S_{\mathfrak p} = K[X_1, ..., X_s]_{(X_1, ..., X_s)}$ 
for $K := \F_p(X_{s+1}, ..., X_r)$. 
In particular, 
we get  $\dim{S_\mathfrak p}=s = \dim{R}$. 
Furthermore, we have $\dim{R/\mathfrak p R}=\dim{R/\m}=0$. 
Then 
it follows from 
\cite[Theorem 23.1]{Matsumura} or 
\cite{stacks-project}*{\href{https://stacks.math.columbia.edu/tag/00R4}{Tag 00R4}} 
that  $\varphi'$ is flat. 
\end{proof}

\begin{lemma}\label{lemma:compare W_nOmega}
Let $\varphi \colon S \to R$ be a flat ring homomorphism of regular $F$-finite Noetherian $\F_p$-algebras. 
Assume that the induced $R$-module homomorphism 
\[
\theta^{e, i} : F^e_*\Omega_S^i \otimes_{S} R \to F^e_*\Omega_R^i
\]
is an isomorphism for every integer $e \geq 0$ and every integer $i \geq 0$. 
Then 
the following hold. 
\begin{enumerate}
    \item[(i)] We have the induced $W_nR$-module isomomorphism 
    \[
    (F_*^eW_{m}S) \otimes_{W_{n}S} W_{n}R \xrightarrow{\simeq} 
    F_*^eW_{m}R
    \]
    for every triple $(e, n, m) \in \Z^3$ satisfying $e \geq 0$ and $n \geq m \geq 1$. 
    \item[(ii)] The induced ring homomorphism $W_nS \to W_nR$ is flat for every integer $n>0$. 
    \item[(iii)] We have the induced $W_nR$-module isomorphism 
    \[
    (F_*^eW_n\Omega_S^i) \otimes_{W_nS} W_nR \xrightarrow{\simeq} F_*^eW_n\Omega_R^i
    \]
    for every triple $(e, n, i) \in \Z^3$ satisfying $e \geq 0, n >0$, and $i\geq 0$. 
\end{enumerate}
\end{lemma}

\begin{proof}
Let us show (i), i.e., the induced homomorphism 
\[
\pi^e_{n, m} : F_*^eW_{m}S \otimes_{W_{n}S} W_{n}R  \to F_*^eW_{m}R
\]
is an isomorprhism. 
For $e \geq 0$ and $n \geq m \geq 1$,  we have the following commutative diagram 
in which each horizontal sequence is exact: 
\begin{equation*} 
    \begin{tikzcd}
     & F_*^m(F^{e}_*S \otimes_{W_{n}S} W_{n}R) \arrow{r}{{\cccred \mathbf{V}}^{m} \otimes {\rm id}} \arrow{d}{F_*^{m}\pi_{n, 1}^e}& 
     F_*^eW_{m+1}S \otimes_{W_{n}S} W_{n}R 
     \arrow{r}{{{\cccred \mathbf{R}}} \otimes {\rm id}} \arrow{d}{\pi_{n, m+1}^e} & F_*^eW_{m}S \otimes_{W_{n}S} W_{n}R \arrow{r} \arrow{d}{\pi_{n, m}^e} & 0 \\
    0 \arrow{r} & F_*^{m+e}R \arrow{r}{{\cccred \mathbf{V}}^{m}} & F_*^eW_{m+1}R \arrow{r}{{\cccred \mathbf{R}}} & F_*^eW_{m}R \arrow{r} & 0. 
    \end{tikzcd}
\end{equation*}
By the snake lemma and induction on $m$, it is enough to show that $\pi^e_{n, 1}$ 
is an isomorphism. 
Note that $\pi^e_{n, 1}$ 
can be written by 
\[
\pi^e_{n, 1}: 
F^{e}_*S \otimes_{W_{n}S} W_{n}R 
\overset{\rho}{\twoheadrightarrow} 
F^{e}_*S \otimes_{S} R\xrightarrow{\theta^{e, 0}, \simeq } F_*^eR,
\]
where the left arrow $\rho$ is the induced surjection given by  
\[
\rho \colon F^{e}_*S \otimes_{W_{n}S} W_{n}R 
\twoheadrightarrow F^{e}_*S \otimes_{W_{n}S} R 
\twoheadrightarrow 
F^{e}_*S \otimes_{S} R. 
\]
In particular, $\pi^e_{n, 1}$ is surjective. 
Taking a localisation of $S$ at a prime ideal, 
we may assume that $F_*^eS \simeq S^{\oplus \ell_e}$ for some integer $\ell_e >0$. 
By the construction of $\rho$, 
we may replace $F_*^eS$ by a direct summand $S$. 
Hence it suffices to show that $\pi^0_{n, 1} : S \otimes_{W_nS} W_nR \to S\otimes_S R$ is an isomorphism. 

Fix an integer $n \geq 1$. 
It is enough to prove that $\pi^0_{n, m}$ is an isomorphism 
for $n  \geq m \geq 1$ by descending induction on $m$. 
It is clear that  $\pi^0_{n, n}$ is an isomorphism. 
Fix an integer $1 \leq m \leq n-1$ and assume that $\pi^0_{n, m+1}$ is an isomorphism. 
Since $\pi^0_{n, 1}$ is surjective and $\pi^0_{n, m+1}$ is an isomorphism, 
the snake lemma implies that $\pi^0_{n, m}$ is an isomorphism. 
Thus (i) holds. 

Let us show (ii).
By assumption, we may assume $n \geq 2$. 
We define ideals $W_n^{\geq 1}S$ and $W_n^{\geq 1}R$ by 
\[
W_n^{\geq 1}S:={\cccred \mathbf{V}}(F_*W_{n-1}S)\qquad \text{and}\qquad 
W_n^{\geq 1}R:={\cccred \mathbf{V}}(F_*W_{n-1}R). 
\]
Then we have $W_nS/W_n^{\geq 1}S \simeq S$ and $W_nR/W_n^{\geq 1}R \simeq R$.
Since we have
\[
W_n^{\geq 1}S \otimes_{W_nS} W_nR \simeq F_*W_{n-1}S \otimes_{W_nS} W_nR \simeq F_*W_{n-1}R \simeq W_n^{\geq 1}R
\]
by (i), we obtain
\[
W_n^{\geq 1}S \otimes_{W_nS} W_nR \simeq W_nS^{\geq 1}S \cdot W_nR = W_n^{\geq 1}(R).
\]
Hence 
the homomorphism
\[
S \simeq W_nS/W_n^{\geq 1}S \to W_nR \otimes_{W_nS} W_nS/W_n^{\geq 1}S \simeq W_nR/W_n^{\geq 1}(R) \simeq R
\]
is flat. 
By  the local criterion for flatness 
(\cite[Theorem 22.3]{Matsumura} or \cite[\href{https://stacks.math.columbia.edu/tag/051C}{Tag 051C}]{stacks-project}), $W_nS \to W_nR$ is flat. 
Thus (ii) holds.

Let us show (iii). 
Given a smooth $\F_p$-algebra $T$, we have the following exact sequence of $W_nT$-modules 
\cite[Ch.\ I, Proposition 3.11]{illusie_de_rham_witt}: 
\[
0\to 
\frac{F_*^{n}\Omega_T^{i-1}}{F_*Z_{n-1}\Omega_T^{i-1}}
\xrightarrow{d{\cccred \mathbf{V}}^{n-1}}
F_*W_{n}\Omega_T^i \xrightarrow{{\cccred \mathbf{V}}} W_{n+1}\Omega_T^i \xrightarrow{F^n} Z_n\Omega_T^i \to 0. 
\]
Note that this exact sequence is  $W_{n+1}T$-linear 
(Lemma \ref{l-dVn-1-linear}). 
Since each of $R$ and $S$ can be written as a direct limit of smooth $\F_p$-algebras, 
we get the corresponding exact sequence for $S$ and $R$ 
(cf.\ \cite[Section 2 and Theorem 2.1]{Shi07}). 
Therefore, we get the following commutative diagram in which each horizontal sequence is exact: 
{\small 
\[
\begin{tikzcd}
0 \arrow[r] & 
f^*\left( F_*^e\frac{F_*^{n}\Omega_S^{i-1}}{F_*Z_{n-1}\Omega_S^{i-1}}\right) \arrow[r, "f^*(d{\cccred \mathbf{V}}^{n-1})"] \arrow[d, "\beta", "\simeq"'] &
f^*(F^{e+1}_*W_{n}\Omega_S^i) \arrow[r, "f^*{\cccred \mathbf{V}}"] \arrow[d, "\alpha_{n}^{e+1, i}", "\simeq"'] &  
f^*(F_*^eW_{n+1}\Omega_S^i) \arrow[r, "f^*(F^n)"]  \arrow[d, "\alpha_{n+1}^{e, i}"]&
f^*(F_*^eZ_n\Omega_S^i)  \arrow[r] \arrow[d, "\gamma", "\simeq"']& 
0 \\ 
0  \arrow[r] & 
F_*^e\frac{F_*^{n}\Omega_R^{i-1}}{F_*Z_{n-1}\Omega_R^{i-1}} \arrow[r, "d{\cccred \mathbf{V}}^{n-1}"] &
F^{e+1}_*W_{n}\Omega_R^i \arrow[r, "{\cccred \mathbf{V}}"] &  
F_*^eW_{n+1}\Omega_R^i \arrow[r, "F^n"] &
F_*^eZ_n\Omega_R^i \arrow[r] &
0, 
\end{tikzcd}
\]
}
where $f^*(-) := (-) \otimes_{W_{n+1}S} W_{n+1}R$. 
Let us prove that $\alpha_n^{e, i}$ is an isomorphism for all $i\geq 0$ and $e \geq 0$ by induction on $n \geq 1$. 
The base case $n=1$ of this induction follows from our assumption. 
In particular, $\beta$ and $\gamma$ are isomorphisms 
(note that the Cartier isomorphism holds for regular $\F_p$-algebras \cite[Proposition 2.5]{Shi07}, and hence $Z_n\Omega_R^i$ and $B_n\Omega_R^i$ can be defined as in the smooth case). 
Fix $n \geq 1$ and assume that $\alpha_n^{e, i}$ is an isomorphism for all $i\geq 0$ and $e \geq 0$ by induction on $n \geq 1$. 
Take integers $i \geq 0$ and $e\geq 0$. 
Then $\alpha_n^{e+1, i}$ is an isomorphism by the induction hypothesis. 
By the 5-lemma, also  $\alpha_{n+1}^{e, i}$ is an isomorphism. 
\end{proof}

\begin{lem}\label{l-dVn-1-linear}
Let $R$ be an $\F_p$-algebra. 
Take integers $i \geq 0$ and $n \geq 1$. 
Then the map
\[
d{\cccred \mathbf{V}}^{n-1}: F^n_*\Omega_R^{i} \to F_*W_n\Omega_R^{i+1}
\]
is a $W_nR$-module homomorphism. 
\end{lem}

\begin{proof}
Clearly, $d{\cccred \mathbf{V}}^{n-1}$ is an additive homomorphism. 
Fix $a \in R$ and $F_*^n \omega \in F_*^n\Omega_R^i$, 
where $\omega \in \Omega_R^i$. 
It suffices to show 
$d{\cccred \mathbf{V}}^{n-1}(a \cdot F^n_*\omega) = [a] \cdot d{\cccred \mathbf{V}}^{n-1}(F^n_*\omega)$. 
The problem is reduced to the case when $\omega = f dg_1 \cdots dg_i$ 
for some $f, g_1, ..., g_i \in R$, 
because $d{\cccred \mathbf{V}}^{n-1}$ is additive and $\omega$ can be written as a finite sum of elements of this form. Set $dg := dg_1 \cdots dg_i$. 
We have 
\begin{eqnarray*}
d{\cccred \mathbf{V}}^{n-1}(a \cdot F_*^e\omega) 
&=&d{\cccred \mathbf{V}}^{n-1}(F_*(a^p\cdot F_*^{n-1}\omega))\\
&=& d(F_*[a^p] \cdot {\cccred \mathbf{V}}^{n-1}(F_*^{n-1}\omega))\\
&=& d(F_*[a^p]) \cdot {\cccred \mathbf{V}}^{n-1}(F_*^{n-1}\omega)) + 
(F_*[a^p]) \cdot d({\cccred \mathbf{V}}^{n-1}(F_*^{n-1}\omega))\\
&=& d(F_*[a^p]) \cdot {\cccred \mathbf{V}}^{n-1}(F_*^{n-1}\omega)) + 
[a] \cdot d({\cccred \mathbf{V}}^{n-1}(F_*^{n}\omega)).
\end{eqnarray*} 
In what follows, we drop $F_*$ and $F_*^{n-1}$. 
Then it is enough to show 
$d([a^p]) \cdot {\cccred \mathbf{V}}^{n-1}(\omega)=0$. 
By 
\[
{\cccred \mathbf{V}}^{n-1}(\omega) = {\cccred \mathbf{V}}^{n-1}(fdg_1\cdots dg_i) = 
({\cccred \mathbf{V}}^{n-1}f)\cdot d({\cccred \mathbf{V}}^{n-1}g_1) \cdots  d({\cccred \mathbf{V}}^{n-1}g_i), 
\]
it holds that  
\[
d([a^p]) \cdot {\cccred \mathbf{V}}^{n-1}(\omega)=(p[a]^{p-1}d[a]) \cdot 
 ({\cccred \mathbf{V}}^{n-1}f) \cdot d({\cccred \mathbf{V}}^{n-1}g_1) \cdots d({\cccred \mathbf{V}}^{n-1}g_i) = 0, 
\]
because $p{\cccred \mathbf{V}}^{n-1}f = F{\cccred \mathbf{V}}^nf=0$ in $W_nR$. 
\end{proof}

\begin{remark}\label{r-FVR-bc}
We use the same notation as in Lemma \ref{lemma:compare W_nOmega}. 
Then the three operations 
\[
F: W_n\Omega_R^i \to F_*W_{n-1}\Omega_R^i, \quad 
{\cccred \mathbf{V}}: F_*W_{n-1}\Omega_R^i \to W_{n}\Omega_R^i, \quad 
{{\cccred \mathbf{R}}}: W_n\Omega_R^i \to W_{n-1}\Omega_R^i, \quad 
\]
are obtained from the corresponding operations for $S$ by taking the base change $(-) \otimes_{W_nS} W_nR$. 
More explicitly, we have the following natural $W_nR$-linear commutative diagrams 
in which all the vertical arrows are the canonical isomorphisms (Lemma \ref{lemma:compare W_nOmega}(iii)): 
    \[
\begin{tikzcd}[column sep=5em]
W_{n}\Omega^i_S \otimes_{W_nS} W_nR  \arrow[r, "F \otimes_{W_nS} W_nR"] \arrow[d, "\simeq"] &F_*W_{n-1}\Omega^i_S\otimes_{W_nS} W_nR\arrow[d, "\simeq"]\\
W_{n}\Omega^i_R \arrow[r, "F"] & F_*W_{n-1}\Omega^i_R
\end{tikzcd}
    \]
    \[
\begin{tikzcd}[column sep=5em]
(F_*W_{n-1}\Omega^i_S) \otimes_{W_nS} W_nR  \arrow[r, "{\cccred \mathbf{V}} \otimes_{W_nS} W_nR"] \arrow[d, "\simeq"] & W_{n}\Omega^i_S \otimes_{W_nS} W_nR\arrow[d, "\simeq"]\\
F_*W_{n-1}\Omega^i_R \arrow[r, "{\cccred \mathbf{V}}"] & W_{n}\Omega^i_R
\end{tikzcd}
    \]   
    \[
\begin{tikzcd}[column sep=5em]
W_{n}\Omega^i_S \otimes_{W_nS} W_nR  \arrow[r, "{\cccred \mathbf{R}} \otimes_{W_nS} W_nR"] \arrow[d, "\simeq"] &W_{n-1}\Omega^i_S\otimes_{W_nS} W_nR\arrow[d, "\simeq"]\\
W_{n}\Omega^i_R \arrow[r, "{\cccred \mathbf{R}}"] & W_{n-1}\Omega^i_R.
\end{tikzcd}
    \]    
This fact can be used to generalise  many exact sequences known for smooth $\F_p$-algebras 
to regular $\F_p$-algebras. 
For example, we have the following factorisation of the $p^m$-multiplication map \cite[Ch.\ I, Proposition 3.4]{illusie_de_rham_witt}: 
\[
p^m : W_n\Omega_S^i 
\overset{{\cccred \mathbf{R}}^m}{\twoheadrightarrow} W_{n-m}\Omega_S^i 
\overset{\underline{p}^m}{\hookrightarrow} W_n\Omega_S^i, 
\]
which induces the $R$-version by taking the base change $(-) \otimes_{W_nS} W_nR$ 
(Lemma \ref{lem:covered be smooth}):  
\[
p^m : W_n\Omega_R^i 
\overset{{\cccred \mathbf{R}}^m}{\twoheadrightarrow} W_{n-m}\Omega_R^i 
\overset{\underline{p}^m}{\hookrightarrow} W_n\Omega_R^i. 
\]
\end{remark}

\begin{remark}\label{r R p-bar F ex}
Let $R$ be a regular $F$-finite integral domain.
For every $\fp \in \Spec{R}$, there exists a flat ring homomorphism $S_\fp \to R_\fp$ as in Lemma \ref{lem:covered be smooth}.
By Remark \ref{r-FVR-bc}, 
we get the decomposition 
\[
p^m \colon 
W_n\Omega_{R_\fp}^i 
\overset{{\cccred \mathbf{R}}^m}{\twoheadrightarrow} W_{n-m}\Omega_{R_\fp}^i 
\overset{\underline{p}^m_\fp}{\hookrightarrow} W_n\Omega_{R_i}^i. 
\]
For the generic point $\eta \in \Spec R$, we have the following commutative diagram, because ${\cccred \mathbf{R}}^m$ is surjective:  
\[
\begin{tikzcd}
W_{n-m}\omega_{R_{\eta}} \arrow[r, "\underline{p}^{m}_{\eta}"] & W_n\omega_{R_{\eta}}\\
W_{n-m}\omega_{R_{\fp}} \arrow[r, "\underline{p}^m_{\fp}"]\arrow[u, hook] & W_n\omega_{R_{\fp}}, \arrow[u, hook]
\end{tikzcd}
\]
where the vertical arrows are the induced injections. 
Then it holds that 
\begin{align*}
    \underline{p}^m_{\eta}(W_{n-m}\omega_{R}) 
&\overset{(\alpha)}{=} 
\underline{p}^m_{\eta}(\bigcap_{\fp \in \Spec R} W_{n-m}\omega_{R_{\fp}})
\subseteq 
\bigcap_{\fp \in \Spec R} 
\underline{p}^m_{\eta}(W_{n-m}\omega_{R_{\fp}}) \\
&= \bigcap_{\fp \in \Spec R} 
\underline{p}^m_{\fp}(W_{n-m}\omega_{R_{\fp}})\subseteq 
\bigcap_{\fp \in \Spec R} W_n\omega_{R_{\fp}} 
\overset{(\beta)}{=}  W_n\omega_R,
\end{align*}
where the equalities $(\alpha)$ and $(\beta)$ can be checked by using the fact that the restriction map $W_n\omega_R(U) \to W_n\omega_R(V)$ is injective for open subsets $U$ and $V$ satisfying $\Spec R \supseteq U \supseteq V \neq \emptyset$. 
Then we get an $R$-module homomorphism  $\underline{p}^m \colon W_{n-m}\omega_{R} \to W_n\omega_{R}$ 
satisfying $\underline{p}^m \circ {\cccred \mathbf{R}}^m = p^m$. 
Moreover, for  $K:=K(R)$ and the non-negative integer $r$ satisfying $[K:K^p] = p^r$, 
the sequence 
\[
0 \to W_m\Omega^r_R \xrightarrow{\underline{p}^n} W_{n+m}\Omega_R^r 
\xrightarrow{F^m}  F_*^m W_n\Omega^r_R \to  0
\]
is exact. 
\end{remark}

\begin{proposition}\label{prop:regular case dualizing complex1}
Let $R$ be a regular $F$-finite integral domain and 
set $K:=K(R)$. 
Let $r$ be the non-negative integer satisfying $[K:K^p] = p^r$. 
Fix an integer $n>0$. 
Let $W_n\Omega_R^{\mydot}$ be the  de Rham-Witt complex of length $n$. 
Set $W_n\omega_R:=W_n\Omega_R^{r}$.
Then 
\begin{enumerate}
    \item[(i)]$W_n\omega_R$ is a dualising $W_nR$-module and 
    \item[(ii)]$W_n\Omega_R^{i}=0$ for every integer $i \geq r+1$.
\end{enumerate}
\end{proposition}

\begin{proof}
We may assume   that $R$ is local \cite[Ch.\ V, Corollary 2.3]{Har66}. 
Take a flat ring homomorphism $\varphi \colon S \to R$ 
from a smooth $\F_p$-algebra $S$ as in Lemma \ref{lem:covered be smooth}. 
We note that
the following properties, corresponding to (i) and (ii), hold on $S$ by \cite{Ekedahl1}*{Theorem 4.1, Lemma 0.6}: 
\begin{enumerate}
    \item[(i)'] $W_n\omega_S$ is a dualising $W_nS$-module. 
    \item[(ii)'] $W_n\Omega_S^{i}=0$ for every integer $i \geq r+1$.
\end{enumerate}
By Lemma \ref{lemma:compare W_nOmega}(iii), we have $W_n\Omega^i_S \otimes_{W_nS} W_nR \xrightarrow{\simeq} W_n\Omega_R^i$ for all integers $n \geq 1$ and $i \geq 0$. 
Therefore, the assertion (ii) follows from (ii)'. 

It suffices to show (i).
Since $W_nS$ is Cohen-Macaulay and $W_n\omega_S$ is a dualising $W_nS$-module, 
the following hold \cite[\href{https://stacks.math.columbia.edu/tag/0A7B}{Tag 0A7B}]{stacks-project}. 
\begin{enumerate}
\item $W_n\omega_S$ is a finitely generated $W_nS$-module. 
\item $W_n\omega_S$ is of finite injective dimension. 
\item The natural homomorphism 
\begin{equation}\label{eq:duality}
W_nS \to R\Hom_{W_nS}(W_n\omega_S,W_n\omega_S)
\end{equation}
is a quasi-isomorphism.
\end{enumerate}
It is enough  
to show that the corresponding properties 
hold for $W_n\omega_R$ \cite[\href{https://stacks.math.columbia.edu/tag/0A7B}{Tag 0A7B}]{stacks-project}. 
By $W_n\omega_R \simeq W_n\omega_S \otimes_{W_nS} W_nR$, $W_n\omega_R$ is a finitely generated $W_nR$-module, i.e., (1) holds for $R$. 
Since $W_nR$ is flat over $W_nS$ \cite[Exercise 7.7]{Matsumura}, 
we get a quasi-isomorphism
\[
W_nR \xrightarrow{\simeq} R\Hom_{W_nS}(W_n\omega_S,W_n\omega_S) \otimes_{W_nS} W_nR \simeq R\Hom_{W_nR}(W_n\omega_R, W_n\omega_R) 
\]
by applying $(-) \otimes_{W_nS} W_nR$ to the quasi-isomorphism (\ref{eq:duality}). 
Thus (3) holds for $R$. 
Set $\mathfrak p := \varphi^{-1}(\m)$ and we get $\mathfrak p R = \m$ by Lemma \ref{lem:covered be smooth}. 
Furthermore, we have
\[
\Ext^i_{W_nS}(S/\mathfrak p,W_n\omega_S) \otimes_{W_nR} W_nS \simeq \Ext^i_{W_nR}(R/\m,W_n\omega_R). 
\]
Since $W_n\omega_S$ has a finite injective dimension, 
there exists $i_0>0$ such that 
\[
\Ext^i_{W_nS}(S/\mathfrak p,W_n\omega_S)=0
\]
for every $i \geq i_0$. 
Therefore, $\Ext^i_{W_nR}(R/\m,W_n\omega_R)=0$ for every $i \geq i_0$. 
By \cite[\href{https://stacks.math.columbia.edu/tag/0AVJ}{Tag 0AVJ}]{stacks-project}, $W_n\omega_R$ is of finite injective dimension, i.e., (2) holds for $R$.
\end{proof}

\begin{proposition}\label{prop:regular case dualizing complex2}
Let $R$ be a regular $F$-finite integral domain and 
set $K:=K(R)$. 
Let $r$ be the non-negative integer satisfying $[K:K^p] = p^r$. 
Let $W_n\Omega_X^{\mydot}$ be the  de Rham-Witt complex of length $n$. 
Set $W_n\omega_R:=W_n\Omega_R^{r}$.
Then there exists a unique set of homomorphisms (called trace maps) 
\[
\{ T^e_{m, n} : F^e_*W_m\omega_R \to W_n\omega_R \,|\, e \geq 0, n \geq m \geq 1 \}  
\]
which satisfies the following. 
\begin{enumerate}
\renewcommand{\labelenumi}{(\roman{enumi})}
    \item[(i)] 
    For all integers $n  \geq m \geq 1$, the $W_nR$-module homomorphism 
    \begin{align*}
\sigma^e_{m, n} :  F_*^eW_m\omega_R &\to \Hom_{W_nR}(F_*^eW_mR, W_n\omega_R)\\
        F_*^e\zeta &\mapsto ( F_*^e\alpha \mapsto T^e_{m, n}(F_*^e(\alpha \zeta))),
    \end{align*}
    induced by the trace map $T^e_{m, n} : F^e_*W_m\omega_R \to W_n\omega_R$, is an isomorphism. 
    Moreover, $\sigma^e_{m, n}$ is $F_*^eW_nR$-linear. 
 \item[(ii)] 
$T^{e'}_{l,m} \circ F^{e'}_*T^e_{m,n} =T^{e+e'}_{l,n}$ 
for all $e,e',l,m,n \in \Z$ satisfying $e,e' \geq 0$ and $n \geq m \geq l \geq 1$. 
\item[(iii)] 
$T^0_{m, n} = \underline{p}^{n-m}$ and 
$\underline{p}^e  \circ T^e_{n, n} = {\cccred \mathbf{V}}^e$ 
for all $e >0$ and $n  \geq m \geq 1$. 
\begin{equation}\label{e1 regular case dualizing complex2}
\begin{tikzcd}
F_*^eW_n\omega_R \arrow[r, "{T^e_{n, n}}"] 
\arrow[rr, bend left, "{\cccred \mathbf{V}}^e"] &W_n\omega_R \arrow[r, hook, "\underline{p}^e"] &W_{n+e}\omega_R
\end{tikzcd}
\end{equation}
\end{enumerate}
\end{proposition}

\begin{proof}
We can directly check the $F_*^eW_nR$-linearity in (i).

\setcounter{step}{0}

\begin{step}\label{s1 regular case dualizing complex2} 
The assertion of 
Proposition \ref{prop:regular case dualizing complex2} holds for the case when $R$ is a finitely generated $\F_p$-algebra. 
\end{step}

\begin{proof}[Proof of Step \ref{s1 regular case dualizing complex2}]
In this case, 
we have $W_n\omega_R = \pi^!\F_p[-r]$ for $\pi: \Spec R \to \Spec \F_p$. 
We then define $T^e_{m, n} \colon F_*^eW_m\omega_R \to W_n\omega_R$ as the usual trace map that induces the Grothendieck duality in (i). 
Then (i) and (ii) hold. 

It is enough to show  (iii). 
Recall that  the Grothendieck dual 
$(-)^* := \Hom_{W_nR}(-, W_n\omega_R)$
of $F, {\cccred \mathbf{V}}, {\cccred \mathbf{R}}, \underline{p}$ 
are ${\cccred \mathbf{V}}, F, \underline{p}, {\cccred \mathbf{R}}$, respectively \cite[Lemma 2.2.4]{Ekedahl1}. 
For example, 
we have the following commutative diagram for $R$: 
\[
\begin{tikzcd}
W_{n-1}\omega_R \arrow[d, "{\sigma^0_{n-1, n}}", "\simeq"'] 
\arrow[r, "\underline p"] &W_n\omega_R \arrow[d, "\text{canonical isom.}", "\simeq"']\\
\Hom_{W_nR}(W_{n-1}R, W_nR) \arrow[r, "{{\cccred \mathbf{R}}^*}"]& 
\Hom_{W_nR}(W_nR, W_n\omega_R). 
\end{tikzcd}
\]
Hence it holds that $T^0_{m, n} = \underline{p}^{n-m}$. 
The other equality $\underline{p}^e  \circ T^e_{n, n} = {\cccred \mathbf{V}}^e$ in (iii)  follows from the fact that 
the diagram (\ref{e1 regular case dualizing complex2}) is 
obtained by applying 
$(-)^* = \Hom_{W_nR}(-, W_n\omega_R)$ to the following commutative diagram: 
\begin{equation*}
\begin{tikzcd}
F_*^eW_nR \arrow[r, leftarrow, "F^e"] 
\arrow[rr, bend left, leftarrow, "F^e"] &W_nR \arrow[r, hook, leftarrow, "{\cccred \mathbf{R}}^e"] &W_{n+e}R. 
\end{tikzcd}
\end{equation*}
This completes the proof of Step \ref{s1 regular case dualizing complex2}. 
\end{proof}

In what follows, we shall repeatedly reduce the problem to the case when $R$ is a finitely generated $\F_p$-algebra.  
If we may assume that $R$ is a local ring, 
then the problem is reduced to the case when ($\star$) holds 
(Lemma \ref{lem:covered be smooth}, Lemma \ref{lemma:compare W_nOmega}). 
\begin{enumerate}
\item[($\star$)] 
There exists a flat ring homomorphism $\varphi : S \to R$ 
 of regular $F$-finite Noetherian domains 
such that $S$ is a finitely generated $\F_p$-algebra, 
$[K(S):K(S)^p]=p^r$, and 
we have the induced isomorphism 
$(F^e_*W_n\Omega_S^i) \otimes_{W_nS} W_nR \xrightarrow{\simeq} F^e_*W_n\Omega_R^i$ 
for all $e \geq 0, n>0$, and $i \geq 0$. 
\end{enumerate}


\begin{step}\label{s2 regular case dualizing complex2} 
If each of $\{T^e_{m, n}\,|\, e \geq 0, n \geq m \geq 1 \}$ and $\{T'^e_{m, n}\,|\, e \geq 0, n \geq m \geq 1 \}$ satisfies {\rm (ii)} and {\rm (iii)}, then $T^e_{m, n} =T'^e_{m, n}$ for all $e\geq 0$ and 
$n \geq m \geq 1$. 
\end{step}

\begin{proof}[Proof of Step \ref{s2 regular case dualizing complex2}]
By (iii), we get $T^0_{m, n} = \underline{p}^{n-m} = T'^0_{m, n}$ and $\underline{p}^e  \circ T^e_{n, n} = {\cccred \mathbf{V}}^e = \underline{p}^e  \circ T'^e_{n, n}$. 
Since $\underline{p}^e$ is injective, we get 
$T^e_{n, n} = T'^e_{n, n}$. 
The general case is reduced to these cases by (ii), because 
\[
T^e_{m,n} \overset{{\rm (ii)}}{=} T^{0}_{m,n} \circ T^e_{n, n} 
= T'^{0}_{m,n} \circ T'^e_{n, n} 
\overset{{\rm (ii)}}{=}
T'^e_{m,n}. 
\]
This completes the proof of Step \ref{s2 regular case dualizing complex2}. 
\end{proof}



\begin{step}\label{s3 regular case dualizing complex2} 
Fix integers $e>0$ and $n>0$. 
Then there exists a $W_nR$-module homomorphism 
$T^e_{n, n} \colon F^e_*W_n\omega_{R} \to W_n\omega_{R}$ such that $\underline{p}^e  \circ T^e_{n, n} = {\cccred \mathbf{V}}^e$. 
\end{step}

\begin{proof}[Proof of Step \ref{s3 regular case dualizing complex2}]
By a similar argument to the one of Remark \ref{r R p-bar F ex}, 
we may assume that $R$ is a local ring. 
In this case, 
there exists a flat ring homomorphism 
$S \to R$ satisfying the properties in $(\star)$. 
By Step \ref{s1 regular case dualizing complex2}, there is an $W_nS$-module homomorphism $T^e_{S, n, n}$ satisfying $\underline{p}^e \circ T^e_{S, n, n} = {\cccred \mathbf{V}}^e$. 
Set $T^e_{n, n} := T^e_{S, n, n} \otimes_{W_nS} W_nR$. 
We then get  $\underline{p}^e \circ T^e_{n, n} = {\cccred \mathbf{V}}^e$, as required. 
This completes the proof of  Step \ref{s3 regular case dualizing complex2}.

\end{proof}

\begin{step}\label{s4 regular case dualizing complex2} 
There exists a set of homomorphisms 
\[
\{ T^e_{m, n} : F^e_*W_m\omega_R \to W_n\omega_R \,|\, e \geq 0, n \geq m \geq 1 \}  
\]
which satisfies {\rm (i)-(iii)}. 
\end{step}

\begin{proof}[Proof of Step \ref{s4 regular case dualizing complex2}]
For integers $e>0$ and $n  \geq m \geq 1$, 
we set $T^0_{m, n} := \underline{p}^{n-m}$ and  
\[
T^e_{m,n} :=T^{0}_{m,n} \circ T^e_{n, n}, 
\]
where 
$T^e_{n, n}$  is defined in Step \ref{s3 regular case dualizing complex2}. 
Then (iii) holds. 

Let us show (ii). 
The required equality can be checked after  replacing $R$ by $R_{\fp}$ for a prime ideal $\fp$ of $R$,  and hence we may assume ($\star$). 
Then  the problem is reduced to the case when $R$ is a finitely generated $\F_p$-algebra. 
Thus (ii) holds by Step \ref{s1 regular case dualizing complex2}. 

Let us show (i). 
Again, we may assume that ($\star$) holds. 
By Step \ref{s2 regular case dualizing complex2}, (ii) and (iii) assure the uniqueness of $T^e_{m, n}$. 
Therefore, we obtain $T^e_{R, m, n} = T^e_{S, m, n} \otimes_{W_nS} W_nR$, where $T^e_{R, m, n}$ and $T^e_{S, m, n}$ denote the trace maps on $R$ and $S$, respectively. 
Hence (i) holds by Step \ref{s1 regular case dualizing complex2}. 
This completes the proof of  Step \ref{s4 regular case dualizing complex2}. 
\end{proof}
Step \ref{s2 regular case dualizing complex2} and 
Step \ref{s4 regular case dualizing complex2} complete the proof of Proposition \ref{prop:regular case dualizing complex2}.
\end{proof}

\begin{proof}[Proof of Theorem \ref{thm:X-exist dualizing complex}] 
The proof consists of three steps.

\setcounter{step}{0}

\begin{step}\label{s1:X-exist dualizing complex}
Theorem \ref{thm:X-exist dualizing complex}(1) and 
Theorem \ref{thm:X-exist dualizing complex}(2) imply  
Theorem \ref{thm:X-exist dualizing complex}(3). 
\end{step}

\begin{proof}[Proof of Step \ref{s1:X-exist dualizing complex}]
Take integers $e,e',l,m,n$ satisfying $e,e' \geq 0$ and $n \geq m \geq l \geq 1$.
Set $f:=\iota^e_{X,m,n}$, $g:=\iota^{e'}_{X,l,m}$, $h:=\iota^{e+e'}_{X,l,n}$, $\rho_f:=\rho^e_{X,m,n}$, $\rho_g:=\rho^{e'}_{X,l,m}$, $\rho_h:=\rho^{e+e'}_{X,l,n}$, $T_f:=T^e_{X,m,n}$, $T_g:=T^{e'}_{X,l,m}$ and $T_h:=T^{e+e'}_{X,l,n}$.
Then 
Theorem \ref{thm:X-exist dualizing complex}(1) and 
Theorem \ref{thm:X-exist dualizing complex}(2) imply  
Theorem \ref{thm:X-exist dualizing complex}(3) 
by the diagram
\begin{equation}\label{eq:diag}
    \begin{tikzcd}
    h_*W_l\omega_X^{\mydot} \arrow[r,"f_*T_g"] \arrow[d,"h_*\rho_g"'] \arrow[dd,"h_*\rho_h"', bend right=90] & f_*W_m\omega_X^{\mydot} \arrow[r,"T_f"] \arrow[d,"f_*\rho_f"'] & W_n\omega_X^{\mydot} \\
    h_*g^{!}W_m\omega_X^{\mydot} \arrow[ru,"{\rm adj}"'] \arrow[d,"h_*g^{!}\rho_f"'] & f_*f^{!}W_n\omega_X^{\mydot} \arrow[ru,"{\rm adj}"'] & \\
    h_*h^{!}W_n\omega_X^{\mydot} \arrow[ru,"{\rm adj}"'] \arrow[rruu,bend right,"{\rm adj}"'].
\end{tikzcd}
\end{equation}
This completes the proof of Step \ref{s1:X-exist dualizing complex}. 
\end{proof}

\begin{step}\label{s2:X-exist dualizing complex}
Let $f \colon Y \to X$ be a morphism of separated schemes which are of finite type over an $F$-finite Noetherian $\F_p$-algebra. 
If the assertion of Theorem \ref{thm:X-exist dualizing complex} holds for $X$, then so does for $Y$.
\end{step}

\begin{proof}[Proof of Step \ref{s2:X-exist dualizing complex}]
Take $W_n\omega_X^{\mydot}$ and $\rho^e_{X,m,n}$ as in Theorem \ref{thm:X-exist dualizing complex} for integers $e,m,n$ satisfying $e \geq 0$ and $n \geq m \geq 1$.
For the induced morphism $f_n \colon W_nY \to W_nX$, 
we set $W_n\omega_Y^{\mydot}:=f_n^{!}W_n\omega_X^{\mydot}$ and $\rho^e_{Y,m,n}:=f_m^{!}\rho^e_{X,m,n}$. 
Then $\{W_n\omega_Y^{\mydot}\}_n$ and $\{\rho^e_{Y,m,n}\}_{e, m, n}$ satisfy 
Theorem \ref{thm:X-exist dualizing complex}(1) and 
Theorem \ref{thm:X-exist dualizing complex}(2)  (cf.\ \cite{Har66}). 
By Step \ref{s1:X-exist dualizing complex}, they also satisfy Theorem \ref{thm:X-exist dualizing complex}(3). 
This completes the proof of Step \ref{s2:X-exist dualizing complex}. 
\end{proof}

\begin{step}\label{s3:X-exist dualizing complex}
The assertion of Theorem \ref{thm:X-exist dualizing complex} holds for the case when $X$ is affine. 
\end{step}

\begin{proof}[Proof of Step \ref{s3:X-exist dualizing complex}]
We have $X=\Spec{R}$, where $R$ is an $F$-finite Noetherian $\F_p$-algebra.
By \cite{Gabber}*{Remark 13.6}, there exists a surjective ring homomorphism $S \to R$ such that $S$ is a regular  $F$-finite $\F_p$-algebra.
Let $S\simeq S_1 \times \cdots \times S_r$, where every $\Spec S_i$ is connected, so that every $S_i$ is a regular domain.
By Step \ref{s2:X-exist dualizing complex}, 
it is enough to show that 
Theorem \ref{thm:X-exist dualizing complex} holds for $\Spec{S}$.
We pick $W_n\omega_{S_i}$ and $T^e_{S_i,m,n}$ as in Proposition \ref{prop:regular case dualizing complex2}, and set
\[
W_n\omega_S:=W_n\omega_{S_1} \otimes \cdots \otimes W_n\omega_{S_r}
\]
and $T^e_{m,n}:=T^e_{S_1,m,n} \otimes \cdots \otimes T^e_{S_r,m,n}$.
Then $W_n\omega_S$ is a dualising module on $W_nS$ 
(Proposition \ref{prop:regular case dualizing complex1}) 
and $\{W_n\omega_S\}$ and $\{T^e_{m,n}\}$  satisfy the conditions (i)-(iii) in Proposition \ref{prop:regular case dualizing complex2}. 
Let $\iota^e_{m,n} \colon \Spec{W_mS} \to \Spec{W_nS}$ be the composition of $F^e$ and the closed immersion $\Spec{W_mS} \to \Spec{W_nS}$.
By 
Proposition \ref{prop:regular case dualizing complex2}(i), we obtain a $W_mS$-module isomorphism
\[
\rho^e_{m,n} \colon W_m\omega_S \xrightarrow{\simeq} (\iota^e_{m,n})^{!}W_nS
\]
such that $(\iota^e_{m,n})_*\rho^e_{m,n}=\sigma^e_{m,n}$ (see Remark \ref{rem:finite-tr-isom}(1)). 

Then Theorem \ref{thm:X-exist dualizing complex}(1) holds for $S$. 
By Step \ref{s1:X-exist dualizing complex}, it suffices to check 
that Theorem \ref{thm:X-exist dualizing complex}(2) holds for $S$. 
Let  $\Psi_{\iota^e_{m,n}}$ be as in Remark \ref{rem:finite-tr-isom}(1). 
By the construction of $\sigma^e_{m,n}$ 
and a commutative diagram similar to (\ref{eq:f-diag}), we have $\Psi_{\iota^e_{m,n}}(T^e_{m,n})=\rho^e_{m,n}$.
Furthermore, by  Proposition \ref{prop:regular case dualizing complex2}(ii) and Remark \ref{rem:finite-tr-isom}(2), we obtain  Theorem \ref{thm:X-exist dualizing complex}(2). 
This completes the proof of Step \ref{s3:X-exist dualizing complex}. 
\end{proof}
Step \ref{s2:X-exist dualizing complex} 
and 
Step \ref{s3:X-exist dualizing complex} complete the proof of 
Theorem \ref{thm:X-exist dualizing complex}. 
\end{proof}

\begin{remark}\label{rem:finite-tr-isom}
Let $f \colon Y \to X$ and $g \colon Z \to Y$ be {finite} morphisms 
of separated Noetherian schemes. 
Set $h:=g \circ f$.
Assume that $X$, $Y$, and $Z$ are Cohen-Macaulay and admit dualising complexes $\omega_X,$ $\omega_Y$, and $\omega_Z$, respectively. 
We identify $\omega_X,$ $\omega_Y$, and $\omega_Z$ with the corresponding coherent sheaves. 
\begin{enumerate}
    \item By \cite[\href{https://stacks.math.columbia.edu/tag/0AWZ}{Tag 0AWZ}]{stacks-project}, we have
    \[
    f_*f^{!}M=\cHom_{\cO_X}(f_*\cO_Y,M)
    \]
    as an $f_*\cO_Y$-module for every coherent $\cO_X$-module $M$.
    Furthermore, for a coherent $\cO_X$-module $M$ and a coherent $\cO_Y$-module $N$, we have
    \[
    \begin{tikzcd}[row sep=0.1cm, column sep=0.3cm]
        f_*\cHom_{\cO_Y}(N,f^{!}M) \arrow[r,equal] & f_*\cHom_{\cO_Y}(N,\cHom_{\cO_X}(f_*\cO_Y,M)) \arrow[r,"\simeq"]  & \cHom_{\cO_X}(f_*N,M) \\
         & (x \mapsto \varphi_x) \arrow[r,mapsto] & (x \mapsto \varphi_x(1)),
    \end{tikzcd}
    \]
    which implies that $f^{!}$ is the right adjoint functor of $f_*$ (cf.~\cite[\href{https://stacks.math.columbia.edu/tag/0AWZ}{Tag 0AWZ}]{stacks-project}). 
Applying this for $N = \omega_Y$ and $M =\omega_X$, we get 
   \[
  \Psi_f \colon \Hom_{\cO_Y}(\omega_Y,f^{!}\omega_X) \xrightarrow{\simeq} \Hom_{\cO_X}(f_*\omega_Y, \omega_X).   
    \]
Specifically, this adjoint bijection is given by 
    \[
    \begin{tikzcd}[row sep=0.1cm, column sep=0.3cm]
 {\rm adj} \colon f_*f^{!}\omega_X \arrow[r,phantom,"="] & f_*\cHom_{\cO_X}(f_*\cO_Y,\omega_X) \arrow[r] & \omega_X \\
        & \varphi \arrow[r,mapsto] & \varphi(1).
    \end{tikzcd}
    \]
In particular, for an $\MO_Y$-module isomorphism $\rho_f \colon \omega_Y \xrightarrow{\simeq} f^!\omega_X$, 
the corresponding $\MO_X$-module homomorphism is given by 
\[   
T_f \colon f_*\omega_Y \xrightarrow{f_*\rho_f, \simeq} f_*f^{!}\omega_X \xrightarrow{\mathrm{adj}} \omega_X.
\]
In this case, we get the following commutative diagram: 
    \begin{equation}\label{eq:f-diag}
        \begin{tikzcd}
            f_*\omega_Y \arrow[r,"T_f"] \arrow[d,"f_*\rho_f"', "\simeq"] & \omega_X \\
            f_*f^{!}\omega_X \arrow[d,equal] \arrow[ru, bend right,"\mathrm{adj}"] \\
            \cHom_{\cO_X}(f_*\cO_Y,\omega_X). \arrow[ruu,bend right,"(\varphi \mapsto \varphi(1))"']
        \end{tikzcd}
    \end{equation}
    \item 
    We  define the adjoint bijections $\Psi_g$ and $\Psi_h$ in the same way as the one of $\Psi_f$. 
    For isomorphisms 
    $\rho_f \colon \omega_Y \xrightarrow{\simeq} f^!\omega_X$ and 
    $\rho_g \colon \omega_Z \xrightarrow{\simeq} g^{!}\omega_Y$, 
    the equation
    \[
    \Psi_h(g^{!}\rho_f \circ \rho_g)=\Psi_f(\rho_f) \circ f_*\Psi_g(\rho_g)
    \]
    holds by a diagram similar to (\ref{eq:diag}). 
    Therefore, for another isomorphism $\rho_h \colon \omega_Z \xrightarrow{\simeq} h^{!}\omega_X$, 
    the equation $g^{!}\rho_f \circ \rho_g=\rho_h$ holds if and only if 
    the equation $\Psi_f(\rho_f) \circ f_*\Psi_g(\rho_g)=\Psi_h(\rho_h)$ holds. 
\end{enumerate}
\end{remark}

\bibliographystyle{skalpha}
\bibliography{bibliography.bib}

\end{document}